UNIVERSITY OF QUEBEC IN MONTREAL

# BRAVE NEW CATEGORICAL SPECTRAL POSITIVE SCHUBERT GEOMETRY

# AND THE CATEGORICAL DUAL AMPLITUHEDRON

## APPLICATION TO THE CATEGORICAL EXISTENCE AND THE DE RHAM VOLUME OF

## THE COHESIVE DIFFERENTIAL SPECTRAL ALGEBRAIC SPACE OF THE DUAL AMPLITUHEDRON

DISSERTATION

SUBMITTED IN PARTIAL SATISFACTION OF THE

REQUIREMENTS FOR THE DEGREE OF

DOCTOR OF PHILOSOPHY

IN

MATHEMATICS

BY

JULIEN DALPAYRAT-GLUTRON

Author's note: *This ArXiv preprint of my doctoral dissertation is intended both to lay the groundwork for a series of papers that build upon one another, paving the way for the future publication of this series in reviewed journals with referee — and to confirm the existence of a first proposed solution to the twelve-year-old problem posed to the mathematical community by Princeton physicists Nima Arkani-Hamed and Jaroslav Trnka, in a landmark 2014 article [13] in the world of theoretical physics, namely that of constructing a dual to the Amplituhedron — a purely mathematical object they invented to encapsulate, within a single mathematical object, the full complexity of Feynman diagrams and the calculation of scattering amplitudes in high-energy particle colliders.*

*Since ArXiv preprints are subject to requirements regarding mentoring (in this case, Prof. François Bergeron, UQAM, Quebec, Canada), non-plagiarism, the authenticity of the research, the novelty and relevance of the preprint to the mathematical community, and the author's ownership of the work, new ideas, and the manuscript, it is therefore particularly relevant to make this document already available on ArXiv — even before the acceptation by the doctoral thesis committee — as it is now accessible and freely available in accordance with ArXiv's rules, while also being protected by ArXiv's copyright policies.*

*Any mathematical errors, typographical errors, or other mistakes are the sole and entire responsibility of the author.*

*Thank you for your understanding, and please exercise caution when using the results of this first proposed solution to the* "Dual Amplituhedron".

## 0.1.Abstract.

This thesis is divided into two parts. The first deals with what is called "*A brave new spectral positive Schubert geometry*" which consists of rewriting the positive real Grassmannian and performing the singular gluing of its positroid varieties in a new way, allowing for a new interpretation via "*Spectral Algebraic Geometry*" on $\mathbb{E}_\infty$-ring spectra and spaces "$\mathcal{G}$-*structured*" in the sense of Jacob Lurie, finding a compact yet holistic formulation via "*intersection complexes*" in the "*theory of perverse sheaves*" with enough perfect perverse sheaves to cover the $\infty$-perverse geometry.

This algebraico-geometric perspective is paired with a synthetic differential-geometric perspective, characterized by the differential cohesion of Bill Lawvere and Urs Schreiber, which addresses the question of formal moduli functors, infinitesimal thickenings, De Rham stacks and their crystalline cohomology, which subsumes the "*modalities of structured geometries*" and the unification of their cohomologies of the underlying concrete topological étale algebraic space and its deformation theory.

The second part uses this rewriting of positive Schubert geometry and develops the "*Tannaka duality for the quasi-coherent stacks*" already present in the first part, as well as "*the cohomology of Brauer-Grothendieck spaces*", to show that, instead of the Grassmannian spectral Deligne-Mumford stack which is autodual in the $\infty$-category of $t$-structured monoïdal prestable $\infty$-categories of modules on $\mathbb{E}_\infty$-ring spectra, the Amplituhedron — an object introduced in 2013 by physicists Nima Arkani-Hamed and Jaroslav Trnka in the context of high-energy particle colliders within the $\mathcal{N}$=4 SYM framework — , viewed as a spectral algebraic space, possesses a concrete dual, which is non-trivial, and has a concrete volume, facts of interest for the duality between the Standard Model of particles and String Theory. This new construction yields the pure mathematical object that theoretical physicists have been awaiting for the past dozen years, known as the "*Dual Amplituhedron*".

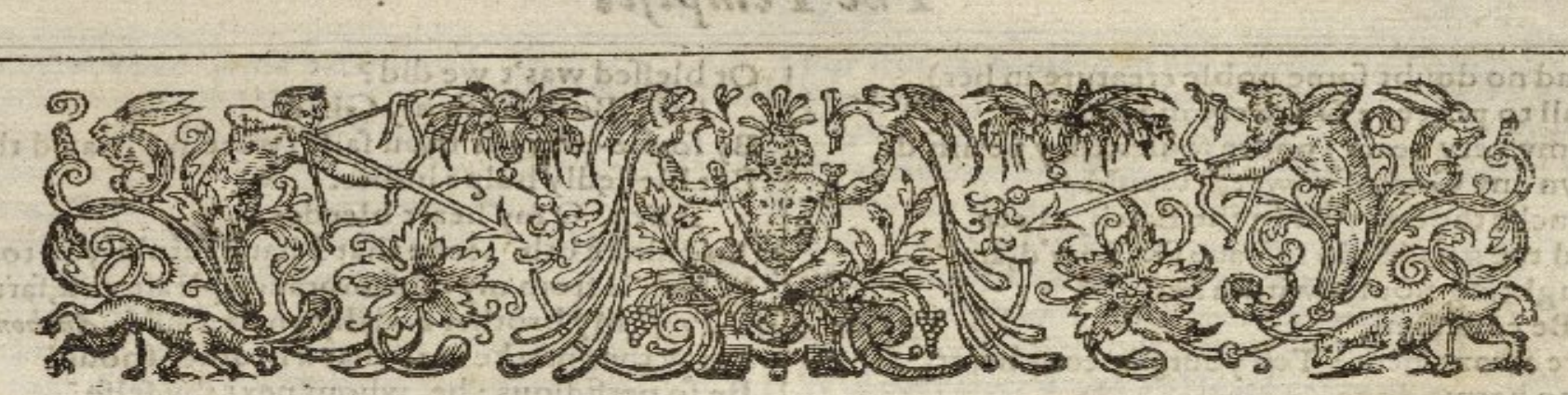

# THE TEMPEST.

## *Actus primus, Scena prima.*

*A tempeſtuous noiſe of Thunder and Lightning heard: Enter a Ship-maſter, and a Boteſwaine.*

*Maſter.*

BOte-ſwaine.

*Boteſ.* Heere Maſter: What cheere?

*Maſt.* Good: Speake to th'Mariners: fall too't, yarely, or we run our ſelues a ground, beſtirre, beſtirre. *Exit.*

*Enter Mariners.*

*Boteſ.* Heigh my hearts, cheerely, cheerely my harts: yare, yare: Take in the toppe-ſale: Tend to th'Maſters whiſtle: Blow till thou burſt thy winde, if roome enough.

*Enter Alonſo, Sebaſtian, Anthonio, Ferdinando, Gonzalo, and others.*

*Alon.* Good Boteſwaine haue care: where's the Maſter? Play the men.

*Boteſ.* I pray now keepe below.

*Anth.* Where is the Maſter, Boſon?

*Boteſ.* Do you not heare him? you marre our labour, Keepe your Cabines: you do aſsiſt the ſtorme.

*Gonz.* Nay, good be patient.

*Boteſ.* When the Sea is: hence, what cares theſe roarers for the name of King? to Cabine; ſilence: trouble vs not.

*Gon.* Good, yet remember whom thou haſt aboord.

*Boteſ.* None that I more loue then my ſelfe. You are a Counſellor, if you can command theſe Elements to ſilence, and worke the peace of the preſent, wee will not hand a rope more, vſe your authoritie: If you cannot, giue thankes you haue liu'd ſo long, and make your ſelfe readie in your Cabine for the miſchance of the houre, if it ſo hap. Cheerely good hearts: out of our way I ſay. *Exit.*

*Gon.* I haue great comfort from this fellow: methinks he hath no drowning marke vpon him, his complexion is perfect Gallowes: ſtand faſt good Fate to his hanging, make the rope of his deſtiny our cable, for our owne doth little aduantage: If he be not borne to bee hang'd, our caſe is miſerable. *Exit.*

*Enter Boteſwaine.*

*Boteſ.* Downe with the top-Maſt: yare, lower, lower, bring her to Try with Maine-courſe. A plague ——

*A cry within.* *Enter Sebaſtian, Anthonio & Gonzalo.*

vpon this howling: they are lowder then the weather, or our office: yet againe? What do you heere? Shal we giue ore and drowne, haue you a minde to ſinke?

*Sebaſ.* A poxe o'your throat, you bawling, blaſphemous incharitable Dog.

*Boteſ.* Worke you then.

*Anth.* Hang cur, hang, you whoreſon inſolent Noyſemaker, we are leſſe afraid to be drownde, then thou art.

*Gonz.* I'le warrant him for drowning, though the Ship were no ſtronger then a Nutt-ſhell, and as leaky as an vnſtanched wench.

*Boteſ.* Lay her a hold, a hold, ſet her two courſes off to Sea againe, lay her off.

*Enter Mariners wet.*

*Mari.* All loſt, to prayers, to prayers, all loſt.

*Boteſ.* What muſt our mouths be cold?

*Gonz.* The King, and Prince, at prayers, let's aſſiſt them, for our caſe is as theirs.

*Sebaſ.* I'am out of patience.

*An.* We are meerly cheated of our liues by drunkards, This wide-chopt-raſcall, would thou mightſt lye drowning the waſhing of ten Tides.

*Gonz.* Hee'l be hang'd yet,
Though euery drop of water ſweare againſt it,
And gape at widſt to glut him. *A confuſed noyſe within.*
Mercy on vs.
We ſplit, we ſplit, Farewell my wife, and children,
Farewell brother: we ſplit, we ſplit, we ſplit.

*Anth.* Let's all ſinke with' King

*Seb.* Let's take leaue of him. *Exit.*

*Gonz.* Now would I giue a thouſand furlongs of Sea, for an Acre of barren ground: Long heath, Browne firrs, any thing; the wills aboue be done, but I would faine dye a dry death. *Exit.*

## Scena Secunda.

*Enter Proſpero and Miranda.*

*Mira.* If by your Art (my deereſt father) you haue
Put the wild waters in this Rore; alay them:
The skye it ſeemes would powre down ſtinking pitch,
But that the Sea, mounting to th' welkins cheeke,
Daſhes the fire out. Oh! I haue ſuffered
With thoſe that I ſaw ſuffer: A braue veſſell

A (Who

*Shakespeare The Tempest (1610)*

If by your art, my dearest father, you have
Put the wild waters in this roar, allay them.
The sky, it seems, would pour down stinking pitch
But that the sea, mounting to th' welkin's cheek,
Dashes the fire out. O, I have suffered
With those that I saw suffer –
A brave vessel
(Who had no doubt some noble creature in her)
Dashed all to pieces. O, the cry did knock
Against my very heart! Poor souls, they perished.
Had I been any god of power, I would
Have sunk the sea within the earth or ere
It should the good ship so have swallowed and
The fraughting souls within her.

(olde english)

*If you have used your magic powers, my dear Father,*
*To make the sea so dangerous, please use them to calm it down.*
*The sky is so black that it looks like tar would rain down,*
*But the waves are going so high,*
*All the way to heaven, that they could cool the boiling tar down.*
*Oh I feel the pain of those that I see suffering!*
*A great (brave) ship*
*(no doubt with loving souls on board), wrecked and destroyed!*
*Oh their cries hurt my heart!*
*Those poor people have all died!*
*If I were a God with any power at all*
*I would have made the sea go beneath the earth*
*Before it could swallow that ship and all the souls aboard.*

*(modern langage)*

*Miranda, Monologue, Acte I, Scène 2*

# Contents

# Chapter 1

# Introduction

## 1.1. Foreword.

Preliminary remark: Category Theory is not "*Abstract Nonsense*". This thesis provides further, non-trivial proof of this fact, since the proposed solution for the dual of the Amplituhedron is derived using purely categorical arguments. Let's see why and how.

Reconstructing a whole from its subparts is one of the oldest problems facing mathematicians, philosophers, archaeologists, geologists, ethnologists or physicists of particles colliders, in all fields of the so-called "*exact*" sciences, as well as in the so-called "*human*" sciences. One could say that it is the strategies used to reconstruct this whole that define the fields of study of all these disciplines and distinguish them from one another.

In mathematics, the Grassmannian stratified into Schubert cells is one such mathematical object, and the "*ramifications*" in the mathematical sense, or the "*plate tectonics*" for geologists, that is the intertwining of these cells, are a subject of study in themselves, which has given rise to much mathematical literature and will continue to do so.

One of the problems lies in the braiding of the very definition of this Grassmannian, which is both an algebraic variety and a combinatorial object, a linear algebra structure and a differential manifold. Its definition itself is multiple and reflects this intertwining of its origins. In general, matrices and "*Young tableaux*" are used for the combinatorics of its stratification into the so-called "*Schubert cells*", and "*Plücker equations*" for its algebraico-geometric definition.

More recently than the 19th century, which saw the birth of the so-called "*Schubert calculus*", the domain of "*cluster algebras*" and the combinatorics of the positivity of the totally non-negative Grassmannian burst onto the scene at the end of the 20th century.

At the same time, the science of reconstructing algebraico-geometric objects was advancing with N. Saavedra Rivano reconstructing a group scheme from the category of its representations ([125]), with P. Gabriel reconstructing a smooth Noetherian scheme from the category of its coherent sheaves ([45]), with A. Rosenberg removing the Noetherian condition ([123, 124]), with A. Joyal and R. Street ([62]) and S. L. Woronowicz ([153, 152]) reconstructing a compact quantum group from a C*-category of comodules over Hilbert spaces in the framework of noncommutative geometry (see also A. Connes [32]), with P. Balmer reconstructing a scheme from the triangulated tensor category of its perfect complexes ([20]), and H. Krause, A. B. Buan, and Ø. Solberg reconstructing a ringed space, that is a "*spectrum*" with a sheaf of rings, from the ideal lattice of thick tensor ideals of an abelian or triangulated tensor category with a tensor identity. The word is out: "*spectrum*". We will come back to this at length.

At the same time, M. Hakim ([53]), A. Joyal ([61]), P. Johnstone ([59, 58]), C. Simpson ([139]), C. Rezk ([116, 117, 118]), B. Toën ([148, 146, 149, 150]), G. Vezzosi ([149, 150]), D. Cisinski ([30]), E. Riehl ([121, 120, 122, 119]) et many others are interested in toposes, in $\infty$-categories, and in $\infty$-toposes.

However, the issue of whether or not the "*recollement*" of the substructures was smooth was not addressed. But derived algebraic geometry and deformation theory emerged through the work of B. Toën, G. Vezzosi, J. Lurie, and others.

In 1982, Bernstein, Beilinson, Gabber, and Deligne had addressed the question of coding the stratification of substructures using perverse sheaves theory.

Finally, J. Lurie synthesizes, develops, and produces encyclopedic knowledge ([102]) on $\infty$-topoi, in particular, Chapter III of his (pre-published) book "*Spectral Algebraic Geometry*" [103] entitled: "*Tannakian Reconstruction and Quasi-Coherent Stacks*" the latest version of which dates from 2018. His work on Chapter 11 is particularly related, that is to say, both: derived from, and: extending, that of B. Toën and B. Antieau on "*Quasi-Coherent Stacks for Spectral Deligne-Mumford Stacks*" and "$\infty$-*Categories of modules over* $\mathbb{E}_\infty$*-rings spectra in an* $\infty$*-category of stable linear* $\infty$*-categories*".

J. Lurie constructed spectral algebraic geometry, which allowed the structure of perverse sheaves to be integrated as a geometric $\mathcal{G}$-structure into algebraic derived geometry dedicated to the study of non-smooth structures, thereby creating what is now called a "*structured space*".

Finally, thanks to Urs Schreiber's differential (infinitesimal) cohesion, thanks to prismatic cohomology (and the infinitesimal site), to A. Kock's synthetic differential geometry, and to M. Bunge's synthetic differential topology — indeed, thanks to all these contributions — the foundations were laid for the emergence of a complete theory of the stratification of Schubert cells and positroid varieties of the TNN Grassmannian.

At the very beginning of this story, a problem arising from theoretical physics emerged, namely the existence of a dual to the physico-mathematical object: the Amplituhedron.

In short, the framework of Deligne-Mumford spectral stacks structured by sheaves of perfect perverse intersection complexes $\mathcal{IC}^\bullet$ on $\infty$-topoi, and categories dualizable in a larger category $\mathsf{Groth}_\infty$ of $\infty$-categories of modules over spectra of $\mathbb{E}_\infty$-rings, was adopted here to solve the problem.

We will explain more precisely what this problem consisted of and how we arrived at the solution in this way via dual coaffine geometric stacks and module geometries in the sense of J. Lurie as perfect perverse structural sheaves arising from the BBDG perverse sheaf geometry encoding the singularities of spaces, stacks and their categories of quasi-coherent sheaves and their quasi-coherent stacks, here respectively 0-affine and 1-affine, allowing for double Tannakian reconstruction (a double decategorification) into a classical étale algebraic space and thus into a topological space equipped with a differential form and a De Rham volume differential form.

To conclude this introduction, we will discuss how this approach and framework can lead to advances in the study of cellular ramifications in Schubert geometry of the positive Grassmannian. This will constitute Part I of this thesis:

*Part I: Rewriting the stratified structure of the positive Grassmannian in a new, compact yet exhaustive and evocative way, enabling the resolution of problems that had previously remained unanswered.*

Part II is then an immediate consequence of this rewriting of the positive Grassmannian:

*Part II: Using this new compact and evocative formulation of the stratified structure of the positive Grassmannian to solve the twelve-year-old problem of the dual of the Amplituhedron, specifically the problem of its existence, its BCFW triangulation, and its volume, all three of which have remained unanswered until now in the context of the* $\mathscr{N} = 4$*SYM planar perturbative Super Yang-Mills theory.*

☆ ☆ ☆

***It should be noted at this point that the relative ingenuity required in a dissertation — the kind that, through profound results and rigorous proofs, wins the jury's approval to award, or withhold, the degree of Doctor of Mathematics — is not formulated here in the usual way.***

***Without wishing to presume to refer to illustrious authors who recount Bourbaki's lectures, where they wondered at what point during A. Grothendieck's four-hour lecture the passage would become difficult, since the title was difficult, and there had to be a moment when it would be hard to understand, and realizing at the end of the four hours that through a well-chosen categorification, the initial intention announced in the title had become, if not trivial, at least easy to understand, since it had become a natural consequence of the entire lecture — here too, in all modesty of course, the same phenomenon occurs.***

***Indeed, once the positive Grassmannian has been written in the correct form (the one concluding Part I), then the dual of the Amplituhedron follows naturally, as a consequence of the relative ingenious formulation in Part I. There is then nothing further to prove, hence the notable absence of lengthy proofs. The double categorification in Part I has already done all the work.***

***The work is exactly building the structure.***

***It is therefore no surprise that Part II is a list of results obtained directly without proofs. There will undoubtedly be exceptions to this rule during the review process, and I will be happy to write the proofs that will nevertheless be necessary for a proper understanding of the manuscript.***

(a humble note from the author)

☆ ☆ ☆

The initial motivation for this work is the construction of a dual to the Amplituhedron, an object from mathematical physics invented in 2013 by N. Arkani-Hamed and J. Trnka to unify a certain approach to scattering amplitudes in the standard model of particles in quantum physics, in the $\mathscr{N} = 4$ Super Yang-Mills ($\mathscr{N} = 4$ SYM) context, planar and perturbative. From the first articles devoted to this subject ([13, 12, 10, 11]), the necessity of the effective and constructive existence of a dual to the Amplituhedron appeared, and since then, numerous results have been obtained in this direction, and beyond, the study itself of the algebraic combinatorics of the underlying objects has made spectacular progress.

In our approach to the subject, and the attempt here to find a solution to the problem of constructing this dual of the Amplituhedron, and even if the conclusion of this work will end up with arguments of purely algebraic combinatorics, the genesis of the argumentation has its source in the type of problem it was chosen to try to solve.

*Is it a problem of the duality of an object whose dual would not yet exist in the mathematical bestiary of existing duals? Or is it a problem of reconstructing a certain dual space with certain combinatorial properties to be respected that come from the starting space?*

In fact, the two problems have to be solved simultaneously.

Within the first problem, many dualities already exist, all more or less related to each other depending on the level of complexity of the object that can be dualized. The invitation is thus extended to follow along this degree of complexity that will occur in the course of this article: the duality of a vector space, the duality of a (co)tangent space, the duality of a presheaf, Isbell duality (or adjunction), the duality of a sheaf, the duality of Serre, Grothendieck, Verdier, the duality of a scheme, the duality of an Artin or Deligne-Mumford stack, the duality "*in*" an $\infty$-category or an $\infty$-topos, and finally: the duality "*of*" an $\infty$-category in a bigger $\infty$-category of certain $\infty$-categories.

The preceding list is a little long but exhaustive for the problem that concerns us and for which an attempt at a solution is proposed in this article. In the course of this list, we begin to see the right phenomenon emerge, the one that is universal in the sense of category theory. What we want to dualize is always a dualizable object, or a complex of dualizable objects, in a (possibly $\infty$-)category that is stable by dualization and contains dualizable objects.

In other words, we take the point of view that A. Grothendieck teaches us: "*Everything is a functor*", and the desired duality must therefore be a functor. Ideally, we look for a functor that commutes in the sense of the other constraints, and therefore of the other functors, of the problem posed. Here, it is then a question of finding the right category of dualizable objects and determining which is the right duality functor that commutes with the structure-functors of the dualized objects.

In the second problem, numerous reconstruction theorems we saw above already exist, with a long list like the previous list of dualities. The right Tannaka reconstruction we had to choose was then the one corresponding to the dualization functor we were going to construct and which commutes with it.

This thesis manuscript indeed attempts to answer the purely mathematical question of the existence and construction of a dual to the Amplituhedron, a question first posed by Nima Arkani-Hamed and his student J. Trnka in 2014 in the following article [13], which marked a milestone in the scientific literature, initially in the literature of particle physics and cosmology, before finding its way into the mathematical literature and joint reflections among physicists, mathematicians, and even philosophers.

This question was restated and developed during his "*11-hour marathon*", so named, at UC Davis on scattering amplitudes in 2019, "*marathon*" to which we will return from time code 7:45 to time code 8:15, because it is ultimately at this precise moment, beyond Nima Arkani-Hamed's justification of the importance of the existence of this "*dual of the Amplituhedron*" through the duality it would bring to string theory and quantum field theory (QFT), that the idea emerges, the phenomenon, the concept, or whatever and however one might name it—in any case, the "*question ultimately left unresolved*" in a satisfactory manner regarding "*the emergence of spacetime and the fundamental forces*" from the vacuum of Anti-DeSitter space.

It is difficult in a mathematics thesis to answer Nima-Arkani-Hamed: "This has, in fact, already been answered, "*but by philosophers*" ! And in particular, this is exactly what the concept of "*Aufhebung*" by F.H. Hegel prescribes in his work "*Science of Logic*" (1807)", according to this term, which is difficult to translate and is generally used as such in the original German language, meaning "(*roughly*)": "*negation in overcoming*", with this notion intrinsic to the phenomenon known as "*the destruction of something at a certain stage to make way for a creation at the next stage in dialectical evolution*", here of the Universe of "*things*" and the "*ideas*" that "*surround*" us.

It is of course understood that all these terms would need to be defined very precisely in a text, whether philosophical or even mathematical, which is not the purpose here. Let us therefore attempt a gradual shift toward the domain that concerns us, mathematics, and thus reformulate the question within the history of mathematics itself, to arrive at a statement of the question that is as modern as possible, and to which we know how to respond.

This word "*Aufhebung*", which is at the origin of "*dialectics*", and thus of the modernity of science's evolution since the 19th century, is nevertheless present as early as B. Lawvere's first texts in his comments on category theory, complementing (in his own way) the definition of A. Grothendieck and J. -L. Verdier of a "*topos* "in 1963, the date of the first definition of a topos, which would become that of a "*Grothendieck (1-)topos*", yielding in B. Lawvere a "*elementary topos* ".

Since the "*elementary topos*" led B. Lawvere to what he considered to be "*improvements*" on the basic concept, he introduced this notion, which was demonstrated (by Giraud) to be even more abstract, as it fell more within the realm of categorical logic, "*improving*", in his view, the "*behavior*" of Grothendieck's topos relative to the nascent categorical logic, which has since been shown to be precisely the conceptual benefits of A. Grothendieck's definition in the sense of liberation through a beneficial abstraction, rather than an encapsulation within a more restrictive framework imposed by B. Lawvere's categorical logic, which is in fact just as "*liberating*"... but equivalent with respect to Giraud's axioms (see [98] for a demonstration of the nuances in the equivalence between "*A. Grothendieck's topos*" and "*B. Lawvere's elementary topos*" with respect to Giraud's axioms).

Let us now be more precise: In a Grothendieck topos, a monomorphism that is also an epimorphism is not necessarily an isomorphism. In other words, the topos, which is a 2-category, is not necessarily "*balanced*". Worse still: topoï constructed from a topological space are very rarely balanced, whereas the 2-category of presheaves on the lattice of open sets of a topological space is very often used as an introductory example (for pedagogical purposes) of what a presheaf is and then of what the topos of sheaves on a topological space is. This topos is not balanced; monomorphisms that are epimorphisms are not isomorphisms, and consequently the double negation is not the affirmation. In other words, the internal logic of this topos, which "*naturally*" serves as a foundational example, is not "*Boolean*" logic—that is, classical logic in the sense upon which our usually trivial everyday understanding of life is based. Indeed, the open set algebra is not a Boolean algebra, since the complement of an open set is not generally an open set. Moreover, a topos may have "*not enough points*", or "*no points at all*". However, these topos are still Heyting algebras, and these topos do indeed have an internal categorical logic. It is called in this case: "*intuitionistic*" or "*constructive*" (the proofs must be constructive, and there is no law of excluded middle nor double negation elimination).

The fact that double negation is not an affirmation leads to this "*negation in overcoming*", this "*Aufhebung*", and thus to dialectics in the epistemological sense that we now attribute to all evolving, historical, and social sciences, and of course to the exact sciences, including mathematics. Not that mathematics had waited for F.H. Hegel to become evolutionary, but the historical framework of the sciences could then be completed and rest firmly on mathematics in constant evolution.

On the other hand, that this "*modality of Being*" — since that is the name given by F.H. Hegel, and later by B. Lawvere, to this concept of "*transcendence through negation*" — could finally be formalized in mathematics, and here in category theory, or better yet: in topos theory, whether Grothendieck or elementary topos, was a real breakthrough in favor of category theory and topos theory as a framework conducive to the formalization of the foundations of mathematics. The formalization of $\infty$-categories by A. Joyal and J. Lurie, to name just a few for this introduction to homotopic theories, followed by the connection with Homotopic Type Theory (HoTT) by Voevoedsky, M. Schulman, E. Riehl, D.-C. Cisinski, A. Kock, M. Bunge — again, to name just a few — constitute the most current phase of what is now called "*synthetic mathematics*", that is, mathematics without external intervention, but arising from its very depths to re-found itself.

Finally, U. Schreiber made the connection with theoretical physics, in a joint advancement of mathematics and physics, of the logical and $\infty$-categorical foundations of $\infty$- topoi, monads, and comonads (which also have their significance in Plato, Leibniz, Husserl, and Deleuze in philosophy, as mathematical terms often have their origins in the history of science or philosophy), arising from the addition of functors linked to the differential cohesion of local $\infty$-topoi, to the calculation of the canonical volume form in positive geometry, and to the redefinition of the spaces within which modern physics evolves.

U. Schreiber thus undertook the translation of gauge fields, cobordism, and Chern-Simons theory into the types of spaces within which relativistic, quantum, and unified physics had to evolve to align with the new types of spaces being studied by derived algebraic geometry. Work will therefore need to continue to integrate into his work the modular spaces of derived $\infty$-stacks, the homotopy of noncommutative spaces, and spectral geometry on the $\mathbb{E}_\infty$-rings.

But for now, the task is already to provide, through mathematics, the possibility of a physical duality between Field Theory and String Theory, which would validate the latter as a theory describing the World — that is, to provide the possibility of a unification of fundamental forces, particularly quantum and gravitational forces, in an Anti-de Sitter space that would thus be empty at the origin of the Big Bang singularity, but whose spacetime should "*emerge*", just as an adjunction of a mathematical (co)monad brings forth a structure for which we have previously been given the construction plan for all "*possible worlds*", which is, both in its function and its denomination, another modality of modal logic in Homotopic Type Theory; these "*possible worlds*" are, in turn, demonstrably "*emergent*" in our particular situation.

*The aim of this thesis is thus to provide a strict formulation* and a "*resolution of the duality of the Amplituhedron*", a phenomenon of duality that is perfectly algebraico-geometric in the sense of A. Grothendieck, but singular, in the sense that it is non-smooth for the theory of intersection and singularities in the sense of the modular deformation theory of B. Toën and J. Lurie, yet sufficiently "*differential and cohesive*" in the sense of B. Lawere and U. Schreiber to give volume to the de Rham space that is the dual constructed here of the Amplituhedron, which is expected by L. Williams, T. Lam, and finally Nima Arkani-Hamed, for his invention of the Amplituhedron, initially combining Schubert's algebraic combinatorics, the volumes of positive Grassmannians, and the motivic theories of generalized cohomology.

The work of J. Lurie in spectral algebraic geometry on Brauer-Grothendieck spaces and the derived Azuyama algebras introduced by B. Toën and B. Antieau within the framework of the stabilization of coconnective prestable $\infty$-categories, such as the categories of derived quasi-coherent coaffine $\infty$-stacks, categorifications/decategorifications, and adjoint $n$-affine Tannakian dualities that allow for "*unfolding the singularities*" of this singular Amplituhedron, will then be joined, in particular through the work of U. Schreiber in his continuation of B. Lawvere's work on the modal logic of spatial, cohesive, and differential $\infty$-toposes, to give a De Rham differential volume form to this dual Amplituhedron.

Indeed, this Amplituhedron is combinatorially stratified in every degree thanks to A. Postnikov for the topological work on étale overlaps enabled by L. Williams following Even-Zohar with the BCFW triangulation, making the positive geometries introduced by T. Lam the appropriate framework for defining a canonical volume for this dual of the Amplituhedron. To take a reconstructible dual that respects its "*geometries*" in the sense of J. Lurie [91], and thus the Schubert geometry of the gluing of $t$-structures of varieties derived from positroid cells, we must functorialize the initial proper topological map defining the Amplituhedron, allowing, between étale topologies of $\infty$-topoï, the respect for pullbacks and pushforwards of the categories of perfect and perverse quasi-coherent sheaves and stacks over $\infty$-Deligne-Mumford stacks emerging from the underlying topological spaces in a Grothendieck context and within the six-functor formalism. These $\infty$-Grothendieck topologies, here on compact, perfect stacks (so that they are dualizable) and perverse sheaves (so that they are stable by the Deligne-Beilinson decomposition theorem, a subvariety constructor), allow here for the preservation of the Lurie modular geometries of the subtopoi, which constitutes here the reformulation of Nima Arkami-Amed's request for respect via duality of Schubert's positroid cells for the expected object of the "*dual Amplituhedron*".

The amplituhedron $\mathscr{A}_{n,k,m}$ is the image, via a projection $\tilde{Z}$ induced by a positive linear map $Z$ which, a priori at the start of the study, is defined mathematically only by the total positivity of its principal minors, of the positive real Grassmannian $\mathrm{Gr}(k,n)^{\geqslant 0}$ , where $k$ and $n$ depend on the number of particles and their helicity, onto the smaller real Grassmannian $\mathrm{Gr}(k,k+m)$, where $n \geqslant k+m$, and where the case of the model chosen by physicists, in quantum collisions and scattering amplitude calculations, corresponds to $m=4$.

The need to construct a geometric and combinatorial dual to this image $\mathscr{A}_{n,k,4}$ arises from the need to model and calculate a certain "*meromorphic volume form*" with poles, namely the "*canonical form*" of N. Arkani-Hamed, T. Lam, and J. Trnka, which depends on this object and whose integral, according to the theory in [9], is equal to the sought-after scattering amplitude embodied by this Amplituhedron $\mathscr{A}_{n,k,4}$.

The fundamental mathematical object of the problem is therefore the Grassmannian $\mathrm{Gr}(k,n,\mathbb{R})$, which is initially real at the start of the study and is a well-known object dating back to the 19th century, closely linked to Schubert's calculus, to the stratification by Schubert cells $\Omega_\sigma$ of the Grassmannian, and to what we will call "*Schubert Geometry*" throughout this thesis in the many instances where this concept appears.

Then, as particle physics invites us to do, it is the positive part $\mathrm{Gr}(k,n,\mathbb{R})^{\geqslant 0}$ (in the sense of Lusztig) of the Grassmannian that is the object of our interest. By carefully identifying, with its associated large étale topos, its small real topos defined on a Grothendieck topology imposed on the semi-algebraic set defined by the positive Grassmannian, we are then reduced to the positive Grassmannian, but embedded in the complex Grassmannian Gr($k$,$n$,$\mathbb{C}$), which allows us to immediately, yet rigorously, return to the heart of algebraic geometry over the algebraically closed field $\mathbb{C}$ of the complex numbers.

This means that instead of having to work within real algebraic geometry, which behaves rather poorly — since the Nullstellensatz is false from the very origins of real theory, resulting in a mathematical literature and results that are scarce and lack true richness or the desired scope — we can then turn to classical algebraic geometry (then derived, then spectral) as a mathematical discipline as vast as imagination and inventiveness can dream of, from A. Grothendieck to the most recent theories concerning $\infty$-topoi, then the structured spaces of J. Lurie's spectral algebraic geometry, and finally an even higher categorification, with as the globalizing structure of the studied framework, the monoidal $\infty$-category of prestable presentable $\infty$-categories of modules over a spectrum of $\mathbb{E}_\infty$-rings.

In this "*large*" higher $\infty$-category, we can then define, in the sense of the tensor duality viewed as a "*categorified inversion*" — which is in fact the only existing duality to which all well-written dualities reduce — a duality of the Amplituhedron, which we will have categorized in the category of quasi-coherent sheaves on a "*functor*": the Deligne-Mumford spectral stack, in spectral algebraic geometry.

A well-chosen étale topology on the positive Grassmannian, consisting of an étale open covering arising from the $\mathcal{BCFW}$ triangulation by Even-Zohar et al. ([37]) and completed by L. Williams et al. for its closure ([39, 38]), will then allow the definition of the "*Amplituhedron functor*", which will ensure that the "*duality functor*" thus defined respects the $\mathcal{BCFW}$-layered structure of the positive Grassmannian. This constitutes a requirement in the resolution of the construction of this duality, which will then be ensured by the proper transfer of the $\mathcal{G}$-structure (in the sense of J. Lurie) given by the category of perverse sheaves on the $\mathcal{BCFW}$-stratified Grassmannian onto the $\mathcal{G}$-structure of the Amplituhedron, viewed as a spectral algebraic space.

It should be noted that this positive Grassmannian and this Amplituhedron, as spectral algebraic spaces, will have all the expected reconstruction properties, since these spectral algebraic spaces will be "*0-affine*", that is, Tannaka duality will allow their reconstruction from their categories of quasi-coherent sheaves. Furthermore, their categories of quasi-coherent sheaves will be "*1-affine*", which implies their compact generation by their canonical structural sheaves and thus their dualizability as categories of modules over the $\mathbb{E}_\infty$-ring spectrum given by the right-derived functor of their global sections functor.

We will thus see that the gluing of Schubert cells, which can exhibit highly non-trivial behavior—in the sense that they can be as ramified as permitted by the complexity of the underlying combinatorics of non-transverse cell intersections—is encoded in categories of perverse sheaves attached, as structures, to the objects "*positive Grassmannian*" and "*Amplituhedron*", viewed as algebraic stacks. This means that it is by successively categorizing the problem several times that these singularities are resolved, yielding, starting from a "*good étale topology*" of the positive Grassmannian, open positroid varieties that are quasi-affine, whose étale $\infty$-topos gives rise to a 0-affine stack ([7, 106]), on which the category of quasi-coherent sheaves is 1-affine ([46, 144]), etc. Each successive categorification increases $n$-affinity ([144]), and thus resolves the singularities of the underlying 1-localic topos, that is, the singularities of the underlying (topological) space itself.

This phenomenon of increasing $n$-affinity with an $n$-categorization should be viewed in light of the preservation of a space's topological information with its higher homotopy groups; the framework of spaces and (1-)Grothendieck toposes thus reaches its limits, since we will see a resolution of the problem of information loss in higher homotopies by moving to $\infty$-topoi, which themselves preserve this homotopic information, making them the appropriate framework—along with Grothendieck $\infty$-topologies on the appropriate $\infty$-étale sites—for resolving the problem at hand.

Returning to the history of Schubert geometry, and perhaps hoping to claim, through the use of $\infty$-topoi and spectral algebraic geometry, an extension of it via the methods used in this thesis, Schubert geometry has been the subject of intense research and numerous publications for a century and a half. Nevertheless, it continues to hold many secrets for mathematicians, particularly due to the universality of the Grassmannian object. One might then observe, from an initial overview of the mathematical literature devoted to the Grassmannian, that historically, it was the combinatorics of cells that was primarily studied in the early days of Schubert calculus in the 19th and 20th centuries, and then that gradually, a shift took place at the beginning of the 21st century toward algebraic geometry, thanks to the discovery of the cluster algebra structure of the coordinate ring of the positroid varieties stratifying the positive Grassmannian (see [136, 50]).

The combinatorics of the Grassmannian, which consists of solving, for example, the problems in question through discoveries and inventions of the most ingenious arrangements of Young tableaux, plabic graphs, decorated permutations, lattices, and many other combinatorial structures of remarkable inventiveness, is very elegant, and the power and beauty of this way of solving problems is spectacular (see [14, 114, 42, 151, 68, 50, 137, 43, 113, 49, 44, 115, 48, 63]).

It should be noted, as a supplement to this initial outline of the study of this object, that the Grassmannian is also, in a more general and thus more flexible sense, a homogeneous space (see [62] Section 5), as the quotient of a compact topological group by a closed subgroup, and that it is, on the other hand, and in a very structuring way this time, a truly beautiful purely algebraico-geometric object in the sense of scheme theory, just as rigid as it can be for combinatorics.

The main body of this thesis will therefore focus on resolving the construction of the dual of the Amplituhedron object, viewed from the perspective of its algebraico-geometric structure, while respecting its structuring in the sense of the dual stratification of Schubert's "*dual*" geometry projected from the positive Grassmannian stratified into positroid varieties.

To reach this conclusion, we will develop in the main body of this thesis, by higher categorification, via the spectral Deligne-Mumford $\infty$-stacks and "*spectral algebraic geometry*", the construction of a duality functor that commutes with, that is, respects, the structure of the dualized objects.

One of the objects of study in spectral algebraic geometry is J. Lurie's extension to "*structured spaces*" ([91]) of "*derived algebraic geometry*" ([148]), whose concepts, in turn, via "*deformation theory*" ([89]), in particular by B. Toën, study "*singular*" phenomena, that is, "*non-smooth*", that appear in classical algebraic geometry, and whose behavior had not yet been the subject of results published by algebraic geometers and category theorists.

The main body of the text therefore consists of the dualization, within the $\infty$-monoidal category $\mathsf{Groth}_\infty$ of prestable $R$-linear $\infty$-categories, where $R$ is a spectrum of $\mathbb{E}_\infty$-rings, $\infty$ -categories equipped with the J. Lurie's tensor product, of a stable $\infty$-category of modules over $\mathbb{E}_\infty$-ring spectra. This $\infty$-category, which is non-trivially invertible—that is, non-autodual—is the stable $\infty$-category of spectral quasi-coherent stacks on the $\infty$-topos of étale sheaves on our topological object to be dualized, the Amplituhedron. The Grothendieck $\infty$-topology used on this $\infty$-topos is the étale topology giving rise to étale cohomology, which will yield, as the second étale cohomology group, the extended Brauer-Grothendieck group of equivalence classes of extended Azumaya algebras ([103]

Section 11.5.3), that is, the equivalence classes of compactly generated, $R$-linear, stable $\infty$-categories that are invertible under the J. Lurie tensor product. The non-triviality of the invertibility, that is, the non-self-duality of the module category under consideration in the case of the positive Grassmannian and the Amplituhedron, is given by the non-triviality of the second étale cohomology group of the positive Grassmannian and the Amplituhedron, by a unique, non-trivial generator, namely the first Chern class of the positive Grassmannian and the Amplituhedron obtained by the algebraic combinatorics of the cohomology of the Grassmannian.

To provide a brief history of the research, we began by exploring the approach involving the reconstruction of a dual to a homogeneous space. However refined the Hopf-Galois theory on fiber functors may be, this solution does not succeed. Nevertheless, as it allows for the restatement and exploitation of the approach by A. Joyal and R. Street in their well-known article ([62]) on braids in compact quantum groups, we reserve the right to a future publication re-examining this perspective, accompanied by the formulations of classical Tannakian reconstructions by N. Saavedra-Rivano for schemes in affine groups and by S. L. Woronowicz in the context of the duality of compact and discrete quantum groups. A discussion will be presented of the Hopf-Galois theory governing the behavior of fiber functors, initially introduced by N. Saavedra-Rivano and later taken up by S. L. Woronowicz, K. H. Ulbrich, P. Schauenburg, and J. Bichon, to reconstruct an algebraic group and a compact quantum group from the rigid tensor category of its finite-dimensional (co)representations and the neutral forgetful functor to the category of finite-dimensional vector and Hilbert spaces, which paved the way for all Tannakian reconstructions of the 21st century.

There is a justification for mentioning this line of work here in the presentation of this thesis, even though there is no corresponding text. This justification lies in the genealogy of the search for the solution itself, which, admittedly, has not yet been successfully concluded (hence the absence of this unfinished text); however, it is a phenomenon—or rather, it is a hidden structure, long known to A. Grothendieck, which has emerged and has therefore been examined and utilized as a foundation for devising the "*path*" leading to the solution developed here, namely: the notion of "*topos*", and this in its various forms: "*0-topos* (a set)", "*1-topos* (a space)", an "*n-topos*" (a topos truncated at the $n^{\text{th}}$ homotopy group), and an "$\infty$-*topos* (in which all homotopy groups are preserved)".

Before the brief overview of the 20th-century history of space reconstructions, a useful comment should be made here regarding Grothendieck's homotopic hypothesis. The nesting mentioned above regarding $n$-toposes:

$$\tau_{\leqslant 0}\mathcal{X} \subset \tau_{\leqslant 1}\mathcal{X} \subset \ldots \tau_{\leqslant n}\mathcal{X} \cdots \subset \tau_{\leqslant \infty}\mathcal{X} = \mathcal{X}$$

is taken for granted here, since the spaces under consideration are sober as "*loci*" from algebraic geometry. There is thus no need to consider spaces with "*pathological*" homotopy. The difference between whether Grothendieck's homotopy hypothesis is considered an axiom or is proven is irrelevant in this thesis.

To resume the development of the solution to the problem posed, in fact, the initial question had turned into a reconstruction problem, whereas numerous reconstructions do indeed exist, mainly and commonly referred to as "*Tannakian*". To cite the main ones in chronological order, the Tannakian reconstruction by N. Saavedra Rivano focused on the tensor category of representations of schemes over affine groups and on fiber functors, P. Gabriel ([45]) and then A. Rosenberg ([123, 124]) reconstructed a smooth scheme via its category of quasi-coherent sheaves, P. Balmer ([20]) also reconstructed a scheme, but from the triangulated tensor category of perfect complexes over that scheme, and finally, A. Joyal and R. Street ([62]) and S. L. Woronowicz ([153, 152]) reconstructed a compact quantum group from a C*-category of comodules within the framework of noncommutative geometry (see also A. Connes [32]) . These major reconstructions during the second half of the 20th century paved the way for the second period in the 21st century, namely the Tannakian reconstructions that subsequently emerged from the world of Grothendieck stacks, particularly in the work of J. Lurie (see also the entire body of work by D. Schappi [133, 126, 134, 135]).

Even though this postponed attempt to reconstruct a homogeneous space has not yet led to a solution for the Amplituhedron, it is nonetheless worth noting that it has its own merit, particularly in that it involves the quantization of energy for scattering amplitudes, which is no trivial matter in the history of quantum mechanics during the 20th century. However, since this quantization of energy is here a consequence of the duality of compact quantum groups to discrete quantum groups, and the fundamental argument underlying this duality is the same as the historical one for the diagonalization of normal operators in noncommutative geometry ([140, 141, 142, 143]), which yields a discrete spectrum, there is no new progress on the subject except for the potential interest of a comprehensive restatement of the Tannakian reconstruction theorems from the early period. In fact, the failure to consider the Grassmannian solely as a homogeneous space stems from the fact that the object—as a compact quotient topological group—for which we seek the dual, lacks sufficient structure to recover Schubert's geometry. It must be considered as "*a compact quotient algebraic group in derived and spectral algebraic geometry*", rather than merely in the harmonic analysis of compact topological groups, which is the subject of the article by A. Joyal and R. Street on noncommutative geometry of braids and compact quantum groups.

At this stage, the progress made on this path toward reconstructing a dual to a homogeneous space—and thus one lacking sufficient structure—nevertheless allowed us to touch upon the essential and to understand what motivated the underlying, profound problem. Indeed, the field of topology is underlying and omnipresent, but particularly in its most advanced form of "*Grothendieck topologies*" well-chosen on a category, with subcanonical Grothendieck topologies on a topos being the fundamental issue underlying the problem.

The main body of this thesis is thus devoted to the algebraico-geometric approach, which is then explored in the context of derived and spectral algebraic geometry, with an étale $\infty$-topology on an étale $\infty$-topos—a field whose objects of study include singular, i.e., non-smooth, varieties, which is precisely our situation here.

Indeed, if we return to the initial problem, but this time considering the Grassmannian as an algebraico-geometric object, the Grassmannian itself is smooth, but the stratification of the positive Grassmannian occurs in positroid varieties that are glued together, but in a non-smooth manner, which must be reflected in the dual to be reconstructed. In other words, the dual of the Amplituhedron must respect the *Schubert Geometry* of the positroid varieties of the totally nonnegative Grassmannian. We then proposed to encode the Schubert geometry of the stratification into positroid varieties using Deligne's perverse intersection sheaves, whose category could then be grafted as an additional structure onto the derived algebraic geometry, giving rise to what J. Lurie calls "*Spectral Algebraic Geometry*".

Grothendieck's revolutionary ideas are omnipresent here: the fact that "*everything is a functor between categories*" and the omnipresence of homological theory ([51]), scheme theory (EGA1 [52] and subsequent volumes of EGA/SGA), topos theory ([16, 58]) and étale cohomology theory ([16, 107]) transformed the problem into a reformulation of the original problem, where, in order not to eliminate higher homotopy and higher cohomology groups, it was the theory of $\infty$-toposes that emerged as the appropriate framework for formulating the problem. More precisely, not only derived algebraic geometry, but further still, J. Lurie's spectral algebraic geometry — which considers Deligne-Mumford spectral stacks within the framework of the $\infty$-topos theory of an étale $\infty$-site—is deployed here along one of the possible avenues opened up: that of the theory of prestable $\infty$-categories of modules over spectra of $\mathbb{E}_\infty$-rings and their dualization in an even larger category, namely the monoidal category — for J. Lurie's tensor product—of prestable $\infty$-categories with a tensor product unit, the $\infty$-category $\mathcal{S}\mathrm{p}$ of spectra of $\mathbb{E}_\infty$-rings.

Here, it was also necessary to find the appropriate (1-)topos (space) underlying our work, namely a simplified covering of the positive Grassmannian, which allowed us to construct a proper Grothendieck étale $\infty$-topology on the $\infty$-étale topos of spectral Deligne-Mumford $\infty$-stacks, and which was given by the BCFW recursion of Even-Zohar et al. ([37]).

By transferring the intersection complexes to the world of $\mathbb{E}_\infty$-ring spectra, and considering J. Lurie's Tannakian reconstruction of spectral algebraic spaces ([93], and Part III of [103]), we concluded by dualizing the category of quasi-coherent stacks over the spectral algebraic space of the Amplituhedron equipped with the étale topology denoted here as "*BCFW étale topology*", proving the non-triviality of the extended Brauer-Grothendieck group (Toën [147], Lurie [103] Part III, Chapters 10 and 11, Antieau-Gepner [6, 7], Chough [27, 28]), consisting of the isomorphism classes of Azumaya algebras ([103] Chapter 11) obtained by considering the endomorphism algebras (and Morita theory [131, 130]) of simple local systems of intersection complexes given by the decomposition theorem for perverse sheaves.

We were then able to "*invert*" the spectral Deligne-Mumford $\infty$-stack, constructed from the usual Amplituhedron viewed as a functor on an $\infty$-topos of sheaves with the BCFW étale topology, into the extended Brauer-Grothendieck cohomology group, which is not trivial, since it is the second cohomology group of the spectral stack taking values in the invertible perverse intersection sheaves, and is equal to the first Chern class of the positive Grassmannian, by isomorphism of the morphism induced by the usual Amplituhedron map at the level of Brauer-Grothendieck cohomology groups.

In other words, thanks to the 0-affinity property, according to the terminology of B. Antieau ([7, 106]), of the spectral stack under consideration, and thanks to the properness of the Amplituhedron functor thus constructed, we have in fact dualized, in the $\infty$-category of prestable $\infty$-categories (so that the $\mathbb{E}_\infty$-ring spectra are coconnective and correspond to cohomologically positive graded differential algebras, and thus to simplicial algebras and effective topological spaces), a certain category of modules over the $\mathbb{E}_\infty$-ring spectra of the global sections of the underlying spectral Deligne-Mumford stack, via the Verdier duality applied to the right-derived functor of the global sections functor. The 1-affinity of the category of quasi-coherent sheaves on the underlying $\infty$-topos—that is, as a category compactly generated by the canonical sheaf—then allows for the Tannakian reconstruction of the dual 1-localic $\infty$-topos and thus the expected dual topological space.

Since the first Chern class of the positive Grassmannian, as a signature of the non-triviality of the inversibility, is non-zero, this mechanism is not self-dual in the case of a non-smooth spectral stack, unlike the case of a category of quasi-coherent sheaves on a smooth scheme, which is self-dual, and thus does not reconstruct any new category of modules, and therefore, by Tannakian duality, no new dual space other than itself. In conclusion, a Tannakian reconstruction argument of the underlying $\infty$-topos from its spectral stack then allowed us to recover the underlying spectral algebraic space, and thus the dual space of the Amplituhedron, which is quite concrete as a (sober) topological space, namely "*the dual Amplituhedron*".

The preservation of the Schubert geometry of the Amplituhedron is then naturally reflected in the dual via the correct transfer of the perverse structure of the $\infty$-topos of the positive Grassmannian to the perverse structure of the $\infty$-topos of the Amplituhedron. This correct transfer is ensured by the semismallness of the Amplituhedron map in the BCFW topology, which is a condition on the dimension of the fibers of the Amplituhedron map within the chosen étale BCFW topology, and the decomposition theorem for perverse sheaves, which preserves the semisimplicity of the endomorphism algebras in the abelian category of perverse sheaves. Since the duality functor then commutes with the structure functor of the $\infty$-toposes of the spectral algebraic spaces of the positive Grassmannian and the Amplituhedron, we obtain the universality of the "*duality functor*", and thus of the duality solution, which was expected in order to claim to have defined a dual to the Amplituhedron.

It should be noted that the dualization of the (co)tangent complex of the initial Amplituhedron map already allowed us to conclude the existence of a dual to the Amplituhedron, but compliance with Schubert geometry was not reflected in the (co)tangent complex in "*only*" derived algebraic geometry. This could only be achieved thanks to the Lurie-Artin representability theorem, applied in conclusion to a functor between $\mathcal{IC}$-structured $\infty$-topoi of étale sheaves, where $\mathcal{IC}$*-structured* means in this case, according to the terminology of J. Lurie's $\mathcal{G}$*-structured spaces*, the grafting of the $\infty$-category of perverse intersection sheaves as a $\infty$-structural sheaf onto the $\infty$-topos in question, which already possesses a canonical structural sheaf $\mathscr{O}$.

Let the dual of a singular algebraico-geometric object be given by the dual of the cotangent vector spaces to its singularities, and let the structure and cohomology of the extended Brauer-Grothendieck group, or more simply those of the extended Picard group, which are related to the category of representations of the fundamental groupoid, and which are also encoded in the structure of the category of perverse sheaves, are non-trivial here, since the gluing of Schubert cells in the positive Grassmannian is "*non-smooth*", are two very "*moral*" and ultimately "*reassuring*" regarding the validity of the solution and the universality of the adopted method of structured dualization.

Opening up new avenues of research, spectral algebraic geometry thus offers a fresh approach to Schubert geometry and other "*well-stratified*" situations in the sense of an algebraico-geometric definition of strata, which could then also, we hope, benefit from these techniques and concepts, just as spectral algebraic geometry is already widely used for the study of elliptic curves and their cohomology (see [86, 97, 99, 101, 35]). Perhaps algebraic combinatorics, through the geometrico-algebraic nature of Schur polynomials, MacDonald polynomials, and many other concepts in algebraic combinatorics, could benefit from this contribution.

We conclude this introduction with the following remark: The formalism — which might at times be described as relatively virtuosic — and the abstract nature of spectral algebraic geometry must not obscure the combinatorial phenomenon, which is itself very concrete — namely, that the constructed dual can be seen, for example, already on a planar polytope, such as a 2-simplex, a triangle, whose constructed dual is indeed the polar dual (or orthopolar dual, to be more precise).

## 1.2.Foreword in theoretical physics, Background, and Aknowdgements.

### 1.2.1 String Theory and Scattering Amplitudes.

#### 1.2.1.1 The Standard Model, Feynman Diagrams and QFT.

Among the unexplained phenomena in theoretical physics, during the second half of the twentieth century, the Standard Model of the quantum physics of particles was confronted with the problem of unification of the gravitational force to the first three forces of interaction acting on matter already united in the same theory, the Quantum Field Theory (QFT).

By the theory of gauge fields, Maxwell's electromagnetic force mathematically realizing a duality between electric field and magnetic field, had been unified, thanks to the gauge Lie group $U(1)$, to the strong nuclear interaction by the gauge Lie group $SU(2)$, then to the electroweak force by the gauge Lie group $SU(3)$.

However, the nonlocal gravitational force, unlike the other three forces of interaction, resisted unification, which gave rise to the named objective: the Unified Theory, or Theory of Everything, or M-Theory.

In the second half of the 20th century, the theory of superstrings emerged to build a broader theoretical framework, in which gravitational force could integrate with other fundamental laws of interaction.

Within Quantum Field Theory thus rised in theoretical physics, the theory $\mathscr{N} = 4$ Super Yang-Mills emerged as one of the most promising theoretical frameworks in the research carried out within the theory IIB of supersymmetric superstrings theory. Among the fundamental concepts in Quantum Field Theory, scattering amplitudes were previously formalised for calculations by the diagrams of R. Feynman introduced in the 1950's.

#### 1.2.1.2 Scattering amplitudes in planar $\mathscr{N} = 4$ SYM.

In Quantum Field Theory, the calculation of scattering amplitudes is an observable that must be performed to ground each new emerging theory [14].

In 2013, N. Arkani-Hamed and J. Trnka invented the concept of Amplituhedron [13, 12, 10, 11], which not only solved the problem of the exponential number of diagrams generated by R. Feynman's method to model the interactions between particles in Standard Model, but also gave a much more real physical meaning than the formalism of R. Feynman allowed.

In particular, this new theory turned out to carry a particularly innovative vision, since it leads even to redefine certain intrinsic laws of particle physics, such as locality and unitarity.

A better understanding of the calculation of scattering amplitudes was then needed at this step, and it emerged to be the volume of the dual of this new mathematical object, the Amplituhedron.

In Quantum Field Theory within theoretical physics, the calculation of scattering amplitudes is an observable that must be performed to ground each new emerging theory [14]. Attempts to unify the four forces, including gravitational force, led to supersymmetric superstrings theory, and within superstrings IIB theory, one of the most relevant is now given by the planar model $\mathcal{N} = 4$ SYM.

#### 1.2.1.3 The dual Amplituhedron whose existence allows computations.

In the initial paper on the Amplituhedron in 2013 [13] and in numerous papers since [12, 10, 11], see also [41, 40, 54, 69, 55, 15], it was suspected the existence of this « *dual Amplituhedron* » to explain the positivity in the apparent behavior, foreshadowing that the calculation of the scattering amplitudes would be reduced to measure the volume of this « *dual Amplituhedron* ».

The BCFW recursion and the associated BCFW triangulation led to the possibility of calculating these scattering amplitudes, on the condition of having the right to do these calculations, that is to say that the « *dual Amplituhedron* » does indeed exist.

### 1.2.2 Backgroung.

Principally, homotopy theory, Lie theory, and topos theory and $\infty$-topos theory are assumed to be known.

This means that stable homotopy theory, Functor Loop and Suspension, Spectrum, and Cohomology Theories and assumed to be known.

This implies that the concept of adjointness is known, that the categories of $\mathbb{E}_\infty$-rings and $\mathbb{E}_\infty$-ring spectra are known, with their homotopies $\pi_k(R)$ where $R$ is an $\mathbb{E}_\infty$-ring, and that $H * (-, \underline{A})$ is a representable functor where $\underline{A}$ is a constant sheaf of abelian groups.

This implies as well that the Grothendieck construction in the framefork of fibered categories is known, with the $\infty$-version of J. Lurie functor of straightning/unstraightnening, allowing the construction of (co)cartesian (co)fibrations, of $\infty$-adjoints, $\infty$-monoidal categories, and so on.

We end this section with the knowledge of Lie theory with respect to formal moduli problems. However, Lie theory will not be of interest for us here.

The comprehension of the cotantangent complex will instead be very useful and we will return on it.

For the spectra categories, the stables $\infty$-categories and the $\mathbb{E}_\infty$-rings, we recommand:[87, 88, 132, 57, 26]. For the cocartesian fibrations, we recommand: [98]. For Serre-Grothendieck-Verdier duality: [100, 96, 83, 110]. For Lurie's tensor product between prestable $\infty$-categories:[100, 103, 24]. For Brauer-Grothendieck cohomology group:[107, 31, 147, 6, 7, 5, 3, 4]. And more generally else, we recommand the references: [59, 111, 47, 121, 120, 65, 64, 85, 94, 95, 154, 25, 27, 84, 22, 134, 135]. And we will find B. Lawvere's paper there: [72, 73, 74, 75, 76, 77, 78, 79, 80, 81, 82]

### 1.2.3 Acknowledgements.

My first thanks naturally go to my thesis advisor, Prof. Hugh R. Thomas, to my University, l'UQAM, to Quebec, and to Canada, all of whom have welcomed me warmly. I am deeply grateful to all four of them for their trust.

I would also like to acknowledge Professors Jacob Lurie and Urs Schreiber for the influence of their work and writings. Each is among the leading figures in their interrelated disciplines, and they taught me everything (without realizing it) about the mathematics contained in this thesis.

I also thank my first two mentors, the Swiss Professor Bernhard Keller of the University of Paris Cité, as well as his close friend, Rached Mneimne, a distinguished professor at the Ecoles Normales Superieures and the founder and driving force behind his publishing house, Editions Calvage et Mounet, Paris. Bernhard showed me the closest, most friendly attention, as Rached always did — the most benevolent and most fruitful for my work, and always available to offer his best advice as a friend and professor to a student. You both accompanied and guided my mathematical life with friendship, warmth, and delight. I really admire you both. Thank you, Dear Rached, and oh, how very much thank you, Dear Bernhard.

Professor François Bergeron was also an exceptional mentor who inspired admiration and humility and an equally guiding figure at UQAM, without whom none of this would ever have come to fruition. I owe so much to François's humanity and professionalism that these few lines are very poor compared to his legendary wisdom and the five years of moral, mathematical, and human support that François lavished upon me with the utmost elegance of one who stands in the posture of the one who knows. A few lines are worth a thousand and more. You know how much I care for you and how infinitely grateful I am to you. Thank you, Dear François.

## 1.3. Typographic glossary.

Throughout this text, we will see that the same underlying objects give rise to increasingly complex objects, for example by complexifying the structure sheaf associated with a topological space. We will see that simple étale schemes become Deligne-Mumford spectral stacks and, to put it more simply, that algebraic spaces endowed with the étale topology become spectral algebraic spaces, simply by attaching a module geometry and a second structure sheaf to them.

All this complexification is visible in the text, even though we have not changed the typography: this was the choice made by Jacob Lurie. Among its many advantages, this allowed us to compile in LaTeX2e on all machines, and we were not hindered by regular changes in typography. On the other hand, the risk was losing track of whether a scheme $\mathscr{X}$ was a simple scheme, a spectral scheme, an algebraic space, a Deligne-Mumford spectral stack, etc., and this, over a span of just twenty pages to more than five hundred pages. The information was contained within the text itself, which had to be reread to ensure, as the discussion progressed, that all the characteristics of the initial scheme had been lost or preserved.

Another choice was to accept working with a typography that became more complex as the meaning of each object became more complex. This was more cumbersome to handle and harder to follow in the text. But the benefit was significant: if an object had a certain typography, then the nature of the object — the category to which it belonged — was unambiguous.

The second choice also followed the evolution of the LaTeX compiler to XeLaTeX, which allows for the integration of every conceivable typeface. I then tried to find a balance between these two choices, which is why I am providing a glossary here for reference in contexts where things might not be entirely clear.

The first typographie is naturally Times New Roman for the text. For the mathematical objects, `\mathcal{}`, `\mathscr{}`, `\mathfrak{}`, `\mathbb{}` are the most common used. We use also `\textfrak{}`,`\textswab{}`, and `\textgoth{}`.

The typeface used for the "O" of Grothendieck is from the Calligra font: $\mathscr{O}$.
And all the derived typefaces are also from the Calligra font: $\mathscr{O}_{\mathfrak{X}}$, $\mathscr{O}_{\mathfrak{X}}$, $\mathscr{O}_{\mathfrak{X}}$, $\mathscr{O}_{\mathfrak{X}^{\star}}$, $\mathscr{O}_{\mathfrak{Y}}$ , $\mathscr{O}_{\mathfrak{Y}}$, $\mathscr{O}_{\mathfrak{Y}}$, etc.

The Calligra font also addresses formal modular problems: $\mathscr{F}$.
And all the derived typefaces are also from the Calligra font: $\mathscr{F}_{\mathfrak{X}}$, $\mathscr{F}_{\mathfrak{Y}}$, etc.

A Deligne-Mumford stack is generally denoted: $\mathfrak{X}$, $\mathfrak{X}_i$, $\mathfrak{X}$, $\mathfrak{Y}$ , $\mathfrak{Y}$, $\mathfrak{Y}$, etc.

The Grassmannian spectral algebraic space endowed with the positroid-etale $\infty$-topology with its module geometry of perfect perverse sheaves is:

$$\mathfrak{G} := \left(\mathfrak{Gr}^{\geqslant 0}_{\text{ét}}, \mathscr{O}_{\mathfrak{Gr}^{\geqslant 0}_{\text{ét}}}, \mathcal{IC}_{\mathfrak{Gr}^{\geqslant 0}_{\text{ét}}}\right)$$

And the spectral Amplituhedron spectral algebraic space endowed with the $\mathcal{BCFW}_{n,4}$-topology, with its module goemetry of perfect perverse sheaves is :

$$\mathfrak{A} := \left(\mathcal{A}^{\text{ét}}_{n,k,4}, \mathscr{O}_{\mathcal{A}^{\text{ét}}_{n,k,4}}, \mathcal{IC}_{\mathcal{A}^{\text{ét}}_{n,k,4}}\right)$$

But we will see that the intersection complex sheaves verifies the transfer of geometries and that the Amplituhedron functor:

$$\mathfrak{Z} : \mathfrak{G} := \left(\mathfrak{Gr}_{\text{ét}}^{\geqslant 0}, \mathscr{O}_{\mathfrak{Gr}_{\text{ét}}^{\geqslant 0}}, \mathcal{IC}_{\mathfrak{Gr}_{\text{ét}}^{\geqslant 0}}\right) \longrightarrow \left(\mathcal{A}_{n,k,4}^{\text{ét}}, \mathscr{O}_{\mathcal{A}_{n,k,4}^{\text{ét}}}, \mathcal{IC}_{\mathcal{A}_{n,k,4}^{\text{ét}}}\right) =: \mathfrak{A} \tag{1.3.1}$$

sends the $\mathcal{IC}$-geometry $\mathcal{IC}_{\mathfrak{Gr}_{\text{ét}}^{\geqslant 0}}$ of $\mathfrak{Gr}_{\text{ét}}^{\geqslant 0}$ to the $\mathcal{IC}$-geometry $\mathcal{IC}_{\mathcal{A}_{n,k,4}^{\text{ét}}}$ of $\mathcal{A}_{n,k,4}^{\text{ét}}$.
Lurie's extended Brauer-Grothendieck group of the dual Amplituhedron is isomorphic to Lurie's extended Brauer-Grothendieck group of the dual TNN Grassmannian since the underlying map $\widetilde{Z}$ is semi-small and proper:

$$H^2_{\text{ét}}({\mathcal{A}_{n,k,4}^{\text{ét}}}^{\star}, \mathcal{IC}^{\times}_{{\mathcal{A}_{n,k,4}^{\text{ét}}}^{\star}}) \simeq H^2_{\text{ét}}({\mathfrak{Gr}_{\text{ét}}^{\geqslant 0}}^{\star}, \mathcal{IC}^{\times}_{{\mathfrak{Gr}_{\text{ét}}^{\geqslant 0}}^{\star}}) \simeq c_1({\mathfrak{Gr}_{\text{ét}}^{\geqslant 0}}^{\star}) \simeq c_1(\mathfrak{Gr}_{\text{ét}}^{\geqslant 0}) \simeq \mathbb{Z}$$

since $\mathfrak{Gr}_{\text{ét}}^{\geqslant 0}$ is autodual. The first Chern class is spanned by one generator, that is the product of the first Chern classes for each BCFW cell in the étale (dual) TNN Grassmannian, since if one of the tiles flips over, then all the tiles flip over due to their overlap along the gluing of the closed immersions with common closed intersections. The dual twist is here given by the swap between the generators **1** and **-1**.

REMARK 1.3.2 : The previous result comes exactly from:

- The duality functor commutes with the Amplituhedron functor
- We constructed exactly one dual.

We have the following assertion:

- The Lurie's $\mathcal{IC}$-geometry of the dual TNN Grassmannian spectral $\mathcal{IC}$-algebraic space is sent on the dual Lurie's $\mathcal{IC}$-geometry of the dual Amplituhedron spectral $\mathcal{IC}$-algebraic space by a transformation of geometries between perfect perverse admissible geometries on the corresponding spectral algebraic spaces.

REMARK 1.3.3 : The previous result means exactly:

- The duality functor respects the Schubert geometry of the TNN Grassmannian since Lurie's geometry of the corresponding dual Amplituhedron spectral algebraic space is the image of Lurie's geometry of the corresponding coaffine dual TNN Grassmannian spectral algebraic space by a transformation of Lurie's geometries at the level of perfect perverse $\mathcal{IC}$-structures.

We will end the text with the volume of the differential cohesive topos:

$$\mathrm{Volume}(\mathbf{H}_{\mathfrak{th}}) = \int \mathrm{d}\Omega^n(\mathbf{H}_{\mathfrak{th}}) = \mathrm{Volume}(\mathbf{H})$$

since the cohomology is insensible to the reduction, which we will apply to the dual Amplutuhedron $\mathscr{A}^{\star}_{n,k,4}$:

$$\mathrm{Volume}(\mathscr{A}^{\star}_{n,k,4}) = \int \mathrm{d}\Omega^n(\mathscr{A}^{\star}_{n,k,4}) = \mathrm{Volume}(\mathscr{A}^{\star}_{n,k,4})$$

On the same path as described in J. Lurie's "*Structured spaces*" [91] and "*Spectral schemes*" [92], we will have associated a geometry to a spectrally Deligne-Mumford stack which suits well with Schubert Geometry. As well as in the case of Zariski topology, where J. Lurie associated $\mathcal{G}_{\text{ZAR}}$, and in the case of étale topology, where he associated $\mathcal{G}_{\text{ÉT}}$, here, we associate a module geometry $\mathcal{G}_{\text{PERV}}$ incarned in the complex intersection sheaf $\mathscr{F}_{\mathfrak{X}} := \mathcal{IC}^{\bullet}$.

For a goal of complete glossary to which we can refer during the lecture, let us recall the statement of the decomposition theorem (see subsection 2.6.6 and subsection 3.2.1) and its consiquencies, in particular that there are enough perfect perverse sheaves to cover the $\infty$-perverse geometry:

Let $f : X \to Y$ be a projective map of quasi-projective varieties. theorem[[66], Thm. 8.4.3] The right derived functor $\mathbb{R} f_*$ of the induced map $f_*$ sends the category of perverse sheaves on $X$ on the category of perverse sheaves on $Y$. Moreover, we have a decomposition into simple objects in the category of perverse sheaves on $Y$ of the pushforward $\mathbb{R} f_* \mathcal{F}$ for any constant sheaf $\underline{\mathbb{C}}_X$ such that:

$$\mathbb{R}\, f_* \underline{\mathbb{C}}_X \simeq \bigoplus_\alpha \mathcal{IC}^\bullet_{\mathcal{L}_{S_\alpha}}[\ell_\alpha]$$

where the $\mathcal{L}_{S_\alpha}$ are irreductible local systems on the closure $\overline{S_\alpha}$ of connected strata $S_\alpha$ of $Y$, and the $[\ell_\alpha] \in \mathbb{Z}$ shift the complexes $\mathcal{IC}^\bullet_{\mathcal{L}_{S_\alpha}} := \mathcal{IC}^\bullet(\overline{S_\alpha}, \mathcal{L}_{S_\alpha})$. theorem

Here, each $S_\alpha$ is a $4k$-dimensional tile of the Amplituhedron tiling, contractible, that is : $\mathring{\Pi}_f$ so is $\overline{\mathring{\Pi}_f} =: \Pi_f$, and we have then each local system $\mathcal{L}_{\Pi_f} \simeq \mathbb{C}[m]$ for a certain shift $m \in \mathbb{N}$ .

So the sum can be written as semi-simple connective DG algebra:

$$\mathbb{R}\, \mathfrak{Z}_* \underline{\mathbb{C}}_{\Pi_f} \simeq \bigoplus_{n\in\mathbb{N}} \mathbb{C}[n]$$

We remark that $\mathbb{R}\, \mathfrak{Z}_* \underline{\mathbb{C}}_{\Pi_f}$ is in degree $\geqslant 0$ since the DG Algebra is connective ($\mathfrak{G}$ and $\mathfrak{A}$ are coaffine so coconnective stacks) and the sum is finite since $\mathfrak{G}$ and $\Pi_f$ for each $f$ have finite étale cohomological dimension.

Following ([103] Chap 11, et Toën, Antieau, et Chang-Yeong Chough), we look for the group $H^2(\mathcal{A}^{\text{ét}}_{n,k,4}, \mathcal{IC}^\times_{\mathcal{A}^{\text{ét}}_{n,k,4}})$.

The inversibles of $\mathcal{IC}^\bullet_{\mathcal{A}^{\text{ét}}_{n,k,4}}$ as $\mathbb{E}_\infty$-ring spectra correspond to perfect perverse non-zero cochains.

The extended Brauer group of the Amplituhedron will then be non-trivial, as soon as there is at least one cell in the Amplituhedron tiling.

Moreover, for an open cell $\mathring{\Pi}_f$, which is an affine variety given by the cluster structure of (the affinization of) its ring of coordinates, and can thereby be written $\mathring{\Pi}_f = \mathbb{C}[x_1, x_2, \ldots, x_n]/\mathcal{I} =: \mathsf{Spec}\ \ \pi_0 R$, we have:

corollary 1.3.4 ([6] Cor 7.13) : *If $R$ is a connective $\mathbb{E}_\infty$-ring spectrum, then homotopy*

*groups of* $\mathbf{Br}(R)$ *are described by:*

$$\pi_k \mathbf{Br}(R) \simeq \begin{cases} H^1_{\text{ét}}(\mathsf{Spec}\ \ \pi_0 R, \mathbb{Z}) \times H^2_{\text{ét}}(\mathsf{Spec}\ \ \pi_0 R, \mathbb{G}_{\text{m}}) & k = 0 \\ H^0_{\text{ét}}(\mathsf{Spec}\ \ \pi_0 R, \mathbb{Z}) \times H^1_{\text{ét}}(\mathsf{Spec}\ \ \pi_0 R, \mathbb{G}_{\text{m}}) & k = 1 \\ \pi_0 R^\times & k = 2 \\ \pi_{k-2} R & k \geqslant 3 \end{cases} \tag{1.3.5}$$

In our situation:

THEOREM 1.3.6 ([103] Construction 11.5.5.3) : *If* $(\mathring{\mathfrak{P}}_f, \mathscr{O}_{\mathring{\mathfrak{P}}_f}, \mathcal{IC}^\bullet_{\mathring{\mathfrak{P}}_f})$ *is the spectral Deligne-Mumford stack induced from* $\mathring{\Pi}_f$*, then the homotopy groups of the* Brauer sheafs $\underline{\mathscr{B}\mathrm{r}}^\dagger_{\mathring{\mathfrak{P}}_f}$ *are described by:*

$$\pi_k \underline{\mathscr{B}\mathrm{r}}^\dagger_{\mathring{\mathfrak{P}}_f} \simeq \begin{cases} 0 & k = 0 \\ \underline{\mathbb{Z}} & k = 1 \\ (\pi_0 \mathcal{IC}^\bullet_{\mathring{\mathfrak{P}}_f})^\times & k = 2 \\ \pi_{k-2} \mathcal{IC}^\bullet_{\mathring{\mathfrak{P}}_f} & k \geqslant 3 \end{cases} \tag{1.3.7}$$

*where* $\underline{\mathbb{Z}}$ *is the constant sheaf associated to* $\mathbb{Z}$.

Thanks to [103] Example 11.5.5.5, for the spectral Deligne-Mumford stack $\mathring{\mathfrak{P}}_f$, we have $\underline{\mathscr{B}\mathrm{r}}^\dagger_{\mathring{\mathfrak{P}}_f} \simeq H^2(\mathring{\mathfrak{P}}_f, \mathcal{IC}^{\bullet,\times}_{\mathring{\mathfrak{P}}_f}) \times (H^1 \mathring{\mathfrak{P}}_f, \underline{\mathbb{Z}})$. And since $\mathring{\mathfrak{P}}_f$ is a quasi-compact, quasi-separated, normal, locally Noetherian spectral algebraic space, then we obtain $\underline{\mathscr{B}\mathrm{r}}^\dagger_{\mathring{\mathfrak{P}}_f} \simeq H^2(\mathring{\mathfrak{P}}_f, \mathcal{IC}^{\bullet,\times}_{\mathring{\mathfrak{P}}_f})$.

That is $\underline{\mathscr{B}\mathrm{r}}^\dagger_{\mathring{\mathfrak{P}}_f} \simeq H^2(\mathring{\mathfrak{P}}_f, \mathcal{IC}^{\bullet,\times}_{\mathring{\mathfrak{P}}_f}) \simeq (\pi_0 \mathcal{IC}^\bullet_{\mathring{\mathfrak{P}}_f})^\times$

But $H^2(\mathring{\mathfrak{P}}_f, \mathcal{IC}^\times) = \mathring{\mathfrak{P}}_f = c_1(\mathring{\mathfrak{P}}_f)$ for $k = 2$, the first Chern class of the corresponding positroid variety, since the étale map gives an isomorphism of the Brauer group.

And this first Chern class of an (open) positroid variety is $\mathbb{Z}$, with one generator **1** since it is an affine variety given by the cluster structure of its (affinization of the) ring of coordinates.

We will remark that the twisted sheaf has **-1** as generator of its first Chern class.

*And, again, if one of the tiles flips over, then all the tiles flip over due to their overlap*

*along the gluing of the closed immersions with common closed intersections of (closed) positroid varieties.*

We can then conclude with what is the heart of the machinery:

**Definition.** *We define the* Amplituhedron geometric morphism $\mathscr{Z}$ *between the* $\infty$*-topos of the TNN Grassmannian* $\mathfrak{Gr}^{\geqslant 0}_{\text{ét}}$ *endowed with the positroid-etale* $\infty$*-topology and the* $\infty$*-topos of the Amplituhedron* $\mathcal{A}^{\text{ét}}_{n,k,4}$ *endowed with the* $\mathcal{BCFW}_{n,4}$*-topology, such that:*

$$\mathscr{Z}_* \ : \ \mathfrak{Gr}^{\geqslant 0}_{\text{ét}} \underset{\longrightarrow}{\overset{\text{left ex. \& left adj.}}{\underset{\perp}{\longleftarrow}}} \mathcal{A}^{\text{ét}}_{n,k,4} \ : \ \mathscr{Z}^* \tag{1.3.8}$$

## 1.4. Main steps and statements of results.

The thesis is subdivided in two parts. The first one deals with what is called "*a new brave positive Schubert geometry*" which consists of the rewriting of the real positive Grassmannian and the singular gluing of its positroid varieties in a new way which allows a more powerful interpretation thanks to spectral algebraic geometry on $\mathbb{E}_\infty$-rings spectra and structured and fractured spaces in the sense of J. Lurie's geometries, finding its compact but meaningful enunciation via intersection complexes in perverse sheaves theory with enough perfect perverse sheaves to cover the $\infty$-perverse geometry.

This algebro-geometric point of view is jumelled with a synthetic differential geometrical point of view with the characteristic of Lawvere's and Schreiber's differential cohesiveness which treats the question of formal moduli functors, infinitesimal thickenings, De Rham stacks and their crystalline cohomology which subsumes all the modalities of the structured geometries and the unification of their singular cohomologies of the underlying concrete topological étale algebraic space and its deformation theories of this initial grassmannian and its positroid varieties.

The second part uses this powerful rewriting to develop the Tannaka duality for quasi-coherent stacks already present in the first part, to show that, instead of the Grassmannian spectral Deligne-Mumford stack which is Verdier-autodual in the $\infty$-category of $t$-structured prestable $\infty$-categories of modules on $\mathbb{E}_\infty$-ring spectra, the Dual Amplituhedron, an object of the physicist N. Arkani-Hamed in high energy particle colliders in N=4 SYM context, has a concrete dual, non trivial, that is of interest for us for duality between standard model and string theory.

Let us examine further what is announced in this little resume.

### 1.4.1 Which brave new geometry of the positive Grassmannian.

A NEW POSITIVE SCHUBERT GEOMETRY

$$\mathfrak{I}(\mathfrak{X})_{\mathrm{dRh}} := \mathfrak{I}\Big(\mathfrak{Gr}^{\geqslant 0}_{\text{ét}}, \mathscr{O}_{\mathfrak{Gr}^{\geqslant 0}_{\text{ét}}}, \mathcal{IC}^{\bullet}_{\mathfrak{Gr}^{\geqslant 0}_{\text{ét}}}\Big)_{\mathrm{dRh}}$$

as a $\mathcal{G}^{\mathrm{perv}}_{\mathrm{perf}}$-structured $(\mathfrak{X}, \mathscr{O}_{\mathfrak{G}}, \mathscr{F}_{\mathfrak{G}})-$module geometry Lurie's formal moduli problem $\mathscr{F}$ for spectral Deligne-Mumford stably $t$-structured $(\mathsf{QStk}(\mathfrak{X})^{\geqslant 0}, \mathsf{QStk}(\mathfrak{X})^{\leqslant 0})$quasi-coherent stacks $\mathsf{QStk}(\mathfrak{X})$ in $\mathsf{Groth}_{\infty} \subseteq \mathsf{QCoh}^{\mathrm{Pst}}_{\mathrm{comp}}(\mathfrak{X})$ on infinitesimally thickened cohesive $\mathfrak{X}_{/*}$ Grothendieck contexts $\mathfrak{R}(\mathfrak{X}) \dashv \mathfrak{I}(\mathfrak{X}) \dashv \&(\mathfrak{X})$ of Tannaka 1-affine perversely twisted $H^2_{\text{ét}}(\mathfrak{X}, \mathcal{IC}^{\bullet\times})$ $\mathcal{BCFW}_{n,4}$ positroid $\mathring{\Pi}_f$ modal de Rham stacks.

We have to define all the following expressions:

1. Lurie's formal moduli problem $\mathscr{F}$
2. $\mathcal{G}^{\mathrm{perv}}_{\mathrm{perf}}$-structure and Lurie's module geometry $(\mathfrak{X}, \mathscr{O}_{\mathfrak{G}}, \mathscr{F}_{\mathfrak{G}})$
3. symmetric monoidal $\infty$-category of spectral deligne-mumford quasi-coherent stacks $\mathsf{QStk}^{\mathrm{Pst}}_{\mathrm{comp}}(\mathfrak{X})$ of complete Grothendieck prestable $\infty$-categories on $\mathfrak{X}$
4. stable triangulated glued $t$-structures $(\mathsf{QStk}(\mathfrak{X})^{\geqslant 0}, \mathsf{QStk}(\mathfrak{X})^{\leqslant 0})$
5. infinitesimal thickening of the fat point $\mathfrak{X}_{/*}$
6. Lawvere's and Schreiber's cohesiveness for $\infty$-topoï
7. Grothendieck adjoint derived context $\mathfrak{R}(\mathfrak{X}) \dashv \mathfrak{I}(\mathfrak{X}) \dashv \&(\mathfrak{X})$
8. Gaitsgory's Tannaka 1-affineness
9. Verdier's autoduality of the perverse category $\mathcal{P}\mathrm{erv}(\mathfrak{X}) \subset \mathsf{QStk}(\mathfrak{X})^{\heartsuit}$

10. twisted étale Brauer-Grothendieck cohomology $H^2_{\text{ét}}(\mathfrak{F}, \mathcal{IC}^{\bullet\times})$
11. fully dualizable $\mathcal{BCFW}_{n,4}$ positroid varieties $\mathring{\Pi}_f$
12. logical homotopy type theoretical modality
13. infinitesimal shape modal de Rham stack $\mathfrak{J}(\mathfrak{F})_{\text{dRh}}$

The explanations come progressively until the whole description of this new object.

### 1.4.2 First part : the new brave positive Schubert geometry.

At this stage, we choose to decline in steps the decomposition of the path until:

$$\mathfrak{I}(\maltese) \;:=\; \mathfrak{I}\big(\mathfrak{Gr}^{\geqslant 0}_{\text{ét}}, \mathscr{O}_{\mathfrak{Gr}^{\geqslant 0}_{\text{ét}}}, \mathcal{IC}^{\bullet}_{\mathfrak{Gr}^{\geqslant 0}_{\text{ét}}}\big)_{\mathbf{dRh}}$$

Main steps (in desorder with respect to the text):

- Rewriting of the Grassmannian in a formal moduli problem
- Existence of the formal deformation
- Existence of the infinitesimal cohesion in a triple modalities of adjunction
- Existence of the De Rham thickening and the crystalline cohomology
- Existence of monadic crystallin sheaves
- Monadicity, modalities, and adjonctions between an $\infty$-topos and its thickening
- Grassmannian and algebraic space in the étale topology
- Grothendieck context and gluing of $t$-structures
- Swap to $\infty$-topos and spectral Deligne-Mumford stacks
- Gluing closed immersions and scallop decomposition
- Infinitesimally thickened spectral algebraic spaces
- Corporeal structured spaces corporeal et Lurie's geometries
- t-Structures, 0-affineness of $(\mathsf{QCoh}^{\geqslant 0}, \mathsf{QCoh}^{\leqslant 0})$
- t-Structures, 1-affineness of $(\mathsf{QStk}^{\geqslant 0}, \mathsf{QStk}^{\leqslant 0})$
- Tannaka duality to recover the initital $\infty$-topos (double decatorification)

- Perverse perfect sheaves and Intersection complexes with singularities
- De Rham Cohomologies, soothness and singularities
- Unified De Rham Cohomologies and six-functor formalism
- $\infty$-Categories of prestable $\infty$-categories
- De Rham volume form and positive canonical form
- Grassmanian and positroid varieties: A new brave positive Schubert geometry

### 1.4.3 Second part : The dual Amplituhedron.

#### 1.4.3.1 Abstract of the second part.

The dual Amplituhedron $\mathfrak{J}(\mathfrak{A}^{\star})$ will be a singular topological space coming from a spectral Deligne-Mumford stack, differentially (infinitesimally) cohesive as a De Rham stack whose infinitesimal thickening will allow to compute its volume.

A longer abstract of this second part would be the following:

Thanks to the theory of spectral algebraic geometry, associated with the theory of perverse sheaves, this thesis consists of an extension of the combinatorial study of Schubert geometry beyond the classical approaches, whether purely combinatorial or mixed with algebraic geometry through cluster algebras. We use a higher stack categorification in order to structure the complexity of the bad-behaved ramified recollements of stratifications by Schubert or positroid cells in the TNN Grassmannian $\mathrm{Gr}(k,n)^{\geqslant 0}$. Thanks to suitable positroid-étale coverings for the underlying topological spaces of our spectral Deligne-Mumford stacks, and involving Lurie's Tannakian reconstruction for geometric stacks and $\mathcal{G}$-structured spaces, the wonderfully fertile point of view of the magic of functors of A. Grothendieck is adopted here. The developments in derived and spectral algebraic geometry give a new path in the description, the encoding, and the resolution of the singularities of Schubert geometry. A condition for the existence of a non-trivial dual $\mathfrak{X}^{\star}$ of a spectral Deligne-Mumford stack $\mathfrak{X}$ is obtained in the case of a coaffine stack, through extended cohomological Brauer-Grothendieck groups and spaces, namely the nonvanishing of the first Chern class of the underlying space. Indeed, the prestable $\infty$-category $\mathsf{QStk}^{\mathrm{PSt}}(\mathfrak{S})$ of quasi-coherent stacks on such a coaffine spectral Deligne-Mumford stack $\mathfrak{S}$ is symmetric monoidal and the stable $\infty$-category $\mathsf{QCoh}(\mathfrak{S})$ of quasi-coherent sheaves on the stack $\mathfrak{S}$ can be dualized in this prestable $\infty$-category $\mathsf{QStk}^{\mathrm{PSt}}(\mathfrak{S})$. Lurie's Tannaka duality then allows us to reconstruct the stack $\mathfrak{S}$, which is a 1-localic $\infty$-topos in case of the initial coaffine spectral stack $\mathfrak{S}$ is constructed from a concrete topological space $X$ underlying a spectral algebraic space $\mathscr{X}$. The definition of a dual space $X^{\star}$ is then obtained for an underlying space $X$ of a spectral algebraic space thanks to a higher categorical monoidal dualization. In the physical case of the Amplituhedron, we then use the theory of perverse sheaves, namely the sheaves $\mathcal{IC}$ of intersection complexes satisfying Deligne's hypothesis, which give us a $\mathcal{G}$-structure in the sense of Lurie for the spectral Deligne-Mumford stacks defined on the TNN Grassmannian and the Amplituhedron $\mathscr{A}_{n,k,4}$ of N. Arkani-Hamed in the planar model $\mathscr{N} = 4$ Super Yang-Mills. By viewing the induced morphism at the level of Grassmannians and by categorifying this induced map $\widetilde{\mathscr{Z}}$ defining the Amplituhedron, which satisfies now the hypothesis of the decomposition theorem for perverse sheaves as a functor between two such coaffine spectral stacks that respects the $\infty$-topoi $\mathcal{G}$-structured by their categories of perverse sheaves, we obtain a proper semi-small geometric functor $\mathfrak{Z}$ between the two spectral Deligne-Mumford stacks $\mathfrak{G}$ and $\mathfrak{A}$ constructed on the $\infty$-topos and endowed with $\mathcal{IC}_{\mathfrak{Gr}^{\geqslant 0}_{\text{ét}}}$ and $\mathcal{IC}_{\mathcal{A}^{\text{ét}}_{n,k,4}}$ as compatible $\mathcal{G}$-geometries. We obtain also a unique non-trivial dual $\mathscr{A}^{\star}_{n,k,4}$ after considerations about module categories over $\mathbb{E}_\infty$-ring spectra and the invertibility of these categories in bigger ones. Azumaya algebras and extended Brauer-Grothendieck spaces are then non-trivial in our case, and are isomorphic since we have a homeomorphism of the $\mathsf{BCFW}$-étale cover of the $\infty$-topos of the TNN Grassmannian $\mathfrak{Gr}^{\geqslant 0}_{\text{ét}}$ and the $\mathcal{BCFW}$-étale cover of the $\infty$-topos of the Amplituhedron $\mathcal{A}^{\text{ét}}_{n,k,4}$. Since the first Chern class of the $\mathsf{BCFW}$-étale TNN Grassmannian $\mathfrak{Gr}^{\geqslant 0}_{\text{ét}}$ is $\mathbb{Z}$, that is the second cohomological étale group $H^2_{\text{ét}}(\mathfrak{Gr}^{\geqslant 0}_{\text{ét}}, \mathcal{IC}^{\times}_{\mathfrak{Gr}^{\geqslant 0}_{\text{ét}}})$ is $\mathbb{Z}$ and is isomorphic to $H^2_{\text{ét}}(\mathcal{A}^{\text{ét}}_{n,k,4}, \mathcal{IC}^{\times}_{\mathcal{A}^{\text{ét}}_{n,k,4}})$, the existence of a (unique) dual $\mathscr{A}^{\star}_{n,k,4}$ of the Amplituhedron is guaranteed. As an essential constraint of the initial problem, the dual positroid stratification of the TNN Grassmannian $\mathrm{Gr}(k,n)^{\geqslant 0}$ is respected by the construction of this dual Amplituhedron $\mathscr{A}^{\star}_{n,k,4}$. The functors of étale spectrum $\mathsf{Spét}(\text{-})$ and quasi-coherent sheaves $\mathsf{QCoh}(\text{-})$ or stacks $\mathsf{QStk}(\text{-})$ on perversely structured spectral Deligne-Mumford stacks, governed by Isbell duality between categories of $\infty$-sheaves on $\infty$-topoi, lead to solve the singularities of the underlying spaces of spectral algebraic spaces, by the increasing the $n$-affineness of the initially bad-behaved ramified recollements of Schubert cells and positroid varieties during the successive categorifications by applying successively the functor $\mathsf{QCoh}(\text{-})$ then the functor $\mathsf{QStk}(\text{-})$ and decategorification by applying twice Tannaka duality. This is therefore a purely $\infty$-topos-theoretic new approach to the study of Schubert's geometry, preventing the conservation of all the necessary structures involved, in particular to construct a unique non-trivial dual Amplituhedron $\mathscr{A}^{\star}_{n,k,4}$.

#### 1.4.3.2 Glossary and results between "old" and "new" algebraic geometry.

In other words, what is an algebraic variety in the scheme-theoretic vocabulary?

Given that Schubert's enumerative geometry, in which the geometry of positroid varieties (1985) is embedded, dates back to Hermann Schubert in 1879 and is over a century old, that it uses the "*old algebraico-geometric vocabulary*", and that Alexander Grothendieck's "*new ideas in algebraic geometry*" — such as schemes and links to homological algebra, but also topoï, étale cohomology, inifinitesimal site, etc. — emerged around 1958 after his Tohōku's paper (1957, [51]), here is an indispensable glossary as well as a set of results enabling translation between what have become two distinct worlds.

To be more precise, here are the two definitions of an algebraic variety in the "*old*" and the "*new* world:

DEFINITION 1.4.1 : *We define an* algebraic variety $\mathscr{X}_{/\Bbbk}$ *as the commun locus of the solutions of a system of polynomial equations over a field* $\Bbbk$*. One says that it is* irreducible *if it is not the union of two smaller proper subvarieties. If it is the case, we call this locus an* algebraic set*. The* Hilbert's Nullstellensatz *then allows us, in the case of* $\Bbbk$ *is the complex numbers, to make a correspondance between ideals of polynomial rings and algebraic sets.*

Although a basic vocabulary is required to fully understand the following definition, let us nevertheless formulate it rigorously:

DEFINITION 1.4.2 : *In the modern scheme vocabulary, an* algebraic variety $\mathscr{X}_{/\Bbbk}$ *over a field* $\Bbbk$ *is an integral, that is irreducible and reduced, scheme over* $\Bbbk$ *whose structure* *morphism* $\mathscr{O}_{\mathscr{X}_{/\Bbbk}}$ *is separated and of finite type.*

We immediately realize, then, that we have some way to go to move from one language to another. And conversely, to be able to recover a "*concrete space*" of the old world, making the reverse journey will also require some effort to find our way back to that "*common locus*" of the solutions of a system of polynomial equations.

In particular, since we are in a summary section, let us enounce the facts that will not be redefined or proved again after the work of A. Postnikov, T. Lam, D. Speyer, and others.

THEOREM 1.4.3 : *We have the following facts:*

- *The Grassmannian* $\mathrm{Gr}^{\mathbb{R}}(k,n)$ *and* $\mathrm{Gr}^{\mathbb{C}}(k,n)$ *are smooth normal Cohen-Macaulay projective varieties*

- *The Grassmannian* $\mathrm{Gr}^{\mathbb{R}}(k,n)$ *is a subvariety of the variety* $\mathrm{Gr}^{\mathbb{C}}(k,n)$

- *The open positroid variety* $\overset{\circ}{\Pi}{}_f^{\mathbb{R}}$ *is a quasi-projective subvariety of the Grassmannian projective variety* $\mathrm{Gr}^{\mathbb{R}}(k,n)$*, and so is a quasi-projective subvariety of the projective variety* $\mathrm{Gr}^{\mathbb{C}}(k,n)$

- *The open positroid variety* $\overset{\circ}{\Pi}{}_f^{\mathbb{C}}$ *is a quasi-projective subvariety of the closed positroid projective variety* $\Pi_f^{\mathbb{C}}$

- *The closed positroid variety* $\Pi_f^{\mathbb{R}} := \overline{\overset{\circ}{\Pi}{}_f^{\mathbb{R}}}$ *, the Zariski closure in* $\mathbb{C}^N$ *(for a certain* $N$ *which depends on the degree of the Plücker coordinates) of the real open positroid variety* $\overset{\circ}{\Pi}{}_f^{\mathbb{R}}$*, is a projective variety*

- *The Grassmannian* $\mathrm{Gr}^{\mathbb{R}}(k,n)$ *and* $\mathrm{Gr}^{\mathbb{C}}(k,n)$ *are reduced noetherian projective schemes*

- *The Grassmannian* $\mathrm{Gr}^{\mathbb{R}}(k,n)$ *is a closed subscheme of the scheme* $\mathrm{Gr}^{\mathbb{C}}(k,n)$

- *The TNN positroid varieties* $\mathrm{Gr}^{\mathbb{R}}_{\geqslant 0}(k,n)$ *and* $\mathrm{Gr}^{\mathbb{C}}_{\geqslant 0}(k,n)$ *are two identical projective varieties in the "*old*" world, and their functor of points in the "*new*" world determine two isomorphic functors. In other words, the complex points of* $\mathrm{Gr}^{\mathbb{R}}_{\geqslant 0}(k,n)$ *are the complex points of* $\mathrm{Gr}^{\mathbb{C}}_{\geqslant 0}(k,n)$*, and they are real points.*

- *The open positroid variety* $\overset{\circ}{\Pi}{}^{\mathbb{R}}_f$ *is a quasi-projective subscheme of the Grassmannian projective scheme* $\mathrm{Gr}^{\mathbb{R}}(k,n)$*, and so is a quasi-projective subscheme of the projective scheme* $\mathrm{Gr}^{\mathbb{C}}(k,n)$

- *The open positroid variety* $\overset{\circ}{\Pi}{}^{\mathbb{R}}_f$ *is a quasi-projective subscheme of the positroid projective scheme* $\Pi^{\mathbb{R}}_f$

- *The open positroid variety* $\overset{\circ}{\Pi}{}^{\mathbb{C}}_f$ *is a quasi-projective subscheme of the positroid projective scheme* $\Pi^{\mathbb{C}}_f$

- *The open positroid variety* $\overset{\circ}{\Pi}{}^{\mathbb{R}}_f$ *in the real closed positroid variety* $\Pi^{\mathbb{R}}_f$*, is a smooth normal Cohen-Macaulay noetherian open quasi-projective subscheme of the projective smooth normal Cohen-Macaulay noetherian scheme* $\Pi^{\mathbb{R}}_f$

#### 1.4.3.3 General prerequisits (first part): The stratified real and complex TNN Grassmannian.

Let us describe the real and the complex Grassmannians and open and closed positroid varieties in a few words on which we will return later to complexify the landscape but, by this way, earn a lot in this complexification. We use here straight now the langage of algebraic varieties coming from the XIX$^{\text{th}}$ century (so without telling it, with the Zariski topology). But we can point right now however that we will have soon to skip to the étale topology and $\infty$-topoï then to derived stacks to get a more powerful object, that is the powerful functor of points of spectral Deligne-Mumford stacks, to inherit structures we are going to attach to the initial object we present just here, the (classical and over-studied for centuries) Grassmannians $\mathrm{Gr}^{\mathbb{R}}_{\geqslant 0}(k,n)$ and $\mathrm{Gr}^{\mathbb{C}}_{\geqslant 0}(k,n)$.

We just incorporate right now the stratification of A. Postnikov in the real case. Indeed, we can show the following disjoint stratification of $\mathrm{Gr}^{\mathbb{R}}_{\geqslant 0}(k,n)$. Let us just recall that classical algebraic geometry in characteristic 0 works for $\mathbb{C}$, the complex numbers, and not for the real numbers $\mathbb{R}$. So an embedding at the right step will be necessary to allow us to inherit the (first) power of classical algebraic geometry, where we mean "*classical*" versus "*spectral*".

THEOREM 1.4.4 ([114]) : *The real TNN Grassmannian,* $\mathrm{Gr}^{\mathbb{R}}_{\geqslant 0}(k,n)$*, the set of* $k$*-planes in* $\mathbb{R}^n$*, can be decomposed in the following cells:*

$$\mathrm{Gr}^{\mathbb{R}}_{\geqslant 0}(k,n) = \coprod_{f \in \mathscr{B}(k,n)} \Pi_f^{>0}$$

*where the* open positroid cell $\Pi_f^{>0}$ *is defined such as :*

$$\Pi_f^{>0} = \{V \in \mathrm{Gr}^{\mathbb{R}}(k,n) \mid P_I(V) = 0 \text{ if } I \notin \mathcal{M} \text{ and } P_I(V) > 0 \text{ if } I \in \mathcal{M}\}$$

*where the* $P_I$ *denotes the real Plücker coordinates indexed by the set* $I$*, in the set* $\mathcal{M}$ *named a* matroid*, of columns of a point* $V$ *in the Grassmannian* $\mathrm{Gr}^{\mathbb{R}}(k,n)$*, and* $\mathscr{B}(k,n)$ *denotes a set, named* bounded affine permutations*, modulo equivalence of the so-called* plabic graphs *indexed by the bounded affine permutations* $f$*.*

We denote the projective (closed) positroid variety $\Pi_f^{\mathbb{R}} := \overline{\Pi_f^{>0}}$, the Zariski closure of the open (quasi-projective) positroid cell $\Pi_f^{>0}$.

#### 1.4.3.4 How to recover classical algebraic geometry from the semialgebraic topology.

We saw at the last subsection that the open positroid cell is a semialgebraic set. Let us now explain how to be allowed to make classical (then derived and spectral) algebraic geometry.

Let us call $M$ a real semialgebraic set. To be able to work over $\mathbb{C}$ with the scheme langage, in particular with structure sheaves, we need a good object $(M, \mathscr{O}_M)$ where $M$ is settled on a classical topological $\mathbb{C}$-space, and $\mathscr{O}_M$ is a structure sheaf of $\mathbb{C}$-rings on it.

The answer to this thorny issue is not trivial, involves recollements of topoï, but is summarized in principally [127] which collects at the beginning of his book all the needed machinery. To be more precise, we willl quote the initial few papers from which Scheiderer takes its machinery, that is the papers of [33] and [67].

Let us recall the facts demonstrated in these papers.

On a real affine variety $V$, we define a first topology given by the total order of $\mathbb{R}$. We call it the strong topology or euclidian topology.

We define a second one, the semialgebraic topology or weak topology (see [33], 1.a) given by the semialgebraic open sets and the finite coverings. We can also more generically define the prime precones $\alpha$ (see [33], 2). In the case of the ring of coordinates $\mathbb{R}[V]$ of a real affine variety, these are the same sets.

And we define a third one, the real étale topology made of étale morphisms of varieties $W \longrightarrow V$ with the coverings made of families $\{W_i \longrightarrow W\}$ of étale morphisms such that the $\{W_i \longrightarrow W\}$ are a surjective family.

Definition 1.4.5 (Real étale topos) : *We define the* real étale topos $V(\mathbb{R})$ *as the topos of sheaves for this topology.*

We define the real spectrum, denoted $\mathsf{Spec}_{\mathbb{R}}\Bbbk$ of a ring $\Bbbk$, as the topological space whose points are the prime precones $\alpha \subset \Bbbk$ (the semialgebraic open sets) and the open sets are generated by the elementary open sets $D_a = \{\alpha \,\text{prime precone} \mid a \notin \alpha\}$.

The real spectrum is a spectral space, which is functorial that is if $f : \Bbbk \longrightarrow \Bbbk'$ is a ring morphism, then $\mathsf{Spec}_{\mathbb{R}} f : \mathsf{Spec}_{\mathbb{R}}\Bbbk' \longrightarrow \mathsf{Spec}_{\mathbb{R}}\Bbbk$ is continuous.

In the case of the ring of coordinates $\mathbb{R}[V]$ of a real affine variety, the real étale topos $V(\mathbb{R})$ endowed with the strong euclidean topology is dense in $\mathsf{Spec}_{\mathbb{R}}\mathbb{R}[V]$, and we have a bijection between the semialgebraic open sets of $V(\mathbb{R})$ and the quasi-compact opens of $\mathsf{Spec}_{\mathbb{R}}\mathbb{R}[V]$.

We define the structural sheaf on $\mathsf{Spec}_{\mathbb{R}}\mathbb{R}[V]$ such that for $\mathscr{U}$ open in $\mathsf{Spec}_{\mathbb{R}}\mathbb{R}[V]$, we have a filtering system $\{(C,\phi)\}$ where $C$ is an étale $\mathbb{R}$-algebra and $\phi$ is a section of the local homeomorphism from $\mathsf{Spec}_{\mathbb{R}}C$ to $\mathsf{Spec}_{\mathbb{R}}D$. The structure sheaf $\mathscr{O}_{\mathsf{Spec}_{\mathbb{R}}\mathbb{R}[V]}$ is the sheaf associated to the presheaf which associates an open $\mathscr{U}$ to the inductive limit of the filtering system.

We can then formulate and enhance [33], Theorem 1 with [127], but also with [67] and [115]:

THEOREM 1.4.6 ([33], [127] Thm 1.3, [67] and [115]) : *The small topos endowed with the Grothendieck topology defined by the semialgebraic sets coming from* $\mathsf{Spec}_{\mathbb{R}}\mathbb{R}[V]$ *is an open topos, complementary to the closed topos made of the big real etale topos. The recollement of these open and closed topoï gives a topos endowed with a Grothendieck topology denoted* b-topology*, where we have homotopy theory (the segment* $[0,1]$ *is semi-algebraic), pullbacks and complete spaces replacing compact spaces of the strong euclidean topology.*

*If a ring of coordinates* $\mathbb{R}[V]$*, as a structure sheaf, is a cluster algebra in one of the open/closed subtopoï, the complementary object in the closed/open corresponding subtopoï has a reversed structure of cluster algebra.*

Conclusion : We earn a topos in which $(M,\mathscr{O}_M)$ is a $\Bbbk$-variety and $\mathscr{O}_M$ is an étale $\Bbbk$-algebra where we have homotopy, pullback, complete spaces, and where the base ring $\Bbbk$ does not matter anymore, since doing algebraic geometry over $\mathbb{C}$ amounts to choose the codomain of the sheaves of the topos with the b-topology, and a presheaf $\mathscr{F}$ with value in $\mathbb{R}$ is sheafified in a sheaf $\widehat{\mathscr{F}}$ with value in $\mathbb{R}$. No matters we choose $\Bbbk = \mathbb{R}$ or $\mathbb{C}$.

REMARK 1.4.7 : For a further application where $\infty$-topoï are involved, we can also remark the following:

THEOREM 1.4.8 ([36], Thm B.13) : *If* $X$ *is a topogical space of finite Krull dimension, then the real etale topos is a hypercomplete* $\infty$*-topos.*

#### 1.4.3.5 General prerequisits (second part): The stratified real and complex TNN Grassmannian.

Let us now return to our general TNN Grassmannian and positroid varieties.

THEOREM 1.4.9 : *The* open positroid variety*, denoted* $\overset{\circ}{\Pi}{}_f^{\mathbb{R}}$*, is defined as the interior of the (closed) projective positroid variety* $\Pi_f^{\mathbb{R}}$*, which has rational singularities. This open positroid variety* $\overset{\circ}{\Pi}{}_f^{\mathbb{R}}$ *is smooth, normal, irreductible, Cohen-Macaulay, open, quasi-projective in the irreductible normal projective variety* $\Pi_f^{\mathbb{R}}$*, where a suitable collection of (real positive) Plücker coordinates are not 0. The ring of coordinates* $\mathbb{C}\left[\overset{\circ}{\Pi}{}_f^{\mathbb{R}}\right]$ *of the quasi-projective open positroid variety* $\overset{\circ}{\Pi}{}_f^{\mathbb{R}}$ *is a cluster algebra.*

Let $f \in \mathscr{B}(k,n)$ be a bounded affine permutation for a plabic graph $G$ representing a positroid cell. Let us remark now that many bounded affine permutations can represent the same positroid cell. We speak then up to the equivalence relation between these bounded affine permutations representing the same positroid cell.

We have that the real open positroid varieties $\overset{\circ}{\Pi}{}_f^{\mathbb{R}}$ of the real TNN Grassmannian are disjoints. Here is the step of the embedding in the complex world: we denote the (complex) closed *positroid variety* $\Pi_f^{\mathbb{C}} := \overline{\overset{\circ}{\Pi}{}_f^{\mathbb{C}}}$, the Zariski closure in $\mathbb{C}^N$ (for a certain $N$ which depends on the degree of the Plücker coordinates) of a real open positroid varieties $\overset{\circ}{\Pi}{}_f^{\mathbb{R}}$.

From now on, we do not specify real or complex positroid variety, being allowed to jump into the complex world thanks to this closure. The varieties will then be now always complex, unless explicit specifications.

They intersect themselves in a union of (complex) positroid varieties $(\Pi_g^{\mathbb{C}} := \overline{\overset{\circ}{\Pi}{}_g^{\mathbb{C}}})_{g\in I}$, where $I$ is a finite set of representatives of bounded affine permutations.

The frontier of a closed positroid variety $\partial\Pi_f := \Pi_f \setminus \overline{\overset{\circ}{\Pi}_f}$ is a union of (closed) positroid varieties $(\Pi_g := \overline{\overset{\circ}{\Pi}_g})_{g\in J}$, where $J$ is a finite subset of $I$ of representatives of bounded affine permutations and a union of (open) positroid varieties $(\overset{\circ}{\Pi}_{g'})_{g'\in J'}$, where $J'$ is a finite subset of $I$.

The union of open positroid varieties $(\mathring{\Pi}_f)_{f\in I}$, for $I$ a finite set, is the positive Grassmannian $\mathrm{Gr}_{\geqslant 0}(k,n)$ and the interior of $\mathrm{Gr}_{\geqslant 0}(k,n)$ is the unique top positroid cell of dimension $k(n-k)$.

THEOREM 1.4.10 (see section 2.7) : *Let $f$ be any bounded affine permutation, up to isomorphism, for a plabic graph $G$ . Then :*

- *the TNN Grassmannian $\mathrm{Gr}_{\geqslant 0}(k,n)$ is a (spectral) algebraic space over $\mathbb{C}$ for the étale topology*
- *the open positroid varieties $\mathring{\Pi}_f$ are quasiprojective subschemes of this TNN Grassmannian étale topos*

We will return on this theorem later.

DEFINITION 1.4.11 (see definition 3.1.1) : *If we fix $Z$ as a $(k+m)\times n$ real matrix in $\mathbb{R}^{(k+m)\times n}$ with all maximal minors are positive, where $k+m\leqslant n$ (and especially here $m=4$ for our purpose), then the matrix $Z$ induces a map:*

$$\widetilde{Z}:\mathrm{Gr}_{k,n}^{\mathbb{R},\geqslant 0}\longrightarrow \mathrm{Gr}_{k,k+m}^{\mathbb{R}}$$

*defined by*

$$\widetilde{Z}(\langle v_1,\ldots,v_k\rangle):=\langle Z(v_1),\ldots Z(v_k)\rangle$$

*where $\langle v_1,\ldots,v_k\rangle$ represents an element of $\mathrm{Gr}_{k,n}^{\mathbb{R},\geqslant 0}$ with $k$ basis vectors. We then define the* (tree) Amplituhedron *$\mathcal{A}_{n,k,m}(Z)$ as the image $\widetilde{Z}(\mathrm{Gr}_{k,n}^{\mathbb{R},\geqslant 0})$ of $\mathrm{Gr}_{k,n}^{\mathbb{R},\geqslant 0}$ inside $\mathrm{Gr}_{k,k+m}$.*

Thanks to the BCFW recursion, we can extract some $4k-$dimensional open positroid cells indexed by $f\in\mathsf{BCFW}(k,n)$, so as to make an étale homéomorphism with a pure $4k-$dimensional étale covering of the Amplituhedron indexed by these $f\in\mathsf{BCFW}(k,n)$ :

$$\mathcal{A}_{n,k,m}(Z)=\coprod_{f\in\mathsf{BCFW}(k,n)}\widetilde{Z}(\Pi_f^{>0})$$

This étale homeomorphism is sent on what we will name a:

$$\mathcal{BCFW}\text{-étale covering of : } \mathcal{A}_{n,k,m}(Z) = \coprod_{f \in \mathcal{BCFW}(k,n)} \widetilde{Z}\,(\Pi_f^{>0})$$

of the Amplituhedron, which we will call its $\mathcal{BCFW}$-topology in the sense of Grothendieck topologies, and we have the following:

THEOREM 1.4.12 (see section 2.7) : *This* BCFW-*indexing and these* $\mathcal{BCFW}$*-indexing and* $\mathcal{BCFW}$*-covering are such that the induced Amplituhedron map between them verifies:*

- *The initial Amplituhedron map* $\widetilde{Z}$ *induces an induced map also denoted* $\widetilde{Z}$ *between the* $4k-$*dimensional open positroid cells of* $\mathsf{BCFW}(k,n)$ *of the TNN Grassmannian to the* $4k-$*dimensional cells of the Amplituhedron* $\mathcal{A}_{n,k,4}(Z)$ *constituting the* $\mathcal{BCFW}$*-topology, which is étale, proper, semi-small, and a local diffeomorphism, being indeed, in particular, a map between tiles of the same* $4k-$*dimension.*

- *This means that, in case of physics, with* $m=4$, $\widetilde{Z}$ *sends pure* $4k-$*dimensional tiles of the TNN Grassmannian* $\mathrm{Gr}_{4,n}^{\mathbb{R},\geqslant 0}$ *to* $4k-$*pure dimensional tiles of the Amplituhedron* $\mathcal{A}_{n,k,4}(Z)$*.*

- *The Amplituhedron* $\mathcal{A}_{n,k,4}(Z)$ *is an algebraic space over* $\mathbb{C}$ *for the* $\mathcal{BCFW}$*-Grothendieck étale topology.*

- *The same phenomenon occurs for any* $m$ *and not only for* $m=4$*: That is the Amplituhedron* $\mathcal{A}_{n,k,m}(Z)$ *is also an algebraic space over* $\mathbb{C}$ *and for which Grothendieck étale* $\infty$*-topology is structured with a Schubert* $\infty$*-geometry given by perfect perverse intersection complexes, giving rise to a "*new brave Schubert $\infty$-geometry for the Amplituhedron*".*

- *The induced Amplituhedron map is injective on the open positroid varieties* $\mathring{\Pi}_f$ *of the TNN Grassmannian* $\mathrm{Gr}_{k,n}^{\mathbb{R},\geqslant 0}$*.*

- *The Amplituhedron* $\mathcal{A}_{n,k,m}(Z)$*, image of the real TNN Grassmannian* $\mathrm{Gr}_{\geqslant 0}^{\mathbb{R}}(k,n)$ *by the induced Amplituhedron map* $\widetilde{Z}$*, is tiled by the images of the open positroid varieties* $\mathring{\Pi}_f$ *of the TNN Grassmannian* $\mathrm{Gr}_{k,n}^{\mathbb{R},\geqslant 0}$ *by* $\widetilde{Z}$*.*

- *The Amplituhedron $\mathcal{A}_{n,k,m}(Z)$ is a (spectral) algebraic space over $\mathbb{C}$ for a certain Grothendieck étale ($\infty$-)topology*
- *The images of (real) open positroid varieties $\mathring{\Pi}_f$ of the TNN Grassmannian $\mathrm{Gr}_{k,n}^{\mathbb{R},\geqslant 0}$ spectral algebraic space are quasiaffine $\infty$-substacks of the $\infty$-étale TNN Grassmannian spectral algebraic space (over $\mathbb{C}$)*
- *Over $\mathbb{C}$, the Amplituhedron map $\widetilde{Z}$ gives rise to a geometric morphism of $\infty$-topoï, the* Amplituhedron functor*, from the TNN Grassmannian spectral algebraic space to the Amplituhedron spectral algebraic space.*
- *This Amplituhedron functor is a geometric morphism which has plenty of adjoints in a six-functor formalism and respects all formal infinitely thickening differential cohesive contexts, which allow us the functoriality of many characteristics of the TNN Grassmannian structured Schubert geometries in the sense of J. Lurie.*

#### 1.4.3.6 Main theorem 1 of Part II: The existence of the Dual Amplituhedron.

Even if the first part has been developped beyond, we can enounce now the first principal theorem of the second part of this thesis, inherited from the compact holistic writing of the positive Grassmannian (section 2.11).

THEOREM 1.4.13 (**Main Theorem 1: The existence of the Dual Amplituhedron**) : *Let $Z$ be a full rank $(k+4)\times n$ real matrix in $\mathbb{R}^{(k+4)\times n}$ with all maximal minors are positive, inducing a map:*

$$\widetilde{Z} : \mathrm{Gr}_{k,n}^{\geqslant 0} \longrightarrow \mathrm{Gr}_{k,k+m}$$

*defined by*

$$\widetilde{Z}\left(\langle v_1, \ldots, v_k\rangle\right) := \langle Z(v_1), \ldots Z(v_k)\rangle$$

*where $\langle v_1, \ldots, v_k\rangle$ represents with $k$ basis vectors, an element of $\mathrm{Gr}_{k,n}^{\geqslant 0}$, the real TNN Grassmannian, defines the so-called* (tree) Amplituhedron $\mathscr{A}_{n,k,4}$ *as the image $\widetilde{Z}\left(\mathrm{Gr}_{k,n}^{\geqslant 0}\right)$* *of* $\mathrm{Gr}_{k,n}^{\geqslant 0}$ *inside* $\mathrm{Gr}_{k,k+m}^{\geqslant 0}$.
*Then we can construct a* unique non-trivial topological dual Amplituhedron $\mathscr{A}_{n,k,4}^{\star}$ *to this (tree) Amplituhedron* $\mathscr{A}_{n,k,4}$ *satisfying the following properties:*

- *the prestable* $\infty$*-category* $\mathsf{QCoh}(\mathscr{A}_{n,k,4}^{\star})$ *of quasi-coherent* $\infty$*-sheaves on* $\mathscr{A}_{n,k,4}^{\star}$ *is the monoidal dual of* $\mathsf{QCoh}(\mathscr{A}_{n,k,4})$ *for Lurie's tensor product of* $\infty$*-categories in the subcategory* $\mathsf{QStk}(\mathscr{A}_{n,k,4})$ *of the* $\infty$*-category of Grothendieck prestable tensor* $\infty$*-categories* $\mathsf{Groth}_{\infty}$*, where the* $\infty$*-subcategory* $\mathsf{QStk}(\mathscr{A}_{n,k,4})$ *is made of quasi-coherent* $\infty$*-stacks on* $\mathscr{A}_{n,k,4}$ *endowed with the étale* $\infty$*-topology*
- *this dual Amplituhedron* $\mathscr{A}_{n,k,4}^{\star}$ *has a canonical* $\infty$*-sheaf of* $\mathbb{E}_{\infty}$*-rings,* $\mathscr{O}_{\mathscr{A}_{n,k,4}}^{\star}$*, which is the Verdier dual to the canonical* $\infty$*-sheaf* $\mathscr{O}_{\mathscr{A}_{n,k,4}}$ *in the sense that* $\mathscr{O}_{\mathscr{A}_{n,k,4}}^{\star} = \mathrm{Hom}_{\mathcal{D}^{+}(\mathsf{Coh}(\mathscr{A}_{n,k,4}))}(\mathscr{O}_{\mathscr{A}_{n,k,4}}, \omega_{\mathscr{A}_{n,k,4}})$ *where* $\omega_{\mathscr{A}_{n,k,4}}$ *is the dualizing complex in the* $\infty$*-category of* $\mathbb{E}_{\infty}$*-spectra.*

#### 1.4.3.7 Main theorem 2 of Part II: The consistency of the $\mathcal{BCFW}$-structure cells of the topological Dual Amplituhedron.

And we can enounce the second, central, and main theorem before concluding the thesis with the following third theorem in the next subsection:

THEOREM 1.4.14 (**Main Theorem 2: The Schubert geometry of the Dual Amplituhedron**) : *The unique non-trivial topological dual Amplituhedron* $\mathscr{A}_{n,k,4}^{\star}$ *constructed in theorem 3.4.7 respects the Schubert cells of the (tree) Amplituhedron* $\mathscr{A}_{n,k,4}$*, that is it satisfies the following properties, such that the dual Amplituhedron admits a tiling in dual* $\mathcal{BCFW}$ *open positroid cells whose Zariski closure overlap themselves along other dual* $\mathcal{BCFW}$ *open positroid cells, that is:*

- $\widetilde{Z}(\Pi_f)^{\star} = \overline{\widetilde{Z}(\mathring{\Pi}_f^{\star})}$
- $\widetilde{Z}(\Pi_f)^{\star} \backslash \widetilde{Z}(\mathring{\Pi}_f)^{\star} = \coprod_{g \in \mathscr{B} \subset \mathscr{B}(k,n)} \widetilde{Z}(\mathring{\Pi}_g)^{\star}$
- *this* $4k$*-dimensional* $\mathcal{BCFW}$ *tiling is a perfect partition* $\mathscr{A}_{n,k,4}^{\star} = \coprod_{f \in \mathscr{B}(k,n)} \widetilde{Z}(\mathring{\Pi}_f)^{\star}$

*Another characterization is that the $\infty$-topos constructed from the topological dual Amplituhedron $\mathscr{A}_{n,k,4}^{\star}$ has a $\mathcal{G}^{\star}$-structure in the sense of Lurie ([91]) where $\mathcal{G}^{\star}$ is the $\infty$-topos spanned by the perfect perverse complexes of $\infty$-sheaves ([103] Chap 9), that is Deligne's intersection complexes of $\infty$-sheaves with respect to the initial Whitney stratification, which are dualizable as compact objects, and preserved as perfect objects by the decomposition theorem for perverse sheaves, as the Amplituhedron map $Z$ is semi-small and proper, with enough perfect perverse sheaves to cover the $\infty$-perverse geometry.*
*This means that the dual Amplituhedron $\mathscr{A}_{n,k,4}^{\star}$ admits a $\mathcal{G}^{\star}$-structure, where $\mathcal{G}^{\star}$ is the $\infty$-topos spanned by the perfect perverse complexes of $\infty$-sheaves in the heart $\mathsf{QCoh}(\mathscr{A}_{n,k,4}^{\star})^{\heartsuit}$ of the $\infty$-category $\mathsf{QCoh}(\mathscr{A}_{n,k,4}^{\star})$, which admits a triangulated t-structure preserved by duality.*
***Endly, a last identical characterization is that the Schubert geometry of the Amplituhedron is stable with respect to taking the dual $\mathscr{A}_{n,k,4}^{\star}$ of the Amplituhedron.***

#### 1.4.3.8 Main theorem 3 of Part II: The volume of the Dual Amplituhedron.

And we can enounce the third theorem concluding the thesis:

THEOREM 1.4.15 (**Main Theorem 3: The volume of the Dual Amplituhedron**) : *The existence of the Dual Amplituhedron and the respect of the inital Schubert geometry of the Amplituhedron allowing to work further in theoretical physics on the Standard Model and Superstring Theory duality, provided we have a volume for this Dual Amplituhedron, here is a preliminary version of it:*

$$\mathscr{V}\mathrm{olume}(\mathscr{A}_{n,k,4}^{\star}) := \int_{\mathscr{A}_{n,k,4}^{\star}} d_{\mathrm{dRh}}^{\underline{\xi}} \mathscr{A}_{n,k,4}^{\star} \tag{1.4.16}$$

*where $d_{\mathrm{dRh}}^{\underline{\xi}}$ denotes the infinitesimal De Rham volume differentiation developped in section 2.8 and $\underline{\xi}$ is the dimension of the dual Amplituhedron $\mathscr{A}_{n,k,4}^{\star}$.*

# Chapter 2

# PART I: Algebraic Geometry and Synthetic Differentiable Geometry

Traditionnaly, in differentiable geometry, we study differentiable manifolds – whether smooth or not – equipped with a metric.

These are Hausdorff separable topological spaces locally isomorphic to $\mathbb{R}^n, n \in \mathbb{N}$. For each point $x$ on a differential manifold $X$, we can define a vector space $T_xX$, its tangent vector space, along with its formal linear dual with respect to the base field $\mathbb{k}$, its cotangent space $T_xX^*$. The scalar product on $T_xX$ can be chosen to be continuous in $x$, which then gives us a Riemannian metric. The classical non-degenerate Euclidian quadratic forms give us certain standard Riemannian metrics; Lorentzian metrics of signature $(1, n)$ on the quadratic form give us Riemannian metrics for general relativity. Intuitively, one can say that a metric is the measure of infinitesimal volumes on a manifold, and theoretically this concept can be translated in equations and formal diagrams through what is called formal geometry. We will return to this notion.

To return to how the choice of continuity is made for the scalar product on the tangent space, we must in fact *glue* together the vector spaces $T_xX$; this is what we call a vector bundle with some extra conditions for a *good* gluing. More formally, the continuous choice of $T_xX$ along the $x$ is then the choice of a *continuous* section $s : X \longrightarrow TX$ such that $p \circ s = \mathrm{Id}$ where $p : TX \longrightarrow X$ is the projection with each fiber : $p^{-1}(x) = T_x X$.

We are then faced with many possible generalizations, and even problems. What happens if the fiber $F_x := T_x X$ have a structure other than that of a vector space $T_x X$? This involves generalizing to sheaves with values in whatever we could think of, vector spaces for vector bundles, but also simply sets, or even a little more: abelian groups, or modules over a ring, modules over an $\mathbb{E}_\infty$-spectrum, but also complexes of whatever one wants, and even sheaf of categories, all of these with the good compatibility conditions? We will then speak of stacks, sheaves on stacks, quasi-coherent sheaves on stacks, and more: quasi-coherent stacks.

The second scenario is a genuine problem which, through a well-chosen formalism, subsumes the precedent generalization: What if the gluing were not so *intuitive*, that is, *topological* with open sets, covers, and intersections of open sets ? For example, what happens if the manifold is not *coherent* as expected, in the usual sense of the term *coherent*? That is to say, for example, that the infinitesimal volumes do not glue as expected? The mathematical term and the world of cohesiveness encompass all theses *formal problems*, but they must also be resolved.

So let us, after a brief digression into philosophical lexicology, which will become relevant later on, and a brief evocation of deformation theory and formal neighborhoods, begin with an overview of the *method* — or rather, to do justice to the significance of this approach to objects and morphisms via category theory and then topos and $\infty$-topos, let us call it: — A. Grothendieck's philosophy of functor fo points, which essentialy asserts that: "*Everything is a functor*".

### 2.0.1 Emmanuel Kant's philosophy of space, pure synthetic judgments and Lawvere's cohesion.

About the concept of "*space*", we are therefore led to examine first its "*nature*", its "*essence*", or "*son essence (in french)*", or "*sein Wesen (in German philosophical terminology)*", although the following one below would be perhaps the best of all denomi-nations for further discussions about: "*the categories of being/not being/nothing*", but also: "*modalities of possibility/existence/necessity*" and more mathematically speaking: "*adjunctions/idempotent (co)monads/algebras of Kleisli and Eilenberg-Moore*" and logically speaking "*categorical logic/homotopy type theory/relation between ∞-topoï, semantics and logic of space*".

Since Emmanuel Kant tells us it is one of the two, (with that of "*time*", on which we return at the end of this text), let us call that of "*space*" at this level an "a priori form of intuition". It does not learn us much mathematically but let us justify this temporary name:

Indeed, in E. Kant's twelve categories of being that allow "pure synthetic judgments", in particular "*Mathematics*", we see the subdivision in four categories:

- quantity (unity, plurality, totality)
- quality (reality, negation, limitation)
- relation (subsistence, causality, community)
- modality (possibility, Dasein, necessity)

The "*Dasein*" seems to allow for a "pure synthetic judgment" (so it allows to do Mathematics) within the category of "*Modality*". We will return to this major concept below.

In summary, "Space" is an "a priori form of intuition". That is what we have to examine here.

And to begin with, we could say that intuitively, a set of points is "*glued*" together in a way to examine, and these points may not behave freely and independently, since they would then not form a "*set of points*". We are led to Lawvere's concept of "cohesion".

### 2.0.2 A first look at infinitesimal neighborhood and deformation theory.

The set of points in a differential (separated Hausdorff) manifold is the set of continuous maps $* \longrightarrow X$. In algebraic geometry, the situation is more complicated between the reduced case and the unreduced case, since a scheme can have non-closed points, but still have generic points. We shall then see the base ring $\Bbbk$ as a point and all $\Bbbk$-algebras as "*extra points*".

Earlier, we have seen that all the geometrical data of a differential manifold $X$ are contained in all the vector bundle $TX$ over it, so what if we add these data to the data of a space? Since the theory of deformation constitutes a field of research in its own right — one that alone warrants entire sections of the text that follows and forms one of the central themes of this entire thesis — we will address this topic in this introductory subsection only briefly and with caution:

If we consider the infinitesimal neighborhood of each point $p \in X$ as a manifold, what if we add to the initial manifold $X$ the data of the "*tangent spaces to each infinitesimal neighborhood manifold* $T_pP$" $\forall p \in X$ (always roughly speaking since this needs a correct definition), where we would define $T_pP$ to be the formal neighborhood of a point $p \in X$, that is (always roughly speaking) the formal completion $T_pP^\wedge$ where $T_pP$ would be more than $\varprojlim \Bbbk[\epsilon]/(\epsilon^n)$ for $n \in \mathbb{N}$ but the universal pullback of the following diagrams, where $R$ is a ring, $I$ a nilpotent ideal of $R$, and $X$ a smooth scheme above $S$?

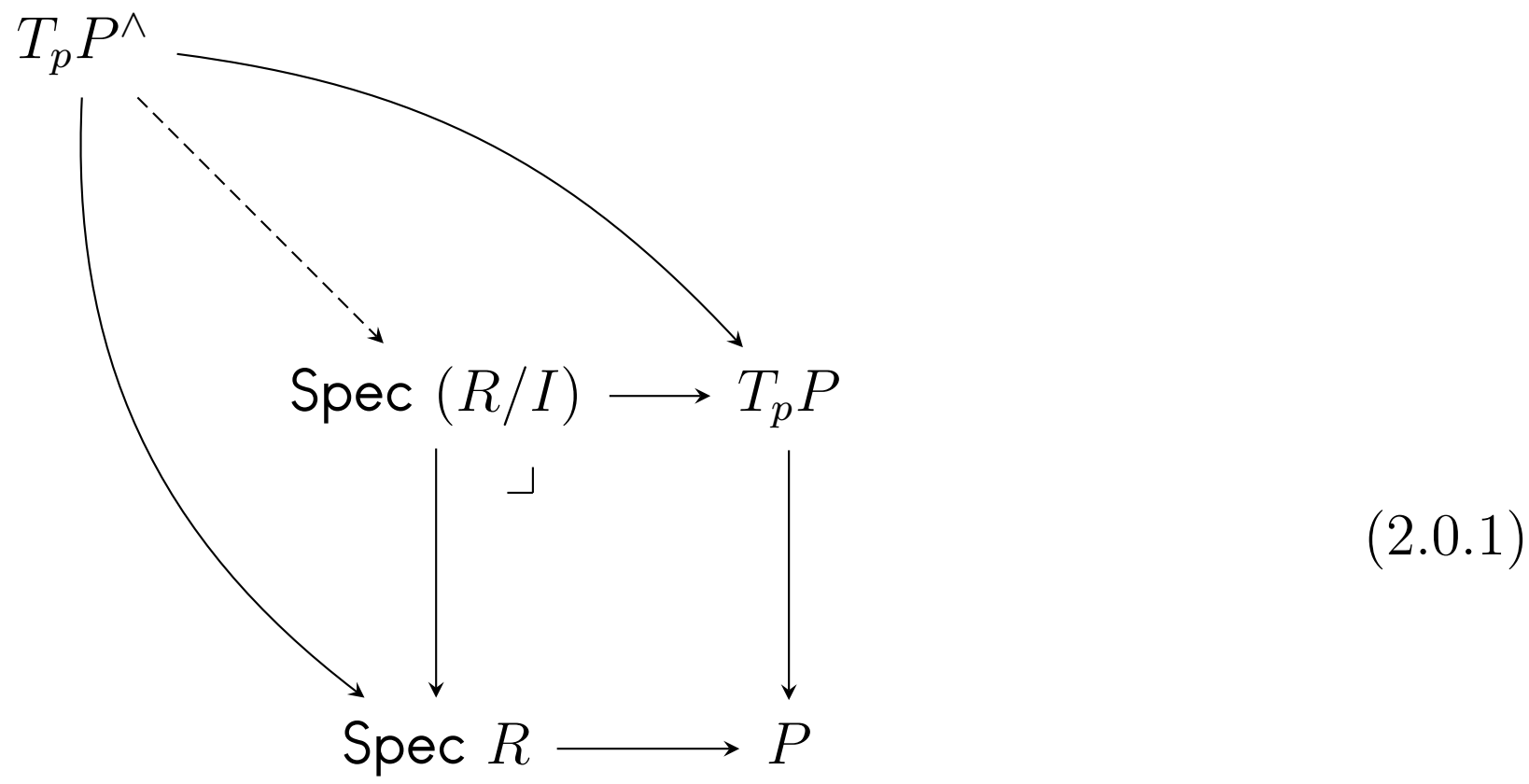

(2.0.1)

While we approach a "*not too bad*" definition, you will remark immediately that this definition is incorrect, and does not correspond to the authentic deformation theory. Moreover it contains the restriction of smoothness that A. Grothendieck stated in his definition of the infinitesimal site, that is just what we want to avoid in the case of Schubert geometry. We will return on this point immediatly below.

To return briefly on our construction of a space equipped with its formal neighbourhoods and their tangent spaces, we learn at least that deformation theory (and later spectral algebraic geometry and structured spaces) and differential cohesion are at the heart of our problems, in particular that is Jacob Lurie, Urs Schreiber and al. are our two primary sources of inspiration, knowledge, personal references, and gratitude, which we humbly acknowledge here.

Let's now address all the issues mentioned above, proceeding with caution as we move forward through the coming stages.

### 2.0.3 Algebraic geometry and the philosophy of functor of points.

Let $\mathcal{C}$ be a locally small category (that is the $\mathrm{Hom}_{\mathcal{C}}(-,-)$ are sets), and let $x$ be an object of the category $\mathcal{C}$. Let $\mathscr{F}$ be any functor from $\mathcal{C}$ to the category $\boldsymbol{\mathcal{S}}\mathbf{et}$ of sets:

$$\mathscr{F} : \mathcal{C} \longrightarrow \boldsymbol{\mathcal{S}}\mathbf{et}$$

The following Yoneda's lemma, where "Nat(-,-)" means : "*natural transformation*", tells us that as complicated a natural transformation from $h_x$ to $\mathscr{F}$ can be, in fact all is under control of the image of the identity morphism: $\mathbf{1}_x : x \longrightarrow x$, as the following tells us:

THEOREM 2.0.2 (Yoneda lemma, [121] Thm 2.2.4) : *Let $\mathcal{C}$ be a locally small category, and let $x$ be an object of the category $\mathcal{C}$. Let $\mathscr{F}$ be any functor from $\mathcal{C}$ to the category $\boldsymbol{\mathcal{S}et}$ of sets. We then have a bijection between the collection of natural transformations from the (representable) functor $h_x := \mathrm{Hom}_{\mathcal{C}}(x, -)$ to the functor $\mathscr{F}$ and the set $\mathscr{F}x$:*

$$\begin{cases} \mathrm{Nat}(h_x, \mathscr{F}) := \mathrm{Nat}(\mathrm{Hom}_{\mathcal{C}}(x, -), \mathscr{F}) & \xrightarrow{\sim} & \mathscr{F}x \\ \alpha : \mathrm{Nat}(h_x, \mathscr{F}) := \mathrm{Hom}_{\mathcal{C}}(x, -) \Rightarrow \mathscr{F} & \longmapsto & \alpha_x(\mathrm{Id}_x) \in \mathcal{S}\mathbf{et} \end{cases} \tag{2.0.3}$$

*The bijection is natural in both $x$ of $\mathcal{C}$ and $\mathscr{F}$ from $\mathcal{C}$ to $\boldsymbol{\mathcal{S}et}$. In particular, we obtain that this collection of natural transformations is a set, even though the category $\mathcal{C}$ is locally small but not necessarily small.*

DEFINITION 2.0.4 : *A functor $\mathscr{F} : \mathcal{C} \longrightarrow \boldsymbol{\mathcal{S}et}$ is said* representable *if exist an object $x \in \mathcal{C}$ and a natural transformation which is an isomorphism of functors $h_x \xrightarrow{\sim} \mathscr{F}$.*

Let us consider now a scheme $\mathscr{X}$ (over a field $\Bbbk$ for example). Here, we use the category $\mathcal{C} := \mathbf{Sch}$, and $x := \mathscr{X}$ as an object of the category of schemes.

DEFINITION 2.0.5 : *The functor $h_{\mathscr{X}} : \mathbf{Sch}^{\mathrm{op}} \longrightarrow \mathcal{S}\mathbf{et}$ is said* the functor of points *of $\mathscr{X}$. For all other schemes $\mathscr{Y}$ in the category* **Sch***, then $h_{\mathscr{X}}(\mathscr{Y}) = \mathrm{Hom}_{\mathbf{Sch}}(\mathscr{Y}, \mathscr{X})$ is* the set of $\mathscr{Y}$-points of $\mathscr{X}$.

REMARK 2.0.6 (Important) : It is in general easier to construct a functor which is geometrically interessant, to take pullbacks and pushouts on it, etc. than to have to do it directly with the scheme, or, worse, to have to pass through $\mathbf{Rings}^{\mathrm{op}}$ (for schemes). Once the geometric functor constructed, then it will be time to ask if this functor is a geometric object, that is if the functor of points is representable by a scheme (or later by an algebraic space, etc.), all of this being well-understood that we identify the scheme (or later the algebraic space) with its functor of points.

REMARK 2.0.7 : Any functor $\mathscr{F} : \mathcal{C} \longrightarrow \boldsymbol{\mathcal{S}}\mathbf{et}$ from a small category $\mathcal{C}$ to the category $\boldsymbol{\mathcal{S}}\mathbf{et}$ of sets is equivalent to an inductive limit of representable functors: $\mathscr{F} \simeq \varinjlim_{x \in \mathcal{C}} h_x$. (This remark will play a role in *accessibility* for topoï and $\infty$-topoï (the $\infty$-Yoneda) and can be skipped for the moment).

### 2.0.4 Cohomology theories of algebraic, analytic and differential varieties.

Note that we have not yet introduced the object of interest — the totally nonnegative Grassmannian — and that at this stage, we do not know whether we should consider it within the framework of algebraic, analytic, or differential geometry. The answer lies in the multitude of equivalence theorems between cohomology theories that were stated during the second half of the 20th century, leading in particular to the one that resolved "*Weil's conjectures*". I speak of "*crystalline cohomology*" and its "*daughters, or derived enhancements*".

We will indeed follow the same path and use the same framework involving these three fields: algebraic geometry with Jacob Lurie's spectral algebraic geometry, analytic geometry with Alexander Grothendieck's infinitesimal site and his equivalence theorems between singular and de Rham cohomologies, and differential geometry with the development of Lawvere's and Schreiber's differential cohesion in the form of infinitesimally thickened cohesive $\infty$-topoï toward the physicists' superpoint $\mathbb{R}^{0|1}$.

However, we will add our own personal touch to all these great theories. I am referring to another major theory that was invented for this purpose and deals with "*singularities as opposed to regularity or smoothness*" — which had always been the central characteristic added in problems to solve until the stabilization of theories and the "*emergence of the necessity of resolving singularities*" — and which "*we also must resolve*" in our case. At least, or at the same time, will we able to begin to open up a field of research for positive Schubert geometry. We are speaking here of "*perverse sheaves and intersection* *cohomology*". Indeed, we will bring into dialogue the structured spaces of spectral algebraic geometry and the structured sites of thickened points of differential cohesion. We will then, hopefully and modestly, be able to open a "*brave new path for the algebraic geometry and combinatorics of positive Schubert geometry*".

### 2.0.5 The Grassmannian and the positroid variety: combinatorial objects, schemes and moduli problems.

#### 2.0.5.1 The Grassmannian and the positroid variety as combinatorial objects.

Let $\mathrm{Gr}(k,n)$ be the set of $k$-dimensional subspaces in $\mathbb{C}^n$.
If $V \in \mathrm{Gr}(k,n)$, this point of the Grassmannian can be represented by a $k \times n$ matrix $M$ of full rank $k$, whose rows span the subspace $V$. We can then associate the *Plücker coordinates*: $(P_I)_{\binom{[n]}{k}}(V) \in \mathbb{P}^{\binom{[n]}{k}}\mathbb{C}$ of maximal minors of $V$ in columns $I \in \binom{[n]}{k}$. These Plücker coordinates satisfie the *Plücker relations*, which are quadratic in the Plücker coordinates. This endowed the Grassmannian $\mathrm{Gr}(k,n)$ with a structure of algebraic projective variety. The totally nonnegative (TNN) Grassmannian $\mathrm{Gr}_{\geqslant 0}(k,n)$ is then defined after G. Lusztig and A. Postnikov by the subset of $V \in \mathrm{Gr}(k,n)$ such as $(P_I)_{\binom{[n]}{k}} \in \mathbb{R}_{\geqslant 0}^{\binom{[n]}{k}}$.
Let then define the *matroid* $\mathcal{M}_V$ of $V$ in $\mathrm{Gr}(k,n)$ by $\mathcal{M}_V = \{I \in \binom{[n]}{k}, (P_I)_{\binom{[n]}{k}} \neq 0\}$.
If $V \in \mathrm{Gr}_{\geqslant 0}(k,n)$, the matroid $\mathcal{M}_V$ is called a *positroid*.
We are then able to define the *open positroid variety* $\overset{\circ}{\Pi}_{\mathcal{M}}$ of a positroid $\mathcal{M}$ as :

$$\overset{\circ}{\Pi}_{\mathcal{M}} = \{V \in \mathrm{Gr}(k,n) \mid P_I(V) = 0 \text{ if } I \notin \mathcal{M} \text{ and } P_I(V) \neq 0 \text{ if } I \in \mathcal{M}\}$$

In terms of ring of coordinates, for a positroid $\mathcal{M}$, we then have :

$$\mathbb{C}[\overset{\circ}{\Pi}_{\mathcal{M}}] = (\mathbb{C}[\mathrm{Gr}(k,n)]/\langle P_I \mid I \notin \mathcal{M}\rangle)_{(P_I \mid I \notin \mathcal{M})}$$

where $A_{(f)}$ for a ring $A$ and $f \in A$, means the localization of the ring $A$ at $f$.

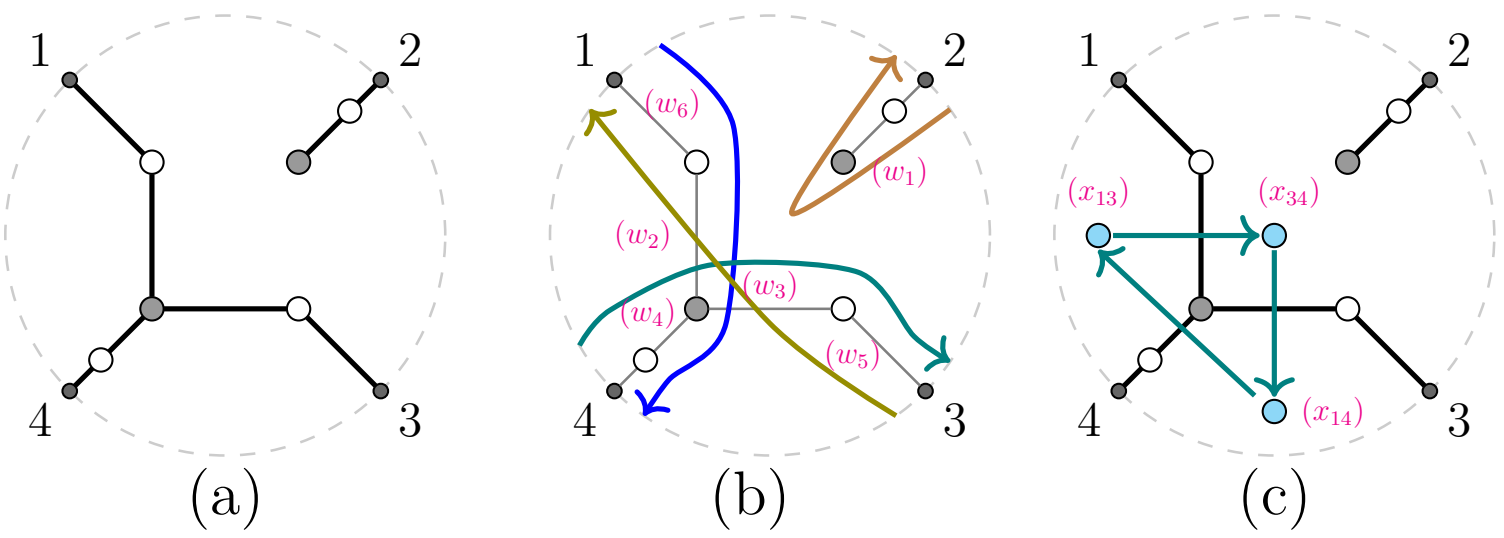

Figure 2.1: (a) One of the representant of the reduced planar bipartite graphs $G$ of the 3-dimensional open positroid variety of the TNN Grassmannian $\mathrm{Gr}_{\geqslant 0}(2,4)$; (b) strands of the strand permutation $f_G$ in $G$ with edge weights $w(e)$, where the unmarked edges have weight 1; (c) The quiver associated with this plabic graph $G$, whose all vertices are all frozen, and all variables $(x_{ij})_{i,j}$ associated in the corresponding cluster algebra are all non mutable frozen variables.

#### 2.0.5.2 The plabic graphs and the open positroid variety $\mathring{\Pi}_f$.

We will here briefly recall the definition of the open positroid variety, constructed thanks to plabic graphs of A. Postnikov [114], with illustrations coming from the example of the open positroid variety of dimension 3 in $\mathrm{Gr}_{\geqslant 0}(2,4)$.

Let then define $G$, a *plabic graph* by a combinatorial object made by a planar bipartite graph embedded in a disk with $n$ black boundary vertices, labeled clockwise by the integers from 1 to $n$, where $n$ is the dimension of the space $\mathbb{C}^n$ in which the $k$-dimensional subspaces $V$ of $\mathrm{Gr}(k,n)$ are embedded.

The number $k$ corresponds to the differences between the number of black vertices and the white vertices of $G$. We are then able to make pairs of a black vertex with a white vertex, leaving eventually black vertices not involved on the boundary of $G$. Each possible bunch of pairs are called an *almost perfect matching* of $G$, and is denoted by $\mathcal{A}$.

We can show that the number of black vertices involved on the boundary of $G$ of a strand permutation is $k$. Let denote $I_{\mathcal{A}} \in \binom{[n]}{k}$, the set of integers involved in an almost perfect matching $\mathcal{A}$. We call then a *strand permutation* in $G$, a path from a black vertex of the boundary to an other, which "turns" right at each black vertex, and "turns" left at each white vertex. We can show that we obtain then a permutation in $\mathfrak{S}_n$ of the boundary vertices.

This plabic graph $G$ is said to be *reduced* if among all plabic graphs with the same strand permutation in $\mathfrak{S}_n$, the number of faces of $G$ is the lowest. Following, this reduced plabic graph $G$ is then said to be *weighted* if a positive real weight $w(e)$ is assigned to each edge $e$ of the plabic graph $G$. The map from the collection of edges of the graph $G$ to $\mathbb{R}_{\geqslant 0}$ is denoted $wt : G \to \mathbb{R}^N_{\geqslant 0}$, where $N$ is the number of edges of $G$.

We then define the *boundary measurement map*, which sends each weighted plabic graph $(G, wt)$ to the collection over each $I \in \binom{[n]}{k}$ of the sums over all possible perfect matchings $\mathcal{A}$ such that $I_{\mathcal{A}} = I$ of the product of weights involved in each perfect matching, such as :

$$\mathcal{M}\text{eas}(G, wt) = (P_I)_{\binom{[n]}{k}} = \left( \sum_{\mathcal{A} | I_{\mathcal{A}} = I} \prod e \in \mathcal{A} w(e) \right)_{\binom{[n]}{k}}$$

We can then show that for a fixed weighted reduced plabic graph $(G, wt)$, the boundary measurement map corresponds to the Plücker coordinates of a certain $V \in \text{Gr}_{\geqslant 0}(k, n)$ which is unique.

By defining an *bounded affine permutation* of a plabic graph, denoted by $f \in \mathscr{B}(k, n)$ from $\mathbb{Z}$ to $\mathbb{Z}$ by the suitable affine linearisation of the strand permutation, and modulo reduction of plabic graphs which are manipulations of vertices leaving invariant the combinatorics of the plabic graph, we then have bijections between plabic graphs, strand permutations, bounded affine permutations, boundary measurement maps, Plücker coordinates, and the $k$-dimensional subspaces $V$ in $\mathbb{C}^n$, viewed as points in $\text{Gr}_{\geqslant 0}(k, n)$.

By fixing the plabic graph $G$ and varying the map of weights on the edges of $G$, we then obtain an other definition of the *open positroid variety*, which is now indexed by the bounded affine permutations $f$ in $\mathscr{B}(k, n)$, and we can show the following disjoint stratification of $\text{Gr}_{\geqslant 0}(k, n)$ :

$$\text{Gr}_{\geqslant 0}(k, n) = \coprod_{f \in \mathscr{B}(k,n)} \Pi_f^{>0}$$

where the *open positroid cell* $\Pi_f^{>0}$ is defined such as :

$$\Pi_f^{>0} = \{V \in \mathrm{Gr}(k,n) \mid P_I(V) = 0 \text{ if } I \notin \mathcal{M} \text{ and } P_I(V) > 0 \text{ if } I \in \mathcal{M}\}$$

with $\Pi_f$ is defined as the Zariski closure of the open positroid cell $\Pi_f := \overline{\Pi_f^{>0}}$, and where we get back the open positroid variety $\mathring{\Pi}_f$ as the open subvariety of $\Pi_f$ where a suitable collection of minors are not 0 (namely the minors of the *Grassmann necklace* defined from the Gale order between labelization of the faces of the plabic graph thanks to the strand permutation).

#### 2.0.5.3 Cluster algebras and rings of coordinates of $\mathbb{C}[\mathring{\Pi}_f]$.

Let $\mathbb{C}$ be the ground field and $(\mathcal{Q}, \underline{x})$ be an ice quiver (with frozen vertices) associated with possibly infinitely many variables $\underline{x} = x_1, x_2, x_3, \dots$ . Formally, a *cluster algebra* can then be defined from this quiver $(\mathcal{Q}, \underline{x})$, as a $\mathbb{C}$-subalgebra of the $\mathbb{C}$-algebra of Laurent polynomials in possibly infinitely many variables $\mathbb{C}[x_1^{\pm 1}, x_2^{\pm 1}, x_3^{\pm 1}, \dots]$, endowed with combinatorial mutation rules on the quiver $\mathcal{Q}$, and on the variables $x_1^{\pm 1}, x_2^{\pm 1}, x_3^{\pm 1}, \dots$. Except those coming from and ending at a frozen vertex, the arrows of the quiver can be modified ("mutated") under these combinatorial rules, and the variables can be mutated under the following rule, at vertex $k$ :

$$x'_k = \frac{\prod\limits_{i\,,\,i \to k} x_i + \prod\limits_{j\,,\,k \to j} x_j}{x_k}$$

Assuming that we can define suitable quivers for the Grassmannian and for the open positroid varieties, like briefly shown in the picture below, where each face correspond to a vertex, and each arrow links two vertex in direction such that the white nodes are on the left of the arrow in the case of the open positroid variety, and each vertex of the Grassmann necklace (along the boundary mostly) is frozen, the two results from [136] and [50] are then of crucial importance :

THEOREM 2.0.8 (Scott 06, [136]) : *Let $k \leqslant n$ be integers. The ring of coordinates $\mathbb{C}[\mathrm{Gr}(k,n)]$ of the Grassmannian $\mathrm{Gr}(k,n)$ is a cluster algebra.*

THEOREM 2.0.9 (Galashin–Lam 21', [50]) : *Let $f \in \mathscr{B}(k,n)$ be a bounded affine permutation for a plabic graph $G$. Then the ring of coordinates $\mathbb{C}[\mathring{\Pi}_f]$ of the open positroid variety $\mathring{\Pi}_f$ is a cluster algebra.*

PROPOSITION 2.0.10 : *If $f$ and $f'$ in $\mathscr{B}(k,n)$ are different bounded affine permutations for a* reduced *plabic graph $G$, then the open positroid varieties $\mathring{\Pi}_f$ and $\mathring{\Pi}_{f'}$ are disjoint.*

Let us call the (closed) *positroid variety*, denoted $\Pi_f := \overline{\mathring{\Pi}_f}$, the Zariski closure of an open positroid varieties $\mathring{\Pi}_f$. Then we have:

PROPOSITION 2.0.11 : *If $f$ and $f'$ in $\mathscr{B}(k,n)$ are different bounded affine permutations for a* reduced *plabic graph $G$, then the intersection of the (closed) positroid varieties $\Pi_f$ and $\Pi_{f'}$ is:*

- *a finite disjoint union of open positroid varieties $\mathring{\Pi}_g$, for a finite collection of $g \in \mathscr{B}(k,n)$*
- *a union of closed positroid subvarieties $\Pi_g$, for a finite collection of $g \in \mathscr{B}(k,n)$*

*In conclusion, the TNN Grassmannian is* "well-stratified" *in the sense of Whitney ([66]).*

Passing to the language of algebraic schemes, we have:

THEOREM 2.0.12 ([50]) : *Let $f \in \mathscr{B}(k,n)$ be a bounded affine permutation for a plabic graph $G$. Then the open positroid variety $\mathring{\Pi}_f$ is a quasi-affine scheme.*

#### 2.0.5.4 The Grassmannian as a scheme and a functor of points.

PROPOSITION 2.0.13 : *We have the following characterizations of a closed immersion of a subscheme $\mathscr{Z}$ in a scheme $\mathscr{X}$ :*

- *A morphism of schemes $f : \mathscr{Z} \hookrightarrow \mathscr{X}$ is a closed immersion*
- *For every open affine $\mathscr{U} = \mathsf{Spec}\ R \subseteq \mathscr{U}$, there exists an ideal $I \subset R$ such that $f^{-1}(\mathscr{U}) = \mathsf{Spec}\ (R/I)$ as schemes over $\mathscr{U}$*

- *There exists an open affine covering* $\mathscr{X} = \bigcup \mathscr{U}_j$, $\mathscr{U}_j = \mathsf{Spec}\ R_j$, *and for each* $j$, *there exists an ideal* $I_j \subset R_j$ *such that* $f^{-1}(\mathscr{U}_j) = \mathsf{Spec}\ (R_j/I_j)$ *as schemes over* $\mathscr{U}_j$
- *There is a quasi-coherent sheaf of ideals* $\mathscr{I}$ *on* $\mathscr{X}$ *such that* $f_*\mathscr{O}_{\mathscr{Z}} \simeq \mathscr{O}_{\mathscr{X}}/\mathscr{I}$ *and* $f$ *is an isomorphism of* $\mathscr{Z}$ *onto the global spectrum of the sheaf* $\mathscr{O}_{\mathscr{X}}/\mathscr{I}$ *over* $\mathscr{X}$ *where the global spectrum is defined by the following isomorphism:*

$$\mathrm{Hom}_{\mathscr{O}_{\mathscr{X}}-\mathrm{alg}}(\mathscr{O}_{\mathscr{X}}/\mathscr{I}, f_*\mathscr{O}_{\mathscr{Z}}) \simeq \mathrm{Hom}_{\mathrm{Sch}/\mathscr{X}}(\mathscr{Z}, \mathsf{GlobalSpec}(\mathscr{O}_{\mathscr{X}}/\mathscr{I}))$$

Let us return to the question of the representability of morphisms. The idea is to define the Grassmannian by its functor of points, and to say that this functor of points is representable by a scheme: the Grassmannian (seen as a scheme), and then to see it as a closed subscheme of the projective space. The same process will be used to define (closed) positroid variety as closed subschemes of the Grassmannian. The question of their gluing will come later, since we will have to swap the topology to the étale topology to obtain an algebraic space.

DEFINITION 2.0.14 : *Let* $\mathscr{F}$ *and* $\mathscr{G}$ *be two non necessarily representable functors:* $\mathscr{F}, \mathscr{G} : \mathbf{Sch}^{\mathrm{op}} \longrightarrow \boldsymbol{\mathcal{S}}\mathbf{et}$. *We say that a natural transformation between* $\mathscr{F}$ *and* $\mathscr{G}$ *is* representable *if and only if for each scheme* $\mathscr{X}$ *(seen as a functor) and each* $g : \mathscr{X} \longrightarrow \mathscr{G}$ *a natural transformation, the functor* $\mathscr{F} \times_{\mathscr{G}} \mathscr{X}$ *is representable (the fiber product is taken pointwise as the codomain of the functors is the category* $\boldsymbol{\mathcal{S}}\mathbf{et}$ *which has pullbacks and then fiber products).*

DEFINITION 2.0.15 : *Let* $\mathcal{P}$ *be a property of morphisms of schemes which is stable by composition with isomorphims on the left and on the right. We say that a representable morphism* $f : \mathscr{F} \longrightarrow \mathscr{G}$ *in the category of presheaves from* $\mathbf{Sch}^{\mathrm{op}}$ *to* $\boldsymbol{\mathcal{S}}\mathbf{et}$ satisfies the property $\mathcal{P}$ *if and only if for each scheme* $\mathscr{X}$ *and each morphism* $G : \mathscr{X} \longrightarrow \mathscr{G}$, *the second projection* $f_{(\mathscr{X})} : \mathscr{F} \times_{\mathscr{G}} \mathscr{X} \longrightarrow \mathscr{X}$ *satisfies the property* $\mathcal{P}$.

THEOREM 2.0.16 : *Let* $\mathscr{F} : \mathbf{Sch}^{\mathrm{op}} \longrightarrow \boldsymbol{\mathcal{S}}\mathbf{et}$ *be a functor such that:*

- $\mathscr{F}$ *is a (Zariski) sheaf*

- *$\mathscr{F}$ admits an open (Zariski) covering by Zariski representable functors (that is each $g_i$ is an open subfunctor of $\mathscr{F}$ (that is: is representable and is an open immersion for $\mathcal{P}$, that is: is an open immersion for the second component) and for each scheme $\mathscr{X}$, and all $g : \mathscr{X} \longrightarrow \mathscr{F}$, the images of $(g_i)_{(\mathscr{X})}$ are an open Zariski covering of the scheme $\mathscr{X}$)*

*then the functor $\mathscr{F}$ is representable.*

Let us now define the Grassmannian functor for all $\Bbbk$ a field, with $k \leqslant n \in \mathbb{N}$:

$$\mathbf{Grass}_{k,n} : \begin{cases} \mathbf{Sch}^{\mathrm{op}} \longrightarrow \boldsymbol{\mathcal{S}}\mathbf{et} \\ \mathsf{Spec}\ \Bbbk \longmapsto \left\{ \begin{array}{r} \mathscr{U} \subseteq \Bbbk^n \mid \mathscr{U} \text{ is a } k\text{-dimensional} \\ \Bbbk\text{-vector subspace of } \Bbbk^n \end{array} \right\} \end{cases}$$

and the Grassmannian functor for all $\mathscr{X}$ a scheme, with $k \leqslant n \in \mathbb{N}$:

$$\mathbf{Grass}_{k,n}^{(\mathscr{X})} : \begin{cases} \mathbf{Sch}^{\mathrm{op}} \longrightarrow \boldsymbol{\mathcal{S}}\mathbf{et} \\ \mathscr{X} \longmapsto \left\{ \begin{array}{l} \mathscr{U} \subseteq \mathscr{O}_{\mathscr{X}}^n \mid \mathscr{O}_{\mathscr{X}}^n/\mathscr{U} \text{ is an } (n-k)\text{-dimensional} \\ \text{locally free } \mathscr{O}_{\mathscr{X}}-\text{module} \end{array} \right\} \end{cases}$$

THEOREM 2.0.17 : *The functor $\mathbf{Grass}_{k,n}^{(\mathscr{X})}$ is representable, and it is by the scheme:*

$$\text{the Grassmannian } \mathrm{Gr}(k,n)$$

*Proof.* The functor $\mathbf{Grass}_{k,n}^{(\mathscr{X})}$ is a (Zariski) sheaf and admits an open Zariski covering by Zariski representable functors, where $I$ are $(n-k)$-subsets of $\{1,\ldots,n\}$:

$$i^I : \mathbf{Grass}_{k,n}^I \hookrightarrow \mathbf{Grass}_{k,n}^{(\mathscr{X})} \tag{2.0.18}$$

which is representable and a closed immersion, given on each scheme $\mathscr{X}$ by:

$$\mathbf{Grass}_{k,n}^I : \begin{cases} \mathbf{Sch}^{\mathrm{op}} \longrightarrow \boldsymbol{\mathcal{S}}\mathbf{et} \\ \mathscr{X} \longmapsto \left\{ \begin{array}{l} \mathscr{U} \in \mathbf{Grass}_{k,n}^{(\mathscr{X})} \text{ such that } \mathscr{O}_{\mathscr{X}}^I \hookrightarrow \mathscr{O}_{\mathscr{X}}^n \twoheadrightarrow \mathscr{O}_{\mathscr{X}}^n/\mathscr{U} \\ \text{is an } (n-k)\text{-dimensional locally free } \mathscr{O}_{\mathscr{X}}-\text{module} \end{array} \right\} \end{cases} \tag{2.0.19}$$

where each $\mathbf{Grass}_{k,n}^I$ is representable by the affine space $\mathbb{A}^{k(n-k)}$. □

#### 2.0.5.5 The closed and open positroid variety of the TNN Grassmannian.

With exactly the same process:

THEOREM 2.0.20 : *We have the following:*

- *The closed positroid variety, denoted* $\Pi_f := \overline{\mathring{\Pi}_f}$*, the Zariski closure of an open positroid variety* $\mathring{\Pi}_f$*, is a closed subscheme of its Grassmannian*
- *The open positroid variety, denoted* $\mathring{\Pi}_f$*, is an open subscheme of its closure, the closed positroid variety, and then an open subscheme of its Grassmannian*

## 2.1.A motivation for Derived Algebraic Geometry.

Most of the material presented here comes from an expository piece by B. Antieau. We follow his notations, which coincide with those of J. Lurie.

Let us pose the first question: What are our affine objects in this context?

So, let us consider the category of abelian groups, denoted $\mathcal{A}\mathrm{b}$, and $\mathscr{D}(\mathcal{A}\mathrm{b})_{\geqslant 0}$ be the connective part of the derived category of the abelian category of abelian groups. The connective part means that we work with complexes $C_*$ such that $H_i(C_*) = 0$ if $i < 0$. Let us endow this derived category with the derived tensor product $\otimes^{\mathbb{L}}$, and with the unit $\mathbf{1}_{\mathcal{A}}$ of the category $\mathcal{A}\mathrm{b}$ seen as a (derived) complex made of the unit $\mathbf{1}_{\mathcal{A}}$ concentrated in degree 0. Then $(\mathscr{D}(\mathcal{A}\mathrm{b})_{\geqslant 0}, \otimes^{\mathbb{L}}, \mathbf{1}_{\mathcal{A}})$ is a symmetric monoidal closed category ([56]). Let us note that Balmer's reconstruction of a scheme from a tensor triangulated category with unit [19] applies to this context.

Let us now look at commutative algebra objects in this symmetric monoidal closed category $(\mathscr{D}(\mathcal{A}\mathrm{b})_{\geqslant 0}, \otimes^{\mathbb{L}}, \mathbf{1}_{\mathcal{A}})$, and replace $\mathscr{D}(\mathcal{A}\mathrm{b})_{\geqslant 0}$ by connective (derived) chain complexes on an abelian category $\mathcal{A}$ with enough injectives endowed with its derived tensor product $\otimes^{\mathbb{L}}$. This leads us to consider *derived algebraic geometry* modeled on the category of *commutative simplicial rings*, or equivalently on the category of *commutative differential graded* $\mathbb{Q}$*-algebras*, or even more on the category of $\mathbb{E}_\infty$*-ring spectra* since we will use J. Lurie's framework in *spectral algebraic geometry.* Let us be more precise:

### 2.1.1 The combinatorial category of simplicial objects.

DEFINITION 2.1.1 ([116] ?) : *The* simplex category $\Delta$ *is the category of non-empty finite ordered sets together with order-preserving maps between them.*

We can visualize the objects of the category $\Delta$ as the collection $\{[n]\}_{n\in\mathbb{N}}$ and its maps as compositions of the following maps:

$$\dots[2] \; \begin{matrix} \xleftarrow{f_0} \\ \xrightarrow{d_0} \\ \xleftarrow{f_1} \\ \xrightarrow{d_1} \\ \xleftarrow{f_2} \end{matrix} \; [1] \; \begin{matrix} \xleftarrow{f_0} \\ \xrightarrow{d_0} \\ \xleftarrow{f_1} \end{matrix} \; [0] \; .$$

where $(f_i)_{i\in\mathbb{N}}$ are the face maps of simplices and $(d_i)_{i\in\mathbb{N}}$ are the degeneracy maps of simplices with compatibility rules (see [116]).

DEFINITION 2.1.2 : *Let* $\mathscr{C}$ *be a category. The* category of simplicial objects*, denoted* $s\mathscr{C}$ *is defined such that* $s\mathscr{C} := \mathfrak{Fun}(\Delta^{\mathrm{op}}, \mathscr{C})$*.* We will consider this category $s\mathscr{C}$ of simplicial objects of the category $\mathscr{C}$ as the category $s\mathscr{C} := [\,\Delta^{\mathrm{op}}\,,\,\mathscr{C}\,]$ of $\mathscr{C}-$valued presheaves on $\Delta$.

EXAMPLE 2.1.3 : If we take $\mathscr{C} = \mathcal{S}$ets, the category of sets, the category of simplicial sets $s\mathcal{S}$ets is a presheaf category, which is generated (Yoneda) by colimits of its representable functors, namely here the $n$-simplices functors $\Delta^n := \mathrm{Hom}_\Delta(-, [n])$ in $s\mathscr{C} := [\,\Delta^{\mathrm{op}}\,,\,\mathscr{C}\,]$. Moreover, we have the Quillen adjunction,which is a homotopically Quillen equivalence of homotopy categories under weak-equivalences, given by the following:

$$\text{(ho) } s\mathcal{S}\text{ets} \underset{|\cdot|_{\mathrm{geom}}}{\overset{\mathrm{Sing}_\bullet(-)}{\leftrightarrows}} \text{(ho) } \mathcal{T}\text{op} \qquad (2.1.4)$$

Since we will evolute in the category of finite CW-complexes in the perspective of the TNN Grassmannian, the Grothendieck homotopy hypothesis does not play a role in this model.

In this diagram, we have:

- the functor $|\,.\,|_{\text{geom}}$ is the *geometric realization*, that is the extension by colimits of the maps sending the $n$-simplex representable functors $\Delta^n$ to the topological $n$-simplex $\Delta^n_{\text{top}} := \{(x_0, \ldots, x_n) \in \mathbb{R}^{n+1} \mid x_i \geqslant 0 \text{ and } \sum x_i = 1\}$.

- the functor $\text{Sing}_\bullet(-)$ is the *singular complex construction* sending a topological space $X$ to the continuous maps from the cosimplicial object $\Delta^\bullet_{\text{top}}$ of $[\,\Delta\,,\mathscr{C}\,]$ to $X$ in the category of topological spaces $\mathcal{T}\text{op}$, that is $\text{Sing}_\bullet(X) := \text{Hom}(\Delta^\bullet_{\text{top}}, X)$ in $s\mathcal{S}\text{ets}$.

The functor $\text{Sing}_\bullet(-)$ is right adjoint to the functor $|\,.\,|_{\text{geom}}$ in a Quillen adjunction of the underlying model (see [56]).
In perspective of singular intersection homology that we see further, let us notice that if we form the abelian group $\mathbb{Z}[\text{Sing}_\bullet(X)]$, and the chain complexes on it: $C_\bullet(\mathbb{Z}[\text{Sing}_\bullet(X)])$, then for all $i \in \mathbb{Z}$, the homology $H_i(\mathbb{Z}[\text{Sing}_\bullet(X)])$ is isomorphic to the singular homology $H_i^{\text{sing}}(X, \mathbb{Z})$.

Example 2.1.5 : More generally, if we take $\mathscr{C} = \mathcal{A}\text{b}$, the category of abelian groups, the category of simplicial abelian groups $s\mathcal{A}\text{b}$ is homotopically Quillen equivalent to the category $\mathscr{D}(\mathcal{A}\text{b})_{\geqslant 0}$ of the connective part of the derived category of the abelian category of abelian groups. The equivalence is given by the Dold-Kan correspondence (or nerve/realization Quillen adjunction for the homology functor):

$$(\text{ho})\ s\mathcal{A}\text{b} \underset{\longrightarrow}{\overset{\longleftarrow}{\top}} (\text{ho})\ \mathscr{D}(\mathcal{A}\text{b})_{\geqslant 0}$$

EXAMPLE 2.1.6 : Finally, in the perspective of derived algebraic geometry, let us consider $\mathscr{C} = \mathcal{CA}\mathrm{lg}_{\Bbbk}$, the category of commutative algebras over the field $\Bbbk$. The category $s\,\mathcal{CA}\mathrm{lg}_{\Bbbk}$ of simplicial commutative algebras over the field $\Bbbk$ is $\mathfrak{Fun}(\Delta^{\mathrm{op}}, \mathcal{CA}\mathrm{lg}_{\Bbbk}) =: [\,\Delta^{\mathrm{op}}\,,\, \mathcal{CA}\mathrm{lg}_{\Bbbk}\,]$, that is the $\mathcal{CA}\mathrm{lg}_{\Bbbk}$-valued presheaves on the simplex category $\Delta$.
This category $s\,\mathcal{CA}\mathrm{lg}_{\Bbbk}$ of simplicial commutative algebras has an underlying simplicial set, and considering the weak equivalences as the weak equivalences of the underlying simplicial sets (or topological spaces thanks to Equation 2.1.4), we get a model category structure which allows us to take *resolutions* to define derived functors. Let us be more precise:
The following Quillen adjunction:

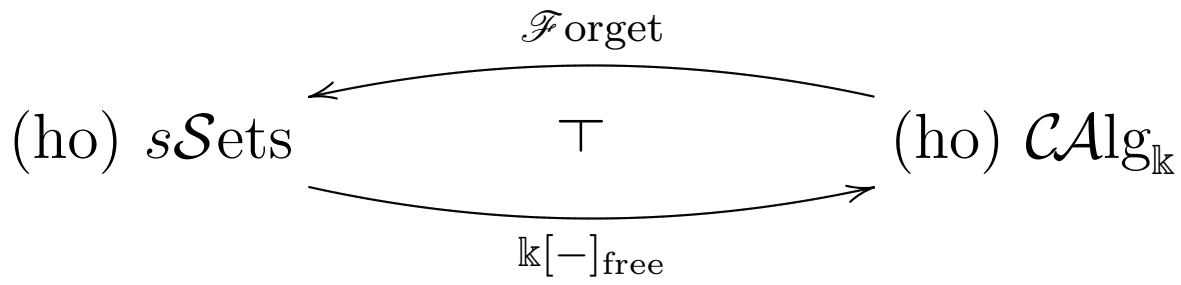

where the functor $\mathscr{F}$orget is the forgetful functor, forgetting the underlying structure of an algebra to get a set, and the functor $\Bbbk[-]_{\mathrm{free}}$ is the free construction $\Bbbk[S]$ of a commutative free $\Bbbk$-algebra from the underlying set $S$.
Then, the unit and counit maps of the adjunction allow us to construct a canonical simplicial commutative algebra $S_\bullet$ from a commutative algebra $R$ such that $S_\bullet \xrightarrow[\sim]{\text{htpy}} R$ is a homotopy equivalence, where $S_\bullet$ is a polynomial ring (in possibly infinitely variables) in each degree. This (canonical) simplicial commutative algebra $S_\bullet$ is, by definition, the "projective" resolution we use to define the following derived functor described below.
Let $\mathcal{CA}\mathrm{lg}_{\Bbbk}^{\mathrm{f.g.}}$ be the category of finitely generated $\Bbbk$-algebras, that is the category of finitely generated polynomial rings over the ground field $\Bbbk$, which underlies (classical) algebraic geometry. Let $F$ be a functor from $\mathcal{CA}\mathrm{lg}_{\Bbbk}^{\mathrm{f.g.}}$ to a category $\mathscr{C}$ with filtered colimits (we will see later the case of an $\infty$-category $\mathscr{C}$).
Let $\mathrm{Ind}(\mathcal{CA}\mathrm{lg}_{\Bbbk}^{\mathrm{f.g.}})$ be the Ind category obtained from $\mathcal{CA}\mathrm{lg}_{\Bbbk}^{\mathrm{f.g.}}$ by adjoining filtered colimits, allowing the number of variables of the polynomial rings to be infinite. This Ind category is a subcategory of $s\,\mathcal{CA}\mathrm{lg}_{\Bbbk}^{\mathrm{f.g.}}$ by considering constant diagrams for the filtering index category. Then we have the following diagram where the functor $\mathbb{L}F$ is a left Kan extension of the functor $F$ defined under the following diagram:

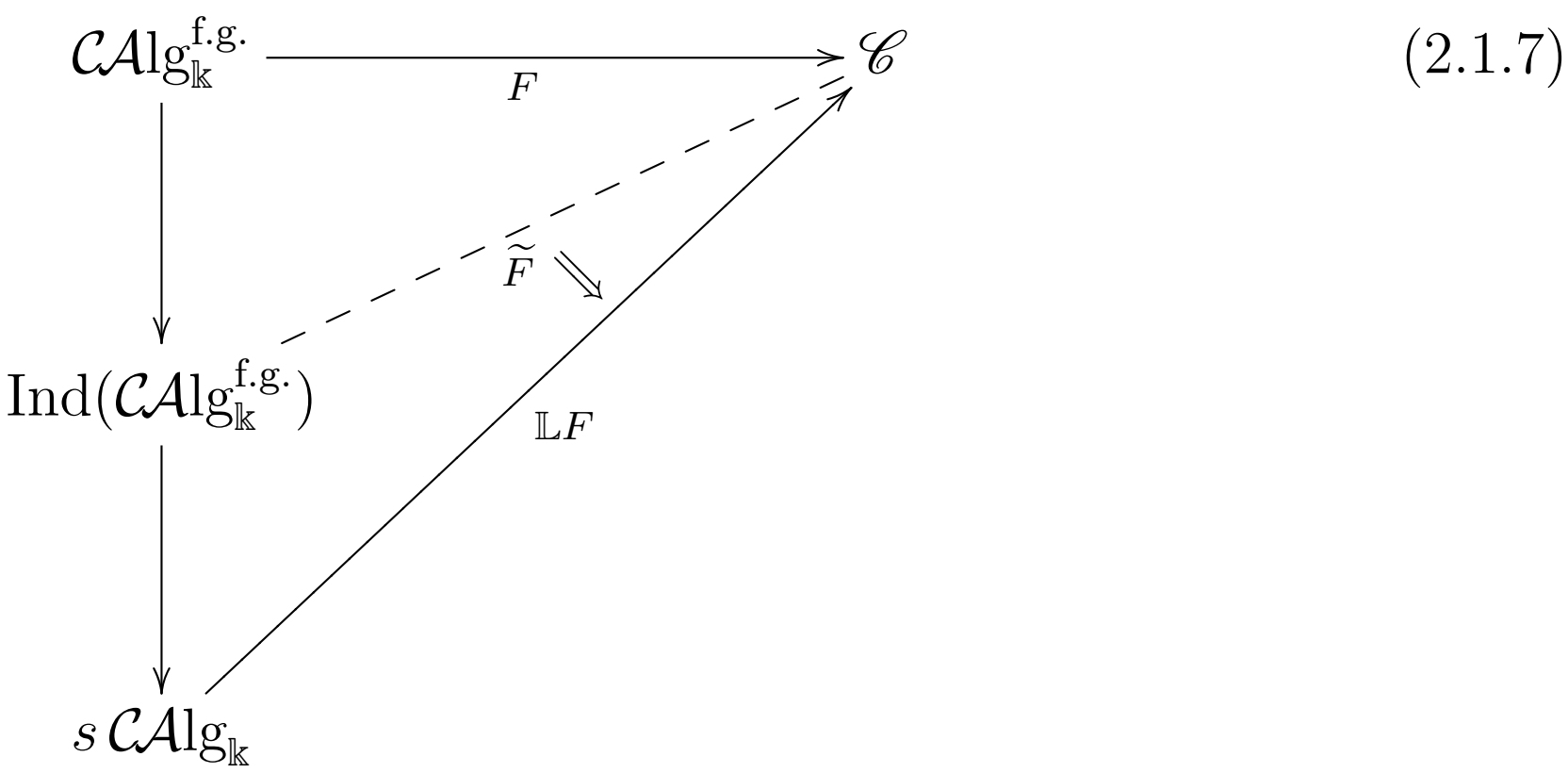

(2.1.7)

Indeed, from the resolution $S_\bullet \xrightarrow[\sim]{\text{htpy}} R$, we define the derived functor $\mathbb{L}F$ from the simplicial object $F(S_\bullet)$ in $s\mathscr{C}$ by taking the homotopy colimit on the simplicial category $\Delta$, which is also the geometric realization:

$$\mathbb{L}F(R) := \underset{\Delta}{\text{ho.colim}}\, F(S_\bullet) = |F(S_\bullet)|_{\text{geom}} \tag{2.1.8}$$

remark 2.1.9 : The derived functor is universal for this property.

remark 2.1.10 : For our further construction, let us remark right now that the simplicial commutative algebras will be *the coordinate rings of our derived affine schemes* in derived algebraic geometry.

definition 2.1.11 : *If* $\Bbbk$, $R$, *and* $S$ *are simplicial commutative rings, then the* Hom *between* $R$ *and* $S$ *in the category of simplicial commutative rings is called the mapping space, denoted* $\mathrm{Map}_{s\,\mathcal{CA}\mathrm{lg}_{\Bbbk}}(R, S)$, *and is a topological space (or an* $\infty$*-groupoid in the* $\infty$*-categories world, or a simplicial set in the simplicial sets world).*

example 2.1.12 : We have:

- $\mathrm{Map}_{s\,\mathcal{CA}\mathrm{lg}_{\Bbbk}}(\Bbbk, S) \xrightarrow[\text{htpy}]{\sim} *$, that is the point (or the contractible space) since $\Bbbk$ is the initial object of $s\,\mathcal{CA}\mathrm{lg}_{\Bbbk}$

- $\mathrm{Map}_{s\,\mathcal{CA}\mathrm{lg}_{\Bbbk}}(\Bbbk[t], S) \xrightarrow[\text{htpy}]{\sim} S$, but with $S$ viewed as a simplicial set (or a topological space) instead of the underlying set of $S$ in the ordinary commutative ring theory.

DEFINITION 2.1.13 : *A simplicial commutative ring $R$ in $s\mathcal{CA}\mathrm{lg}_{\Bbbk}$ is said to be* discrete *if $\pi_i R = 0$ for all $i > 0$.*

REMARK 2.1.14 : A simplicial commutative ring has an underlying simplicial set, that is a topological space, so it makes sense to talk about its higher homotopy groups.

REMARK 2.1.15 : We have the following adjunction:

$$s\mathcal{CA}\mathrm{lg}_{\Bbbk} \underset{\pi_0}{\overset{\text{ful. faith.}}{\underset{\longrightarrow}{\overset{\hookleftarrow}{\top}}}} s\mathcal{CA}\mathrm{lg}_{\Bbbk}^{\text{discrete}} \simeq s\mathcal{CA}\mathrm{lg}_{\pi_0 \Bbbk} \tag{2.1.16}$$

So the classical theory of commutative algebras, that is exactly $s\mathcal{CA}\mathrm{lg}_{\pi_0 \Bbbk}$ for some $\Bbbk$, is fully faithfully included in the theory of simplicial commutative algebras.

DEFINITION 2.1.17 : *A simplicial commutative $\Bbbk$-algebra $R$ in $s\mathcal{CA}\mathrm{lg}_{\Bbbk}$ is said to be* of locally finite presentation (l.f.p.) *if it is a compact object in $s\mathcal{CA}\mathrm{lg}_{\Bbbk}$, that is the functor:*

$$\mathrm{Map}_{s\mathcal{CA}\mathrm{lg}_{\Bbbk}}(R, -) : s\mathcal{CA}\mathrm{lg}_{\Bbbk} \longrightarrow \mathcal{S}\mathrm{paces} =: \mathcal{S}$$

*commutes with filtered colimits.*

### 2.1.2 The cotangent complex.

#### 2.1.2.1 Introduction:Generalized derived manifold.

We follow a presentation of generalized derived manifolds made by J. Lurie in 2017. Let $\mathscr{X}$ be a smooth manifold (in the differentiable sense). It is equipped with a canonical sheaf $\mathscr{O}_{\mathscr{X}}$ of classical rings which verifies a certain bunch of properties ($\mathscr{O}_{\mathscr{X},x}$ is a local ring, etc.).

Let $\mathscr{X}_0, \mathscr{X}_1$ be two other smooth manifolds (always in the differentiable sense). Consider

the pullback square:

$$\begin{array}{ccc} \mathscr{X}_1 \times_{\mathscr{X}} \mathscr{X}_0 & \to & \mathscr{X}_0 \\ \downarrow & \lrcorner & \downarrow f \\ \mathscr{X}_1 & \xrightarrow[g]{} & \mathscr{X} \end{array} \tag{2.1.18}$$

where $f$ and $g$ are smooth maps and $f$ is a submersion. Then $\mathscr{X}_1 \times_{\mathscr{X}} \mathscr{X}_0$ is a smooth manifold.

We say that the category of smooth manifolds with smooth maps has pullbacks.

Consider now that the manifold $\mathscr{X}$ is no more smooth. Then the pullback, if it exists (at least in the category $\mathcal{T}$op of topological spaces) is no more in the category of smooth manifolds. In which category does this pullback can happen ? We have to construct it.

DEFINITION 2.1.19 (Generalized derived manifold) : *Let us consider the functor:*

$$\mathscr{O}_{\mathscr{X}} : \left\{ \begin{array}{ccc} \{\text{open subsets} \mathscr{U} \text{ of} \mathscr{X}\}^{\text{op}} \times \{\text{smooth manifolds}\} & \longrightarrow & \boldsymbol{\mathcal{S}}\text{paces} \\ (\mathscr{U}, \mathcal{M}) & \longmapsto & \mathscr{O}_{\mathscr{X}}^{\mathcal{M}}(\mathscr{U}) \end{array} \right\} \tag{2.1.20}$$

*veryfying the three following axioms:*

- *($\mathscr{O}_{\mathscr{X}}^{\mathcal{M}}$ is a sheaf) If $\mathcal{M}$ is fixed, the map $\{\mathscr{U} \longmapsto \mathscr{O}_{\mathscr{X}}^{\mathcal{M}}(\mathscr{U})\}$ is a homotopoy sheaf, that is if $\mathscr{U}$ is covered by open sets $\mathscr{U}_\alpha$, then we have an isomorphism:*

$$\mathscr{O}_{\mathscr{X}}^{\mathcal{M}}(\mathscr{U}) \xrightarrow{\sim} (\text{ho}) \lim_{\mathscr{V} \subseteq \mathscr{U}} . \mathscr{O}_{\mathscr{X}}^{\mathcal{M}}(\mathscr{V}) \tag{2.1.21}$$

- *($\mathscr{O}_{\mathscr{X}}$ is local) If $\mathcal{M}$ has an open cover by $\mathcal{M}_\alpha \subseteq \mathcal{M}$, then we have a surjection of sheaves:*

$$\coprod_{\alpha} \mathscr{O}_{\mathscr{X}}^{\mathcal{M}_\alpha} \xrightarrow{\text{surj. of sheaves}} \mathscr{O}_{\mathscr{X}}^{\mathcal{M}} \tag{2.1.22}$$

- *(pullbacks) If*

$$\begin{array}{ccc} \mathcal{M}_1 \times_{\mathcal{M}} \mathcal{M}_0 & \longrightarrow & \mathcal{M}_0 \\ \downarrow & \lrcorner & \downarrow f \\ \mathcal{M}_1 & \longrightarrow & \mathcal{M} \end{array} \tag{2.1.23}$$

*where $f$ is a submersion, is a pullback square, then*

$$\begin{array}{ccc} \mathscr{O}_{\mathscr{X}}^{\mathcal{M}_0 \times_{\mathcal{M}} \mathcal{M}_1}(\mathscr{U}) & \longrightarrow & \mathscr{O}_{\mathscr{X}}^{\mathcal{M}_0}(\mathscr{U}) \\ \downarrow & \lrcorner & \downarrow f \\ \mathscr{O}_{\mathscr{X}}^{\mathcal{M}_1}(\mathscr{U}) & \longrightarrow & \mathscr{O}_{\mathscr{X}}^{\mathcal{M}}(\mathscr{U}) \end{array} \tag{2.1.24}$$

*is a homotopy pullback square*

*We say that $(\mathscr{X}, \mathscr{O}_{\mathscr{X}})$ is a* generalized derived manifold*.*

DEFINITION 2.1.25 : *Let $(\mathscr{X}, \mathscr{O}_{\mathscr{X}})$ and $(\mathscr{Y}, \mathscr{O}_{\mathscr{Y}})$ be two such generalized derived manifolds. Let $f_0 : X \longrightarrow Y$ be a continuous function, and $f_0^* \mathscr{O}_{\mathscr{Y}} \longrightarrow \mathscr{O}_{\mathscr{X}}$ be a map of sheaves, then we have an analogue of loacl maps: If $\mathcal{N} \hookrightarrow \mathcal{M}$ is an immersion of smooth*

*manifolds, then:*

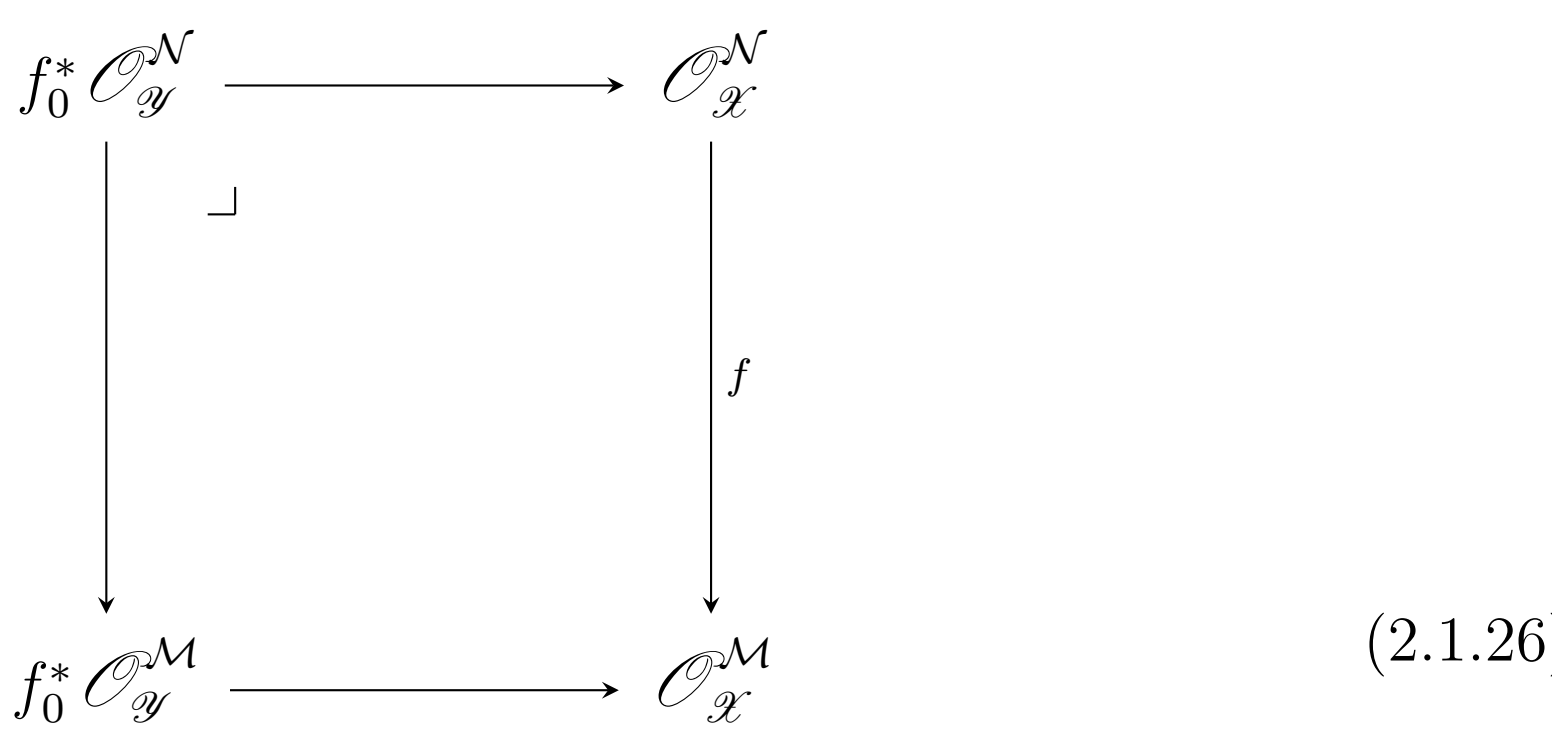

$$\begin{array}{ccc} f_0^* \mathscr{O}_{\mathscr{Y}}^{\mathcal{N}} & \longrightarrow & \mathscr{O}_{\mathscr{X}}^{\mathcal{N}} \\ \downarrow & & \downarrow f \\ f_0^* \mathscr{O}_{\mathscr{Y}}^{\mathcal{M}} & \longrightarrow & \mathscr{O}_{\mathscr{X}}^{\mathcal{M}} \end{array} \tag{2.1.26}$$

*is a homotopy pullback square This defines the morphism between two generalized derived manifold, making a category which is an* $\infty$*-category.*

THEOREM 2.1.27 : *The category of generalized derived maniford has pullbacks, and the category of smooth manifolds is fully faithfully imbedded in it.*

Of course, in general, even if $(\mathscr{X}, \mathscr{O}_{\mathscr{X}})$, $(\mathscr{Y}, \mathscr{O}_{\mathscr{Y}})$ and $(\mathscr{Z}, \mathscr{O}_{\mathscr{Z}})$ are smooth, $(\mathscr{X} \times_{\mathscr{Z}} \mathscr{Y}, \mathscr{O}_{\mathscr{X}} \times_{\mathscr{O}_{\mathscr{Z}}} \mathscr{O}_{\mathscr{Y}})$ is not a smooth manifold (else if there is a submersion between two of the smooth manifolds). The underlying space can ever been particularly "*awful*". However, we can right now see on an example, the gain of this operation:

EXAMPLE 2.1.28 : Consider the inclusion $\{0\} \xhookrightarrow{i} \mathbb{R}$ and let the following pullback

diagram be:

$$\begin{array}{ccc} (\mathscr{X}, \mathscr{O}_{\mathscr{X}}) & \xrightarrow{i} & \{0\} \\ \downarrow i & \searrow h & \downarrow \\ \{0\} & \longrightarrow & \mathbb{R} \end{array} \tag{2.1.29}$$

Then $\mid X \mid=\mid \{0\} \mid \times_{\mathbb{R}} \mid \{0\} \mid=\mid \{0\} \mid= \{0\}$, set -theorically. Moreover,

$$\begin{array}{ccc} h^* \mathscr{O}^{\mathbb{R}}_{\mathbb{R}}) & \longrightarrow & \mathscr{O}_{\{0\}} \\ \downarrow & & \downarrow \\ \mathscr{O}^{i}_{\{0\}} & \longrightarrow & \mathscr{O}^{\mathbb{R}}_{\mathbb{X}} \end{array} \tag{2.1.30}$$

that is:

$$\begin{array}{ccc} A & \longrightarrow & \mathbb{R} \\ \downarrow & & \downarrow \\ \mathbb{R} & \longrightarrow & \mathscr{O}^{\mathbb{R}}_{\mathbb{X}} \end{array} \tag{2.1.31}$$

where $A := \{$ Rings of germs of smooth maps $\mathbb{R} \to \mathbb{R}$ at $0$ $\}$. That is:

$$\mathscr{O}^{\mathbb{R}}_{\mathbb{X}} = \mathbb{R} \otimes^{\mathbb{L}}_{A} \mathbb{R} \tag{2.1.32}$$

We have a short exact sequence:

$$0 \longrightarrow A \longrightarrow A \longrightarrow \mathbb{R} \longrightarrow 0 \tag{2.1.33}$$

that is:

$$0 \longrightarrow A \longrightarrow A \longrightarrow 0 \tag{2.1.34}$$

is quasi-sisomorphic to $\mathbb{R}$, so:

$$0 \longrightarrow \mathbb{R} \longrightarrow \mathbb{R} \longrightarrow 0 \tag{2.1.35}$$

is a resolution of $\mathscr{O}_{\mathscr{X}}^{\mathbb{R}} = \mathbb{R} \otimes_A^{\mathbb{L}} \mathbb{R}$.

We have then obtained that the pullback has two not null degrees, in degree 0: $\mathbb{R}$, but also in degree 1: $\mathbb{R}$, and the spectrum of the coconnective DGAlgebra $\mathbb{R} \oplus \mathbb{R}$ is then the connective generalized derived manifold:

$$(\mathscr{X}, \mathscr{O}_{\mathscr{X}}) = (\mathsf{Spec}\,(\mathbb{R} \oplus \mathbb{R}), \mathbb{R} \otimes_A^{\mathbb{L}} \mathbb{R}) \tag{2.1.36}$$

DEFINITION 2.1.37 : *A generalized derived manifold is said to be* quasi-smooth *if locally on* $\mathscr{X}$*, there is a pullback square:*

$$\begin{array}{ccc} (\mathscr{X}, \mathscr{O}_{\mathscr{X}}) & \longrightarrow & (\mathscr{X}_0, \mathscr{O}_{\mathscr{X}_0}) \\ \downarrow & \lrcorner & \downarrow \\ (\mathscr{X}_1, \mathscr{O}_{\mathscr{X}_1}) & \longrightarrow & (\mathscr{X}, \mathscr{O}_{\mathscr{X}}) \end{array} \tag{2.1.38}$$

Let us now look at the tangent space of a generalized derived manifold:

If $(\mathscr{X}, \mathscr{O}_{\mathscr{X}})$ is a generalized derived manifold, let us look at:

$$\boldsymbol{\mathcal{S}}\mathbf{hv}_{\boldsymbol{\mathcal{S}}}(\mathscr{X})_{/\mathscr{O}_{\mathscr{X}}} := \left\{ \begin{array}{l} \text{Sheaves } \mathscr{O}_{\mathscr{X}}' \text{ satisfying the axioms of generalized derived manifolds} \\ \text{with a local map } \mathscr{O}_{\mathscr{X}}' \longrightarrow \mathscr{O}_{\mathscr{X}} \end{array} \right\} \tag{2.1.39}$$

and the algebraic version, which is a 1-approximation, that is in fact a sort of tangent space:

$$\boldsymbol{\mathcal{S}}\mathbf{hv}_{\boldsymbol{\mathcal{S}}}(\mathscr{X})^{\mathrm{alg}}_{/\mathscr{O}_{\mathscr{X}}} := \left\{ \begin{array}{l} \text{Sheaves of simplicial commutative } \mathbb{R}\text{-algebras} \\ \text{with a local map } A \longrightarrow \mathscr{O}^{\mathbb{R}}_{\mathscr{X}} \end{array} \right\} \tag{2.1.40}$$

Let us take the spectrum of each of these slice $\infty$-categories:

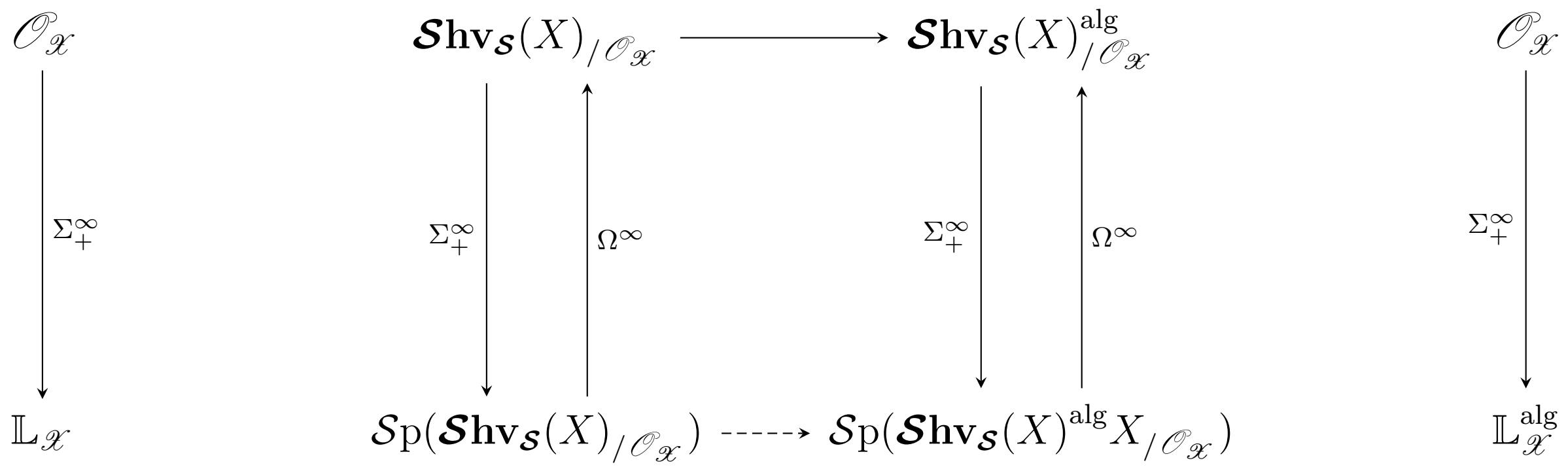

$$\tag{2.1.41}$$

where $\mathcal{S}\mathrm{p}(\boldsymbol{\mathcal{S}}\mathbf{hv}_{\boldsymbol{\mathcal{S}}}(\mathscr{X})^{\mathrm{alg}}_{/\mathscr{O}_{\mathscr{X}}}) =$ Sheaves of $\mathscr{O}^{\mathbb{R}}_{\mathscr{X}}$-modules

We recognize on the left and on the right the future cotangent complex. The one on the left $\mathbb{L}_{\mathscr{X}}$ is the spectrum of the stable $\infty$-category $\boldsymbol{\mathcal{S}}\mathbf{hv}_{\boldsymbol{\mathcal{S}}}(X)_{/\mathscr{O}_{\mathscr{X}}}$. The one on the right $\mathbb{L}^{\mathrm{alg}}_{\mathscr{X}}$ is the spectrum of the stable $\infty$-category $\boldsymbol{\mathcal{S}}\mathbf{hv}_{\boldsymbol{\mathcal{S}}}(X)^{\mathrm{alg}}_{/\mathscr{O}_{\mathscr{X}}}$.

THEOREM 2.1.42 : *The bottom dashed arrow is an isomorphism, allowing to speak of the* cotangent complex *identical to the first-order algebraic approximation, the* algebraic cotangent complex $\mathbb{L}_{\mathscr{X}}$*, which we will identify as* $\mathbb{L}_{\mathscr{X}}$*.*

REMARK 2.1.43 : If $\mathscr{X}$ is smooth, $\mathcal{S}\mathrm{p}(\boldsymbol{\mathcal{S}}\mathbf{hv}_{\boldsymbol{\mathcal{S}}}(\mathscr{X})^{\mathrm{alg}}_{/\mathscr{O}_{\mathscr{X}}}) =$ Sheaves of $\mathscr{O}^{\mathbb{R}}_{\mathscr{X}}$-modules would be the derived category of the category of $\mathbb{R}$-modules.

After having given a motiveation for generalized derived manifold, and a first smell of what the cotangent complex would be, let us now be more precise in the construction of this cotangent complex.

#### 2.1.2.2 The first definition of the cotangent complex.

Let us recall the following adjunction giving a projective resolution of a commutative algebra in the category of simplicial commutative algebras:

$$\text{(ho) } s\mathcal{S}\text{ets} \underset{\mathcal{G}:=\Bbbk[-]_{\text{free}}}{\overset{\mathcal{F}:=\mathscr{F}\text{orget}}{\rightleftarrows}} \text{(ho) } \mathcal{CA}\text{lg}_{\Bbbk} \qquad (\top)$$

where $\mathcal{G}(S) := \Bbbk_{\text{free}}[S]$, the free commutative algebra on $S$. This adjunction gives us unit and counit natural transformations: $\mathcal{F} \circ \mathcal{G} \to \mathbf{1}$ and $\mathbf{1} \to \mathcal{F} \circ \mathcal{G}$, which allow us to construct the following resolution $S_{\bullet} \longrightarrow R$ for each $R \in \mathcal{CA}\text{lg}_{\Bbbk}$ such that:

$$\dots\ \Bbbk_{\text{free}}[\Bbbk_{\text{free}}[\Bbbk_{\text{free}}[S]]] \begin{array}{c} \xleftarrow{f_0} \\ \xrightarrow{d_0} \\ \xleftarrow{f_1} \\ \xrightarrow{d_1} \\ \xleftarrow{f_2} \end{array} \Bbbk_{\text{free}}[\Bbbk_{\text{free}}[S]] \begin{array}{c} \xleftarrow{f_0} \\ \xrightarrow{d_0} \\ \xleftarrow{f_1} \end{array} \Bbbk_{\text{free}}[S]$$

whose homology gives a homotopy equivalence : $S_{\bullet} \xrightarrow[\text{htpy}]{\sim} R$, where $S_{\bullet}$ is a polynomial ring in each degree.

DEFINITION 2.1.44 : *Let $\Bbbk \to R$ be two commutative rings in $\mathcal{CA}\text{lg}_{\Bbbk}$ such that the following commutative diagram is a resolution of $R$ by $S_{\bullet}$ :*

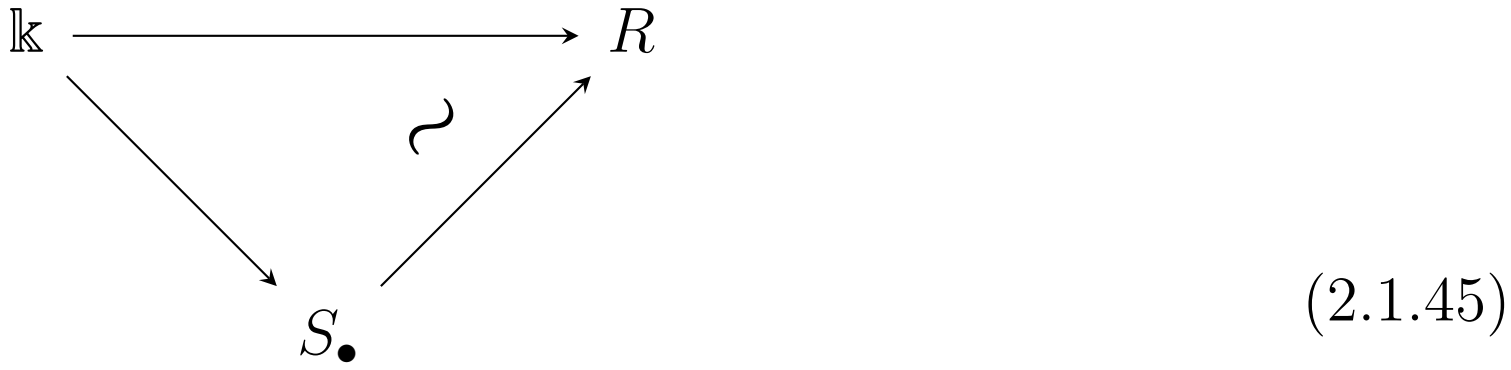

$$(2.1.45)$$

*Let us define the* cotangent complex $L_{R/\Bbbk}$ *as the simplicial* $R$*-module in* $s\mathsf{Mod}_R$ :

$$L_{R/\Bbbk} := \Omega^1_{S_\bullet/\Bbbk} \otimes_{S_\bullet} R$$

*where* $\Omega^1_{S_\bullet/\Bbbk}$ *is the Kähler* 1*-differentials degreewise of* $S_\bullet$.

REMARK 2.1.46 : The cotangent complex $L_{R/\Bbbk}$ is a simplicial $R$-module in $s\mathsf{Mod}_R$ by base change over the simplicial $R$-module $S_\bullet$ , and then an element of the connective derived category $\mathscr{D}(R)_{\geqslant 0}$ of $R$ by Dold-Kan correspondence.

REMARK 2.1.47 : This first definition does not exhibit the fact that the definition of the cotangent complex does not depend on the resolution $S_\bullet \longrightarrow R$. The second definition does.

#### 2.1.2.3 The second definition of the cotangent complex.

Let us consider the definition of the derived functor of the functor $\Omega^1_{-/\Bbbk}$ seen in example 2.1.6, where the double arrow means the Kan functor extension:

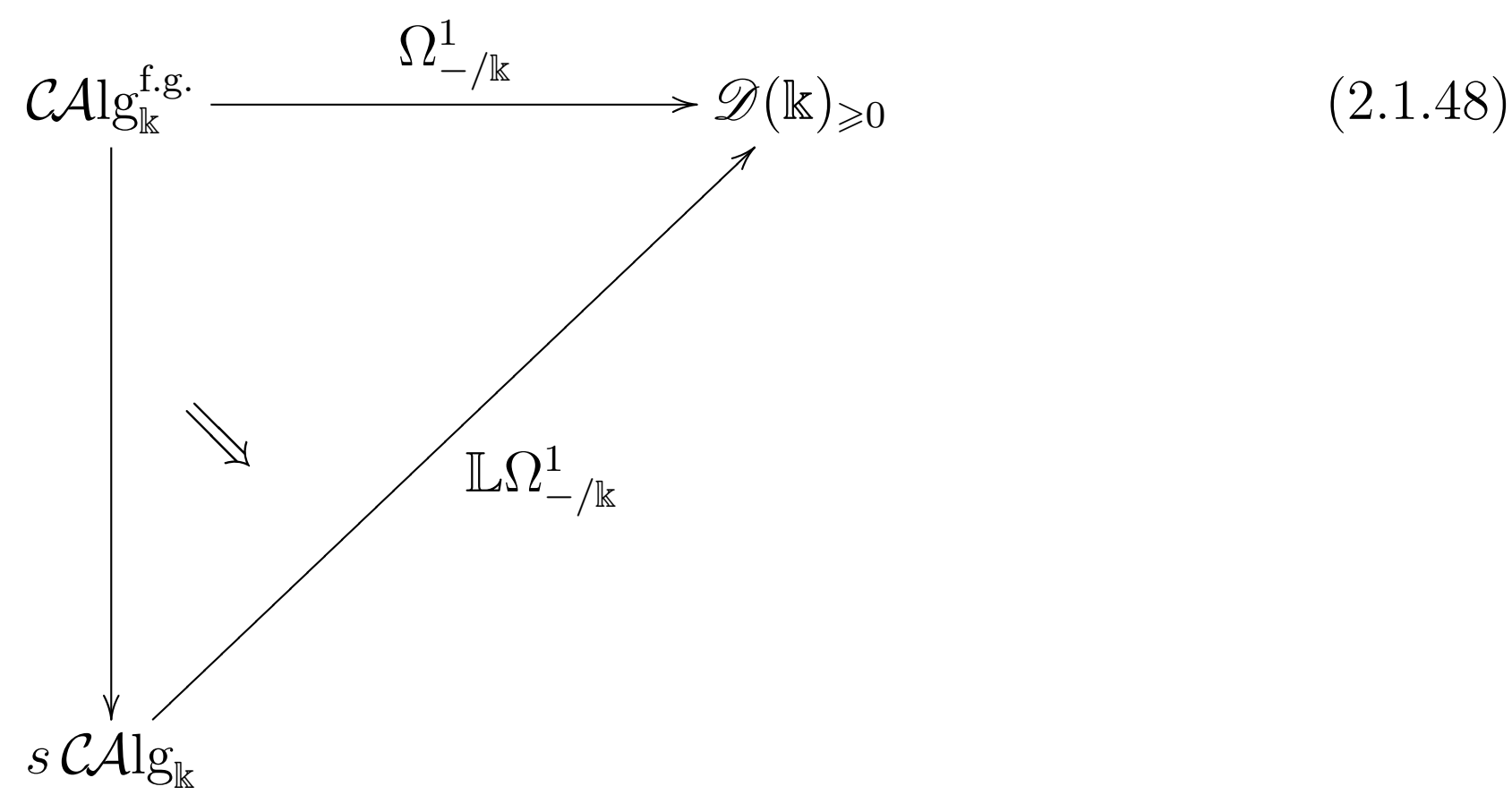

(2.1.48)

DEFINITION 2.1.49 : *The* cotangent complex $L_{R/\Bbbk}$ *is the simplicial* $R$*-module in* $s\mathsf{Mod}_R$ *defined such that:*

$$L_{R/\Bbbk} := \mathbb{L}\Omega^1_{R/\Bbbk}$$

REMARK 2.1.50 : This second definition exhibits the fact that the definition of the cotangent complex does not depend on the resolution $S_\bullet \longrightarrow R$. On the other hand, it is not immediate that the functor $\Omega^1_{-/\Bbbk}$ sends $\mathcal{CA}\mathrm{lg}_{\Bbbk}^{\mathrm{f.g.}}$ onto $\mathscr{D}(\Bbbk)_{\geqslant 0}$ but rather onto $\mathscr{D}(R)_{\geqslant 0}$. But it is the case. So the best way to define the cotangent complex is by the third definition, which involves a universal property.

#### 2.1.2.4 The third definition of the cotangent complex.

Let us consider a chain of maps: $\Bbbk \to R \to S$ in the category of simplicial commutative rings. Let us denote $\mathrm{s}\mathcal{CA}\mathrm{lg}_{\Bbbk}//S$ the comma category of simplicial commutative rings equipped with a fixed map to $S$. Then we have the following:

DEFINITION 2.1.51 : *The* cotangent complex *is the universal object $L_{R/\Bbbk}$ that gives, for each $M$ an $S$-module in $\mathscr{D}(S)_{\geqslant 0}$, an isomorphism between the following mapping topological spaces:*

$$\mathrm{Hom}_R(L_{R/\Bbbk}, M) \xrightarrow{\sim} \mathrm{Hom}_{\mathrm{s}\mathcal{CA}\mathrm{lg}_{\Bbbk}//S}(R, S \oplus M)$$

*where $S \oplus M$ is seen as $S \oplus M \longrightarrow S$ in $s\mathcal{CA}\mathrm{lg}_{\Bbbk}//S$.*

#### 2.1.2.5 Properties of the cotangent complex.

We have the following properties:

PROPOSITION 2.1.52 :

- *We have an isomorphism: $\pi_0 L_{R/\Bbbk} \simeq \Omega^1_{R/\Bbbk}$*
- *(transitivity) The following chain of maps:*

$$S \otimes^{\mathbb{L}}_R L_{R/\Bbbk} \to L_{S/\Bbbk} \to L_{S/R}$$

*is an exact triangle in the derived category $\mathscr{D}(\Bbbk)_{\geqslant 0}$*

- *If we have the pushout diagram:*

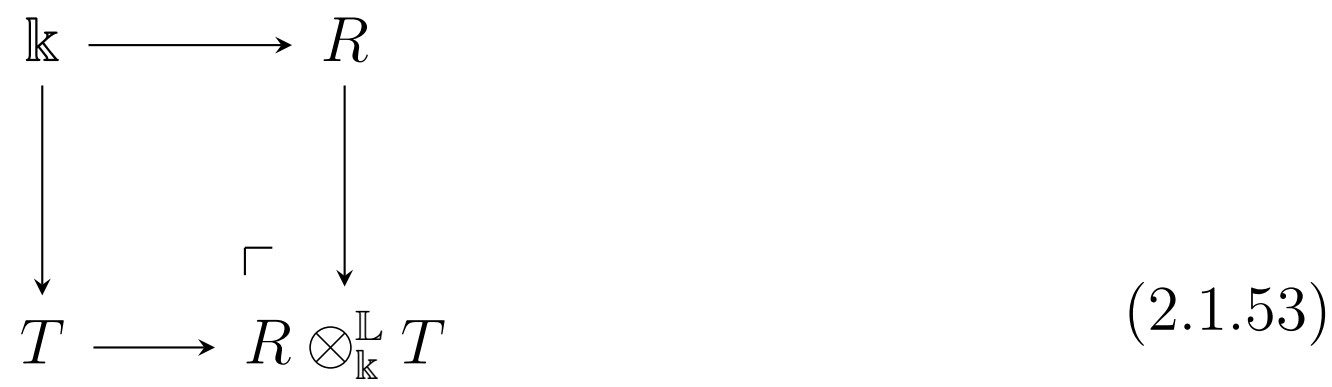

$$(2.1.53)$$

*then:*

$$T \otimes_{\Bbbk} L_{R/\Bbbk} \xrightarrow{\sim} L_{R \otimes^{\mathbb{L}}_{\Bbbk} T\,/\,T}$$

#### 2.1.2.6 Derived De Rham cohomology.

Let us consider the following De Rham functor denoted $\mathrm{dRh}\Omega_{-/\Bbbk}$ defined such that:

$$\begin{array}{ccl} \mathcal{CA}\mathrm{lg}^{\mathrm{f.g.}}_{\Bbbk} & \longrightarrow & \mathscr{D}(\Bbbk)^{\mathrm{b}}_{\geqslant 0} \\ R & \longmapsto & \mathrm{dRh}\Omega_{R/\Bbbk} := (R \to \Omega^1_{R/\Bbbk} \to \Omega^2_{R/\Bbbk} \to \dots) \end{array}$$

where $\Omega^i_{R/\Bbbk} = \bigwedge^i \Omega^1_{R/\Bbbk}$. Since $R$ is finitely generated, the De Rham complex is bounded. Let us consider the derived functor of the De Rham functor denoted $\mathbb{L}\mathrm{dRh}\Omega_{-/\Bbbk}$ :

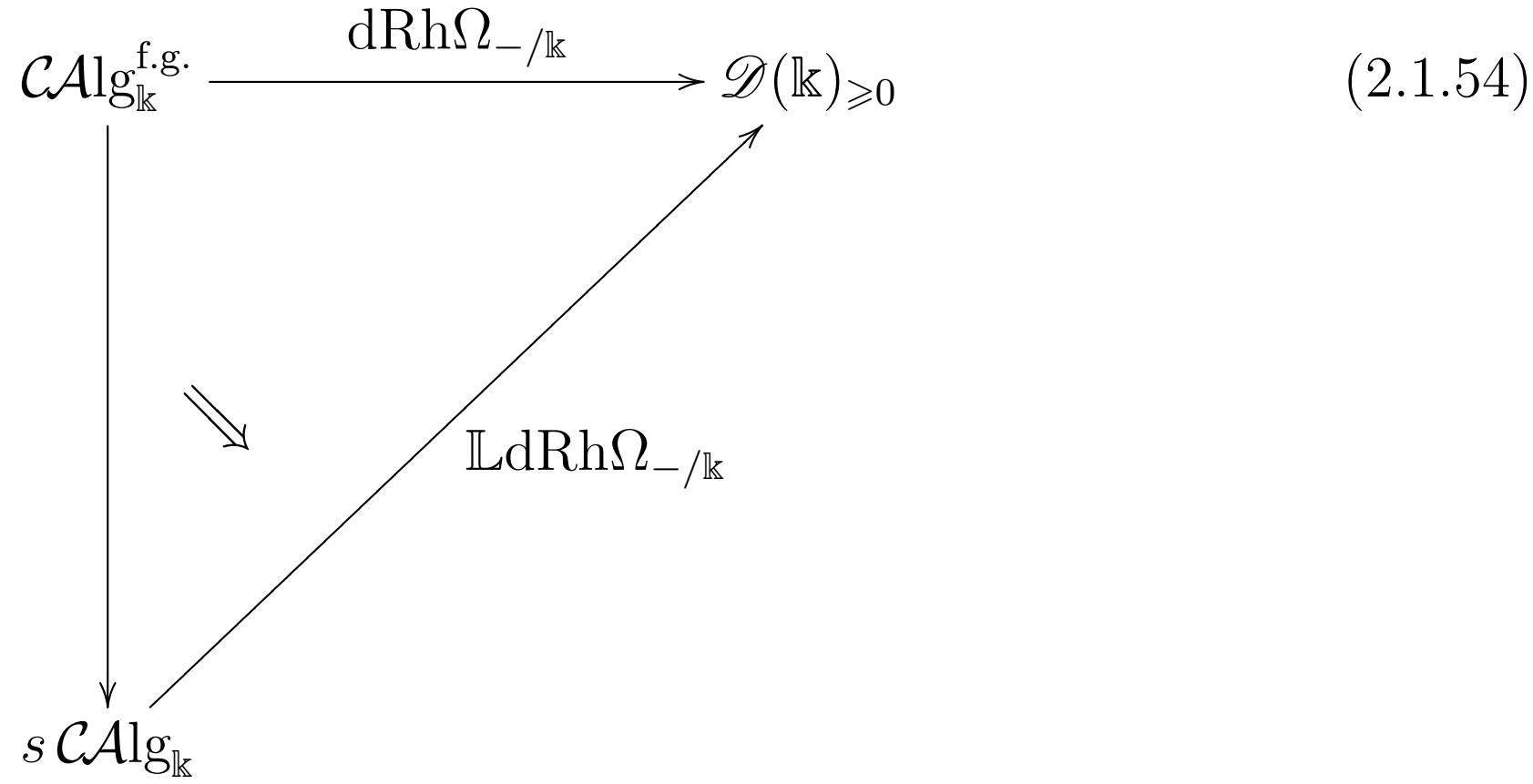

$$(2.1.54)$$

We know that De Rham cohomology behaves well for smooth varieties. So we have the following properties:

PROPOSITION 2.1.55 :

- *The De Rham cohomology of an affine space is trivial, that is:* $\mathrm{dRh}\Omega_{R/\Bbbk} \simeq k$
- *The derived De Rham cohomology of a simplicial ring $S$ is also trivial, that is:* $\mathbb{L}\mathrm{dRh}\Omega_{S/\Bbbk} \simeq k$

But looking at the Hodge canonical filtration of the De Rham complex in each degree $i$:

$$\mathrm{gr}^i_{\mathrm{H}}\, \mathrm{dRh}\Omega_{R/\Bbbk} \simeq \Omega^i_{R/\Bbbk}[-i]$$

we notice that these gradings are compatible with colimits. Thus the definition by a Kan extension of the derived De Rham cohomology allows us to remember the (Hodge canonical) filtration, and we write the following:

$$\mathrm{gr}^i_{\mathrm{H}}\, \mathbb{L}\mathrm{dRh}\Omega_{-/\Bbbk} \xrightarrow{\sim} \mathbb{L} \bigwedge^i L_{R/\Bbbk}\, [-i]$$

where $L_{R/\Bbbk}$ is the cotangent complex seen in $\mathscr{D}(\Bbbk)_{\geqslant 0}$. We can then state:

DEFINITION 2.1.56 : *The total completion of the derived De Rham cohomology functor $L_{R/\Bbbk}$ is the simplicial $R$-module in $\mathsf{sMod}_R$ defined such that:*

$$\widehat{\mathbb{L}\mathrm{dRh}\Omega_{-/\Bbbk}} := \operatorname*{Lim}_i \Big( \mathbb{L}\mathrm{dRh}\Omega_{-/\Bbbk} \;/\; \mathrm{gr}^i_{\mathrm{H}}\, \mathbb{L}\, \mathrm{dRh}\Omega_{-/\Bbbk} \Big)$$

#### 2.1.2.7 Bhatt-Grothendieck-Hartshorne theorem.

Let $X_{/\mathbb{C}}$ be a (not necessarily smooth) scheme of finite type over $\mathbb{C}$. Then:

THEOREM 2.1.57 ([23]) : *We have an isomorphism:*

$$\mathbb{R}\Gamma\Big( X\,,\, \widehat{\mathbb{L}\mathrm{dRh}\Omega_{\mathscr{O}_X/\Bbbk}} \Big) \xrightarrow{\sim} \mathbb{R}\Gamma_{\mathrm{sing.}}\Big( X(\mathbb{C})\,,\, \mathbb{C} \Big) \tag{2.1.58}$$

REMARK 2.1.59 : This theorem allows us to embed a singular scheme $X$ in a larger ambient smooth scheme.

## 2.2.Introduction to Spectral Algebraic Geometry.

### 2.2.1 Open and closed immersions in $\infty$-topoï.

Before diving into the gluing of closed immersions—which is the central focus of the problems to come, since our subject is the pathological topology of the Grassmannian and the gluing of its positroid varieties (open or closed), their gluings (smooth or singular), etc., which will justify the shift from Zariski topology to étale topology, followed by the transition to intersection sheaves and perverse geometry, let us begin to establish the definitions of topology in $\infty$-topoï we will use, especially since it is never certain in the case of $\infty$-topoï that we are indeed dealing with underlying topological spaces, or else with their connective or coconnective homotopy, in duality with their coconnective or connective dual differential graded algebras, in positive or negative cohomology degrees.

We will have to return to these topics regularly, since the goal is indeed to arrive at the end of this text with a description of the Grassmannian and its positroid varieties as a well-defined, concrete underlying topological space, which will moreover be differentially cohesive, etc. In short, it will have all the expected topological and differential properties, but since we are proceeding via spectral Deligne-Mumford $\infty$-stacks, it is not superfluous to begin getting accustomed to our approach to topology, and in particular, at first glance: closed immersions.

First, we recall the remark of [90] (between Remark 4.26 and Definition 4.27) which I quote wihout modification, since it is perfectly adapted to our discussion between the underlying topologies and the morphisms of schemes, stacks or of ($\infty$-)topos : " *If $X$ is a topological space, then there is a bijection from the set of closed subsets of $X$ to the set of open subsets of $X$, given by $Y \longmapsto X - Y$. In algebraic geometry, the situation is more subtle: every closed subscheme $Y \subseteq X$ has an open complement $U = X - Y$. However, this construction is not bijective: a closed subset of $X$ generally admits many different scheme structures. However, we can recover a bijective correspondence by restricting our attention to reduced closed subschemes of $X$.*"

However, we can be more precise like this, and J. Lurie is indeed cautious in many loci of his work. Particularly I will quote the following ones:

DEFINITION 2.2.1 ([98] Def 7.3.2.6) : *If* $\mathfrak{X}$ *is an* $\infty$*-topos and* $\mathfrak{U} \in \mathrm{Sub}(\mathbf{1}_{\mathfrak{X}})$ *where* $\mathrm{Sub}(\mathbf{1}_{\mathfrak{X}})$ *denote the partially ordered set of equivalence classes of* $(-1)$*-truncated objects of* $\mathfrak{X}$*, which is independent of the choice of a final object* $\mathbf{1}_{\mathfrak{X}} \in \mathfrak{X}$*, up to isomorphism, then we refer to* $\mathfrak{X}/\mathfrak{U}$ *as the* closed subtopos of $\mathfrak{X}$ complementary to $\mathfrak{U}$*. More generally, we will say that a geometric morphism* $\pi : \mathfrak{Y} \longrightarrow \mathfrak{X}$ *is a* closed immersion *if there exists* $\mathfrak{U} \in \mathrm{Sub}(\mathbf{1}_{\mathfrak{X}})$ *such that* $\pi_*$ *induces an equivalence of* $\infty$*-categories from* $\mathfrak{Y}$ *to* $\mathfrak{X}/\mathfrak{U}$*.*

PROPOSITION 2.2.2 ([98] Prop 7.3.2.7) : *Let* $\mathfrak{X}$ *be an* $\infty$*-topos, and let* $\mathfrak{U} \in \mathrm{Sub}(\mathbf{1}_{\mathfrak{X}})$*. Then the closed immersion* $\pi : \mathfrak{X}/\mathfrak{U} \longrightarrow \mathfrak{X}$ *induces an isomorphism of partially ordered sets from* $\mathrm{Sub}(\mathbf{1}_{\mathfrak{X}})$ *to* $\{\mathfrak{V} \in \mathrm{Sub}(\mathbf{1}_{\mathfrak{X}}) \mid \mathfrak{U} \subseteq \mathfrak{V}\}$*.*

Another way to define a closed immersion is the following:

DEFINITION 2.2.3 ([103] Intro of §2.5.2) : *Let* $\mathfrak{X}$ *be an* $\infty$*-topos, and let* $\mathfrak{U} \in \mathfrak{X}$ *be a* $(-1)$*-truncated object. We let* $\mathfrak{X}/\mathfrak{U}$ *denote the full subcategory of* $\mathfrak{X}$ *spanned by those objects* $\mathscr{X} \in \mathfrak{X}$ *for which the projection map* $\mathscr{X} \times \mathfrak{U} \longrightarrow \mathfrak{U}$ *is an equivalence. Then* $\mathfrak{X}/\mathfrak{U}$ *is itself an* $\infty$*-topos, and the inclusion functor* $i_* : \mathfrak{X}/\mathfrak{U} \longrightarrow \mathfrak{X}$ *is a geometric morphism of* $\infty$*-topoï. A geometric morphism* $f_* : \mathfrak{Y} \longrightarrow \mathfrak{X}$ *is called a* closed immersion *if it factors as a composition:*

$$\mathfrak{Y} \xrightarrow{g_*} (\mathfrak{X}/\mathfrak{U}) \xrightarrow{i_*} \mathfrak{X} \tag{2.2.4}$$

*where* $\mathfrak{U}$ *is a* $(-1)$*-truncated object of* $\mathfrak{X}$ *and* $g_*$ *is an equivalence.*

And this is the good way to define an open immersion, even if we do not have yet define what is a nonconnective spectral Deligne-Mumford stack. That is why we will come back later on this definition, but let us state it right away, just to complete the picture we begin to draw:

DEFINITION 2.2.5 ([103] Def 1.6.7.2) : *Suppose that* $j : \mathsf{U} \longrightarrow \mathsf{X} = (\mathfrak{X}, \mathscr{O}_{\mathfrak{X}})$ *is a map of nonconnective spectral Deligne-Mumford stacks. We will say that* $j$ *is* open immersion *if it factors as a composition:*

$$\mathsf{U} \xrightarrow{j'} (\mathfrak{X}_{/\mathfrak{U}}, \mathscr{O}_{\mathfrak{X}_{/\mathfrak{U}}}) \xrightarrow{j''} (\mathfrak{X}, \mathscr{O}_{\mathfrak{X}}) \tag{2.2.6}$$

*where $j'$ is an equivalence and $j''$ is the étale morphism associated to a $(-1)$-truncated object $\mathfrak{U} \in \mathfrak{X}$. In this case, we will also say that $\mathsf{U}$ is an* open substack *of $\mathsf{X}$.*

### 2.2.2 Zariski schemes and Deligne-Mumford stacks.

DEFINITION 2.2.7 : *A map $R \to S$ between two (ordinary) commutative rings is said to be an* étale map *if the ring $S$, seen as an $R$-module, is finitely presented and flat, and if it is also flat if it is seen as an $S \otimes_R S$-module.*

DEFINITION 2.2.8 : *Following [103] and stable homotopy theory, the category of ordinary commutative rings will be denoted $\mathcal{CA}\mathrm{lg}^{\heartsuit}_{\mathbb{Z}}$, and if we only consider the étale morphisms, this defines the subcategory $\mathcal{CA}\mathrm{lg}^{\heartsuit,\text{ét}}_{\mathbb{Z}}$.*

DEFINITION 2.2.9 :

- *A simplicial commutative ring $\Bbbk \to R$ in $s\mathcal{CA}\mathrm{lg}_{\Bbbk}$ is said to be* formally étale *if the cotangent complex $L_{R/\Bbbk} \simeq 0$.*
- *A simplicial commutative ring $\Bbbk \to R$ in $s\mathcal{CA}\mathrm{lg}_{\Bbbk}$ is said to be* étale *if it is formally étale and of locally finite presentation (see definition 2.1.17).*

PROPOSITION 2.2.10 : *The following properties are equivalent:*

- *The map $\Bbbk \to R$ is étale*
- *The induced map: $\pi_0\Bbbk \to \pi_0 R$ between ordinary commutative rings is étale and the induced map $\pi_i\Bbbk \otimes_{\pi_0\Bbbk} \pi_0 R \xrightarrow{\sim} \pi_i R$ is an isomorphism for all $i > 0$.*

THEOREM 2.2.11 ([103]) : *The following $\infty$-categories are equivalent:*

$$s\mathcal{CA}\mathrm{lg}^{\text{ét}}_{\Bbbk} \xrightarrow{\sim} s\mathcal{CA}\mathrm{lg}^{\text{ét}}_{\pi_0\Bbbk}$$

*where $s\mathcal{CA}\mathrm{lg}^{\text{ét}}_{\Bbbk}$ denotes the subcategory of $s\mathcal{CA}\mathrm{lg}_{\Bbbk}$ only considering the étale maps.*

remark 2.2.12 : This is a deep theorem saying that the small étale site of the simplicial commutative ring $\Bbbk$ is equivalent to the small étale site of the ordinary commutative ring $\pi_0\Bbbk$. This comes from the fact that a simplicial commutative ring $\Bbbk$ is the (inverse) limit of its truncations $\tau_{\leqslant n}\Bbbk$ where $\pi_i(\tau_{\leqslant n}\Bbbk) \simeq \pi_i\Bbbk$ if $0 \leqslant i \leqslant n$ and 0 otherwise, and in $\pi_i(\tau_{\leqslant n}\Bbbk)$, any element $x \in \pi_i(\tau_{\leqslant n}\Bbbk)$ for $i > 0$ is nilpotent ($x^{n+1} = 0$).

definition 2.2.13 (Zariski affine scheme) : *We define the category* **AffSch** *of* (Zariski) affine schemes *as the opposite of the category* $\mathcal{CA}\mathrm{lg}_{\mathbb{Z}}^{\heartsuit} \simeq \mathcal{CA}\mathrm{lg}_{\pi_0\Bbbk}$ *of (ordinary) commutative rings:*

$$\mathbf{AffSch} \simeq \mathcal{CA}\mathrm{lg}_{\mathbb{Z}}^{\heartsuit,\mathrm{op}} \simeq \mathbf{Ring}^{\mathrm{op}}$$

*and if $R$ is a commutative ring, then the corresponding affine scheme is denoted* Spec $R$. *For each commutative rings $R$ and $S$ in $\mathcal{CA}\mathrm{lg}_{\mathbb{Z}}^{\heartsuit}$, the* Hom *between the affine schemes is a set, which is:*

$$\mathrm{Hom}_{\mathcal{CA}\mathrm{lg}_{\mathbb{Z}}^{\heartsuit}}(R, S) = \mathrm{Hom}_{\mathbf{AffSch}}(\mathsf{Spec}\, S, \mathsf{Spec}\, R) \quad \in \mathcal{S}_{\leqslant 0},\ \textit{that is a}\ \text{set}$$

*The chosen Grothendieck topology is the Zariski topology, that is $\{\coprod R_i \to S\}_{i\in I}$ is a Zariski cover if for each $i \in I$, the map $\{R_i \hookrightarrow S\}$ is a (Zariski) open immersion.*

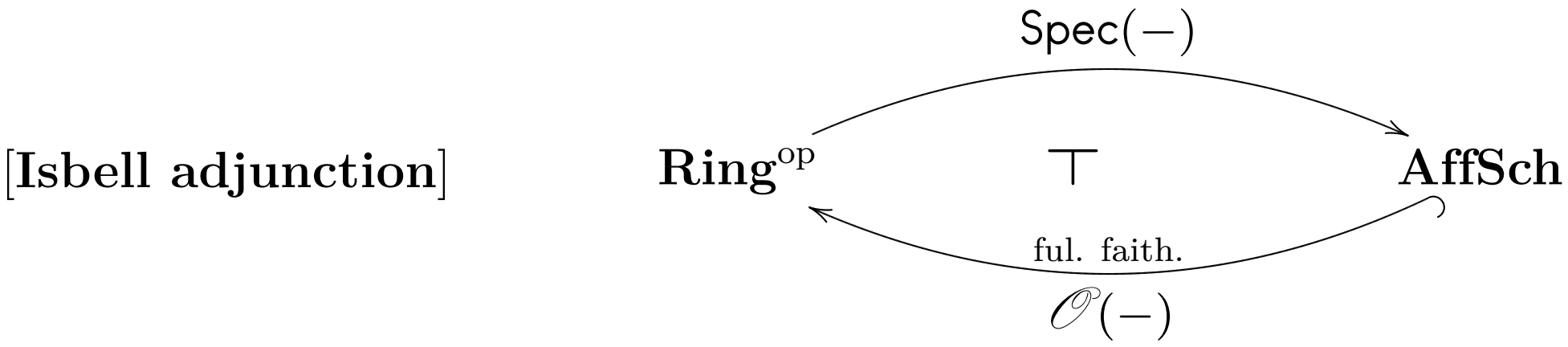

definition 2.2.14 (Affine derived scheme) : *We define the category* **AffDerSch** *of* affine derived schemes *as the opposite of the category $s\,\mathcal{CA}\mathrm{lg}_{\Bbbk}$ of simplicial commutative rings:*

$$\mathbf{AffDerSch} \simeq s\,\mathcal{CA}\mathrm{lg}_{\Bbbk}^{\mathrm{op}}$$

*and if $R$ is a simplicial commutative ring, then the corresponding affine derived scheme is denoted* Spec $R$. *For each simplicial commutative rings $R$ and $S$ in $s\,\mathcal{CA}\mathrm{lg}_{\Bbbk}$, the mapping space between the affine derived schemes then is:*

$$\mathrm{Map}_{s\,\mathcal{C}\mathrm{Alg}_{\Bbbk}}(R,S) = \mathrm{Map}_{\mathbf{AffDerSch}}(\mathsf{Spec}\,S, \mathsf{Spec}\,R) \quad \in \mathcal{S}\mathrm{paces} =: \mathcal{S}$$

*The chosen Grothendieck $\infty$-topology is the étale topology, that is $\{R \to S_i\}_{i\in I}$ is an étale cover if $\{\pi_0 R \to \pi_0 S_i\}_{i\in I}$ is an étale cover in the classical sense.*

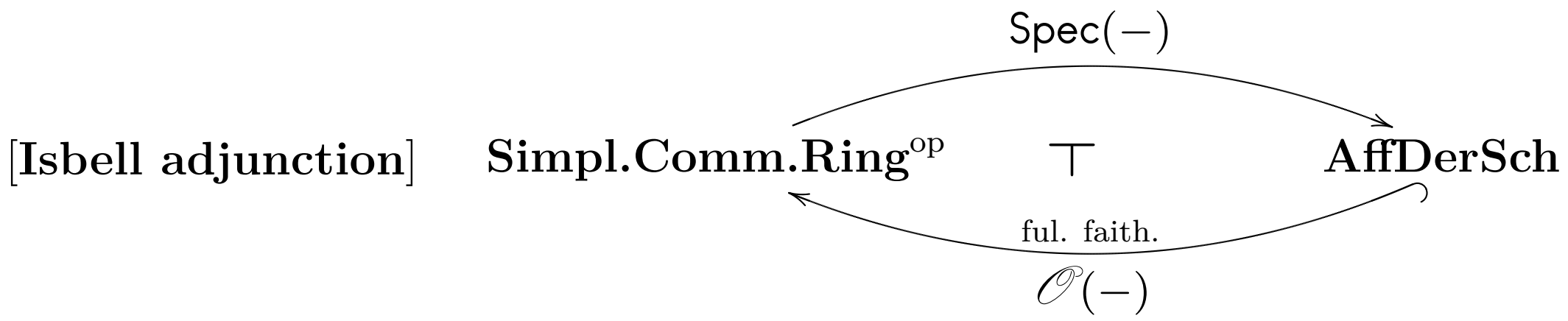

DEFINITION 2.2.15 : *A morphism of schemes over* $\mathsf{Spec}(\mathbb{Z})$*, denoted* $R_{/\mathbb{Z}} \xrightarrow{\text{ét}} S_{/\mathbb{Z}}$ *is said to be an* étale morphism of schemes *if, locally, for the chosen Zariski topology, that is for the Zariski cover:* $\{\coprod R_i \to S\}_{i\in I}$ *where for each* $i \in I$*, the map* $\{R_i \hookrightarrow S\}$ *is a (Zariski) open immersion, that is each open immersion*

REMARK 2.2.16 : An open immersion as map of schemes is étale because it is locally an isomorphism.

DEFINITION 2.2.17 : *A* quasi-finite morphism of schemes *is defined to be of finite type, and each point of the source scheme is isolated in its fiber.*

LEMMA 2.2.18 : *A smooth map which is also quasi-finite is étale.*

LEMMA 2.2.19 : *A smooth map* $R \to S$ *of (ordinary) commutative rings such that the Khäler differentials vanish, that is* $\Omega^1_{S/R} = 0$*, is étale.*

Important note: an étale morphism of schemes is smooth, and for the BCFW étale topology, the immersion of each open positroid variety in the BCFW tiling of the Amplituhedron is étale since it is smooth and a quasi-finite morphism of schemes, that is it is of finite-type and each point of the open positroid variety is isolated in its fiber.

DEFINITION 2.2.20 (étale topology $\infty$-sheaf condition) : *The following compatibility condition of sheaves is called* étale topology $\infty$-sheaf condition with respect to the étale topology *on an $\infty$-category $\mathscr{C}$, or, for short,* étale $\infty$-sheaf condition *(on $\mathscr{C}$):*

- *If* $R \xrightarrow[\text{faithfully}]{\text{étale}} S$ *is a faithfully étale map in the category $\mathscr{C}$, and if $\mathscr{F}$ is a functor* *in* $\mathfrak{Fun}^{\infty}(\mathscr{C}, \mathcal{S})$, then the following map is an isomorphism:

$$\mathscr{F}(R) \xrightarrow{\sim} \operatorname{Tot} \mathscr{F}(S^{\otimes_R \bullet+1})$$

where $\operatorname{Tot} \mathscr{F}(S^{\otimes_R \bullet+1})$ is the totalization of the Amitsur complex:

$$\operatorname{Tot}\left(\ldots \mathscr{F}(S \otimes_R S \otimes_R S) \begin{smallmatrix} \longleftarrow \\ \longrightarrow \\ \longleftarrow \\ \longrightarrow \\ \longleftarrow \end{smallmatrix} \mathscr{F}(S \otimes_R S) \begin{smallmatrix} \longleftarrow \\ \longrightarrow \\ \longleftarrow \end{smallmatrix} \mathscr{F}(S)\right)$$

DEFINITION 2.2.21 (étale $\infty$-sheaves) : *The category of* étale $\infty$-sheaves of spaces on $\mathcal{CA}\mathrm{lg}_{\mathbb{Z}}^{\heartsuit}$ *is the full subcategory of the* $\infty$*-category of functors* $\mathfrak{Fun}^{\infty}(\mathcal{CA}\mathrm{lg}_{\mathbb{Z}}^{\heartsuit}, \mathcal{S})$, denoted $\boldsymbol{\mathcal{S}}\mathbf{hv}^{\text{ét}}$, and satisfying the étale topology $\infty$-sheaf condition for $\mathcal{CA}\mathrm{lg}_{\mathbb{Z}}^{\heartsuit}$.

REMARK 2.2.22 : This means that if we have an étale cover $\coprod R_i \to S$ defining an étale topology of the scheme $S$, we can recover $\mathscr{F}(S)$ from each $\mathscr{F}(R_i)$. This is the $\infty$-*version of the classical étale sheaf condition.*

DEFINITION 2.2.23 : *We define the subcategory, denoted* $\boldsymbol{\mathcal{S}}\mathbf{hv}^{\text{ét},\heartsuit}$, *of discrete-valued sheaves in* $\boldsymbol{\mathcal{S}}\mathbf{hv}^{\text{ét}}$, *the full subcategory of sheaves taking values in discrete spaces, that is* 0*-truncated spaces.*

DEFINITION 2.2.24 : *We call* Spectrum*, denoted* $\mathsf{Spec}$*, the Yoneda embedding functor* $\mathcal{CA}\mathrm{lg}_{\mathbb{Z}}^{\heartsuit} \longrightarrow \boldsymbol{\mathcal{S}}\mathbf{hv}^{\text{ét}}$, *which is a fully faithful embedding.*

REMARK 2.2.25 : This spectrum functor $\mathsf{Spec}$ factors through the subcategory $\boldsymbol{\mathcal{S}}\mathbf{hv}^{\text{ét},\heartsuit}$, and we can verify that if $R$ and $S$ are ordinary commutative rings (which are then in $\mathcal{CA}\mathrm{lg}_{\mathbb{Z}}^{\heartsuit}$), the evaluation in $S$ of the sheaf $\mathsf{Spec}\, R$ gives $(\mathsf{Spec}\, R)(S) \simeq \operatorname{Hom}_{\mathcal{CA}\mathrm{lg}_{\mathbb{Z}}^{\heartsuit}}(R, S)$ which is indeed the set (or 0-truncated space) of morphisms from $R$ to $S$ in the 1-category of (ordinary) commutative rings.

We have by definition of the category of affine schemes, denoted **AffSch**, that it is the opposite of the category $\boldsymbol{\mathcal{CA}}\mathbf{lg}_{\mathbb{Z}}^{\heartsuit}$. We want now to describe the following fully faithful embeddings of subcategories:

$$\mathbf{AffSch} \hookrightarrow \mathbf{Scheme} \hookrightarrow \mathcal{S}\mathbf{hv}^{\mathbf{\acute{e}t}}$$

REMARK 2.2.26 : The category $\mathcal{S}\mathbf{hv}^{\mathbf{\acute{e}t}}$ of étale sheaves has limits and colimits, so in particular pushouts and pullbacks.

DEFINITION 2.2.27 : *In the category $\mathcal{S}\mathbf{hv}^{\mathbf{\acute{e}t}}$ of étale sheaves, a morphism $f : \mathscr{F} \longrightarrow \mathscr{G}$ is* representable *if in the following (homotopy) pullback square:*

$$\begin{array}{ccc} P & \longrightarrow & \mathsf{Spec} R \\ \downarrow & \lrcorner & \downarrow \\ \mathscr{F} & \xrightarrow{f} & \mathscr{G} \end{array} \tag{2.2.28}$$

*the étale sheaf $P$ is isomorphic to a coproduct $\coprod X_i$ with $X_i$ is a quasiaffine scheme, that is $X_i \subseteq \mathsf{Spec}\, T_i$ is open and quasicompact in $\mathsf{Spec}\, T_i$ where $\mathsf{Spec}\, T_i$ is an affine scheme. In other words, the fibers of the map $f : \mathscr{F} \longrightarrow \mathscr{G}$ are unions of quasiaffine schemes.*

REMARK 2.2.29 : We do not consider in this definition that the representability can be made of infinite covering of affine schemes, since the infinite union of affine schemes is not the spectrum of any ring.

PROPOSITION 2.2.30 : *The essential image of the Yoneda embedding in the point of view of functor of points:*

$$\begin{array}{ccc} \mathbf{Scheme} & \xhookrightarrow{\text{fully faithful}} & \mathcal{S}\mathbf{hv}^{\mathbf{\acute{e}t}} \\ X & \longmapsto & \big(R \mapsto X(R)\big) \end{array}$$

*consists of the sheaves $\mathscr{F}$ in $\mathcal{S}\mathbf{hv}^{\mathbf{\acute{e}t}}$ such that there exists:*

$$\coprod \mathsf{Spec}\ S_i \xrightarrow[\text{surj. (on alg. closed points)}]{\text{representable}} \mathscr{F}$$

*and each* $\mathsf{Spec}\ S_i \hookrightarrow \mathscr{F}$ *is an open immersion, in the sense that in the following pullback square:*

$$\begin{array}{ccc} \coprod X_i & \longrightarrow & \mathsf{Spec}\, R \\ \downarrow & \lrcorner & \downarrow \\ \mathsf{Spec}\, S_i & \xrightarrow{f} & \mathscr{F} \end{array} \tag{2.2.31}$$

*each* $X_i$ *is quasiaffine with* $X_i \hookrightarrow \mathsf{Spec}\ R$ *is an open immersion (the image of each* $X_i$ *is isomorphic to a subset of* $\mathsf{Spec}\ R$*).*

DEFINITION 2.2.32 : *From the point of view of functor of points, we define a* scheme *as a sheaf that we can cover with affine schemes, such that the intersections are also affine schemes. The schemes form the category* **Scheme** *of schemes.*

REMARK 2.2.33 : In the previous definition of the essential image of the category of schemes in the category of étale sheaves, we used the Grothendieck topology, which uses a Zariski open cover with Zariski open immersions. We define this way the category of schemes, which is a certain class of étale sheaves. This suggests the use of other Grothendieck topologies, in particular, the étale topology. For example, the following class of étale sheaves will be the one that interests us in the case of the Amplituhedron functor.

DEFINITION 2.2.34 (Deligne-Mumford stack) : *A sheaf* $\mathscr{F}$ *in* $\boldsymbol{\mathcal{S}}\mathbf{hv}^{\text{ét}}$ *is a* Deligne-Mumford stack *if there exists a surjective étale map* $\coprod \mathsf{Spec}\ S_i \to \mathscr{F}$*, that is there exists an étale cover of the sheaf* $\mathscr{F}$ *in* $\boldsymbol{\mathcal{S}}\mathbf{hv}^{\text{ét}}$*.*

REMARK 2.2.35 : No condition is asked for the affineness or the representability of the diagonal.

REMARK 2.2.36 : When we consider a scheme as a functor, and we evaluate a scheme on a ring, we get a set. But if we evaluate a Deligne-Mumford stack on a ring, we obtain a space, which, in particular, has a higher homotopy. Considering a space as an $\infty$-category, and the category of sheaves on an $\infty$-category with values in an $\infty$-category, we can then reach higher homotopy in the framework of $\infty$-toposes. For our purpose, we will consider $\infty$-categories of sheaves with values in the $\infty$-category of simplicial commutative algebras or better: in the $\infty$-category of $\mathbb{E}_\infty$-rings spectra, which allows us to take a dual $\mathbb{E}_\infty$-rings spectrum that encodes the Schubert geometry of the Grassmannian and of the Amplituhedron. That is why we consider the framework of $\infty$-toposes, étale covers, and Deligne-Mumford stacks as the most relevant one for the purpose we assigned us, namely construct a dual Amplituhedron that respects Schubert geometry of the (dual) TNN Grassmannian.

### 2.2.3 Spectral Deligne-Mumford stacks.

#### 2.2.3.1 Definition via étale spectra and sheaves of $\mathbb{E}_\infty$-rings spectra.

DEFINITION 2.2.37 : *Let $\mathscr{F}$ be a Deligne-Mumford stack. We define the full subcategory $\mathscr{F}^{\textbf{ét}}$ of the slice (comma) category $\mathcal{S}\mathbf{hv}^{\textbf{ét}}/\mathscr{F}$, such that if $\mathscr{G}$ is an object of $\mathscr{F}^{\textbf{ét}}$, there exists a natural transformation $\mathscr{G} \xrightarrow{\text{ét}} \mathscr{F}$ which is étale (and then representable).*

REMARK 2.2.38 : This full subcategory $\mathscr{F}^{\textbf{ét}}$ is small, that can be turned into a Grothendieck $\infty$-category, by making it an $\infty$-site with the étale topology, where an étale cover of any object in $\mathscr{F}^{\textbf{ét}}$ is a family of surjective étale maps (composition of étale maps is étale).

We can then consider the subcategory of sheaves $\mathbf{Shv}(\mathscr{F}^{\textbf{ét}})$ in the presheaves category $\mathfrak{Fun}((\mathscr{F}^{\textbf{ét}})^{\text{op}}, \mathcal{S})$, satisfying the *étale $\infty$-sheaf condition.* The category $\mathscr{F}^{\textbf{ét}}$ being equipped in an $\infty$-site, we then have a prototypal $\infty$-topos, and we then have the following:

PROPOSITION 2.2.39 ([103]) : *The $\infty$-category $\mathcal{S}\mathbf{hv}(\mathscr{F}^{\textbf{ét}})$ associated to a Deligne-Mumford stack $\mathscr{F}$ in $\mathcal{S}\mathbf{hv}^{\textbf{ét}}$ is an $\infty$-topos.*

Let us now consider a point $\mathsf{Spec}\, R \longrightarrow \mathscr{F}$ in the Deligne-Mumford stack $\mathscr{F}$ with $R$ in $\mathcal{CA}\mathrm{lg}_{\mathbb{Z}}^{\heartsuit}$. We can equip $\mathbf{Shv}(\mathscr{F}^{\textbf{ét}})$ with a sheaf of rings $\mathscr{O}$ such that $\mathscr{O}(\mathsf{Spec}\, R) = R$. The couple $\left(\boldsymbol{\mathcal{S}}\mathbf{hv}(\mathscr{F}^{\text{ét}}), \mathscr{O}\right)$ is the prototypal model of spectral Deligne-Mumford stacks.

We now want to characterize Deligne-Mumford stacks among $\infty$-topoi instead of in sheaf theory and enlarge them to their spectral version.

DEFINITION 2.2.40 (Presentable $\infty$-category: [98]) : *A* presentable $\infty$-category $\mathcal{X}$ *is such that there exists a small* $\infty$*-category* $\mathscr{C}$ *with an (accessible) functor called* localization *between the category of* $\mathcal{S}$*-valued presheaves on* $\mathscr{C}$ *and the* $\infty$*-category* $\mathcal{X}$ *with the following adjunction:*

$$\mathfrak{Fun}(\mathscr{C}^{\text{op}}, \mathcal{S}) \underset{\text{localization}}{\overset{\text{ful. faith.}}{\underset{\longrightarrow}{\overset{\hookleftarrow}{\top}}}} \mathcal{X} \tag{2.2.41}$$

REMARK 2.2.42 : The localization functor $\mathfrak{Fun}(\mathscr{C}^{\text{op}}, \mathcal{S}) \longrightarrow \mathcal{X}$ is right exact since it is a left adjoint. The right adjoint to the localization functor is the fully faithful functor $\mathcal{X} \xhookrightarrow{\text{ful. faith.}} \mathfrak{Fun}(\mathscr{C}^{\text{op}}, \mathcal{S})$.

DEFINITION 2.2.43 ($\infty$-topos: [98]) : *An* $\infty$-topos *is a presentable* $\infty$*-category* $\mathfrak{X}$ *such that the localization functor is left exact:*

$$\mathfrak{Fun}(\mathscr{C}^{\text{op}}, \mathcal{S}) \underset{\text{loc. left ex.}}{\overset{\text{ful. faith.}}{\underset{\longrightarrow}{\overset{\hookleftarrow}{\top}}}} \mathfrak{X} \tag{2.2.44}$$

REMARK 2.2.45 : The localization functor $\mathfrak{Fun}(\mathscr{C}^{\text{op}}, \mathcal{S}) \xrightarrow[\text{loc. left ex.}]{} \mathfrak{X}$ is right exact since it is a left adjoint. Still here in the definition of an $\infty$-topos, this functor is asked to be exact (as a classical sheafification in the 1-categories world), that is it is also asked to be left exact. It preserves finite limits (left exactness) and all colimits (right exactness from a category of presheaves on a small category).

REMARK 2.2.46 : If we consider the category $\mathscr{C}$ is a site, then the category $\mathfrak{Fun}(\mathscr{C}^{\text{op}}, \mathcal{S})$ is an $\infty$-category of $\mathcal{S}$-valued $\infty$-presheaves. The $\infty$-topos $\mathfrak{X}$ is then a full subcategory of an $\infty$-category of presheaves that satisfy a sheaf condition. The localization is then a sheafification. Therefore, in this case, it is indeed exact.

REMARK 2.2.47 : Even if all 1-topoi come from a localization of a category of presheaves, we will notice that not all $\infty$-topoi come from the localization of a certain category of presheaves.

DEFINITION 2.2.48 : *Let* $R \xrightarrow{f} S$ *be a map of* $\mathcal{CA}\mathrm{lg}$*, that is a map of* $\mathbb{E}_\infty$*-ring spectra. The map* $f$ *is said to be* flat *if the map:* $\pi_0 R \to \pi_0 S$ *is flat and the induced map:* $\pi_n R \otimes_{\pi_0 R} \pi_0 S \xrightarrow{\sim} \pi_n S$ *is an isomorphism for all* $n \in \mathbb{N}$*.*

DEFINITION 2.2.49 : *Let* $R \xrightarrow{f} S$ *be a map in* $\mathcal{CA}\mathrm{lg}$*. The map* $f$ *is said to be* étale *if it is flat and the induced map:* $\pi_0 R \to \pi_0 S$ *is étale.*

DEFINITION 2.2.50 : *If* $R$ *is fixed in* $\mathcal{CA}\mathrm{lg}$*, then we can define the (slice) category* $\mathcal{CA}\mathrm{lg}_R^{\text{ét}}$ *of étale maps* $R \xrightarrow[\text{ét}]{f} S$ *in* $\mathcal{CA}\mathrm{lg}$*.*

REMARK 2.2.51 : This is an essentially small $\infty$-category and the composition law between étale maps is fulfilled:

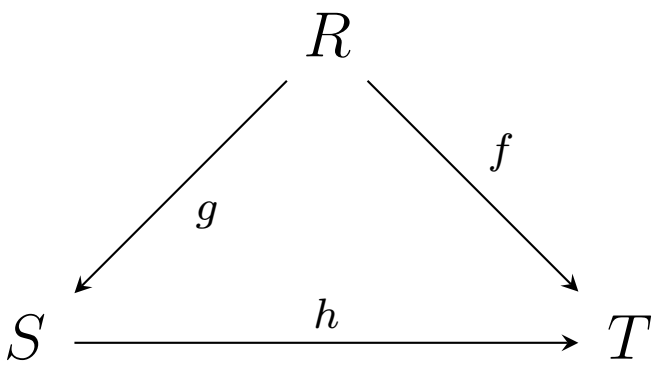

since the map $h$ is étale as soon as the maps $f$ and $g$ are étale.

DEFINITION 2.2.52 (étale spectrum) : *If* $R \in \mathcal{CA}\mathrm{lg}$*, that is* $R$ *is a* $\mathbb{E}_\infty$*-ring spectrum, and if we consider the* $\infty$*-subcategory* $\mathcal{CA}\mathrm{lg}_R^{\text{ét}}$*, then we can define the category of étale sheaves:* $\boldsymbol{\mathcal{S}}\mathbf{hv}_{\mathcal{S}}^{\textbf{ét}}(\mathcal{CA}\mathrm{lg}_R^{\text{ét}})$ *of spaces, which we denote:* Spét $R$*, and call it the* étale spectrum*. It is the subcategory of the presheaves category* $\mathfrak{Fun}((\mathcal{CA}\mathrm{lg}_R^{\text{ét}})^{\mathrm{op}}, \mathcal{S})$ satisfying the étale $\infty$-sheaf condition (see definition 2.2.20):

$$\mathsf{Sp\acute{e}t}\, R := \boldsymbol{\mathcal{S}}\mathbf{hv}_{\mathcal{S}}^{\textbf{ét}}(\mathcal{CA}\mathrm{lg}_R^{\text{ét}}) \tag{2.2.53}$$

PROPOSITION 2.2.54 ([103]) : *The étale spectrum is a localization of the category of* $\mathcal{S}$*-valued presheaves* $\mathfrak{Fun}((\mathcal{CA}\mathrm{lg}_R^{\text{ét}})^{\mathrm{op}}, \mathcal{S})$*. The* $\infty$*-category* Spét $R$ *is then an* $\infty$*-topos.*

REMARK 2.2.55 : The étale spectrum is not yet a spectral Deligne-Mumford stack, since we need to attach a sheaf of $\mathbb{E}_\infty$-rings spectra. We say here that the $\infty$-topos $\mathsf{Sp\acute{e}t}\, R$ play for our spectral Deligne-Mumford stacks the model role of an affine scheme in the theory of classical algebraic geometry.

We need now to define what is a sheaf on an $\infty$-topos.

DEFINITION 2.2.56 : *Let $\mathscr{C}$ be an $\infty$-category, and $\mathfrak{X}$ an $\infty$-topos. The* category of sheaves on the $\infty$-topos $\mathfrak{X}$*, denoted $\mathcal{S}\mathbf{hv}_{\mathscr{C}}(\mathfrak{X})$, is the category of $\mathscr{C}$-valued presheaves on $\mathfrak{X}$ which preserves limits on $\mathfrak{X}^{\mathrm{op}}$ to $\mathscr{C}$. We write:*

$$\mathcal{S}\mathbf{hv}_{\mathscr{C}}(\mathfrak{X}) := \mathfrak{Fun}^{\lim}(\mathfrak{X}^{\mathrm{op}}, \mathscr{C}) \tag{2.2.57}$$

REMARK 2.2.58 : The $\infty$-topos $\mathfrak{X}$ has all colimits, so the opposite $\infty$-category $\mathfrak{X}^{\mathrm{op}}$ (which is not necessarily an $\infty$-topos) has all limits. Thus this makes sense to ask for the functors to preserve limits on $\mathfrak{X}^{\mathrm{op}}$.

REMARK 2.2.59 : It is important to choose right now to define an $\infty$-sheaf by the Equation 2.7.7 containing a condition on preserving limits for $\mathscr{C}$-valued presheaves instead of the classical condition via Čech complexes. Indeed, Čech complexes and the totalization functor for a functor $f : U \to V$ in $\mathfrak{X}$ is

$$\check{C}(f) := (\dots \begin{matrix}\longrightarrow\\ \longrightarrow\\ \longrightarrow\end{matrix} U \times_V U \rightrightarrows U.) \tag{2.2.60}$$

We then have $f$ an *effective epimorphism* if $|\check{C}(f)|_{\mathrm{geom}} \simeq V$, that is: $\operatorname*{colim}_{\Delta^{\mathrm{op}}} \check{C}(f) \simeq V$. Then, since we work in the opposite $\infty$-category $\mathfrak{X}^{\mathrm{op}}$, we have to verify the descent condition for limits, which means exactly that the descent condition we must impose to a functor in $\mathfrak{Fun}(\mathfrak{X}^{\mathrm{op}}, \mathscr{C})$ in Equation 2.7.7 to be a sheaf must be to preserve limits. However, the condition involving preserving certain limits will have the virtue to be more generic later when we use the equivalence of this notion to the one of $\mathcal{M}$od-valued sheaves, and the equivalence with defining them by triples $(\mathfrak{X}, \mathscr{O}_{\mathfrak{X}}, \mathscr{F}_{\mathfrak{X}})$ (see subsection 2.3.6) and therefore speak of $\mathcal{G}$-geometries in the sense of J. Lurie.

Now, if we take the $\infty$-category of spectra $\mathsf{Sp}$ for the $\infty$-category $\mathscr{C}$, we obtain the category of sheaves of spectra on $\mathfrak{X}$, denoted: $\mathcal{S}\mathbf{hv}_{\mathsf{Sp}}(\mathfrak{X})$. Moreover, if we take for $\mathscr{C}$ the $\infty$-category of $\mathbb{E}_\infty$-rings spectra $\mathcal{CA}\mathrm{lg}$, we have the category of sheaves of $\mathbb{E}_\infty$-rings spectra $\mathcal{S}\mathbf{hv}_{\mathcal{CA}\mathrm{lg}}(\mathfrak{X})$ on the $\infty$-topos $\mathfrak{X}$.

DEFINITION 2.2.61 : *We call a* spectrally ringed $\infty$-topos*, a pair* $(\mathfrak{X}, \mathscr{O}_{\mathfrak{X}})$*, where* $\mathfrak{X}$ *is an* $\infty$*-topos, and* $\mathscr{O}_{\mathfrak{X}} \in \mathcal{S}\mathbf{hv}_{\mathcal{CA}\mathrm{lg}}(\mathfrak{X})$ *is a sheaf of* $\mathbb{E}_\infty$*-rings spectra on* $\mathfrak{X}$*.*

DEFINITION 2.2.62 : *Let* $R$ *be an* $\mathbb{E}_\infty$*-rings spectrum in* $\mathcal{CA}\mathrm{lg}$*, and let* $\mathfrak{X} := \mathsf{Sp\acute{e}t}\, R$ *be the* $\infty$*-topos of the étale spectrum of* $R$*. We define the* canonical sheaf $\mathscr{O}_R$ *in* $\mathcal{S}\mathbf{hv}_{\mathsf{Sp}}(\mathfrak{X})$ *as the functor of points:* $\mathscr{O}_R(S) = S$ *for all* $S \in \mathsf{Sp\acute{e}t}\, R$*.*

REMARK 2.2.63 : We have the embedding:

$$\mathcal{S}\mathbf{hv}_{\mathsf{Sp}}(\mathsf{Sp\acute{e}t}\, R) := \mathfrak{Fun}^{\lim}((\mathsf{Sp\acute{e}t}\, R)^{\mathrm{op}}, \mathsf{Sp}) \longhookrightarrow \mathfrak{Fun}((\mathsf{Sp\acute{e}t}\, R)^{\mathrm{op}}, \mathsf{Sp})$$

So the canonical sheaf $\mathscr{O}_R$ in $\mathcal{S}\mathbf{hv}_{\mathsf{Sp}}(\mathfrak{X})$ is just the forgetful functor of the étale topology $\infty$-sheaf condition for all sheaves $S$ of $\mathbb{E}_\infty$-rings spectra with $S \in \mathfrak{X}$, where:

$$\mathfrak{X} := \mathsf{Sp\acute{e}t}\, R := \mathcal{S}\mathbf{hv}_{S}^{\mathbf{\acute{e}t}}(\mathcal{CA}\mathrm{lg}_R^{\acute{e}t})$$

LEMMA 2.2.64 : *The canonical sheaf* $\mathscr{O}_R$ *is not only in* $\mathcal{S}\mathbf{hv}_{\mathsf{Sp}}(\mathfrak{X})$ *but better in* $\mathcal{S}\mathbf{hv}_{\mathcal{CA}\mathrm{lg}}(\mathfrak{X})$*.*

PROPOSITION 2.2.65 : *The pair* $(\mathfrak{X}, \mathscr{O}_{\mathfrak{X}}) := (\mathsf{Sp\acute{e}t}\, R, \mathscr{O}_R)$ *is a spectrally ringed* $\infty$*-topos.*

DEFINITION 2.2.66 : *Let us say that the terminal object* $*$ *of an* $\infty$*-topos* $\mathfrak{X}$ *has an* effective cover *if there exists a collection of objects* $\mathfrak{U}_i \in \mathfrak{X}$ *and an effective epimorphism from the disjoint union of these objects of* $\mathfrak{X}$ : $\coprod \mathfrak{U}_i \to *$*.*

DEFINITION 2.2.67 : *Let us define the spectrally ringed topos* $\mathfrak{X}_i$ *such that:*

$$\mathfrak{X}_i := (\mathfrak{X}_{/\mathfrak{U}_i}, \mathscr{O}_{\mathfrak{X}}\big|_{\mathfrak{X}_{/\mathfrak{U}_i}})$$

*where* $\mathfrak{X}_{/\mathfrak{U}_i}$ *is the slice* $\infty$*-topos* $\mathfrak{X}$ *by* $\mathfrak{U}_i$ *endowed with the canonical sheaf* $\mathscr{O}_{\mathfrak{X}_{/\mathfrak{U}_i}}$ *which is the restriction* $\mathscr{O}_{\mathfrak{X}}\big|_{\mathfrak{X}_{/\mathfrak{U}_i}}$ *to* $\mathfrak{X}_{/\mathfrak{U}_i}$ *of the canonical sheaf* $\mathscr{O}_{\mathfrak{X}}$*.*

REMARK 2.2.68 : If $\mathfrak{X}$ is an $\infty$-topos, then $\mathfrak{X}_{/\mathfrak{U}_i}$ is an $\infty$-topos.

DEFINITION 2.2.69 (Nonconnective Spectral Deligne-Mumford stack:[103] Def. 1.4.4.2 ) : *A* nonconnective spectral Deligne-Mumford stack *is a spectrally ringed $\infty$-topos $\mathfrak{X} = (\mathcal{X}, \mathcal{O}_{\mathcal{X}})$ for which the terminal object $*$ of $\mathcal{X}$ has an effective cover $\{\coprod \mathcal{U}_i \to *\}$, and such that each $\mathcal{U}_i \in \mathcal{X}$ of the effective cover is such that $\mathfrak{X}_i := (\mathcal{X}_{/\mathcal{U}_i}, \mathcal{O}_{\mathcal{X}}\big|_{\mathcal{X}_{/\mathcal{U}_i}})$ is equivalent to some $(\mathsf{Sp\acute{e}t}\, R_i, \mathcal{O}_{R_i})$ for an $R_i \in \mathcal{CA}\mathrm{lg}$.*

REMARK 2.2.70 : Each $\mathfrak{X}_i := (\mathcal{X}_{/\mathcal{U}_i}, \mathcal{O}_{\mathcal{X}}\big|_{\mathcal{X}_{/\mathcal{U}_i}}) \simeq (\mathsf{Sp\acute{e}t}\, R_i, \mathcal{O}_{R_i})$ play the role of some "affine piece".

DEFINITION 2.2.71 (Spectral Deligne-Mumford stack: [103], Def. 1.4.4.2 ) : *A* spectral Deligne-Mumford stack *is a nonconnective spectral Deligne-Mumford stack such that the negative homotopy groups of $\mathcal{O}_{\mathcal{X}}$ are zero: $\pi_i(\mathcal{O}_{\mathcal{X}}) = 0$ for $i < 0$.*

REMARK 2.2.72 : This means that in the definition 2.2.71, we can choose the $R_i$ of definition 2.2.69 to be connective $\mathbb{E}_\infty$-ring spectra.

REMARK 2.2.73 : We can pass from a nonconnective spectral Deligne-Mumford stack to a (connective) spectral Deligne-Mumford stack by taking the connective cover of the sheaf of $\mathbb{E}_\infty$-rings spectra.

REMARK 2.2.74 : If $(\mathcal{X}, \mathcal{O}_{\mathcal{X}})$ is a spectrally ringed $\infty$-topos, then $(\mathcal{X}^\heartsuit, \pi_0 \mathcal{O}_{\mathcal{X}})$ is an ordinary ringed topos.

REMARK 2.2.75 : We have equivalences :

$$\mathsf{Sp\acute{e}t}\, R \xrightarrow{\sim} \mathsf{Sp\acute{e}t}\, \tau_{\geqslant 0} R \xleftarrow{\sim} \mathsf{Sp\acute{e}t}\, \pi_0 R \tag{2.2.76}$$

since the underlying étale algebras are the same:

$$\mathcal{CA}\mathrm{lg}_R^{\text{ét}} \xleftarrow{\sim} \mathcal{CA}\mathrm{lg}_{\tau_{\geqslant 0} R}^{\text{ét}} \xrightarrow{\sim} \mathcal{CA}\mathrm{lg}_{\pi_0 R}^{\text{ét}} \tag{2.2.77}$$

The categories of presheaves on these étale algebras are the same, so the topologies are the same, so the $\infty$-toposes are the same. But even if the $\infty$-toposes are the same, their associated sheaves of $\mathbb{E}_\infty$-rings spectra differ, so the spectral Deligne-Mumford stacks differ.

#### 2.2.3.2 Quasicoherent sheaves on a spectral Deligne-Mumford stack.

DEFINITION 2.2.78 : *Let $\mathfrak{X}$ be an $\infty$-topos. We denote* $\mathsf{Mod}_{\mathscr{O}_{\mathfrak{X}}} := \mathsf{Mod}_{\mathscr{O}_{\mathfrak{X}}}(\boldsymbol{\mathcal{S}}\mathbf{hv}_{\mathsf{Sp}}(\mathfrak{X}))$ *the category of modules on* $\mathscr{O}_{\mathfrak{X}}$ *in the category* $\boldsymbol{\mathcal{S}}\mathbf{hv}_{\mathsf{Sp}}(\mathfrak{X})$ *of sheaves of spectra on* $\mathfrak{X}$.

DEFINITION 2.2.79 : *Let us define the full subcategory of* quasicoherent sheaves*, denoted* $\mathsf{QCoh}(\mathfrak{X})$*, in* $\mathsf{Mod}_{\mathscr{O}_{\mathfrak{X}}}$*, such that, if* $\mathscr{F} \in \mathsf{QCoh}(\mathfrak{X})$ *and* $f : U \longrightarrow V$ *a morphism between affine pieces in an atlas of covers of* $\mathfrak{X}$*, that is* $f : U := \mathsf{Sp\acute{e}t}\, R \longrightarrow \mathsf{Sp\acute{e}t}\, S =: V$ *for some* $R$ *and* $S$ *in* $\mathcal{CA}\mathrm{lg}$*, then:*

$$\mathscr{F}(V) \otimes_{\mathscr{O}(V)} \mathscr{O}(U) \simeq \mathscr{F}(U) \tag{2.2.80}$$

REMARK 2.2.81 : The sheaf $\mathscr{O}_{\mathfrak{X}}$ is quasi-coherent, like all objects constructed as colimits of $\mathscr{O}_{\mathfrak{X}}$.

REMARK 2.2.82 : The subcategory $\mathsf{QCoh}(\mathfrak{X})$ of quasicoherent sheaves on $\mathfrak{X}$ is presentable and closed under colimits since the inclusion of the subcategory $\mathsf{QCoh}(\mathfrak{X})$ in $\mathsf{Mod}_{\mathscr{O}_{\mathfrak{X}}}$ is exact and admits a right adjoint named the coherator.

PROPOSITION 2.2.83 : *The $\infty$-category* $\mathsf{QCoh}(\mathsf{Sp\acute{e}t}\, R)$ *of quasicoherent sheaves on the étale spectrum of* $R \in \mathcal{CA}\mathrm{lg}$ *and the $\infty$-category* $\mathsf{Mod}_R$ *of $R$-module spectra are equivalent:*

$$\mathsf{QCoh}(\mathsf{Sp\acute{e}t}\, R) \quad \simeq \quad \mathsf{Mod}_R$$

PROPOSITION 2.2.84 ([106]Def. 3.16) : *Let $\mathfrak{X}$ be an $\infty$-topos. The following $\infty$-categories are equivalent:*

$$\mathsf{QCoh}(\mathfrak{X}) \xrightarrow{\sim} \lim_{\mathsf{Sp\acute{e}t}\, R \xrightarrow{\text{étale}} \mathfrak{X}} \mathsf{Mod}_R \tag{2.2.85}$$

*where the limit is taken in the category of étale maps.*

REMARK 2.2.86 : This means that giving a quasicoherent sheaf of spectra on an $\infty$-topos $\mathfrak{X}$ is the same as giving an $\infty$-category of $R_i$-module spectra on each affine piece of the atlas $\{ (\mathsf{Sp\acute{e}t}\, R_i, \mathscr{O}_{R_i}) \}$ for $R_i \in \mathcal{CA}\mathrm{lg}$.

Let us recall the following proposition before the following remark below:

PROPOSITION 2.2.87 : *If* $\mathfrak{X}$ *is an ordinary scheme* $\mathscr{X}$*, then the homotopy category* $\mathsf{ho}(\mathsf{QCoh}(\mathscr{X}))$ *is the classical derived category* $\mathscr{D}(\mathscr{X}) := \mathscr{D}(\mathsf{QCoh}(\mathscr{X}))$ *of the category of quasicoherent sheaves on* $\mathscr{X}$*.*

REMARK 2.2.88 : In order to construct a spectral Deligne-Mumford stack, it is easier to define the underlying $\infty$-topos, which is generally given, than to build a sheaf on it that fulfills the local sheaf conditions of definition 2.2.69 and definition 2.2.71. In other words, we can reverse the problem by saying that it is by choosing a suitable sheaf associated with this $\infty$-topos in the perspective that we have defined that the structure of this spectral Deligne-Mumford stack will reveal all its power and will allow us to obtain deeper results than remaining on the surface of the phenomena. (For deeper considerations, see Lurie, "*Structured spaces*" [91])

#### 2.2.3.3 Quasi-affineness.

DEFINITION 2.2.89 ([93] Def. 4.4.2, [145] Def. 2.2.4, [108] Def. 2.1.11) : *A nonconnective spectral Deligne-Mumford stack* $\mathfrak{X}$ *is* 0-affine *if the following map is an equivalence:*

$$\Gamma \,:\, \mathsf{QCoh}(\mathfrak{X}) \xrightarrow{\sim} \mathsf{Mod}_{\Gamma(\mathfrak{X},\mathscr{O}_{\mathfrak{X}})} \tag{2.2.90}$$

THEOREM 2.2.91 (Serre) : *Let* $(\mathscr{X}, \mathscr{O}_{\mathscr{X}})$ *be a (classical) scheme. Then the following are equivalent:*

- *the scheme* $\mathscr{X}$ *is affine* ($\mathscr{X} = \mathsf{Spec}\ R$*, for an (ordinary) commutative ring* $R$*)*
- *the cohomology functors* $H^i(\mathscr{X}, \mathscr{F}) = 0$*, for all* $i > 0$ *and all* $\mathscr{F} \in \mathsf{QCoh}(\mathscr{X})$
- *the cohomology functors* $H^i(\mathscr{X}, \mathscr{O}_{\mathscr{X}}) = 0$*, for all* $i > 0$ *and* $\mathsf{QCoh}(\mathscr{X}) \xrightarrow{\sim} \mathsf{Mod}_{H^0(\mathscr{X},\mathscr{O}_{\mathscr{X}})}$ *is an equivalence.*

REMARK 2.2.92 : If the higher cohomologies of the structure sheaf $\mathscr{X}$ vanish and $\mathscr{X}$ is 0-affine as a nonconnective spectral Deligne-Mumford stack then the scheme $\mathscr{X}$ is affine.

REMARK 2.2.93 : If the scheme $\mathscr{X}$ is quasi-affine, that is $\mathscr{X}$ is an open subscheme in a scheme $\mathscr{Y}$ $=\mathsf{Spec}\ R$, for an (ordinary) commutative ring $R$ and $\mathscr{X}$ is quasi-compact, the cohomology functors $H^i(\mathscr{X},\mathscr{F})$ do not vanish anymore for all $i>0$ and all $\mathscr{F}\in\mathsf{QCoh}(\mathscr{X})$, but we have the following:

THEOREM 2.2.94 (0-affineness: [93]) : *Let $(\mathscr{X},\mathscr{O}_{\mathscr{X}})$ be a (classical) scheme. Then the following are equivalent:*

- *the scheme $\mathscr{X}$ is quasi-affine*
- *the derived category of the stable $\infty$-category $\mathsf{QCoh}(\ \mathscr{X}\ )$ is equivalent to the following category of sheaves of $\mathbb{E}_\infty$-rings:*

$$\mathscr{D}^+(\mathsf{QCoh}(\mathscr{X}))\xrightarrow{\sim}\mathsf{Mod}_{\mathbb{R}\Gamma(\mathscr{X},\mathscr{O}_{\mathscr{X}})}\tag{2.2.95}$$

  *where $\mathbb{R}\Gamma$ is the right derived functor of global sections functor. So $\mathbb{R}\Gamma(\mathscr{X},\mathscr{O}_{\mathscr{X}})$ is an $\mathbb{E}_\infty$-ring which is coconnective (it has only negative homotopy groups), and taking its global sections gives a functor which is an equivalence.*

*We say then that a scheme is 0-*affine *and in this case, the global sections should be also discrete.*

LEMMA 2.2.96 : *The stack $(\mathfrak{X},\mathscr{O}_{\mathfrak{X}})$ is 0-affine if and only if :*

- *the sheaf $\mathscr{O}_{\mathfrak{X}}$ is compact in* $\mathsf{QCoh}(\mathfrak{X})$ *(it commutes to colimits)*
- *the sheaf $\mathscr{O}_{\mathfrak{X}}$ generates* $\mathsf{QCoh}(\mathfrak{X})$ *(that is $\Gamma$ is conservative in Equation 2.5.2)*

*Proof.* We have (Equation 2.5.2) : $\Gamma(-)=\mathsf{Map}(\mathscr{O}_{\mathfrak{X}},-)$ is the mapping spectrum. So if $\Gamma$ is an equivalence, then $\mathscr{O}_{\mathfrak{X}}$ has to be compact. And it is then also conservative. Schwede-Shipley ([132]) implies $\mathscr{O}_{\mathfrak{X}}$ is a compact generator via Morita theory since the ring $\mathsf{End}(\mathscr{O}_{\mathfrak{X}})$ is exactly the ring of global sections. □

### 2.2.4 Quasi-Coherent Sheaves on a Functor.

Let $\mathfrak{X} = (\mathfrak{X}, \mathcal{O}_{\mathfrak{X}})$ be a spectral Deligne-Mumford stack and let us denote $h_{\mathfrak{X}} := \text{Morph}_{\infty\mathcal{T}\text{op}^{\text{str. Hen}}_{\mathcal{CA}\text{lg}}}(\mathsf{Sp\acute{e}t}(-), \mathfrak{X})$ the fully faithful Yoneda embedding of the category of spectral Deligne-Mumford stacks in the category of functors $\mathfrak{Fun}(\mathcal{CA}\text{lg}^{\text{cn,op}}, \mathsf{Sp})$ :

$$\begin{array}{ccc} \mathbf{SpDMStk} & \xhookrightarrow[\text{Yoneda}]{\text{fully faithful}} & \mathfrak{Fun}(\mathcal{CA}\text{lg}^{\text{cn,op}}, \mathsf{Sp}) \\ \mathfrak{X} = (\mathfrak{X}, \mathcal{O}_{\mathfrak{X}}) & \longmapsto & h_{\mathfrak{X}} := \text{Morph}_{\infty\mathcal{T}\text{op}^{\text{str. Hen}}_{\mathcal{CA}\text{lg}}}(\mathsf{Sp\acute{e}t}(-), \mathfrak{X}) \end{array} \tag{2.2.97}$$

PROPOSITION 2.2.98 ([103], Prop. 6.2.3.4) : *(1) For every functor $\mathfrak{X} : \mathcal{CA}\text{lg}^{\text{cn}} \longrightarrow \mathsf{Sp}$, the category $\mathsf{QCoh}(\mathfrak{X})$ is a stable $\infty$-category which admits small colimits.*
*(2) If $\varphi : \mathfrak{X} \to \mathfrak{Y}$ is a natural transformation in $\mathfrak{Fun}(\mathcal{CA}\text{lg}^{\text{cn,op}}, \mathsf{Sp})$, then the pullback functor $\varphi^* : \mathsf{QCoh}(\mathfrak{Y}) \to \mathsf{QCoh}(\mathfrak{X})$ preserves small colimits, and, in particular, this pullback functor is exact.*

Since in the $\infty$-category of presheaves $\mathfrak{Fun}(\mathcal{CA}\text{lg}^{\text{cn,op}}, \text{Sp})$ (and in the subcategory of spectral Deligne-Mumford stacks), all the affine objects are small colimits of representable functors, that is colimits of $\{h_{\mathfrak{X}_i}\}_{i\in\mathcal{D}}$ for a small diagram $\mathcal{D}$, we can define quasi-coherent sheaves on all presheaves functors $h_{\mathfrak{X}}$ in $\mathfrak{Fun}(\mathcal{CA}\text{lg}^{\text{cn,op}}, \text{Sp})$ and extend the definition to all functors in $\mathfrak{Fun}(\mathcal{CA}\text{lg}^{\text{cn,op}}, \text{Sp})$ by right Kan extension of the Yoneda embedding from affine spectral Deligne-Mumford stacks to the category of all presheaves. This is summarized in the following:

PROPOSITION 2.2.99 ([103], Prop. 6.2.4.1) : *Let $\mathfrak{X} = (\mathfrak{X}, \mathcal{O}_{\mathfrak{X}})$ be a spectral Deligne-Mumford stack and let $h_{\mathfrak{X}} := \text{Morph}_{\infty\mathcal{T}\text{op}^{\text{str. Hen}}_{\mathcal{CA}\text{lg}}}(\mathsf{Sp\acute{e}t}(-), \mathfrak{X})$ denote the fully faithful Yoneda embedding of the category of spectral Deligne-Mumford stacks in the category of functors $\mathfrak{Fun}(\mathcal{CA}\text{lg}^{\text{cn,op}}, \mathsf{Sp})$. Then there is a canonical equivalence of $\infty$-categories $\mathsf{QCoh}(\mathfrak{X}) \simeq \mathsf{QCoh}(h_{\mathfrak{X}})$, where $\mathsf{QCoh}(\mathfrak{X})$ is the category of quasi-coherent sheaves on $\mathfrak{X}$ and $\mathsf{QCoh}(h_{\mathfrak{X}})$ is the essential image of $\mathsf{QCoh}(\mathfrak{X})$ by the right Kan extension of the fully faithfull Yoneda embedding.*

### 2.2.5 Proper morphisms and direct image functor.

In our purpose, we will also need a result about the direct image (or pushforward) functor $\varphi_* : \mathsf{QCoh}(\mathfrak{X}) \to \mathsf{QCoh}(\mathfrak{Y})$ for two spectral Deligne-Mumford stacks $\mathfrak{X}$ and $\mathfrak{Y}$ .

Let us recall that if $A$ is an $\mathbb{E}_\infty$-ring and $M$ an $A$-module, then:

- $M$ is a *compact* object of $\mathrm{Mod}_A$ if mapping out of it preserves filtered colimits.
- $M$ is a *perfect* object of $\mathrm{Mod}_A$ if it is dualizable as a chain complex. As a chain complex of $\mathscr{O}_{\mathfrak{X}}$-modules, it is called perfect if it is locally quasi-isomorphic to a bounded complex of free $\mathscr{O}_{\mathfrak{X}}$-modules of finite type.

REMARK 2.2.100 : The distinction between *perfect* objects and *almost perfect* objects, made in Lurie's work [103], is not necessay here, since the spectral algebraic spaces involved here are of finite cohomological dimension.

DEFINITION 2.2.101 ([103], Def. 5.1.2.1) : *Let $\varphi : \mathfrak{X} \to \mathfrak{Y}$ be a map between two spectral Deligne-Mumford stacks. We say that $\varphi$ is* proper *if it is quasi-compact, locally of finite type, and universally closed.*

REMARK 2.2.102 : We bring attention to the fact that being a proper morphism between spectral Deligne-Mumford stacks uniquely depends of 0-truncations. So for spectral algebraic spaces and spectral schemes, we only look at the properness of $(\mathfrak{X}, \pi_0\mathscr{O}_{\mathfrak{X}}) \to (\mathfrak{Y}, \pi_0\mathscr{O}_{\mathfrak{Y}})$, assuming that the diagonal of the proper morphism $\varphi : \mathfrak{X} \to \mathfrak{Y}$ is closed.

We have the easy case for our purpose here which we reuse later:

EXAMPLE 2.2.103 ([103], Ex. 5.1.2.4) : A closed immersion of spectral Deligne-Mumford stacks is proper.

And we have:

PROPOSITION 2.2.104 ([103], Prop. 5.6.0.2) : *Let $\varphi : \mathfrak{X} \to \mathfrak{Y}$ be a morphism of spectral Deligne-Mumford stacks of finite cohomological dimensions. Suppose $\varphi$* *is proper and locally almost of finite presentation. Then the pushforward functor* $\varphi_* : \mathsf{QCoh}(\mathfrak{X}) \to \mathsf{QCoh}(\mathfrak{Y})$ *carries perfect (dualizable) objects of* $\mathsf{QCoh}(\mathfrak{X})$ *to perfect (dualizable) objects of* $\mathsf{QCoh}(\mathfrak{Y})$*.*

REMARK 2.2.105 : In our case, the perverse sheaves will be of finite cohomological dimension, and they will be perfect objects as finite extensions of simple objects in the categories of perverse sheaves on the $\infty$-topoi of the TNN Grassmannian and of the Amplituhedron. We will also note that subsubsection 2.6.6.1 there are enough perfect perverse sheaves to cover the $\infty$-perverse geometry.

The main proposition summarizing this section is the following:

PROPOSITION 2.2.106 ([103] Prop. 5.6.1.1) : *Let* $f : \mathfrak{X} \to \mathfrak{Y}$ *be a morphism of spectral Deligne-Mumford stacks which is finite and of finite presentation. Then an object* $\mathscr{F} \in \mathsf{QCoh}(\mathfrak{X})$ *is (compact) perfect (dualizable) if and only if its direct image* $f_*\mathscr{F} \in \mathsf{QCoh}(\mathfrak{Y})$ *is (compact) perfect (dualizable).*

### 2.2.6 Tannaka duality and reconstruction of a coaffine geometric stack.

DEFINITION 2.2.107 : *Let* $\mathfrak{X}$ *and* $\mathfrak{Y}$ *be spectral Deligne-Mumford stacks. A map is said to be* affine *if, for any map* $\mathsf{Sp\acute{e}t} R \to \mathfrak{X}$*, the pullback* $\mathsf{Sp\acute{e}t}\ R \times_{\mathfrak{X}} \mathfrak{Y}$ *is affine, that is of the form* $\mathsf{Sp\acute{e}t}\ A$ *for some* $\mathbb{E}_\infty$*-ring* $A$*.*

DEFINITION 2.2.108 ([103] Def. 9.1.0.1 and Def. 9.3.0.1) : *We call a* geometric stack *a functor* $\mathfrak{X} : \mathcal{CA}\mathrm{lg}^{\mathrm{cn}} \longrightarrow \widehat{\mathcal{S}}$ *which satisfies the following conditions:*

- *The functor* $\mathfrak{X}$ *satisfies descent with respect to the fpqc topology*
- *The diagonal map* $\delta : \mathfrak{X} \longrightarrow \mathfrak{X} \times \mathfrak{X}$ *is affine*
- *There exists a connective* $\mathbb{E}_\infty$*-ring* $A$ *for which there exists a faithfully flat morphism* $f : \mathsf{Spec}\ \ A \to \mathfrak{X}$

REMARK 2.2.109 : Even if the definition is made for a functor which takes values in the $\infty$-category $\widehat{\mathcal{S}}$ of spaces which are not necessarily small, a set-theoretic technicality allows us to replace the $\infty$-category $\widehat{\mathcal{S}}$ by the $\infty$-category $\mathcal{S}$, since for every connective $\mathbb{E}_\infty$-ring $R$, the space $\mathfrak{X}(R)$ will be essentially small.

EXAMPLE 2.2.110 ([103] Cor. 9.1.4.6) : Any quasi-compact, quasi-separated spectral algebraic space is a geometric stack when identified with its functor of points.

PROPOSITION 2.2.111 ([103] Cor. 9.1.3.1) : *Let* $\mathfrak{X} : \mathcal{CA}\mathrm{lg}^{\mathrm{cn}} \longrightarrow \mathcal{S}$ *be a geometric stack. Then* $\mathsf{QCoh}(\mathfrak{X})^{\mathrm{cn}}$ *is a complete Grothendieck prestable* $\infty$*-category. Moreover, the inclusion* $\mathsf{QCoh}(\mathfrak{X})^{\mathrm{cn}} \hookrightarrow \mathsf{QCoh}(\mathfrak{X})$ *exhibits* $\mathsf{QCoh}(\mathfrak{X})$ *as a stabilization of* $\mathsf{QCoh}(\mathfrak{X})^{\mathrm{cn}}$.

PROPOSITION 2.2.112 ([103] Prop. 9.1.5.3) : *Let* $\mathfrak{X}$ *be a geometric stack. The following conditions are equivalent:*

- *The geometric stack* $\mathfrak{X}$ *has finite cohomological dimension*
- *The structure sheaf* $\mathscr{O}_{\mathfrak{X}} \in \mathsf{QCoh}(\mathfrak{X})$ *is compact*
- *The global sections functor* $\Gamma : \mathsf{QCoh}(\mathfrak{X}) \longrightarrow \mathsf{Sp}$ *commutes with small colimits*
- *Every perfect (dualizable) object* $\mathscr{F} \in \mathsf{QCoh}(\mathfrak{X})$ *is compact*
- *An object* $\mathscr{F} \in \mathsf{QCoh}(\mathfrak{X})$ *is compact if and only if it is perfect*

THEOREM 2.2.113 (Tannaka Duality for Geometric Stacks, [103] Thm. 9.3.0.3) : *Let* $\mathfrak{X}, \mathfrak{Y} : \mathcal{CA}\mathrm{lg}^{\mathrm{cn}} \longrightarrow \mathcal{S}$ *be functors, and suppose that* $\mathfrak{X}$ *is a geometric stack. Then the construction*

$$(f : \mathfrak{Y} \longrightarrow \mathfrak{X}) \longmapsto (f^* : \mathsf{QCoh}(\mathfrak{X}) \longrightarrow \mathsf{QCoh}(\mathfrak{Y}))$$

*determines a fully faithful embedding :*

$$\mathrm{Hom}_{\mathcal{F}\mathrm{un}(\mathcal{CA}\mathrm{lg}^{\mathrm{cn}},\mathcal{S})}(\mathfrak{Y}, \mathfrak{X}) \longrightarrow \mathcal{F}\mathrm{un}^{\otimes}(\mathsf{QCoh}(\mathfrak{X}), \mathsf{QCoh}(\mathfrak{Y})) \tag{2.2.114}$$

*whose essential image is spanned by those symmetric monoidal functors* $F : \mathsf{QCoh}(\mathfrak{X}) \longrightarrow \mathsf{QCoh}(\mathfrak{Y})$ *which satisfy:*

- *The functor* $F$ *preserves small colimits (denoted* $F_{\mathrm{sc}}$*)*

- *The functor $F$ carries connective objects of* $\mathsf{QCoh}(\mathfrak{X})$ *to connective objects of* $\mathsf{QCoh}(\mathfrak{Y})$ *(denoted $F_{\mathrm{co}}$)*

- *The functor $F$ carries flat objects of* $\mathsf{QCoh}(\mathfrak{X})$ *to flat objects of* $\mathsf{QCoh}(\mathfrak{Y})$ *(denoted $F_{\mathrm{fl}}$)*

REMARK 2.2.115 : For further reading on the dual Amplituhedron, it would suffice to state the Tannaka duality for algebraic spaces. But in order to respect Schubert geometry, we need to reconstruct the associated geometries (the intersection sheaves $\mathcal{IC}$). So we need a Tannaka duality for geometric stacks between spectral Deligne-Mumford stacks for tensor $\infty$-categories of functors and quasi-coherent sheaves on these functors.

The following diagram illustrates the Tannaka reconstruction :

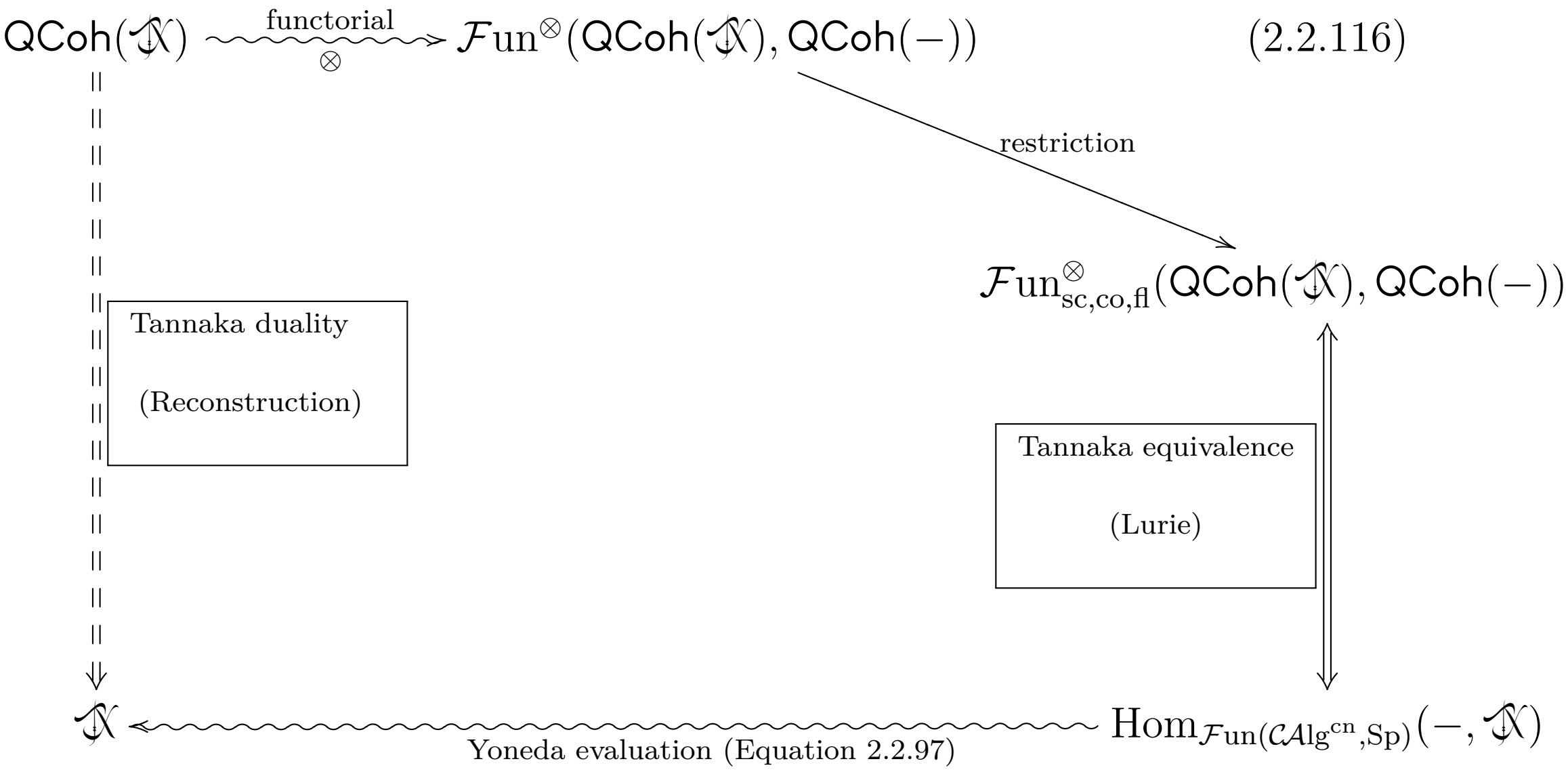

where:

- The map $\mathsf{QCoh}(\mathfrak{X}) \xrightarrow{\text{functorial}} \mathcal{F}\mathrm{un}^{\otimes}(\mathsf{QCoh}(\mathfrak{X}), \mathsf{QCoh}(-))$ is functorial in $\mathfrak{X}$.

- The theorem 2.2.113 invites to invoke the restriction functor with respect to the essential image in Equation 2.2.114:

$$\mathcal{F}\text{un}^{\otimes}(\mathsf{QCoh}(\mathfrak{X}), \mathsf{QCoh}(-)) \xrightarrow{\text{restriction}} \mathcal{F}\text{un}^{\otimes}_{\text{sc,co,fl}}(\mathsf{QCoh}(\mathfrak{X}), \mathsf{QCoh}(-))$$

- The fully faithful embedding of Equation 2.2.114 in theorem 2.2.113 :

$$\text{Hom}_{\mathcal{F}\text{un}(\mathcal{C}\mathcal{A}\text{lg}^{\text{cn}}, \mathcal{S})}(\mathfrak{Y}, \mathfrak{X}) \longmapsto \mathcal{F}\text{un}^{\otimes}(\mathsf{QCoh}(\mathfrak{X}), \mathsf{QCoh}(\mathfrak{Y}))$$

  composed after the fully faithful Yoneda embedding:

$$h_{\mathfrak{X}} : \mathfrak{Y} \longmapsto \text{Hom}_{\mathcal{F}\text{un}(\mathcal{C}\mathcal{A}\text{lg}^{\text{cn}}, \mathcal{S})}(\mathfrak{Y}, \mathfrak{X})$$

  gives the *Tannaka equivalence* (theorem 2.2.113)

- The last functor at the bottom is the Yoneda evaluation (Equation 2.2.97)

REMARK 2.2.116 : For a further analysis of this phenomenon of *up and down shift*, that is:

- the *categorification* : $\mathfrak{X} \longmapsto \mathsf{QCoh}(\mathfrak{X})$
- the *decategorification* : $\mathsf{QCoh}(\mathfrak{X}) \longmapsto \mathfrak{X}$

we notice that we can pass from a geometric stack to the category of quasi-coherent sheaves on it, with a way back by Tannaka reconstruction without ambiguity. In the case which concerns us here, the 0-affineness and the 1-affineness of the geometric stacks will play a central role, which we see in the next **??** below.

## 2.3. Lurie's geometries and structured spaces.

Let us freely use the deep construction of J. Lurie's [91] and [93]. No consideration gets out of his vision of structured spaces.

### 2.3.1 Preliminaries on Ind and Pro objects and presheaves categories.

DEFINITION 2.3.1 : *We define the* ind-objects *of a category $\mathcal{C}$ as the formal filtered colimit of objects of $\mathcal{C}$, that is the colimit is taken in the category of presheaves of $\mathcal{C}$ with values in the category $\mathcal{S}$ of spaces, seen as the free cocompletion of $\mathcal{C}$. The category of ind-objects of $\mathcal{C}$ is written Ind($\mathcal{C}$).*

REMARK 2.3.2 : One definition is to define the objects of $\mathrm{Ind}(\mathcal{C})$ to be diagrams $\mathcal{F} : D \to \mathcal{C}$ where $D$ is a small filtered category. Ind-objects need not belong to the category $\mathcal{C}$. An ordinary object of $\mathcal{C}$ is the corresponding diagram $1 \to \mathcal{C}$, and the morphisms between $\mathcal{F} : D \to \mathcal{C}$ and $\mathcal{G} : E \to \mathcal{C}$ should be such that:

- the embedding $\mathcal{C} \longrightarrow \mathrm{Ind}(\mathcal{C})$ is full and faithful
- each diagram $\mathcal{F} : D \to \mathcal{C}$ is the colimit of itself in $\mathrm{Ind}(\mathcal{C})$
- the objects of $\mathcal{C}$ are compact in $\mathrm{Ind}(\mathcal{C})$

Then $\mathrm{Ind}(\mathcal{C})(\mathcal{F}, \mathcal{G}) = \lim\limits_{d \in D} \operatorname*{colim}\limits_{e \in E} \mathrm{Hom}_{\mathcal{C}}(\mathcal{F}d, \mathcal{G}e)$

By the co-Yoneda lemma, every presheaf $X \in \mathrm{PSh}(\mathcal{C}, \mathcal{S})$ is a colimit over representable presheaves, that is there is a functor $F : D \to \mathcal{C}$ (with the diagram $D$ possibly large) such that:

$$X \simeq \operatorname*{colim}_{d \in D} y(F(d))$$

where $y$ is the Yoneda embedding.

PROPOSITION 2.3.3 : *The category* $\mathrm{Ind}(\mathcal{C}) \subseteq \mathrm{PSh}(\mathcal{C}, \mathcal{S})$ *is the full subcategory of the presheaf category* $\mathrm{PSh}(\mathcal{C}, \mathcal{S}) =: [\mathcal{C}^{\mathrm{op}}, \mathcal{S}]$ *on those presheaves functors which are filtered colimits of representables, i.e. those for which*

$$X \simeq \operatorname*{colim}_{d \in D} y(F(d))$$

*for $D$ some filtered category.*

REMARK 2.3.4 : Given that the category of presheaves $\mathrm{PSh}(\mathcal{C},\mathcal{S}) =: [\mathcal{C}^{\mathrm{op}},\mathcal{S}]$ is the *free cocompletion* of $\mathcal{C}$, then $\mathrm{Ind}(\mathcal{C})$ is the *free cocompletion under filtered colimits*.

DEFINITION 2.3.5 : *The category* $\mathrm{Pro}(\mathcal{C})$ *of pro-objects of* $\mathcal{C}$ *is defined such that* $\mathrm{Pro}(\mathcal{C}) := \mathrm{Ind}(\mathcal{C}^{\mathrm{op}})^{\mathrm{op}}$.

REMARK 2.3.6 : Given that the category of presheaves $\mathrm{PSh}(\mathcal{C},\mathcal{S}) =: [\mathcal{C}^{\mathrm{op}},\mathcal{S}]$ is the *free cocompletion* of $\mathcal{C}$, then $\mathrm{Pro}(\mathcal{C})$ is the *free completion under cofiltered limits*, that is equivalently *filtered colimits of presheaves*.

PROPOSITION 2.3.7 : *If the category* $\mathcal{C}$ *admits finite limits, then* $\mathrm{Pro}(\mathcal{C}) = \mathcal{F}\mathrm{un}^{\mathrm{lex}}(\mathcal{C},\mathcal{S})^{\mathrm{op}}$.

PROPOSITION 2.3.8 : *If the category* $\mathcal{C}$ *is small, then the* $\infty$*-category* $\mathrm{Ind}(\mathcal{C})$ *admits small filtered colimits, and the* $\infty$*-category* $\mathrm{Pro}(\mathcal{C})$ *admits small filtered limits.*

DEFINITION 2.3.9 ([98], Def. 5.2.6.1) : *Let* $\mathcal{C}$ *and* $\mathcal{D}$ *be two* $\infty$*-categories. We let* $\mathcal{F}\mathrm{un}^{\mathrm{L}}(\mathcal{C},\mathcal{D})$ *denote the full subcategory of* $\mathcal{F}\mathrm{un}(\mathcal{C},\mathcal{D})$ *spanned by those functors which admit right adjoints, and* $\mathcal{F}\mathrm{un}^{\mathrm{R}}(\mathcal{C},\mathcal{D}) \subseteq \mathcal{F}\mathrm{un}(\mathcal{C},\mathcal{D})$ *the full subcategory spanned by those functors which admit left adjoints.*

PROPOSITION 2.3.10 ([98], Prop. 5.2.6.2) : *We have canonical equivalences:*

$$\mathcal{F}\mathrm{un}^{\mathrm{L}}(\mathcal{C},\mathcal{D}) \simeq \mathcal{F}\mathrm{un}^{\mathrm{L}}(\mathcal{D},\mathcal{C})^{\mathrm{op}} \simeq \mathcal{F}\mathrm{un}^{\mathrm{R}}(\mathcal{D}^{\mathrm{op}},\mathcal{C}^{\mathrm{op}}) \simeq \mathcal{F}\mathrm{un}^{\mathrm{R}}(\mathcal{C},\mathcal{D})$$

### 2.3.2 Sheaves on an $\infty$-category with values in an $\infty$-category.

Let us recall the following:

DEFINITION 2.3.11 (see Equation 2.7.7) : *Let* $\mathscr{C}$ *be an* $\infty$*-category, and* $\mathfrak{X}$ *an* $\infty$*-topos. The* category of sheaves on the $\infty$-topos $\mathfrak{X}$, *denoted* $\boldsymbol{\mathcal{S}}\mathbf{hv}_{\mathscr{C}}(\mathfrak{X})$, *is the category of* $\mathscr{C}$*-valued presheaves on* $\mathfrak{X}$ *which preserves limits on* $\mathfrak{X}^{\mathrm{op}}$ *to* $\mathscr{C}$*. We write:*

$$\boldsymbol{\mathcal{S}}\mathbf{hv}_{\mathscr{C}}(\mathfrak{X}) := \mathfrak{Fun}^{\lim}(\mathfrak{X}^{\mathrm{op}},\mathscr{C}) \tag{2.3.12}$$

But we can state the following:

PROPOSITION 2.3.13 ([103], Prop. 1.3.1.7) : *Let* $\mathcal{T}$ *be a small* $\infty$*-category equipped with a Grothendieck topology. Let* $j : \mathcal{T} \longrightarrow \mathfrak{Fun}(\mathcal{T}^{\mathrm{op}},\mathcal{S})$ *denote the Yoneda embedding and* $L : \mathfrak{Fun}(\mathcal{T}^{\mathrm{op}}, \mathcal{S}) \longrightarrow \boldsymbol{\mathcal{S}}\mathbf{hv}_{\mathcal{S}}(\mathcal{T})$ *a left adjoint to this inclusion* $j$*. Let* $\mathcal{C}$ *be an* $\infty$*-category which admits small colimits. Then a* $\mathcal{C}$*-valued sheaf on* $\mathcal{T}$ *is equivalent to a* $\mathcal{C}$*-valued sheaf on* $\boldsymbol{\mathcal{S}}\mathbf{hv}_{\mathcal{S}}(\mathcal{T})$*, that is composition with* $L \circ j$ *induces an equivalence of* $\mathcal{C}$*-valued sheaves* $\boldsymbol{\mathcal{S}}\mathbf{hv}_{\mathcal{C}}(\boldsymbol{\mathcal{S}}\mathbf{hv}_{\mathcal{S}}(\mathcal{T})) \simeq \boldsymbol{\mathcal{S}}\mathbf{hv}_{\mathcal{C}}(\mathcal{T})$.

We can then state the following:

DEFINITION 2.3.14 : *Let* $\mathcal{C}$ *be an* $\infty$*-category and let* $\mathfrak{X}$ *be an* $\infty$*-topos. A* $\mathcal{C}$-valued sheaf on $\mathfrak{X}$ *is a functor* $\mathfrak{X}^{\mathrm{op}} \to \mathcal{C}$ *which preserves limits. We denote* $\boldsymbol{\mathcal{S}}\mathbf{hv}_{\mathcal{C}}(\mathfrak{X})$ *this full sub-*$\infty$*-category of* $\mathcal{F}\mathrm{un}(\mathfrak{X}^{\mathrm{op}}, \mathcal{C})$.

REMARK 2.3.15 : The usual sheaf condition of compatibility on a Grothendieck ($\infty$-)site is replaced by the preservation of limits. It is a simple rewriting of the same data. (see remark 2.2.59).

REMARK 2.3.16 ([91],Rem. 1.1.6) : Let $\mathcal{G}$ be a small $\infty$-category which admits finite limits, and let $y : \mathcal{G} \to \mathrm{Pro}(\mathcal{G})$ denote the Yoneda embedding. Let $\mathfrak{X}$ be an $\infty$-topos. Then we have the following equivalences:

$$\mathcal{F}\mathrm{un}^{\mathrm{R}}(\mathrm{Pro}(\mathcal{G}), \mathfrak{X}) \simeq \mathcal{F}\mathrm{un}^{\mathrm{lex}}(\mathcal{G}, \mathfrak{X}) \simeq \mathcal{F}\mathrm{un}^{\mathrm{R}}(\mathfrak{X}, \mathrm{Ind}(\mathcal{G}^{\mathrm{op}})) \simeq \boldsymbol{\mathcal{S}}\mathbf{hv}_{\mathrm{Ind}(\mathcal{G}^{\mathrm{op}})}(\mathfrak{X})$$

REMARK 2.3.17 : The equivalences above mean that we can view objects of these equivalent $\infty$-categories from two points of view. On the one hand, they are sheaves on the $\infty$-topos $\mathfrak{X}$ with values in $\mathrm{Ind}(\mathcal{G}^{\mathrm{op}})$, that is we send the "*opens*" of $\mathfrak{X}$ to colimits of (maybe for example compact) objects of $\mathcal{G}$ (if $\mathcal{G}$ is compactly generated). On the other hand, these objects "*are*" (embedded as) objects of $\mathfrak{X}$ with the additional structure given by $\mathcal{G}$.

### 2.3.3 Admissibility and geometries.

DEFINITION 2.3.18 ([91], Def. 1.2.1) : *Let* $\mathcal{G}$ *be an* $\infty$*-category. An* admissibility structure *on* $\mathcal{G}$ *consists of the following data:*

- *A subcategory* $\mathcal{G}^{\mathrm{ad}} \subseteq \mathcal{G}$*, containing every object of* $\mathcal{G}$*, and whose morphisms belonging to* $\mathcal{G}^{\mathrm{ad}}$ *will be called* admissible *morphisms in* $\mathcal{G}$.

- *A Grothendieck topology on $\mathcal{G}$, which is* generated *by admissible morphisms in the following sense: let $\mathcal{G}^{(0)}_{/X} \subseteq \mathcal{G}_{/X}$ be a covering sieve on an object $X \in \mathcal{G}$. Then $\mathcal{G}^{(0)}_{/X}$ contains a covering sieve which is generated by a collection of admissible morphisms $\{U_\alpha \to X\}$.*

*Moreover, these data are required to satisfy the following conditions:*

*i) Let $f : U \to X$ be an admissible morphism in $\mathcal{G}$, and $g : X' \to X$ any morphims. Then there exists a pullback diagram*

$$\begin{array}{ccc} U' & \longrightarrow & U \\ \downarrow{\scriptstyle f'} & \lrcorner & \downarrow{\scriptstyle f} \\ X' & \xrightarrow{g} & X \end{array} \tag{2.3.19}$$

*where $f'$ is admissible.*

*ii) Suppose given a commutative diagram:*

$$\begin{array}{ccccc} & & Y & & \\ & \overset{f}{\nearrow} & & \overset{g}{\searrow} & \\ X & & \xrightarrow{h} & & Z \end{array} \tag{2.3.20}$$

*in $\mathcal{G}$, where $g$ and $h$ are admissible. Then $f$ is admissible.*

*iii) Every retract of an admissible morphim of $\mathcal{G}$ is admissible.*

REMARK 2.3.21 : Every admissibility structure $\mathcal{G}^{\mathrm{ad}} \subseteq \mathcal{G}$ determines a Grothendieck topology on $\mathcal{G}^{\mathrm{ad}}$, namely the collection of all sieves on objects $X \in \mathcal{G}^{\mathrm{ad}}$ which generate covering sieves in $\mathcal{G}$. The condition that every covering sieve in $\mathcal{G}$ be generated by admissible morphisms guarantees that the topology on $\mathcal{G}$ is generated by the topology on $\mathcal{G}^{\mathrm{ad}}$. However, the definition definition 2.3.18 says more than the Grothendieck topology on $\mathcal{G}$ induces a Grothendieck topology on $\mathcal{G}^{\mathrm{ad}}$, since a condition would have to be added.

DEFINITION 2.3.22 ([91], Def. 1.2.5) : *A* geometry *consists of the following data:*

- *An essentially small $\infty$-category $\mathcal{G}$ which admits finite limits and is idempotent complete.*

- *An admissibility structure on $\mathcal{G}$.*

### 2.3.4 Transformations and structures.

DEFINITION 2.3.23 ([91], Def. 1.2.6) : *A* transformation of geometries $f : \mathcal{G} \to \mathcal{G}'$ *is such that:*

- *the functor $f$ preserve finite limits*

- *the functor $f$ carries admissible morphisms of $\mathcal{G}$ to admissible morphisms of $\mathcal{G}'$.*

- *for every admissible cover $\{U_\alpha \to X\}_\alpha$ of an object $X \in \mathcal{G}$, the collection of morphisms $\{f(U_\alpha) \to f(X)\}_\alpha$ is an admissible cover of $f(X) \in \mathcal{G}'$.*

DEFINITION 2.3.24 ([91], Def. 1.2.8) : *Let $\mathcal{G}$ be a geometry and $\mathfrak{X}$ an $\infty$-topos. A $\mathcal{G}$-*structure *on $\mathfrak{X}$ is a left exact functor $\mathscr{O} : \mathcal{G} \to \mathfrak{X}$ with the following property: for every collection of admissible morphisms $\{U_\alpha \to X\}_\alpha$ in $\mathcal{G}$ which generates a covering sieve on $X$, the induced map $\coprod_\alpha \mathscr{O}(U_\alpha) \to \mathscr{O}(X)$ is an effective epimorphism in $\mathfrak{X}$. Let us denote $\mathcal{S}\mathrm{tr}_{\mathcal{G}}(\mathfrak{X})$ the full subcategory of $\mathcal{F}\mathrm{un}(\mathcal{G}, \mathfrak{X})$ spanned by the $\mathcal{G}$-structures on $\mathfrak{X}$.*

DEFINITION 2.3.25 ([91], Def. 1.2.8) : *Given a pair of $\mathcal{G}$-structures $\mathscr{O}, \mathscr{O}' : \mathcal{G} \to \mathfrak{X}$ on $\mathfrak{X}$, we say that a natural transformation $\varphi : \mathscr{O} \to \mathscr{O}'$ is a* local transformation of

$\mathcal{G}$-structures *if, for every admissible morphism $U \to X$ in $\mathcal{G}$, the induced diagram:*

$$\begin{array}{ccc} \mathscr{O}(U) & \to & \mathscr{O}'(U) \\ \downarrow & \lrcorner & \downarrow \\ \mathscr{O}(X) & \to & \mathscr{O}'(X) \end{array} \tag{2.3.26}$$

*is a pullback square in $\mathfrak{X}$. Let us denote $\mathcal{S}\mathrm{tr}^{\mathrm{loc}}_{\mathcal{G}}(\mathfrak{X})$ the subcategory of $\mathcal{S}\mathrm{tr}_{\mathcal{G}}(\mathfrak{X})$ spanned by the local transformations of $\mathcal{G}$-structures.*

DEFINITION 2.3.27 ([91], Def. 1.2.10) : *A geometry $\mathcal{G}$ is* discrete *if:*

- *the admissible morphisms in $\mathcal{G}$ are the equivalences*
- *the Grothendieck topology on $\mathcal{G}$ is trivial: that is, a sieve $\mathcal{G}^{(0)}_{/X} \subseteq \mathcal{G}_{/X}$ on an object $X \in \mathcal{G}$ is a covering sieve if and only if $\mathcal{G}^{(0)}_{/X} = \mathcal{G}_{/X}$.*

REMARK 2.3.28 : If $\mathcal{G}$ is an essentially small $\infty$-category that admits finite limits, then the two conditions of the previous definition endow $\mathcal{G}$ with the structure of a discrete geometry.

If $\mathcal{G}$ is a discrete geometry and $\mathfrak{X}$ is an $\infty$-topos, then we have the equivalences:

$$\mathcal{S}\mathrm{tr}^{\mathrm{loc}}_{\mathcal{G}}(\mathfrak{X}) = \mathcal{S}\mathrm{tr}_{\mathcal{G}}(\mathfrak{X}) = \mathcal{F}\mathrm{un}^{\mathrm{lex}}(\mathcal{G}, \mathfrak{X}) \simeq \boldsymbol{\mathcal{S}}\mathbf{hv}_{\mathrm{Ind}(\mathcal{G}^{\mathrm{op}})}(\mathfrak{X}) \tag{2.3.29}$$

where $\mathrm{Ind}(\mathcal{C})$ denotes the Ind category of a category $\mathcal{C}$. These equivalences mean exactly that the theory of $\mathcal{G}$-structures is equivalent to the theory of $\mathrm{Ind}(\mathcal{G}^{\mathrm{op}})$-valued sheaves.

### 2.3.5 Sober spaces and recovering of a space from its category of sheaves.

If a space $X$ is sober, that is if every irreducible closed subset of $X$ has a unique generic point, we can recover $X$ from the $\infty$-topos defined by the $\infty$-category $\boldsymbol{\mathcal{S}}\mathbf{hv}_{\mathcal{S}}(X)$ of sheaves of spaces on $X$ thanks to the following process. The points $x \in X$ can be identified with isomorphism classes of geometric morphisms $x^* : \mathcal{S}\mathbf{hv}_{\mathcal{S}}(X) \to \mathcal{S}$. The open subsets of $X$ are the subobjects of the unit (terminal) object $\mathbf{1}_X \in \mathcal{S}\mathbf{hv}_{\mathcal{S}}(X)$, where $\mathbf{1}_X \simeq \mathscr{O}_X$, the canonical sheaf on $X$. The space $X$ and the $\infty$-topos $\mathcal{S}\mathbf{hv}_{\mathcal{S}}(X)$ are then canonically interchangeable. The theory of $\infty$-topoi is then a generalization of the theory of classical (sober) topological spaces.

remark 2.3.30 : In particular, when we look at Whitney stratifications of the space $X$, and even more of the TNN Grassmannian, we have a characterization of each stratum in terms of sheaves on it, and we will be able to recover the strata from the local systems on them, as well as their duals.

Every $\infty$-topos $\mathfrak{X}$ represents a functor from topological spaces to $\infty$-categories:

$$\begin{cases} \mathcal{T}\mathrm{op}_{\mathrm{sober}} & \longrightarrow & \infty\text{-Cat} \\ X & \longmapsto & \mathcal{F}\mathrm{un}^*(\mathfrak{X}, \mathcal{S}\mathbf{hv}_{\mathcal{S}}(X)) \end{cases}$$

where $\mathcal{F}\mathrm{un}^*(\mathcal{X}, \mathcal{Y})$ between two topoi $\mathcal{X}$ and $\mathcal{Y}$ denotes the geometric morphisms of topoi (that is the pair of adjoint $\infty$-functors $(f_*, f^*)$ between topoi $\mathcal{X} \underset{f_*}{\overset{f^*}{\rightleftarrows}} \mathcal{Y}$ where the inverse image functor $f^*$ is a left exact $\infty$-functor). In general, these functors are not representable in the category of spaces (not all the $\infty$-topoi come from an $\infty$-category of sheaves on a topological space $X$ endowed with a Grothendieck topology, that is, an $\infty$-site on it).

### 2.3.6 Module geometries.

definition 2.3.31 : *Let $\mathfrak{X}$ be an $\infty$-topos and let $\mathscr{O}_{\mathfrak{X}} \in \mathcal{S}\mathbf{hv}_{\mathcal{CA}\mathbf{lg}}(\mathfrak{X})$ be a sheaf of $\mathbb{E}_\infty$-rings on $\mathfrak{X}$. Recall that $\mathscr{O}_{\mathfrak{X}}$ can be identified with a commutative algebra object of the symmetric monoidal $\infty$-category $\mathcal{S}\mathbf{hv}_{\mathcal{S}\mathrm{p}}(\mathfrak{X})$ of sheaves of spectra on $\mathfrak{X}$ (see lemma 2.2.96). We let $\mathcal{M}\mathbf{od}_{\mathscr{O}_{\mathfrak{X}}}$ denote the $\infty$-category $\mathcal{M}\mathbf{od}_{\mathscr{O}_{\mathfrak{X}}}(\mathcal{S}\mathbf{hv}_{\mathcal{S}\mathrm{p}}(\mathscr{O}_{\mathfrak{X}}))$ of $\mathscr{O}_{\mathfrak{X}}$-module objects of $\mathcal{S}\mathbf{hv}_{\mathcal{S}\mathrm{p}}(\mathfrak{X})$. Then $\mathcal{M}\mathbf{od}_{\mathscr{O}_{\mathfrak{X}}}$ can be regarded as a symmetric monoidal*

*$\infty$-category with respect to the relative tensor product $- \otimes_{\mathscr{O}_{\mathfrak{X}}} -$ (see [100], Appendix A.3.4.4). We will refer to the objects of $\boldsymbol{\mathcal{M}}\mathbf{od}_{\mathscr{O}_{\mathfrak{X}}}$ as* sheaves of $\mathscr{O}_{\mathfrak{X}}$-modules on $\mathfrak{X}$.

Let $(\mathfrak{X}, \mathscr{O}_{\mathfrak{X}})$ be a spectrally ringed topos. Let us consider the $\infty$-category, denoted $\boldsymbol{\mathcal{M}}\mathbf{od}$ of pairs $(A, M)$ where $A$ is an $\mathbb{E}_\infty$-ring and $M$ is an $A$-module spectrum. Thanks to [93], Remark 2.2.1, we have an isomorphism of categories:

$$\boldsymbol{\mathcal{S}}\mathbf{hv}_{\boldsymbol{\mathcal{M}}\mathbf{od}}(\mathfrak{X}) \simeq \boldsymbol{\mathcal{M}}\mathbf{od}(\boldsymbol{\mathcal{S}}\mathbf{hv}_{\mathcal{S}\mathrm{p}}(\mathfrak{X})) \tag{2.3.32}$$

We can now identify a $\boldsymbol{\mathcal{M}}\mathbf{od}$-valued sheaf on $\mathfrak{X}$ with a pair $(\mathscr{O}_{\mathfrak{X}}, \mathscr{F}_{\mathfrak{X}})$, where $\mathscr{O}_{\mathfrak{X}}$ is a sheaf of $\mathbb{E}_\infty$-rings on $\mathscr{O}_{\mathfrak{X}}$ and $\mathscr{F}_{\mathfrak{X}}$ is a $\mathscr{O}_{\mathfrak{X}}$-module. The $\infty$-category $(\mathfrak{X}, \boldsymbol{\mathcal{M}}\mathbf{od}_{\mathscr{O}_{\mathfrak{X}}})$, that is the spectrally ringed topos $(\mathfrak{X}, \boldsymbol{\mathcal{M}}\mathbf{od}_{\mathscr{O}_{\mathfrak{X}}}(\boldsymbol{\mathcal{S}}\mathbf{hv}_{\mathcal{S}\mathrm{p}}(\mathscr{O}_{\mathfrak{X}})))$ of $\mathfrak{X}$ endowed with the $\mathscr{O}_{\mathfrak{X}}$-module objects of $\boldsymbol{\mathcal{S}}\mathbf{hv}_{\mathcal{S}\mathrm{p}}(\mathfrak{X})$ can now be identified with triples $(\mathfrak{X}, \mathscr{O}_{\mathfrak{X}}, \mathscr{F}_{\mathfrak{X}})$, where $\mathfrak{X}$ is an $\infty$-topos, $\mathscr{O}_{\mathfrak{X}}$ is a sheaf of $\mathbb{E}_\infty$-rings on $\mathscr{O}_{\mathfrak{X}}$ and $\mathscr{F}_{\mathfrak{X}}$ is a $\mathscr{O}_{\mathfrak{X}}$-module.

THEOREM 2.3.33 : *Considering an $\infty$-subcategory of $\infty$-topoï $\mathfrak{X}$ with one geometry $\mathcal{G}$ (called a "*Lurie's geometry*") endowing then a $\mathcal{G}$-*structure *on each $\mathfrak{X}$, that is a left exact functor $\mathscr{O} : \mathcal{G} \to \mathfrak{X}$, that defines an admissibility structure made of a well-defined subcategory $\mathcal{G}^{\mathrm{ad}}$ of admissible morphisms and a factorization system, is then exactly the same as considering an $\infty$-category of triples $(\mathfrak{X}, \mathscr{O}_{\mathfrak{X}}, \mathscr{F}_{\mathfrak{X}})$, where $\mathfrak{X}$ is an $\infty$-topos, $\mathscr{O}_{\mathfrak{X}}$ is a sheaf of $\mathbb{E}_\infty$-rings on $\mathscr{O}_{\mathfrak{X}}$ and $\mathscr{F}_{\mathfrak{X}}$ is a $\mathscr{O}_{\mathfrak{X}}$-module.*

We have two immediate applications given by J. Lurie of this theorem. The first one concerns the Zariski structure, and the second one the étale structure. In the former case, $\mathcal{G}^{\mathrm{ad}}$ of admissible morphisms is the $\infty$-category of affine charts modelled on $R[x^{-1}]$, which give the Zarisky topology and the Zariski spectrum. The later one concerns the étale topology, wehe admissible morphisms $\mathcal{G}^{\mathrm{ad}}$ are made of étale morphisms, leading to the étale topology of Deligne-Mumford stacks, and to the étale spectrum. We refer to [91] and [92] for a complete treatment of these $\mathcal{G}$-schemes, which are generalizations of Grothendieck's Zariski and étale topologies, spectra and structures of $\mathcal{G}$-schemes.

### 2.3.7 A remainder of our goal for a new positive Schubert geometry.

However, we are concerned here by a third geometry : the one where the sheaves of rings $\mathscr{F}_{\mathfrak{X}}$ of $\mathscr{O}_{\mathfrak{X}}$-modules are the intersection complexes $\mathcal{IC}^{\bullet}$ which are (perfect and then compact and dualizable) perverse $\mathscr{O}_{\mathfrak{X}}$-modules, where there are enough perfect perverse sheaves subsubsection 2.6.6.1 to cover the $\infty$-perverse geometry.

We can then right now state what will be our playgroung : a new positive geometry:

$$\mathfrak{I}(\mathfrak{X})_{\mathrm{dRh}} := \mathfrak{I}\Big(\mathfrak{Gr}^{\geqslant 0}_{\text{ét}}, \mathscr{O}_{\mathfrak{Gr}^{\geqslant 0}_{\text{ét}}}, \mathcal{IC}^{\bullet}_{\mathfrak{Gr}^{\geqslant 0}_{\text{ét}}}\Big)_{/\,\mathrm{dRh}}$$

Indeed, we explain now a little this resume et why it will occupy the remainder of this paper:

The problem of Schubert intersections can be reformulated in a Lurie's formal moduli problem $\mathscr{F}$ since we can adjoin a $\mathcal{G}^{\mathrm{perv}}_{\mathrm{perf}}$-structure and a Lurie's module geometry $(\mathscr{O}_{\mathfrak{G}}, \mathscr{F}_{\mathfrak{G}})$ made of perverse intersection complexes in a symmetric monoidal $\infty$-category of quasi-coherent spectral deligne-Mumford stacks $\mathsf{QStk}^{\mathrm{Pst}}_{\mathrm{comp}}(\mathfrak{X})$ in complete Grothendieck prestable $\infty$-categories on $\mathfrak{X}$ as stable triangulated glued $t$-structures $(\mathsf{QStk}(\mathfrak{X})^{\geqslant 0}, \mathsf{QStk}(\mathfrak{X})^{\leqslant 0})$, whose points are infinitesimally thickened in fat points following Lawvere's and Schreiber's cohesiveness for topoï in a triple adjoint derived Grothendieck context of modalities $\mathfrak{R}(\mathfrak{X}) \dashv \mathfrak{I}(\mathfrak{X}) \dashv \&(\mathfrak{X})$ where Gaitsgory's Tannaka 1- affineness and Verdier's autoduality of the perverse category $\mathcal{P}\mathrm{erv}(\mathfrak{X}) \subset \mathsf{QStk}(\mathfrak{X})^{\heartsuit}$ in the twisted étale Brauer-Grothendieck cohomology $H^2_{\text{ét}}(\mathfrak{X}, \mathcal{IC}^{\bullet\times})$, fully dualizable in the context of our $\mathcal{BCFW}_{n,4}$ positroid vaieties $\mathring{\Pi}_f$ lead to a logical homotopy theoretical inifinitesimal shape modal De Rham space $\mathfrak{I}(\mathfrak{X})_{\mathrm{dRh}}$ which has a dual and a volume, which was expected from the beginning.

Let us now explain this with more details.

### 2.3.8 Nilcompleteness, integrability and infinitesimally cohesiveness.

DEFINITION 2.3.34 ([96], Def. 0.0.4) : *We say that a commutative ring is a* Grothendieck ring *if it is Noetherian, and for every prime ideal* $\mathfrak{p} \subseteq A$*, the map* $A_{\mathfrak{p}} \to \hat{A}$ *is geomet-* *rically regular, where* $\hat{A}$ *denotes the completion of* $A_{\mathfrak{p}}$ *with respect to its maximal ideal.*

DEFINITION 2.3.35 ([96], Def. 2.1.1) : *Let* $X : \mathcal{CA}\mathrm{lg}^{\mathrm{cn}} \to \mathcal{S}$ *be a functor. We will say that* $X$ *is* cohesive *if the following condition is satisfied:*

- *For every pullback diagram*

$$\begin{array}{ccc} A' & \longrightarrow & A \\ \downarrow & \lrcorner & \downarrow \\ B' & \longrightarrow & B \end{array} \tag{2.3.36}$$

  *in* $\mathcal{CA}\mathrm{lg}^{\mathrm{cn}}$ *for which the maps* $\pi_0 A \to \pi_0 B$ *and* $\pi_0 A' \to \pi_0 B'$ *are surjective, the induced diagram*

$$\begin{array}{ccc} X(A') & \longrightarrow & X(A) \\ \downarrow & \lrcorner & \downarrow \\ X(B') & \longrightarrow & X(B) \end{array} \tag{2.3.37}$$

  *is a pullback square in* $\mathcal{S}$.

DEFINITION 2.3.38 ([96], Def. 2.1.3) : *Let* $X : \mathcal{CA}\mathrm{lg}^{\mathrm{cn}} \to \mathcal{S}$ *be a functor. We will say that* $X$ *is* nilcomplete *if, for every connective* $\mathbb{E}_\infty$*-ring* $R$*, the canonical map* $X(R) \longrightarrow \varprojlim X(\tau_{\leqslant n} R)$ *is a homotopy equivalence.*

DEFINITION 2.3.39 ([96], Def. 2.1.5) : *Let* $X : \mathcal{CA}\mathrm{lg}^{\mathrm{cn}} \to \mathcal{S}$ *be a functor. We will say that* $X$ *is* integrable *if the following condition is satisfied:*

- *Let* $R$ *be a local Noetherian* $\mathbb{E}_\infty$*-ring which is complete with respect to its maximal idel* $\mathfrak{m} \subseteq \pi_0 R$*. Then the inclusion of functors* $\mathsf{Spf}\, R \hookrightarrow \mathsf{Spec}^{\mathrm{f}}\, R$ *induces a homotopy equivalence:*

$$X(R) \simeq \mathrm{Map}_{\mathcal{F}\mathrm{un}(\mathcal{CA}\mathrm{lg}^{\mathrm{cn}},\mathcal{S})}(\mathsf{Spec}^{\mathrm{f}}\, R, X) \longrightarrow \mathrm{Map}_{\mathcal{F}\mathrm{un}(\mathcal{CA}\mathrm{lg}^{\mathrm{cn}},\mathcal{S})}(\mathsf{Spf}\, R, X)$$

DEFINITION 2.3.40 ([96], Def. 2.1.9) : *Let* $X : \mathcal{CA}\mathrm{lg}^{\mathrm{cn}} \to \mathcal{S}$ *be a functor. We will say that* $X$ *is* infinitesimally cohesive *if the following condition is satisfied:*

- *For every pullback diagram*

$$\begin{array}{ccc} A' & \longrightarrow & A \\ \downarrow & \lrcorner & \downarrow \\ B' & \longrightarrow & B \end{array} \tag{2.3.41}$$

  *in* $\mathcal{CA}\mathrm{lg}^{\mathrm{cn}}$ *for which the maps* $\pi_0 A \to \pi_0 B$ *and* $\pi_0 A' \to \pi_0 B'$ *are surjections whose kernels are nilpotent ideals in* $\pi_0 A$ *and* $\pi_0 B'$*, respectively. Then the induced diagram*

$$\begin{array}{ccc} X(A') & \longrightarrow & X(A) \\ \downarrow & \lrcorner & \downarrow \\ X(B') & \longrightarrow & X(B) \end{array} \tag{2.3.42}$$

  *is a pullback square in* $\mathcal{S}$.

### 2.3.9 The cotangent complex of Artin–Luries's criteria.

The cotangent complex involved in Artin–Luries's criteria has already been treated in subsection 2.1.2.

### 2.3.10 Artin–Lurie's representability theorem.

We come then to the central theorem which gives us that the formal moduli problems which are going to be described to reformulate Schubert geometry, are representable by spectral Deligne-Mumford stacks. Moreover, the duality of (certain) Deligne-Mumford stacks is inherited from this theorem of representability.

THEOREM 2.3.43 (Artin–Lurie's theorem [96],Thm. 3.2.1) : *Let* $X : \mathcal{CA}\mathrm{lg}^{\mathrm{cn}} \to \mathcal{S}$ *be a functor, and suppose we are given a natural transformation* $X \to \mathsf{Sp\acute{e}t}^{\mathrm{f}} R$*, where* $R$ *is a Noetherian* $\mathbb{E}_\infty$*-ring and* $\pi_0 R$ *is a Grothendieck ring. Let* $n \geqslant 0$*. Then* $X$ *is representable by a spectral Deligne-Mumford* $n$*-stack* $\mathfrak{X}$*, which is locally almost of finite presentation over* $R$*, if and only if the following conditions are satisfied:*

*1) For every discrete commutative ring* $A$*, the space* $X(A)$ *is* $n$*-truncated*

*2) The functor* $X$ *is a sheaf for the étale topology*

*3) The functor* $X$ *is nilcomplete, infinitesimally cohesive, and integrable.*

*4) The natural transformation* $X \to \mathsf{Sp\acute{e}t}^{\mathrm{f}} R$ *admits a connective cotangent complex denoted* $L_{X/\mathsf{Sp\acute{e}t}^{\mathrm{f}}}$

*5) The natural transformation* $f$ *is locally almost of finite presentation*

## 2.4. GLUING CLOSED IMMERSIONS AND REPRESENTABILITY.

### 2.4.1 Prerequisites for closed immersions.

Let us recall that in the categories of schemes, stacks, etc. we have pullbacks et fiber products but in general no pushouts.

We will remedy this situation in a case that is exactly the same as that encountered with positroid varieties, which will allow us to join them along their boundaries. Indeed, the BCFW stratification of the positive Grassmannian is carried out using open positroid varieties.

However, as soon as we take their closure, there are boundary overlaps. This is where closed immersions come into play and allow us to obtain the representability of the glued space by a structured spectral algebraic space .

To do this, we have three nested theorems. First, we will gain a simple pushout (a pullback in the opposite category) which is structured by a $\mathcal{G}$-geometry, then a pushout of $\mathcal{G}$-structured $\infty$-topoï, and then a pushout of spectral Deligne-Mumford stacks.

It is of great importance for our goal that the underlying topological space we obtain for the spectral algebraic space which represents the pushout is given by the pushout of underlying spaces in the category of topological spaces and the structure sheaf by the fiber product over the diagonal.

For the definition of a closed immersion, we refer to the section 2.4.

For the definition of geometries, we refer to definition 2.3.22.

It is noteworthy that the good behaviour of the pushouts is obtained if the (pre)geometry is unramified. We define this in the following subsection.

Let us denote $\mathcal{T}\mathrm{op} := {}^{\mathrm{L}}\mathcal{T}\mathrm{op}^{\mathrm{op}}$ where ${}^{\mathrm{L}}\mathcal{T}\mathrm{op}$ denotes the subcategory of presheaves $\widehat{\mathcal{C}\mathrm{at}_\infty}$ on $\mathcal{C}\mathrm{at}_\infty$ whose objects are $\infty$-topoï and morphisms are geometric morphisms of $\infty$-topoï $f^* : \mathfrak{X} \longrightarrow \mathfrak{Y}$. The category $\mathcal{T}\mathrm{op}$ is then the $\infty$-category of $\infty$-topoï. If $\mathcal{G}$ is a geometry, $\mathcal{T}\mathrm{op}(\mathcal{G})$ is the category whose objects are pairs $(\mathfrak{X}, \mathscr{O}_{\mathfrak{X}})$ where $\mathfrak{X}$ is an $\infty$-topos and $\mathscr{O}_{\mathfrak{X}}$ is a $\mathcal{G}$-structure.

If $R$ is an $\mathbb{E}_\infty$-ring, we denote $\mathsf{Spec}^{\text{ét}} R$ the nonconnective spectral Deligne-Mumford stack $(\mathfrak{X}, \mathscr{O}_{\mathfrak{X}})$, where $\mathfrak{X} \simeq \boldsymbol{\mathcal{S}}\mathbf{hv}_{\mathcal{S}\mathrm{p}}(\boldsymbol{\mathcal{CA}}\mathbf{lg}_R^{\text{ét}^{\mathrm{op}}})$ is the $\infty$-topos of étale sheaves over $R$, and $\mathscr{O}_{\mathfrak{X}}$ is the sheaf of $\mathbb{E}_\infty$-rings on $\mathfrak{X}$ given by $\mathscr{O}_{\mathfrak{X}}(R') = R'$.

#### 2.4.1.1 Unramification.

DEFINITION 2.4.1 ([90] Def. 1.3) : *Let $\mathcal{G}$ be a geometry. We will say that $\mathcal{G}$ is* unramified *if the following condition is satisfied:*

*For every morphism $f : \mathfrak{X} \longrightarrow \mathfrak{Y}$ in $\mathcal{G}$ ($\mathcal{G}$ is an $\infty$-topos), and every object $\mathfrak{Z} \in \mathcal{G}$, the diagram*

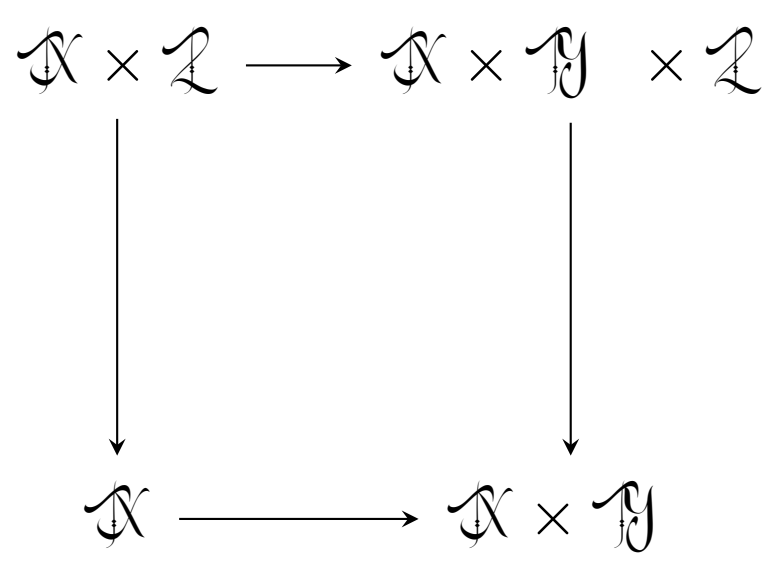

*induces a pullback square:*

$$\begin{array}{ccc}
\mathfrak{X}_{\mathcal{X}\times\mathcal{Z}} & \longrightarrow & \mathfrak{X}_{\mathcal{X}\times\mathcal{Y}\times\mathcal{Z}} \\
\downarrow & \lrcorner & \downarrow \\
\mathfrak{X}_{\mathcal{X}} & \longrightarrow & \mathfrak{X}_{\mathcal{X}\times\mathcal{Y}}
\end{array}$$

*in the category $\mathcal{T}op$*

EXAMPLE 2.4.2 :

- If the geometry $\mathcal{G}$ is discrete, then for every $\mathcal{X}_i \in \mathcal{G}$, the $\infty$-topos $\mathfrak{X}_{\mathcal{X}_i}$ is equivalent to the category $\mathcal{S}$ of spaces. Then a discrete geometry is unramified.

- If $\mathcal{G}$ is a geometry, and the composite functor $\mathcal{G} \xrightarrow{\mathsf{Spec}^{\mathcal{G}}} \mathcal{T}\mathrm{op}(\mathcal{G}) \longrightarrow \mathcal{T}\mathrm{op}$ preserves finite products, then the geometry $\mathcal{G}$ is unramified.

### 2.4.2 Three nested theorems of gluing.

We will now state the three nested theorems and corrolaries one after another:

THEOREM 2.4.3 ([90] Th. 1.6) : *Let $\mathcal{G}$ be an unramified geometry, and suppose we are given morphisms $f : (\mathfrak{Y}, \mathscr{O}_{\mathfrak{Y}}) \longrightarrow (\mathfrak{X}, \mathscr{O}_{\mathfrak{X}})$ and $g : (\mathfrak{X}', \mathscr{O}_{\mathfrak{X}'}) \longrightarrow (\mathfrak{X}, \mathscr{O}_{\mathfrak{X}})$ in the category $\mathcal{T}\mathrm{op}(\mathcal{G})$. Assume that $f$ induces an aeffective epimorphism $f^* \mathscr{O}_{\mathfrak{X}} \longrightarrow \mathscr{O}_{\mathfrak{Y}}$. Then:*

- *There exists a pullback square:*

$$\begin{array}{ccc} (\mathfrak{Y}', \mathscr{O}_{\mathfrak{Y}'}) & \longrightarrow & (\mathfrak{X}', \mathscr{O}_{\mathfrak{X}'}) \\ \downarrow & \lrcorner & \downarrow \\ (\mathfrak{Y}, \mathscr{O}_{\mathfrak{Y}}) & \longrightarrow & (\mathfrak{X}, \mathscr{O}_{\mathfrak{X}}) \end{array}$$

  *in the category $\mathcal{T}\mathrm{op}(\mathcal{G})$*

- *The map $f'$ induces an effective epimorphism $f'^* \mathscr{O}\mathfrak{X}' \longrightarrow \mathscr{O}\mathfrak{Y}'$ in $\mathrm{Str}^{\mathrm{loc}}_{\mathcal{G}}(\mathfrak{Y}')$*

- *The underlying diagram of $\infty$-topoï*

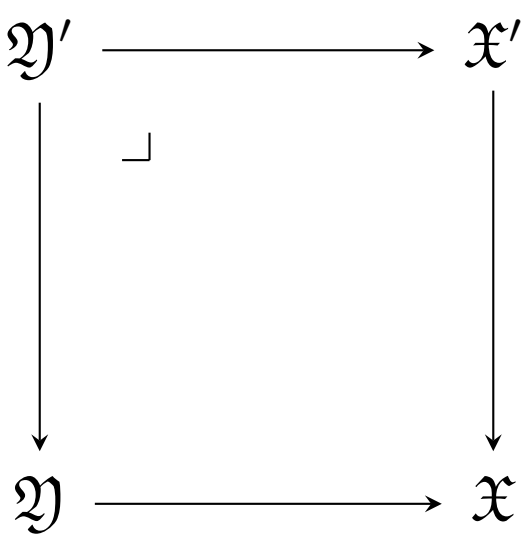

*is a pullback square in the category* $\mathcal{T}\mathrm{op}$*.*

- *For every geometric morphism of* $\infty$*-topoï* $h^* : \mathfrak{Y}' \longrightarrow \mathfrak{Z}$*, the diagram:*

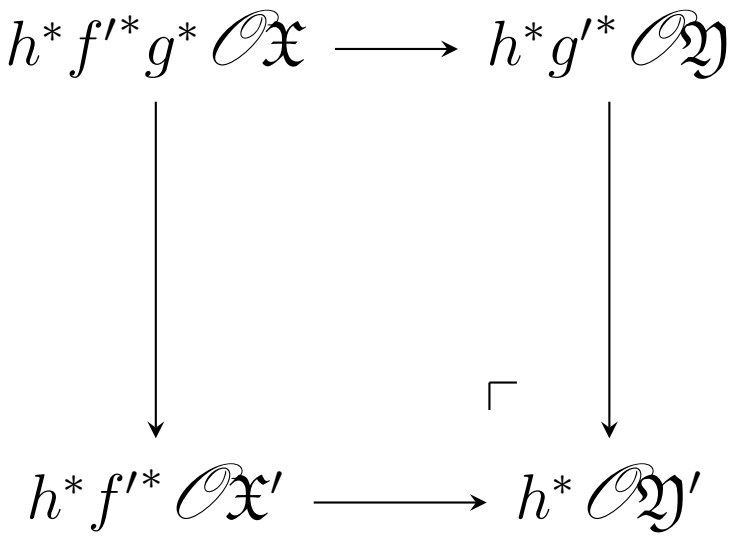

*is a pushout square in* $\mathrm{Str}^{\mathrm{loc}}_{\mathcal{G}}(\mathfrak{Z})$

corollary 2.4.4 (Fiber products in ${}^{\mathrm{L}}\mathcal{T}\mathrm{op}(\mathcal{G})^{\mathrm{op}}$, [90] Cor. 1.7) : *Let* $\mathcal{G}$ *be a geometry and suppose we are given maps:*

$$(\mathfrak{Y}, \mathscr{O}_{\mathfrak{Y}}) \leftarrow (\mathfrak{X}, \mathscr{O}_{\mathfrak{X}}) \rightarrow (\mathfrak{X}', \mathscr{O}_{\mathfrak{X}'})$$

*in the category* ${}^{\mathrm{L}}\mathcal{T}\mathrm{op}(\mathcal{G})^{\mathrm{op}}$*. Assume that* $\mathcal{G}$ *is unramified and that* $f$ *is a closed immersion. Then there exists a pullback square:*

$$\begin{array}{ccc} (\mathfrak{Y}', \mathscr{O}_{\mathfrak{Y}'}) & \longrightarrow & (\mathfrak{X}', \mathscr{O}_{\mathfrak{X}'}) \\ \downarrow & \lrcorner & \downarrow \\ (\mathfrak{Y}, \mathscr{O}_{\mathfrak{Y}}) & \longrightarrow & (\mathfrak{X}, \mathscr{O}_{\mathfrak{X}}) \end{array}$$

*in the category* ${}^{\mathrm{L}}\mathcal{T}\mathrm{op}(\mathcal{G})^{\mathrm{op}}$*. Moreover,* $f'$ *is a closed immersion and the diagram of* $\infty$*-topoï*

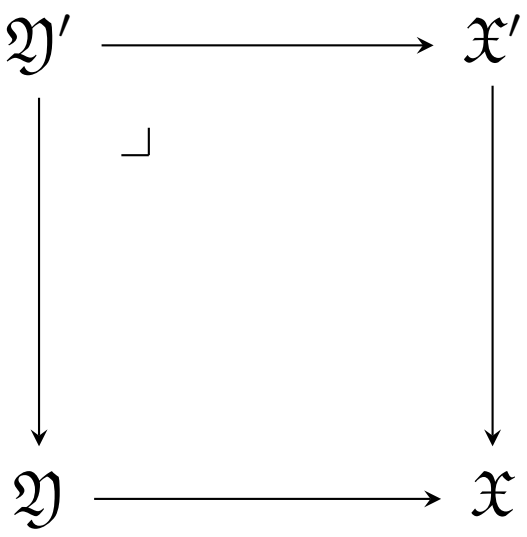

*is a pullback square in* ${}^{\mathrm{L}}\mathcal{T}\mathrm{op}^{\mathrm{op}}$

REMARK 2.4.5 : As an analogous of Prop 4.1 of [90], we have to remember to prove that the $\mathcal{IC}^{\bullet}$-geometry for spectral Deligne-Mumford stacks is unramified.

Now, in the category of schemes, if we have a diagram

$$\mathcal{X}_0 \xleftarrow{f} \mathcal{X}_{01} \xhookrightarrow{g} \mathcal{X}_1$$

where $f$ and $g$ are closed immersions, then there exists a pushout diagram:

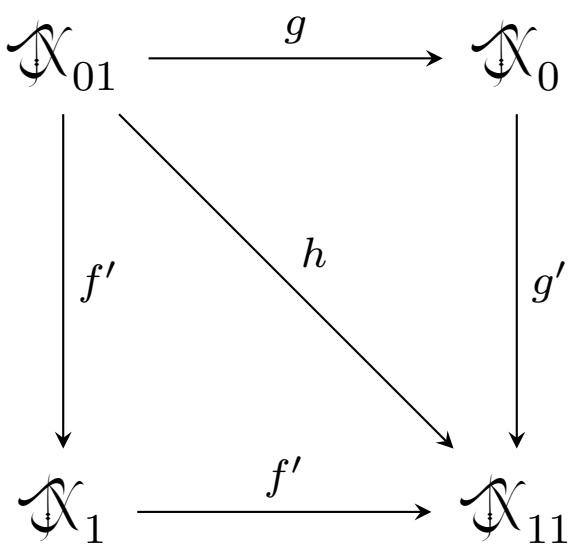

in the category of schemes. We glue along the common closed subscheme $\mathscr{X}_{01}$ and the structure sheaf $\mathscr{O}_{\mathscr{X}_{11}}$ can be described by the pushout $g'_* \mathscr{O}_{\mathscr{X}_0} \coprod_{h_* \mathscr{O}_{\mathscr{X}_{01}}} f'_* \mathscr{O}_{\mathscr{X}_1}$. And the category can even be larger than that of schemes.

THEOREM 2.4.6 (pushouts of $\mathcal{G}$-structured $\infty$-topoï, [90] Thm. 5.1) : *Let $\mathcal{G}$ be a geometry and suppose we are given morphisms:*

$$(\mathfrak{X}_0, \mathscr{O}_{\mathfrak{X}_0}) \leftarrow (\mathfrak{X}_{01}, \mathscr{O}_{\mathfrak{X}_{01}}) \rightarrow (\mathfrak{X}_1, \mathscr{O}_{\mathfrak{X}_1})$$

*in $\mathcal{T}\mathrm{op}(\mathcal{G})$ which induce closed immersions of $\infty$-topoï $\mathfrak{X}_0 \hookleftarrow \mathfrak{X}_{01} \hookrightarrow \mathfrak{X}_1$. Let $\mathfrak{X}_{11}$ be the pushout $\mathfrak{X}_0 \coprod_{\mathfrak{X}_{01}} \mathfrak{X}_1$ in $\mathcal{T}\mathrm{op}$, so that we have a commutative diagram of geometric morphisms:*

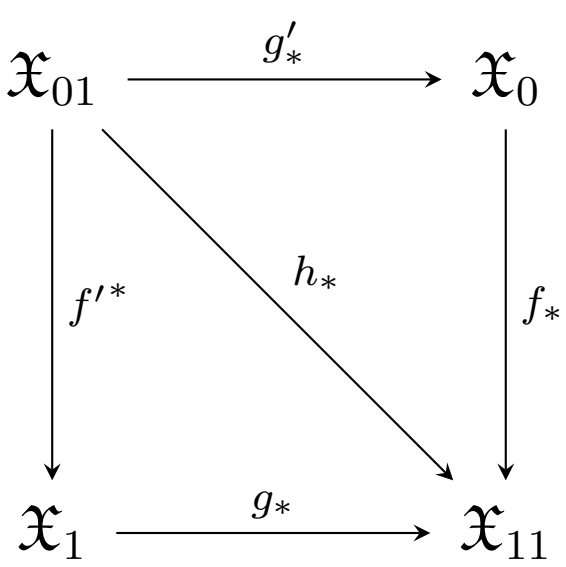

*Let* $\mathscr{O}_{\mathfrak{X}_{11}} : \mathcal{G} \longrightarrow \mathfrak{X}_{11}$ *be the functor given by:*

$$\mathscr{O}_{\mathfrak{X}_{11}}(\mathfrak{X}) = f_* \mathscr{O}_{\mathfrak{X}_0}(\mathfrak{X}) \times_{h_* \mathscr{O}_{\mathfrak{X}_{01}}} g_* \mathscr{O}_{\mathfrak{X}_1}(\mathfrak{X})$$

*then* $\mathscr{O}_{\mathfrak{X}_{11}}$ *is a* $\mathcal{G}$*-structure on* $\mathfrak{X}_{11}$*, and the diagram:*

$$\begin{array}{ccc} (\mathfrak{X}_{01}, \mathscr{O}_{\mathfrak{X}_{01}}) & \longrightarrow & (\mathfrak{X}_0, \mathscr{O}_{\mathfrak{X}_0}) \\ \downarrow & & \downarrow \\ (\mathfrak{X}_1, \mathscr{O}_{\mathfrak{X}_1}) & \longrightarrow & (\mathfrak{X}_{11}, \mathscr{O}_{\mathfrak{X}_{11}}) \end{array}$$

*is a pushout square in* $\mathcal{T}\mathrm{op}(\mathcal{G})$*.*

THEOREM 2.4.7 (pushouts of spectral DM stacks, [90] Thm. 6.1) : *Suppose we are given morphisms:*

$$\mathfrak{X}_0 \longleftarrow \mathfrak{X}_{01} \hookrightarrow \mathfrak{X}_1$$

*which are closed immersions of spectral Deligne-Mumford stacks. Then:*

- *There exists a pushout diagram:*

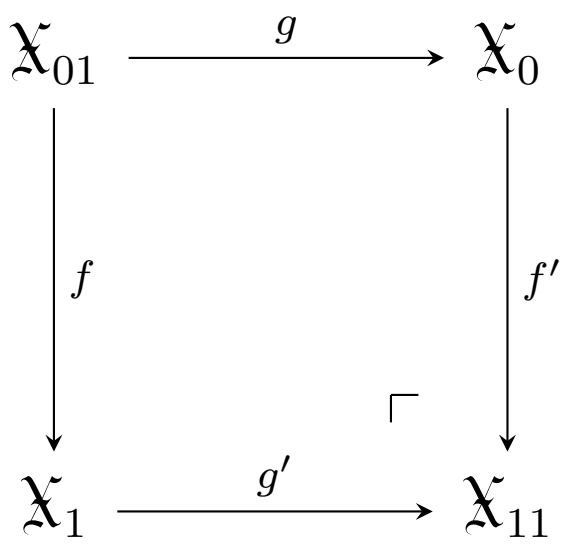

*in* $\mathbf{RingTop}_{\text{ét}}$

- *The image of the diagram in* $\mathcal{T}\mathrm{op}$ *is a pushout diagram of* $\infty$*-topoï*
- *The maps* $f'$ *and* $g'$ *are closed immersions*
- *The pushout* $\mathfrak{X}_{11}$ *is a spectral Deligne-Mumsford stack*
- *If* $\mathfrak{X}_0$ *and* $\mathfrak{X}_1$ *are affine, then* $\mathfrak{X}_{11}$ *is also affine*

## 2.5. Duality of coaffine spectral Deligne-Mumford stacks.

### 2.5.1 The dual of a coaffine geometric stack: First step.

DEFINITION 2.5.1 ([46] Section. 0.0.1) : *Let us call a* prestack*, a functor of* $\infty$*-categories, without more conditions like those of definition 2.2.108:*

$$\mathfrak{X} : \mathcal{CA}\mathrm{lg}^{\mathrm{cn}} \longrightarrow \mathcal{S}$$

We have a sheaf of categories $\mathcal{C}$ over a prestack $\mathfrak{X}$, denoted $\mathsf{ShvCat}(\mathfrak{X})$, which is a functorial assignement for every affine spectral scheme $\mathfrak{Y} := \mathsf{Spec}\ R$, where $R$ is an $\mathbb{E}_\infty$-ring spectrum, such that there exists a faithfully flat map $\mathfrak{Y} := \mathsf{Spec}\ R \longrightarrow \mathfrak{X}$, of the category $\Gamma(\mathfrak{X}, \mathcal{C})$ of its global sections over $\mathfrak{X}$, if we take $\mathcal{C}$ in $\mathsf{ShvCat}(\mathfrak{X})$ itself:

$$\Gamma(\mathfrak{X},-):\begin{cases}\mathsf{ShvCat}(\mathfrak{X}) & \longrightarrow \quad \mathsf{QCoh}(\mathfrak{X})\text{-Mod}\\ \quad \mathcal{C} & \longmapsto \quad \Gamma(\mathfrak{X},\mathcal{C})\end{cases}$$

It will be naturally acted on by $\mathsf{QCoh}(\mathfrak{X})$. Thus, we obtain an object $\Gamma^{\mathrm{enh}}(\mathfrak{X},\mathcal{C}) \in \mathsf{QCoh}(\mathfrak{X}) - \mathrm{Mod}$, and we obtain a functor:

$$\Gamma^{\mathrm{enh}}_{\mathfrak{X}} := \Gamma^{\mathrm{enh}}(\mathfrak{X},-) : \mathsf{ShvCat}(\mathfrak{X}) \longrightarrow \mathsf{QCoh}(\mathfrak{X}) - \mathrm{Mod}$$

which we can denote, on each affine spectral scheme $\mathfrak{Y}$ with $\mathfrak{Y} \to \mathfrak{X}$ :

$$(\mathfrak{Y} := \mathsf{Spec}\ R \to \mathfrak{X}) \in \mathcal{CA}\mathrm{lg}^{\mathrm{cn}}/\mathfrak{X} \rightsquigarrow \Gamma^{\mathrm{enh}}_{\mathfrak{X}}(\mathcal{C}) \in \mathrm{Mod}_{\mathsf{QCoh}(\mathfrak{X})} \tag{2.5.2}$$

where we denote $\mathrm{Mod}_{\mathsf{QCoh}(\mathfrak{X})}$ the $\infty$-category of $\mathsf{QCoh}(\mathfrak{X})$-modules categories.

The assignement Equation 2.5.2 must be functorial in $\mathfrak{Y}$ in the sense that for a map $f : \mathfrak{Y}_1 \to \mathfrak{Y}_2$, we must be given an isomorphism in $\mathrm{Mod}_{\mathsf{QCoh}(\mathfrak{Y}_1)}$ :

$$\mathsf{QCoh}(\mathfrak{Y}_1) \otimes_{\mathsf{QCoh}(\mathfrak{Y}_2)} \Gamma(\mathfrak{Y}_2,\mathcal{C}) \xrightarrow{\sim} \Gamma(\mathfrak{Y}_1,\mathcal{C})$$

together with a homotopy-coherent system of compatibilities for composition of morphisms.

The goal of the study here is then the connection between the $\infty$-categories $\mathsf{ShvCat}(\mathfrak{X})$ and $\mathrm{Mod}_{\mathsf{QCoh}(\mathfrak{X})}$:

$$\mathsf{ShvCat}(\mathfrak{X}) \longrightarrow \mathrm{Mod}_{\mathsf{QCoh}(\mathfrak{X})} \tag{2.5.3}$$

If $\mathfrak{X}$ is an affine spectral scheme $\mathfrak{X} = \mathsf{Spec}\ R$ for $R$ an $\mathbb{E}_\infty$-ring, these categories are tautologically equivalent and we say that $\mathfrak{X}$ is 1-affine.

More generally, we say:

Definition 2.5.4 : *The prestack $\mathfrak{X}$ is* 1-affine *if the $\infty$-categories* $\mathsf{ShvCat}(\mathfrak{X})$ *and* $\mathrm{Mod}_{\mathsf{QCoh}(\mathfrak{X})}$ *are equivalent:*

$$\mathsf{ShvCat}(\mathfrak{X}) \xrightarrow{\sim} \mathrm{Mod}_{\mathsf{QCoh}(\mathfrak{X})} \tag{2.5.5}$$

Definition 2.5.6 : *Let $\mathfrak{X}$ be a prestack. We say that $\mathfrak{X}$ is* passable *if*

- *The diagonal morphism* $\delta : \mathfrak{X} \longrightarrow \mathfrak{X} \times \mathfrak{X}$ *is schematic, quasi-affine and quasi-compact*

- *The sheaf* $\mathscr{O}_{\mathfrak{X}} \in \mathsf{QCoh}(\mathfrak{X})$ *is (perfect, dualizable, that is) compact*

- *The category* $\mathsf{QCoh}(\mathfrak{X})$ *is dualizable in the category* $\mathsf{Groth}_\infty$

PROPOSITION 2.5.7 ([46] Prop. 5.1.3) : *Let* $\mathfrak{X}$ *be a passable prestack. Then the functor*

$$\Gamma^{\text{enh}}_{\mathfrak{X}} := \mathsf{ShvCat}(\mathfrak{X}) \longrightarrow \mathsf{Mod}_{\mathsf{QCoh}(\mathfrak{X})} \tag{2.5.8}$$

*is fully faithful.*

THEOREM 2.5.9 ([46] Thm. 5.1.5) : *A coconnective locally of finite type quasi-compact geometric stack is passable, and then is 1-affine.*

Now, let us restrict the domain $\mathsf{ShvCat}(\mathfrak{X})$ of $\Gamma^{\text{enh}}_{\mathfrak{X}}$ to $\mathsf{QStk}^{\text{PSt}}(\mathfrak{X})$ where we denote:

- $\mathsf{Groth}_\infty$ : the $\infty$-category of Grothendieck prestable $\infty$-categories, which admits a symmetric monoidal structure ( see [103] Appendix C, and section C.4 for the monoidal structure)

- $\mathsf{Mod}_{\mathsf{QCoh}(\mathfrak{X})^{\text{cn}}}(\mathsf{Groth}_\infty)$ : the associated $\infty$-category of modules where $\mathsf{QCoh}(\mathfrak{X})^{\text{cn}}$ is indeed a commutative algebra object of the $\infty$-category $\mathsf{Groth}_\infty$

- $\mathsf{QStk}^{\text{PSt}}(\mathfrak{X})$ : the $\infty$-category of prestable quasi-coherent stacks $\mathfrak{C}$ on $\mathfrak{X}$, for $\mathfrak{X}$ a spectral Deligne-Mumford stack (see [103] Definition 10.1.2.1)

THEOREM 2.5.10 ([103],Thm. 10.2.0.2) : *Let* $\mathfrak{X}$ *be a quasi-compact, quasi-separated spectral algebraic space. Then the following functor constructed in [103] 10.1.7.1:*

$$\mathsf{QCoh}(\mathfrak{X}\,,\bullet) : \mathsf{QStk}^{\text{PSt}}(\mathfrak{X}) \xrightarrow{\sim} \mathsf{Mod}_{\mathsf{QCoh}(\mathfrak{X})^{\text{cn}}}(\mathsf{Groth}_\infty) \subsetneq \mathsf{Mod}_{\mathsf{QCoh}(\mathfrak{X})} \tag{2.5.11}$$

*can be viewed as a* "global sections functor" *and is an equivalence of categories for these hypotheses on* $\mathfrak{X}$. *This is* 1-affineness *for quasi-compact, quasi-separated spectral algebraic spaces.*

Let $(\mathscr{X}, \mathscr{O}_{\mathscr{X}})$ be a (classical) quasi-affine scheme, such that the category $\mathsf{QCoh}(\mathscr{X})$ is compactly generated by $\mathscr{O}_{\mathscr{X}}$, then we had defined that $\mathscr{X}$ is 0-affine, or weakly 0-affine ([46]), and the global sections functor:

$$\Gamma(\mathscr{X}, -) : \mathsf{QCoh}(\mathscr{X}) \xrightarrow{\sim} \mathsf{Mod}_{\mathbb{R}\Gamma(\mathscr{X}, \mathscr{O}_{\mathscr{X}})}$$

is an equivalence.

Now, if we start from an étale covering of $\mathfrak{X}$ a spectral Deligne-Mumford stack by quasi-affine schemes $\coprod_{f \in \mathscr{B}(k,n)} \mathscr{X}_f \longrightarrow \mathfrak{X}$, where $\mathscr{X}_f \hookrightarrow \mathfrak{X}$ is a quasi-affine classical scheme, then, recalling theorem 2.2.94 for $(\mathscr{X}_f, \mathscr{O}_f)$ a classical scheme which is quasi-affine, the derived category of the stable $\infty$-category $\mathsf{QCoh}(\mathscr{X})$ is equivalent to the following category of sheaves of $\mathbb{E}_\infty$-rings:

$$\mathscr{D}^+(\mathsf{QCoh}(\mathscr{X}_f)) \xrightarrow{\sim} \mathsf{Mod}_{\mathbb{R}\Gamma(\mathscr{X}_f, \mathscr{O}_f)} \tag{2.5.12}$$

where $\mathbb{R}\Gamma$ is the right derived functor of global sections functor. So $\mathbb{R}\Gamma(\mathscr{X}_f, \mathscr{O}_f)$ is an $\mathbb{E}_\infty$-ring which is coconnective (it has only negative homotopy groups), and taking its global sections gives a functor which is an equivalence.

The quasi-affine scheme $\mathscr{X}_f$ is then 0-$\mathsf{affine}$, since in this case, the global sections are also discrete.

But recalling lemma 2.2.96, the geometric stack $(\mathfrak{X}, \mathscr{O}_{\mathfrak{X}})$ is then 0-affine since the sheaf $\mathscr{O}_{\mathfrak{X}}$ is compact in $\mathsf{QCoh}(\mathfrak{X})$ and this sheaf $\mathscr{O}_{\mathfrak{X}}$ generates $\mathsf{QCoh}(\mathfrak{X})$.

But the conditions of the definition 2.5.6 being fulfilled in our case of a geometric stack $\mathfrak{X}$ with diagonal morphism which is schematic, affine, and quasi-compact, which is indeed the case here of a prestable quasi-coherent stack $\mathfrak{Z}$ on $\mathfrak{X}$, the sheaf $\mathscr{O}_{\mathfrak{X}}$ being compact in $\mathsf{QCoh}(\mathfrak{X})$ , the theorem 2.5.9 tells us that the 0-affine stack $\mathfrak{Z}$ is also 1-affine.

To get the dual of a coaffine geometric stack $\mathfrak{X}$, in particular a coaffine Deligne-Mumford stack, amounts to the following steps:

- Take the category of quasi-coherent sheaves $\mathsf{QCoh}(\mathfrak{X})$ on the coaffine geometric stack $\mathfrak{X}$ with values in the category $\mathsf{Groth}_\infty$

- Identify by 0-affineness the category $\mathsf{QCoh}(\mathfrak{X})$ with the category $\mathfrak{C} := \mathsf{Mod}_{\Gamma(\mathfrak{X},\mathscr{O}_{\mathfrak{X}})}$
- By 1-affineness (Equation 2.5.5), dualize the category $\mathfrak{C} := \mathsf{Mod}_{\Gamma(\mathfrak{X},\mathscr{O}_{\mathfrak{X}})}$ in the category $\mathsf{Mod}_{\mathsf{QCoh}(\mathfrak{X})^{\mathrm{cn}}}(\mathsf{Groth}_{\infty})$ to get an element $\mathfrak{C}^{\star}$ (see subsection 2.5.2):

$$\mathfrak{C}^{\star} \in \mathsf{QStk}^{\mathrm{PSt}}(\mathfrak{X}) \simeq \mathsf{Mod}_{\mathsf{QCoh}(\mathfrak{X})^{\mathrm{cn}}}(\mathsf{Groth}_{\infty}) \subset \mathsf{Mod}_{\mathsf{QCoh}(\mathfrak{X})} \tag{2.5.13}$$

- The category $\mathfrak{C}^{\star}$ is in fact an object of $\mathsf{Mod}_{\mathsf{QCoh}(\mathfrak{X})}$ which is a quasi-coherent sheaf for the étale topology $\infty$-sheaf condition, as a quasi-coherent stack on $\mathfrak{X}$, so it is an object of $\mathsf{QCoh}(\mathfrak{X})$
- Use Tannaka reconstruction (Equation 2.2.116) to reconstruct from $\mathfrak{C}^{\star}$ in the category $\mathsf{QCoh}(\mathfrak{X})$, what we define to be the dual coaffine geometric stack $\mathfrak{X}^{\star}$ of $\mathfrak{X}$.

definition 2.5.14 (Dual coaffine geometric stack) : *We define the* dual coaffine geometric stack $\mathfrak{X}^{\star}$ *of a coaffine geometric stack* $\mathfrak{X}$ *the object obtained from the preceding construction.*

At this point, it remains then to explain how we dualize the category $\mathsf{Mod}_{\Gamma(\mathfrak{X},\mathscr{O}_{\mathfrak{X}})}$ in the category $\mathsf{Mod}_{\mathsf{QCoh}(\mathfrak{X})^{\mathrm{cn}}}(\mathsf{Groth}_{\infty})$ thanks to Serre-Grothendieck-Verdier duality and dualizing complexes (next subsection 2.5.2), and we have to check that we indeed reconstruct a stack which is Deligne-Mumford, coaffine and geometric (see subsection 2.5.3). And finally, we have to check the non-triviality of the dualization of these stacks in the spectral world, thanks to the cohomology of their Brauer-Grothendieck groups (see subsection 2.5.6).

### 2.5.2 The Verdier duality functor: Second step.

Let $\mathfrak{X}$ be an affine spectral Deligne-Mumford stack. Thanks to lemma 2.2.96, we know that $\mathscr{O}_{\mathfrak{X}}$ is compact and generates the stable monoidal $\infty$-category $\left(\mathsf{QCoh}(\mathfrak{X}), \otimes, \mathbf{1}\right)$, where $\mathbf{1}$ is a terminal object of $\mathfrak{X}$, that is $\mathbf{1} \simeq \mathscr{O}_{\mathfrak{X}}$.

In the spirit of the classical definition of a dual, we would like to write something like the dual of $\mathfrak{X}$ in the category of spectral Deligne-Mumford stack is given by " $\mathfrak{X}^{\star} := \mathrm{Hom}_{\mathrm{SpDMStk}}(\mathfrak{X}, \mathscr{O}_{\mathfrak{X}})$ ". But for the moment, this definition has no sense. Indeed, $\mathscr{O}_{\mathfrak{X}}$ is not a spectral Deligne-Mumford stack, since it is just a sheaf defined in definition 2.2.69. But Tannaka duality and Grothendieck-Serre-Verdier duality will solve the problem. Let us recall theorem 2.2.94 which says that if $\mathscr{X}$ is an 0-affine (classical) scheme, then the derived category of the stable $\infty$-category $\mathsf{QCoh}(\mathscr{X})$ is equivalent to the following category of sheaves of $\mathbb{E}_\infty$-rings:

$$\mathsf{QCoh}(\mathscr{X}) \xrightarrow{\sim} \mathsf{Mod}_{\mathbb{R}\Gamma(\mathscr{X}, \mathscr{O}_{\mathscr{X}})} \tag{2.5.15}$$

where $\mathbb{R}\Gamma$ is the right derived functor of the global sections functor. We also have that $\mathbb{R}\Gamma(\mathscr{X}, \mathscr{O}_{\mathscr{X}})$ is an $\mathbb{E}_\infty$-ring which is coconnective (it has only negative homotopy groups and then nonnegative cohomological groups), and taking its global sections gives a functor which is an equivalence.

But we have more than this equivalence if the 0-affine classical scheme is a perfect stack. Indeed, we recall, for $\mathfrak{X}$ a quasi-compact, quasi-separated spectral algebraic space, that the following functor constructed in [103] 10.1.7.1 (see theorem 2.5.10):

$$\mathsf{QCoh}(\mathfrak{X}, \bullet) : \mathsf{QStk}^{\mathrm{PSt}}(\mathfrak{X}) \xrightarrow{\sim} \mathsf{Mod}_{\mathsf{QCoh}(\mathfrak{X})^{\mathrm{cn}}}(\mathsf{Groth}_\infty) \tag{2.5.16}$$

can be viewed as a "*global sections functor*" and is an equivalence of categories for these hypotheses on $\mathfrak{X}$.

In particular, if $\mathfrak{X}$ is a coaffine geometric spectral Deligne-Mumford stack, then we have this equivalence of $\infty$-categories.

The monoidal structure $\mathsf{QStk}^{\mathrm{PSt}}(\mathfrak{X})$ then allows us to take the monoidal dual $\mathsf{QCoh}(\mathfrak{X})^\vee$ of $\mathsf{QCoh}(\mathfrak{X})$ in $\mathsf{QStk}^{\mathrm{PSt}}(\mathfrak{X})$.

Tannaka duality described in Equation 2.2.116 and theorem 2.2.113 allows us to recover the coaffine spectral Deligne-Mumford stack $\mathfrak{X}$ from $\mathsf{QCoh}(\mathfrak{X})$. We can then define the dual $\mathfrak{X}^{\star}$ of $\mathfrak{X}$ as follows:

DEFINITION 2.5.17 : *The dual $\mathfrak{X}^{\star}$ of a coaffine geometric spectral Deligne-Mumford stack $\mathfrak{X}$ is the affine spectral Deligne-Mumford stack recovered by J. Lurie's Tannaka* *duality for geometric stacks from the stable $\infty$-category* $\mathsf{QCoh}(\mathfrak{X})^\vee$ *such that:*

$$\left(\mathsf{QCoh}(\mathfrak{X}^{\star}), \otimes, \mathscr{O}_{\mathfrak{X}^{\star}}\right) \simeq \left(\mathsf{QCoh}(\mathfrak{X}), \otimes, \mathscr{O}_{\mathfrak{X}}\right)^{\vee}. \tag{2.5.18}$$

But this definition is not constructive. See definition 2.5.14 for the construction.

Consider now a coconnective spectral scheme $\mathscr{X}$ and the category of sheaves $\mathsf{Mod}_{\mathbb{R}\Gamma(\mathscr{X},\mathscr{O}_{\mathscr{X}})}$. We want to take the dual of $\mathbb{R}\Gamma(\mathscr{X}, \mathscr{O}_{\mathscr{X}})$. This is exactly the situation of Grothendieck-Serre-Verdier duality.

DEFINITION 2.5.19 ([109],Def. 3.1) : *Let* $\mathscr{X}$ *be a noetherian separated scheme. A* dualizing complex *of* $\mathscr{X}$ *is an object* $\omega_{\mathscr{X}} \in \mathcal{D}^+(\mathsf{Coh}(\mathscr{X}))$ *so that the functor:*

$$\begin{array}{ccc} \mathcal{D}^+(\mathsf{Coh}(\mathscr{X})) & \longrightarrow & \mathcal{D}^+(\mathsf{Coh}(\mathscr{X})) \\ \mathcal{F} & \longmapsto & \mathcal{H}\mathrm{om}_{\mathcal{D}^+(\mathsf{Coh}(\mathscr{X}))}(\mathcal{F}, \omega_{\mathscr{X}}) \end{array}$$

*gives an equivalence*

$$\mathcal{D}^+(\mathsf{Coh}(\mathscr{X}))^{\mathrm{op}} \xrightarrow{\sim} \mathcal{D}^+(\mathsf{Coh}(\mathscr{X}))^{\mathrm{op}}$$

THEOREM 2.5.20 : *Let* $f : X \to Y$ *be a continuous map between locally compact topological spaces. There exists a right adjoint* $f^! : \mathcal{D}^+(\mathsf{Coh}(Y)) \longrightarrow \mathcal{D}^+(\mathsf{Coh}(X))$ *such that:*

$$\mathcal{H}\mathrm{om}_{\mathcal{D}^+(\mathsf{Coh}(Y))}(\mathbb{R}f_!\mathcal{F}, \mathcal{G}) \simeq \mathcal{H}\mathrm{om}_{\mathcal{D}^+(\mathsf{Coh}(X))}(\mathcal{F}, f^!\mathcal{G})$$

DEFINITION 2.5.21 : *We define the dual of* $\mathbb{R}\Gamma(\mathscr{X}, \mathscr{O}_{\mathscr{X}})$ *such that:*

$$\mathbb{R}\Gamma(\mathscr{X}, \mathscr{O}_{\mathscr{X}})^{\star} = \mathcal{H}\mathrm{om}_{\mathcal{D}^+(\mathsf{Coh}(\mathscr{X}))}(\mathbb{R}\Gamma(\mathscr{X}, \mathscr{O}_{\mathscr{X}}), \omega_{\mathscr{X}})$$

DEFINITION 2.5.22 : *In the case of a coaffine geometric spectral stack* $\mathfrak{X}$*, we define the dual* $\mathfrak{X}^{\star}$ *from the category defined in definition 2.5.14 :* $\mathfrak{F}^{\star}$ *in* $\mathsf{QStk}^{\mathrm{PSt}}(\mathfrak{X}) \simeq \mathsf{Mod}_{\mathsf{QCoh}(\mathfrak{X})^{\mathrm{cn}}}(\mathsf{Groth}_\infty) \subset \mathsf{Mod}_{\mathsf{QCoh}(\mathfrak{X})}$ *by the coaffine geometric spectral stack recovered by Tannaka duality from the following equivalence of $\infty$-categories:*

$$\mathfrak{F}^{\star} \underset{\mathrm{htpy}}{\simeq} \mathsf{Mod}_{\Gamma(\mathfrak{X},\mathscr{O}_{\mathfrak{X}})^{\star}} \tag{2.5.23}$$

*with*

$$\mathsf{Mod}_{\Gamma(\mathfrak{X},\mathscr{O}_{\mathfrak{X}})} \otimes \mathsf{Mod}_{\Gamma(\mathfrak{X},\mathscr{O}_{\mathfrak{X}})^{\star}} \underset{\text{htpy}}{\simeq} \mathbf{1}_{\mathsf{Groth}_{\infty}} \underset{\text{htpy}}{\simeq} \mathsf{Mod}_{\mathrm{Sp}} \tag{2.5.24}$$

*and*

$$\Gamma(\mathfrak{X}, \mathscr{O}_{\mathfrak{X}})^{\star} = \mathcal{H}\mathrm{om}_{\mathcal{D}^{+}(\mathsf{Coh}(\mathfrak{X}))}(\Gamma(\mathfrak{X}, \mathscr{O}_{\mathfrak{X}}), \omega_{\mathfrak{X}}) \tag{2.5.25}$$

In the construction of definition 2.5.14, we have then dualized the category $\mathsf{Mod}_{\Gamma(\mathfrak{X},\mathscr{O}_{\mathfrak{X}})}$ in the category $\mathsf{Mod}_{\mathsf{QCoh}(\mathfrak{X})^{\mathrm{cn}}}(\mathsf{Groth}_{\infty})$ thanks to Serre-Grothendieck-Verdier duality and dualizing complexes

### 2.5.3 Artin-Lurie's dual representation: Third step.

#### 2.5.3.1 The dual of a simplicial object.

Let us recall the construction of the opposite of an $\infty$-category given by C. Rezk ([116], 6.16), which is the same as the definition given by D-C. Cisinski ([30], 1.5.7) and J. Lurie ([98],???). Let $\rho : \Delta \to \Delta$ be the functor defined as the identity on objects and by the formula

$$\rho(f)(i) = n - f(m-i)$$

for any map $f : [m] \longrightarrow [n]$, for all $0 \leqslant i \leqslant m$. The functor $\rho$ is a non-trivial involution.

REMARK 2.5.26 : This involution is the functor which "reverses the ordering" of the totally-ordered sets [n]. We get a totally ordered set "[n] *with the order of its elements reversed*", which isn't actually an object of $\Delta$, but is uniquely isomorphic to [n], via the function $x \mapsto n - x$. On the other hand, the following definition agrees with the world of $\infty$-categories:

DEFINITION 2.5.27 : *Let $\mathscr{C}$ be a category. We define the* dual of a simplicial object *$X$ in $\mathscr{C}$, denoted $X^{\star} : \Delta^{\mathrm{op}} \longrightarrow \mathscr{C}$ as the composite functor $X^{\star} := X \circ \rho$. Given a morphism of simplicial objects $X : S \longrightarrow T$, we get, by functoriality, a morphism $X^{\star} : S^{\star} \longrightarrow T^{\star}$. The duals of simplicial objects form an $\infty$-category denoted $s\mathscr{C}^{\star}$. If $\mathscr{C}$ is an $\infty$- category, the opposite category is denoted $\mathscr{C}^{\mathrm{op}}$.*

EXAMPLE 2.5.28 : We define the **dual of a simplicial set** $X$, denoted $X^{\star}$ as : $\Delta^{\mathrm{op}} \longrightarrow \mathcal{S}\mathrm{ets}$. The duals of simplicial sets form an $\infty$-category denoted $s\mathcal{S}\mathrm{ets}^{\star}$.

REMARK 2.5.29 : We have that $(\Delta^n)^{\star} = \Delta^n$, while the horns $(\Lambda^n_j)^{\star} = \Lambda^n_{n-j}$, so that the opposite of an inner horn is another inner horn. Thus the dual of an $\infty$-category is an $\infty$-category.

REMARK 2.5.30 : The covariant functor *"dual"* : $(-)^{\star}$ commutes with the nerve functor $\mathcal{N}$, that is: $\mathcal{N}(\mathscr{C}^{\star}) \simeq \mathcal{N}(\mathscr{C})^{\star}$. As a consequence, the Quillen adjunctions of example 2.1.3, example 2.1.5 and example 2.1.6 are valid is the dual world.

### 2.5.4 Left Kan extension, Yoneda functor, presheaves and Isbell adjunction.

Let $\mathfrak{X} := (\mathfrak{X}, \mathscr{O}_{\mathfrak{X}})$ be a spectrally ringed $\infty$-topos. Let us now define a first version of a dual spectrally ringed $\infty$-topos.

DEFINITION 2.5.31 (Dual of a spectrally ringed $\infty$-topos) : *The dual $\mathfrak{X}^{\star}$ of $\mathfrak{X}$ is defined such that $\mathfrak{X}^{\star} := \mathfrak{X}^{\mathrm{op}}$ as in subsubsection 2.5.3.1, that is, the opposite of the $\infty$-category $\mathfrak{X}$.*

To validate this dual, we have to verify that the structure sheaf of the dual $\mathfrak{X}^{\star}$ is the dual, in the sense of subsubsection 2.5.3.1 for $\infty$-categories, of the structure sheaf $\mathscr{O}_{\mathfrak{X}}$. Indeed, $\mathfrak{X}$ and $\mathscr{O}_{\mathfrak{X}}$ are $\infty$-categories to which the construction of the simplicial dual of subsubsection 2.5.3.1 applies.

But we have the following:

THEOREM 2.5.32 : *Let $\mathscr{C}$ and $\mathscr{D}$ be two $\infty$-categories. The following two duals coincide:*

$$\boldsymbol{\mathcal{S}}\mathbf{hv}_{\mathscr{C}^{\mathrm{op}}}(\mathscr{D}) \simeq \boldsymbol{\mathcal{S}}\mathbf{hv}_{\mathscr{C}}(\mathscr{D}^{\mathrm{op}}) \tag{2.5.33}$$

*Proof.* We have an Isbell $\infty$-adjunction:

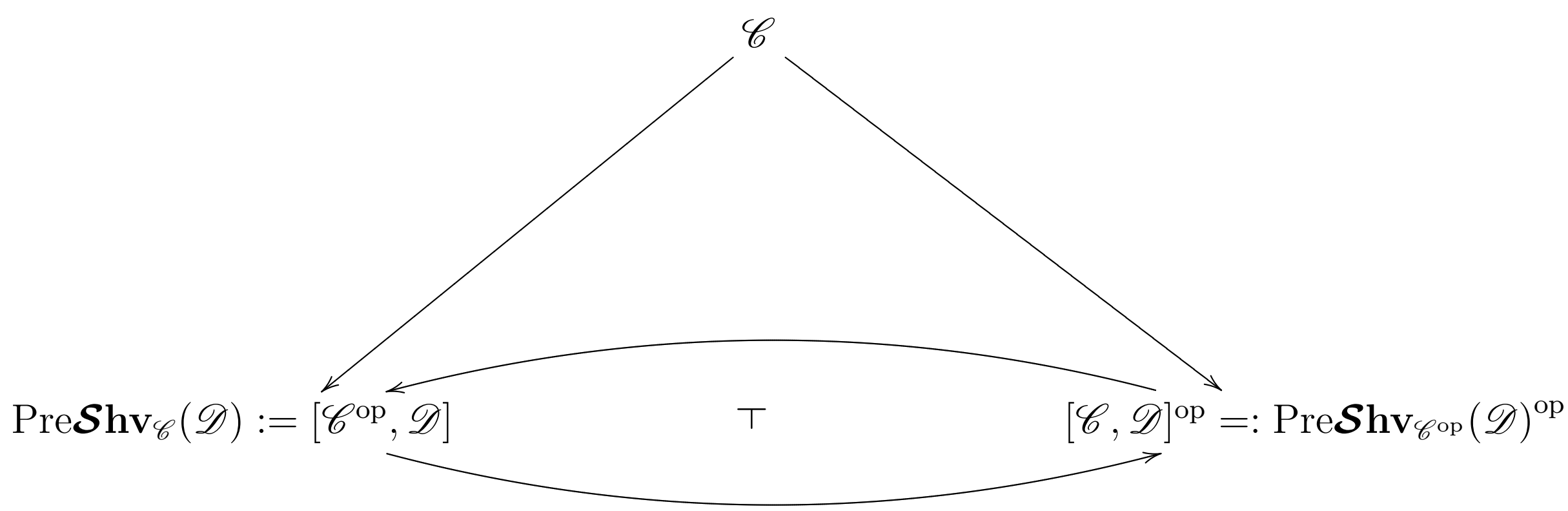

But in the case of sheaves instead of presheaves, this $\infty$-adjunction is an equivalence. □

We then have the following equivalence which ensures the compatibility of the definition of the dual of $\mathfrak{X}$:

$$\mathscr{O}_{\mathfrak{X}^{\star}} \simeq \mathscr{O}_{\mathfrak{X}}^{\mathrm{op}} \tag{2.5.34}$$

Let us denote now indifferently:

$$\mathscr{O}_{\mathfrak{X}}^{\star} := \mathscr{O}_{\mathfrak{X}}^{\mathrm{op}}$$

### 2.5.5 Artin-Lurie's representation theorem and duality.

The question of the suitable writing of the existence of a solution for the dual is reformulated here by saying that in the classical setting, the datum of an affine scheme $\mathscr{X}$ is exactly the datum of the functor of points $\mathcal{F}_{\mathscr{X}} : R \to \mathrm{Hom}(\mathsf{Spec}\, R, \mathscr{X})$ where $R$ is a classical ring in $\mathcal{CA}\mathrm{lg}_{\mathbb{Z}}^{\heartsuit}$. The scheme $\mathscr{X}$ is determined by the functor $\mathcal{F}_{\mathscr{X}}$ up to canonical isomorphism.

In the spectral setting, we know that a spectral affine scheme $\mathfrak{X}$ is exactly the datum of the functor of points $\mathcal{F}_{\mathfrak{X}} : R \mapsto \mathrm{Hom}(\mathsf{Sp\acute{e}t}\, R, X)$ where $R$ is an $\mathbb{E}_{\infty}$-ring in $\mathcal{CA}\mathrm{lg}$. The construction $\mathfrak{X} \mapsto \mathcal{F}_{\mathfrak{X}}$ determines a fully faithful embedding from the $\infty$-category of spectral schemes to the $\infty$-category $\mathcal{F}\mathrm{un}(\mathcal{CA}\mathrm{lg}^{\mathrm{cn}}, \mathcal{S})$.

Looking here for a dual coaffine spectral Deligne-Mumford stack $\mathfrak{X}^{\star}$ to $\mathfrak{X}$ is then exactly looking for a spectral Deligne-Mumford stack that represents the functor $\mathcal{F}^{\star}$ :

$\mathcal{CA}\mathrm{lg}^{\mathrm{cn}} \to \mathcal{S}^{\mathrm{op}}$ which is the composition of the functor $\mathcal{F}_{\mathfrak{X}} : \mathcal{CA}\mathrm{lg}^{\mathrm{cn}} \to \mathcal{S} : R \mapsto \mathrm{Hom}(\mathsf{Sp\acute{e}t}\, R, \mathfrak{X})$ and the dualizing functor $\mathcal{S} \xrightarrow{i} \mathcal{S}^{\mathrm{op}}$ of subsubsection 2.5.3.1, or equivalently the composition of the functor $\mathcal{F}^{\star}_{\mathfrak{X}} : \mathcal{CA}\mathrm{lg}^{\mathrm{cn}} \to \mathcal{S}^{\mathrm{op}} : R \mapsto \mathrm{Hom}((\mathsf{Sp\acute{e}t}\, R)^{\star}, \mathfrak{X})$ and the dualizing functor $\mathcal{S}^{\mathrm{op}} \xrightarrow{i} \mathcal{S}$.

One way is to prove that the conditions of Artin–Lurie's theorem are fulfilled and that the functor $\mathscr{F}^{\star} : \mathcal{CA}\mathrm{lg}^{\mathrm{cn}} \to \mathcal{S}^{\mathrm{op}}$ dualizing the functor yet to define $: \mathscr{F} \mathcal{CA}\mathrm{lg}^{\mathrm{cn}} \to \mathcal{S}^{\mathrm{op}}$ is well-represented by a coaffine geometric stack which will be the expected dual coaffine geometric stack.

But, among these conditions, each of the conditions 1), 2), 3), 4) 5) must be examined in détail. The first one is:

The fourth non-trivial item is the following:
4)The natural transformation $X \to \mathsf{Sp\acute{e}t}^{\mathrm{f}}\, R$ admits a connective cotangent complex denoted $L_{X/\mathsf{Sp\acute{e}t}^{\mathrm{f}}\, R}$.

That is, what we need here, is that the functor $\mathfrak{X}^{\star} : \mathcal{CA}\mathrm{lg}^{\mathrm{cn}} \to \mathcal{S}$ involved in Artin–Lurie's representation theorem admits a connective cotangent complex $L_{\mathfrak{X}^{\star}/\mathsf{Sp\acute{e}t}^{\mathrm{f}}\, R}$ for each $\mathsf{Sp\acute{e}t}^{\mathrm{f}}\, R$.

But, following subsection 2.1.1, it suffices to take the dual (coconnective) simplicial resolution given by Verdier duality (Equation 2.5.25) from the coconnective dual $\mathfrak{X}^{\star}$ given by reconstruction from $\mathsf{QCoh}(\mathfrak{X})^{\mathrm{cn}}$.

In fact, more easily, $\mathfrak{X}^{\star}$ is the dual of a (coaffine geometric) perfect stack $\mathfrak{X}$, so it is perfect, geometric, and coaffine, since in particular it is of finite cohomological dimension thanks to proposition 2.2.112.

Conclusion : The definition 2.5.14 of the dual coaffine geometric stack $\mathfrak{X}^{\star}$ of a coaffine geometric stack $\mathfrak{X}$ is then fully justified.

remark 2.5.35 : In the case of the dual Amplituhedron, the spectral Deligne-Mumford stacks that represent these functors will be respectively the duals $\mathfrak{G}^{\star}$ and $\mathfrak{A}^{\star}$.

### 2.5.6 The extended Brauer-Grothendieck group: Fourth step.

We will follow the Chapter 11 of [103]. We can also use the foundational [147], as well as [6, 7, 27, 28] for a deep dive, and [31] for a survey.

In definition 2.5.14, we conditioned the existence of a dual to the coaffine geometric stack $\mathfrak{X}$ to the capacity of taking a dual $\mathfrak{X}^{\star}$ in Equation 2.5.13, as well as the existence of a monoidal dual to $\left(\mathsf{QCoh}(\mathfrak{X}), \otimes, \mathscr{O}_{\mathfrak{X}}\right)$ in definition 2.5.17, and finally we concluded with Verdier's duality in Equation 2.5.24, which guaranted the existence of this dual.

But the subject is more general and can be recalled there in a few words:

To cite [103], Chapter 11, page 957: "*One way to guarantee that a (compactly generated) R-linear $\infty$-category $\mathcal{C}$ be smooth and proper is to assume that it is* invertible: *that is, that there exists another R-linear $\infty$-category $\mathcal{C}^{-1}$* such that $\mathcal{C} \otimes_R \mathcal{C}^{-1} \simeq \mathsf{Mod}_R$. *The collection of all equivalences classes of (compactly generated) invertible R-linear $\infty$-categories can be organized into an abelian group $\mathrm{Br}^{\dagger}(R)$, which we will refer to as the* extended Brauer group of *$R$.* [...] *In particular, we will show that every invertible R-linear $\infty$-category is locally equivalent to $\mathsf{Mod}_R$ with respect to the étale topology (Theorem 11.5.5.1) and that $\mathrm{Br}^{\dagger}(R)$ depends only on the commutative ring $\pi_0 R$ (Proposition 11.5.5.6).*"

Let us be inspired by this extract, and let us take $R := \mathsf{Groth}_{\infty}$. Let us define the tensor product of categories and the monoidal "*inversion*":

$$\mathfrak{X} \otimes_{\mathsf{Groth}_{\infty}} \mathfrak{X}^{\star} \underset{\text{htpy}}{\simeq} \mathbf{1}_{\mathsf{Groth}_{\infty}} \underset{\text{htpy}}{\simeq} \mathsf{Mod}_{\mathrm{Sp}} \tag{2.5.36}$$

recalling that the "*inversion*" must occur in the smaller category $\mathsf{Mod}_{\mathsf{QCoh}(\mathfrak{X})^{\mathrm{cn}}}(\mathsf{Groth}_{\infty})$.

DEFINITION 2.5.37 ([103] Def. 11.5.3.1) : *Let $R$ be a connectivve $\mathbb{E}_{\infty}$-ring and let $A$ be an $\mathbb{E}_1$-algebra over $R$. We say that $A$ is a* derived Azumaya algebra over *$R$ if it satisfies the following conditions:*

- *The algebra $A$ is a compact generator of* $\mathsf{Mod}_R$

- *The left and right actions of $A$ on itself induce an equivalence $A \otimes A^{\mathrm{op}} \to \mathrm{End}_R(A)$*

PROPOSITION 2.5.38 ([103] Prop.11.5.3.3) : *We have:*

- *A derived Azumaya algebra is connective*
- *A derived Azumaya algebra is a locally free of finite rank as an $R$-module*

If we denote $\mathrm{RMod}_A$, the $\infty$-category of right $A$-modules, and $\mathrm{LinCat}_R^{\mathrm{St}}$, the $\infty$-category of $R$-linear stable $\infty$-category, we have:

PROPOSITION 2.5.39 ([103] Prop. 11.5.3.4) : *Let $R$ be a connective $\mathbb{E}_\infty$-ring et let $A \in \mathcal{A}\mathrm{lg}_R$. Then $A$ is a derived Azumaya algebra over $R$ if and only if the stable $R$-linear $\infty$-category $\mathrm{RMod}_A$ is an invertible object of $\mathrm{LinCat}_R^{\mathrm{St}}$.*

REMARK 2.5.40 : We highlight the importance of this proposition in the existence of an invertible category of modules in the problem that concerns us. This is exactly the point around which we turn, and the condition we require.

Let us state now Theorem 11.5.5.1 of [103]:

THEOREM 2.5.41 ([103] Thm. 11.5.5.1) : *Let $\mathfrak{X}$ be a spectral Deligne-Mumford stack, and let $u \in \mathrm{Br}^\dagger(\mathfrak{X})$. Then there exists an étale surjection $f : \mathfrak{U} \to \mathfrak{X}$ such that $f^*u = 0$ in $\mathrm{Br}^\dagger(\mathfrak{U})$.*

For this purpose of the existence of an invertible category of modules, let us quote an example of Chapter 11 of [103] which will be of importance for us:

PROPOSITION 2.5.42 ([103] Ex. 11.5.5.5) : *Let $\mathfrak{X} := (\mathfrak{X}, \mathcal{O}_{\mathfrak{X}})$ be a spectral Deligne-Mumford stack, and suppose that the structure sheaf $\mathcal{O}_{\mathfrak{X}}$ is 0-truncated. Then we can identify $\mathcal{O}_{\mathfrak{X}}$ with a commutative ring object of the topos of discrete objects of $\mathfrak{X}$. Let $\mathcal{O}_{\mathfrak{X}}^\times$ denote its group of units. After several considerations for which we refer to [103], we obtain, for a quasi-compact, quasi-separated spectral algebraic space which is normal and locally Noetherian, an isomorphism:*

$$\mathrm{Br}^\dagger(\mathfrak{X}) \simeq H^2(\mathfrak{X}\,; \mathcal{O}_{\mathfrak{X}}^\times) \tag{2.5.43}$$

And let us state Proposition 11.5.5.6 of [103]:

PROPOSITION 2.5.44 ([103] Thm. 11.5.5.6) : *Let $R$ be a connective $\mathbb{E}_\infty$-ring. Then the canonical map $\mathrm{Br}^\dagger(R) \to \mathrm{Br}^\dagger(\pi_0 R)$ is an isomorphism.*

Let us follow the study of Chapter 11 of [103], and let us state the following:

DEFINITION 2.5.45 ([103] Def. 11.5.7.1) : *Let $\mathfrak{X} : \mathcal{CA}\mathrm{lg}^{\mathrm{cn}} \to \mathcal{S}$ be a functor and let $\mathsf{QStk}^{\mathrm{PSt}}(\mathfrak{X})$ denote the $\infty$-category of prestable quasi-coherent stacks on $\mathfrak{X}$. We let $\mathscr{B}\mathrm{r}(\mathfrak{X})$ denote the full subcategory of $\mathsf{QStk}^{\mathrm{PSt}}(\mathfrak{X})^{\simeq}$ spanned by the invertible quasi-coherent stacks $\mathcal{C} \in \mathsf{QStk}^{\mathrm{PSt}}(\mathfrak{X})$ for which $\mathcal{C}$ and $\mathcal{C}^{-1}$ are compactly generated. We will refer to $\mathscr{B}\mathrm{r}(\mathfrak{X})$ as the* Brauer space *of $\mathfrak{X}$. We let $\mathrm{Br}(\mathfrak{X})$ denote the set $\pi_0\mathscr{B}\mathrm{r}(\mathfrak{X})$, which we will refer to as the* Brauer group *of $\mathfrak{X}$.*

REMARK 2.5.46 ([103] Rem. 11.5.7.2) : Let $\mathfrak{X} : \mathcal{CA}\mathrm{lg}^{\mathrm{cn}} \to \mathcal{S}$ be a functor. Then the Brauer space $\mathscr{B}\mathrm{r}(\mathfrak{X})$ admits the structure of a symmetric monoidal $\infty$-category in which every object and every morphism is invertible. It is then a $\mathbb{E}_\infty$-space, and the set $\mathrm{Br}(\mathfrak{X}) = \pi_0\mathscr{B}r(\mathfrak{X})$ is an abelian group.

Now, let $\mathfrak{X} : \mathcal{CA}\mathrm{lg}^{\mathrm{cn}} \to \mathcal{S}$ be a functor and let $\mathscr{A} \in \mathcal{A}\mathrm{lg}(\mathsf{QCoh}(\mathfrak{X})^{\mathrm{cn}})$ be a connective Azumaya algebra. Then the construction:

$$(\eta \in \mathfrak{X}(R)) \longmapsto (\mathrm{RMod}^{\mathrm{cn}}_{\mathscr{A}_\eta} \in \mathrm{LinCat}^{\mathrm{St}}_R)$$

determines a compactly generated prestable quasi-coherent stack on $\mathfrak{X}$. It is an invertible object of $\mathrm{LinCat}^{\mathrm{St}}_R$ for an analoguous argument of that of proposition 2.5.39, whose inverse is compactly generated, and can therefore be identified with a point of the Brauer space $\mathscr{B}\mathrm{r}(\mathfrak{X})$. We let $[\mathscr{A}] \in \mathrm{Br}(\mathfrak{X})$ denote the equivalence class of this quasi-coherent stack. We will refer to $[\mathscr{A}]$ as the Brauer class of $\mathscr{A}$. Then, we can state the principal results of this section, which give us the condition for a coaffine geometric stack not to be autodual, and will be used further in the case of the Amplituhedron.

PROPOSITION 2.5.47 ([103] Prop. 11.5.7.9) : *Let $\mathfrak{X} : \mathcal{CA}\mathrm{lg}^{\mathrm{cn}} \to \mathcal{S}$ be a functor and let $\mathscr{A} \in \mathcal{A}\mathrm{lg}(\mathsf{QCoh}(\mathfrak{X}))$. The following conditions are equivalent:*

- *The object* $\mathscr{A}$ *is a connective Azumaya algebra and the Brauer class* $[\mathscr{A}] \in Br(\mathfrak{X})$ *vanishes.*
- *There exists an equivalence* $\mathscr{A} \simeq \mathrm{End}(\mathscr{E}) = \mathscr{E} \otimes \mathscr{E}^{\vee}$*, where* $\mathscr{E} \in Vect(\mathfrak{X})$ *is a vector bundle on* $\mathfrak{X}$ *of rank* $> 0$*.*

THEOREM 2.5.48 ([103] Thm. 11.5.7.11) : *Let* $\mathfrak{X}$ *be a spectral Deligne-Mumford stack, and let* $u \in Br(\mathfrak{X})$*. Then there exists an étale surjection* $f : \mathfrak{U} \longrightarrow \mathfrak{X}$ *such that* $f^*u = 0$ *in* $\mathrm{Br}^{\dagger}(\mathfrak{U})$*.*

COROLLARY 2.5.49 ([103] Cor. 11.5.7.12) : *Let* $\mathfrak{X}$ *be a spectral Deligne-Mumford stack, and let* $\mathscr{A} \in \mathcal{A}\mathrm{lg}(\mathsf{QCoh}(\mathfrak{X}))$ *be a connective Azumaya algebra on* $\mathfrak{X}$*. Then there exists an étale surjection* $f : \mathfrak{U} \longrightarrow \mathfrak{X}$ *and an equivalence* $f^*\mathscr{A} \simeq \mathrm{End}(\mathscr{E})$*, where* $\mathscr{E} \in Vect(\mathfrak{X})$ *is a vector bundle on* $\mathfrak{X}$ *of rank* $> 0$*.*

Let us conclude this section with the example 11.5.7.15 of [103] in case of a 0-truncation, which will be of importance for us:

PROPOSITION 2.5.50 ([103] Ex. 11.5.7.15) : *Let* $\mathfrak{X} := (\mathfrak{X}, \mathscr{O}_{\mathfrak{X}})$ *be a spectral Deligne-Mumford stack, and suppose that the structure sheaf* $\mathscr{O}_{\mathfrak{X}}$ *is 0-truncated. Then we can identify* $\mathscr{O}_{\mathfrak{X}}$ *with a commutative ring object of the topos of discrete objects of* $\mathfrak{X}$*. Let* $\mathscr{O}_{\mathfrak{X}}^{\times}$ *denote its group of units. After several considerations for which we refer to [103], we obtain, for a quasi-compact, quasi-separated spectral algebraic space which is normal and locally Noetherian, an isomorphism:*

$$\mathrm{Br}(\mathfrak{X}) \simeq H^2(\mathfrak{X}\,;\mathscr{O}_{\mathfrak{X}}^{\times}) \tag{2.5.51}$$

*We then have an isomorphism between the Brauer group of* $\mathfrak{X}$ *and the cohomological Brauer group of* $\mathfrak{X}$*:*

$$\mathrm{Br}(\mathfrak{X}) \simeq \mathrm{Br}^{\dagger}(\mathfrak{X}) \tag{2.5.52}$$

with its corollary:

COROLLARY 2.5.53 ([103] Cor. 11.5.7.16) : *Let* $R$ *be a connective* $\mathbb{E}_{\infty}$*-ring. Then the canonical map:*

$$\mathrm{Br}(R) \longrightarrow \mathrm{Br}(\pi_0 R) \tag{2.5.54}$$

*is an isomorphism of abelian groups.*

Conclusion: This equivalence between the spectral world of Brauer groups and the classical world of Brauer groups of 0-truncation of $\mathbb{E}_\infty$-rings will be of importance for us in the reconstruction of a concrete topological space for the dual Amplituhedron.

## 2.6. Decomposition theorem and perverse geometry on spectral DM stacks.

We assume that notions of Whitney stratification and of perverse sheaves are known (see for example [66, 1, 8]).

We just recall that an open subset $\mathcal{U}$ is called *admissible* if it is included in a Whitney strata and that the codimension of the $\Sigma_{n-1}$-stratum in $\Sigma_n = X$ is $\geqslant 2$, which is the case as long as we are in $\mathbb{C}$-projective geometry.

### 2.6.1 Resolutions and complexes of intersection.

Definition 2.6.1 : *Let $\mathscr{X}$ be an algebraic scheme. If $(I^\bullet, \mathrm{d}^\bullet)$ is an injective, flabby, or fine* resolution *of a complex $K^\bullet$, then $\mathcal{H}^i(K^\bullet) \simeq \mathcal{H}^i(I^\bullet) := \mathrm{Ker}(\mathrm{d}^i)/\mathrm{Im}(\mathrm{d}^{i-1}))$ is the cohomology of the complex $K^\bullet$ with injective, flabby, fine components. We define then the* cohomology groups *$H^\bullet(\mathscr{X}, K^\bullet)$ of $K^\bullet$ to be the cohomology groups of the complex of global sections $\Gamma(\mathscr{X}, I^\bullet)$ of $\mathscr{X}$ with coefficients in $I^\bullet$. It is the* cohomology of the scheme *$\mathscr{X}$ with coefficients in the complex of sheaves $I^\bullet$.*

Remark 2.6.2 : We will come back later on the various cohomology theories we can apply to the scheme $\mathscr{X}$ or to a spectral stack in general. Among various problems are the followings: the problem of the coefficients which are considered, the problem of smoothness of the scheme $\mathscr{X}$ or of the geometric stack $\mathfrak{X}$ or of the spectral Deligne-Mumford stack or algebraic space $\mathfrak{X}$ involved, and then the problem of singulariries, the kind of geometric object we look at: is it a classical topological space, a Kan complex, a CW-complex, an $\mathbb{E}_\infty$-ring spectrum, a classical scheme, a spectral scheme, a topos, an $\infty$-topos, a cohesive topos, a differentiel (infinitesimal) topos, or even more a $\mathcal{G}$-structured corporeal Deligne-Mumford $\infty$-stack... The types of objects and the types of cohomology theories will play a central role. But allow us to remain somewhat vague for the moment.

REMARK 2.6.3 : For the definition of the cohomology of a single sheaf $\mathcal{I}$, we identify the sheaf with the complex of sheaves $I^\bullet$ concentrated in degree 0, that is with 0 in each degree except $\mathcal{I}$ in degree 0.

DEFINITION 2.6.4 : *Let $Z$ be a closed subscheme of the scheme $\mathscr{X}$. A* quasi-isomorphism $K^\bullet \xrightarrow[\simeq]{\text{quasi-iso}} L^\bullet$ *induces isomorphisms on the cohomology,* $\forall i \in \mathbb{Z}, H^i(U, K^\bullet) \simeq H^i(U, L^\bullet)$ *of any open subscheme $U \subset Z$ and these isomorphisms are compatible with the maps induced by inclusions and with Mayer Vietoris sequences.*

DEFINITION 2.6.5 : *Let $\mathscr{X}$ be an algebraic scheme. The* bounded, bounded below, bounded above, unbounded derived category*, denoted* $\mathfrak{D}^{\mathrm{b}}(\mathscr{X})$ , $\mathfrak{D}^{+}(\mathscr{X})$ , $\mathfrak{D}^{-}(\mathscr{X})$, $\mathfrak{D}(\mathscr{X})$, *are the categories of bounded, bounded below, bounded above, unbounded complexes of quasi-coherent sheaves* $\mathsf{QCoh}(\mathscr{X})$ *on the scheme $\mathscr{X}$ to which we formally invert quasi-isomorphisms (it is a localization) after taking homotopic equivalence classes.*

REMARK 2.6.6 : We will pay attention to the fact that the constant sheaf is not an injective object in the category of quasi-coherent sheaves. This means that we can not compute its cohomology without finding first an injective resolution. We do not even then have right now yet the cohomology of the point with coefficients in the constant sheaf. This justifies the introduction of local systems in the next paragraph.

DEFINITION 2.6.7 : *Let* $K^\bullet \xrightarrow{\text{inj.}} I^\bullet$ *be an injective resolution, which is then a* quasi-iso morphism*, we define and denote* $\mathbb{R}\Gamma(K^\bullet) := \Gamma(Z, I^\bullet)$ *the* (right) derived functor *of the functor* $\Gamma$ *of global sections.*

REMARK 2.6.8 : In the perspective of a resolution by (injective) local systems in the next paragraph, then for the typography of the Decomposition theorem, let us remark that if $f : Z \longrightarrow W$ is a continuous map between algebraic varieties, and $I^\bullet$ is a complex of injective sheaves on $Z$, then the pushforward complex $f_*(I^\bullet)$ satisfies:

$$H^i(U, f_*(I^\bullet)) \simeq H^i(f^{-1}(U), I^\bullet)$$

However, if $I^\bullet$ is no more injective, the isomorphism fails, and we have to replace $I^\bullet$ first by an injective resolution before passing to the derived category. This will be what we do:

DEFINITION 2.6.9 : *We call* (derived) direct image functor *and denote* $\mathbb{R}f_*(-)$ *such that for each* $K^\bullet$*, complex of sheaves on a closed subscheme* $Z$ *of a scheme* $\mathscr{X}$*, we have* $\mathbb{R}f_*K^\bullet := \Gamma(Z, K^\bullet) = \Gamma(Z, I^\bullet)$ *where* $K^\bullet \xrightarrow{\text{inj.}} I^\bullet$ *is an injective resolution of* $K$*.*

REMARK 2.6.10 : If we map $Z$ to the point: $f : Z \longrightarrow \{*\}$, then $\mathbb{R}f_! \simeq \mathbb{R}f_*K^\bullet = \mathbb{R}\Gamma(K^\bullet)$ and $\mathbb{R}^i f_* K^\bullet = H^i(Z, K^\bullet)$.

By a similar process, we define the (derived) direct image functor and denote $\mathbb{R}f_!(-)$if $f$ is a proper map, since $\mathbb{R}f_! \simeq \mathbb{R}f_*$, and thanks to the replacement by an injective resolution made in the same way, and so on for $f^*$ and $f^!$ (in our special case, a section dedicated to six-functor formalism would not have brought any extra informations else the ones that are in the section subsection 2.10.1 or scattered throughout in the text).

### 2.6.2 Constructible sheaves, Local systems.

DEFINITION 2.6.11 : *Let* $Z$ *be a complex algebraic variety. We say that a subset* $V \subset Z$ *is* constructible *if it is obtained from a finite sequence of unions, intersections or complements of algebraic subvarieties of* $Z$*. A* local system *is a locally constant sheaf, denoted* $\mathcal{L}_Z$ *with finite dimensional stalks. A complex of sheaves* $K^\bullet$ *is* constructible *or with* constructible cohomology sheaves *if there exists a decomposition* $Z = \coprod_\alpha Z_\alpha$ *into finitely many constructible subsets such that each of the cohomology sheaves* $\mathcal{H}^i(K^\bullet)$ *is locally constant, so a local system* $\mathcal{L}_{Z_\alpha}$ *along each* $Z_\alpha$ *with finite dimensional stalks.*

REMARK 2.6.12 : This implies that the limit in the right hand part of the following:

$$\mathcal{H}^i_x(K^\bullet) := \varinjlim H^i(U_x, K) \tag{2.6.13}$$

is attained by any "*regular*" neighborhood $U_x$ of the point $x$. For example, one may embed (locally) $Z$ into a manifold and take $U_x := Z \cap B_\epsilon(x)$ to be the intersection of $Z$ with a sufficiently small ball centered at $x$. It also implies that $H^i(Z, K)$ is finite dimensional.

DEFINITION 2.6.14 : *Let $U \subset Z$ be a nonsingular Zariski open subscheme of the algebraic closed subscheme $Z$ of the scheme $\mathscr{X}$ and let $\mathcal{L}_U$ be a local system on $U$: it is a locally constant sheaf on $U$ with finite dimensional stalks. The* intersection complex *$IC_Z(\mathcal{L})$ is a complex of sheaves on $Z$, which extends the complex $\mathcal{L}_U[\dim Z]$ on $U$, uniquely determined up to quasi-isomorphism by the following support conditions, where for a local system $\mathcal{L}_U$ on a nonsingular Zariski open subset $U$, there exists a constructible complex of sheaves $IC(\mathcal{L}) \in \mathfrak{D}_Y$, unique up to canonical isomorphism in $\mathfrak{D}_Y$ such that $IC(\mathcal{L})_{|U} \simeq \mathcal{L}_U$:*

- $dim_{\mathbb{C}}\{y \in Y | \mathcal{H}^i_y(IC(\mathcal{L})) \neq 0\} < -i, \text{ if } i > -n$
- $dim_{\mathbb{C}}\{y \in Y | \mathcal{H}^i_{c,y}(IC(\mathcal{L})) \neq 0\} < -i, \text{ if } i > -n$ *(compact support)*

*where, for any complex $S$ of sheaves, $\mathcal{H}^i_{c,y}(S) = \varprojlim H^i_c(U_y, S)$ is the local compactly supported cohomology at $x$. The intersection complex $IC_Z(\mathcal{L})$ is sometimes called the* middle extension *of $\mathcal{L}$. Its shifted cohomology is the intersection cohomology $IH^{n+\bullet}(Y, \mathcal{L}) := H^\bullet(Y, IC(\mathcal{L}))$. More precisely:*
*The cohomology of this intersection complex is named the* intersection cohomology*, defined and denoted $IH^i(Z, \mathcal{L}_Z) := H^{i-\dim Z}(Z, \mathcal{L}_Z)$. The intersection complexes of local systems supported on open sets of subvarieties of $Z$ are the building blocks of the abelian subcategory of* perverses sheaves *in the category of constructible complexes of sheaves, verifying some condition on the support of its cohomology sheaves and that of its cohomology with compact support sheaves. The category of perverse sheaves is Artinian, that is the iteration of extensions is of finitely many simple perverse sheaves. The* simple perverse sheaves *are the intersection complexes $IC_Y(\mathcal{L}_U)$ of* irreductible subvarieties *$Y \subset Z$ and* simple local systems *$\mathcal{L}_U$, defined on a nonsingular Zariski open subset $U \subset Y$. The constructible sheaves, thought as complexes* *concentrated in degree 0, form an abelian subcategory of the constructible derived category $\mathfrak{D}_Z$. An object $K$ in $\mathfrak{D}_Z$ is isomorphic to an object of this subcategory if and only if $\mathcal{H}^i(K)$ if $i \neq 0$.*

REMARK 2.6.15 : If $U$ is a nonempty, nonsingular and pure dimensional open subset of $Z$ on which all the cohomology sheaves $\mathcal{H}^i(K)$ are local systems, then the restriction to $U$ of $\mathcal{H}^m(K)$ and $\mathcal{H}^{m-\dim Z}(K)[\dim Z]$ coincide. In general, the two differ. For example, $\mathcal{H}^0(\mathbb{R}f_*\mathbb{Q}_X[2]) = IC_Y(R^1) \oplus T_\Sigma$, and not only $IC_Y(R^1)$ if $f : X \longrightarrow Y$ is a projective map with connected fibers from a smooth surface $X$ onto a smooth curve $Y$, with $\Sigma \subseteq Y$ is the finite set of critical values and $U = Y \backslash \Sigma$ its complement, and $R^1 = (\mathbb{R}^1 f_*\mathbb{Q}_X)_{|U}$ is the local system on $U$ with stalk the first cohomology of the typical fiber. We then have indeed a non canonical isomorphism

$$\mathbb{R}f_*\mathbb{Q}_X[2] \simeq \mathbb{Q}_Y[2] \oplus (IC_Y(R^1) \oplus T_\Sigma) \oplus \mathbb{Q}_Y$$

where $T_\Sigma$ is a skycreaper sheaf over $\Sigma$ with stalks $T_s = H_2(f^{-1}(s))/ < [f^{-1}(s)] >$ at every $s \in \Sigma$ and $\mathcal{H}^m(K)$ and $\mathcal{H}^{m-\dim Z}(K)[\dim Z]$ do not coincide.

### 2.6.3 Statement of the decomposition theorem.

THEOREM 2.6.16 ((Decomposition Theorem) De Cataldo-Migliori, Thm 1.14) : *Let $f : X \longrightarrow Y$ be a proper map. There is a noncanonical isomorphism in the constructible derived category $\mathfrak{D}_Y$ of sheaves on $Y$:*

$$\mathbb{R}f_*IC_X \simeq \bigoplus_i \mathcal{H}^i(\mathbb{R}f_*IC_X)[-i] \tag{2.6.17}$$

*Furthermore, the perverse sheaves $\mathcal{H}^i(\mathbb{R}f_*IC_X)$ are semi-simple, i.e. there is a decomposition into finitely many disjoint locally closed and nonsingular varieties $Y = \coprod S_\beta$ and a canonical decomposition into a direct sum of intesection complexes of semisimple local systems*

$$\mathcal{H}(\mathbb{R}f_*IC_X) = \bigoplus_i IC_{\overline{S_\beta}}(\mathcal{L}_\beta)[-i] \tag{2.6.18}$$

*The first Chern class of a line bundle $\eta$ on $X$ yields, for $i \geqslant 0$, maps:*

- $\eta^i : \mathbb{R}f_*IC_X \longrightarrow \mathbb{R}f_*IC_X[2i]$
- $\eta^i : \mathcal{H}^{-i}(\mathbb{R}f_*IC_X) \longrightarrow \mathcal{H}^i(\mathbb{R}f_*IC_X)$

THEOREM 2.6.19 ((Relative Hard Lefschetz Theorem), De Cataldo-Migliori Thm 1.15) : *Assume that $f$ is projective and $\eta$ the first Chern class of a line bundle on $X$ whose restriction to each fiber is ample. Then we have isomorphisms:*

$$\eta^i : \mathcal{H}^{-i}(\mathbb{R}f_*IC_X) \xrightarrow{\simeq} \mathcal{H}^i(\mathbb{R}f_*IC_X)$$

As corollary of the Decomposition Theorem, we have:

COROLLARY 2.6.20 ((The Global and Local Invariant Cycle Theorems)) : *Let $f : X \longrightarrow Y$ be a proper map. Let $U \subseteq Y$ be a Zariski open subset on which the sheaf $\mathbb{R}f_*IC_X$ is locally constant.*

- *The natural restriction map*

$$IH^i(X) \longrightarrow H^0(U, \mathbb{R}f_*IC_X) \quad \text{is surjective} \tag{2.6.21}$$

- *Let $u \in U$ and $B_u \subseteq U$ be the intersection with a sufficiently small Euclidian ball centered at $u$. Then the natural restriction/retraction map:*

$$H^i(f^{-1}(u), IC_X) = H^i(f^{-1}(B_u), IC_X) \longrightarrow H^0(B_u, \mathbb{R}f_*IC_X) \quad \text{is surjective} \tag{2.6.22}$$

DEFINITION 2.6.23 : *We call* stratification *for $f$ a decomposition of $Y$ into finitely many locally closed nonsingular subsets such that $f^{-1}(S_k) \longrightarrow S_k$ is a topologically locally trivial fibration. The subsets are called the* strata*.*

### 2.6.4 Semismall maps.

PROPOSITION 2.6.24 : *Let $X$ be a connected nonsingular $n$-dimensional variety, and $f : X \longrightarrow Y$ be a proper surjective map of varieties. Let $Y = \coprod_{k=0}^n S_k$ be a stratification for $f$. Let $y_k \in \S_k$ and set $d_k := \dim f^{-1}(y_k) = \dim f^{-1}(S_k) - \dim S_k$. The following are equivalent:*

- *$f_*\mathbb{Q}_X[n]$ is a perverse sheaf on $Y$*
- $\dim X \times_Y X \leqslant n$
- $\dim S_k + 2d_k \leqslant \dim X$ *for* $k = 0, \ldots, n$

DEFINITION 2.6.25 : *A proper and surjective map $f$ satisfying one of the equivalent properties above is said to be* semi-small

REMARK 2.6.26 : A semismall map $f : X \longrightarrow Y$ must be finite over an open dense stratum in $Y$. Hence, semismall maps are generically finite. The converse is not true.

REMARK 2.6.27 : If the third condition is stronger (the inequality is strict), then the ma is said to be small. In this case, the perverse sheaf $f_*\mathbb{Q}_X[n]$ staisfies support and cosupport conditions for intersection homology. Hence, if $Y_0 \subseteq Y$ denotes a nonsingular dense open subset over which $f$ is a covering, then we have that $f_*\mathbb{Q}_X[n] = IC_Y\mathcal{L}$ where $\mathcal{L} = f_*\mathbb{Q}_{X|Y_0}$

In the case of a semismall map, it is possible to identify precisely the strata $S_k$ contributing to the Decomposition theorem with a non trivial summand $IC_{\overline{S_k}}(\mathcal{L}_k)$ as well as the local system $\mathcal{L}_k$ which turn out to have finite monodromy.

DEFINITION 2.6.28 : *Let $X, Y, S_k$ and $d_k$ defined as earlier. A stratum $S_k$ is said to be* relevant *if* $\dim S_k + 2d_k = \dim X$. *Since* $\dim X = \dim Y$, *a relevant stratum has even dimension. We denote $I_{\mathrm{rel}}$ the set of relevant strata.*

One can prove the following theorem, where $\mathrm{Irr}(\pi_1(S_k))$ denotes the set of irreducible representations of $\pi_1(S_k, y_k)$, and $\mathcal{L}_k = \bigoplus_{\mathrm{Irr}(\pi_1(S_k))} \mathcal{L}_k^\chi$ the isotypical decomposition of $\mathcal{L}_k$. The local systems $\mathcal{L}_k^\chi$ are the tensor product of the irreducible local system $\mathcal{L}^\chi$ associated with the representation $\chi$ with the vector space $V_k^\chi$, whos dimension is the multiplicity of the representation $\chi$ in $\rho_k$. We cal then write:

THEOREM 2.6.29 : *There is a canonical isomorphism in the category of perverse sheaves*

*on $Y$ such that:*

$$f_*\mathbb{Q}_X[n] = \bigoplus_{k\in I_{\text{rel}}} IC_{\overline{S_k}}(\mathcal{L}_k) = \bigoplus_{\substack{k \in I_{\text{rel}} \\ \chi \in \text{Irr}(\pi_1(S_k))}} IC_{\overline{S_k}}(\mathcal{L}_k^\chi) \simeq \bigoplus_{\substack{k \in I_{\text{rel}} \\ \chi \in \text{Irr}(\pi_1(S_k))}} IC_{\overline{S_k}}(\mathcal{L}^\chi) \otimes V_k^\chi \tag{2.6.30}$$

COROLLARY 2.6.31 : *Since the derived category $\mathfrak{D}_Y$ for $Y$ like above, the $\mathbb{Q}$-vector space $\text{End}_{\mathfrak{D}_Y}(f_*\mathbb{Q}_X[n])$ is endowed with an algebra structure, and we can show that if $f$ is a semismall map, this algebra is semi-simple:*

$$\text{End}_{\mathfrak{D}_Y}(f_*\mathbb{Q}_X[n]) \simeq \bigoplus_{\substack{k \in I_{\text{rel}} \\ \chi \in \text{Irr}(\pi_1(S_k))}} \text{End}(V_k^\chi) \tag{2.6.32}$$

### 2.6.5 $\mathbb{E}_\infty$-ring structure of intersection complexes.

Let $X$ be a topological space. Let us denote $\mathbf{E}_\mathbb{Q}$ the category of singular cochains $C^\bullet(X,\mathbb{Q})$ equipped with an $\mathbb{E}_\infty$-structure.

DEFINITION 2.6.33 ([26], Def. 2.13) : *We define the* derived global sections functor*:*

$$\mathbb{R}_{\mathbf{E}_\mathbb{Q}}\Gamma(X,-) : \text{Shv}(\mathbf{E}_\mathbb{Q}) \longrightarrow \text{Shv}(\mathbf{E}_\mathbb{Q})$$

*where $\Gamma(X,-)$ is the global sections functor.*

THEOREM 2.6.34 ([26], Thm 2.16) : *The assignment $X \mapsto \mathbb{R}_{\mathbf{E}_\mathbb{Q}}\Gamma(X,\underline{\mathbb{Q}}_X)$, where $\underline{\mathbb{Q}}_X$ denotes the constant sheaf on X, defines a contravariant functor from the category of Hausdorff paracompact and locally contractible topological spaces to the category of $\mathbb{E}_\infty$-rings. This functor is naturally quasi-isomorphic to the functor $C^\bullet(-,\mathbb{Q})$ defined by the complex of singular cochains with its $\mathbb{E}_\infty$-structure.*

THEOREM 2.6.35 ([26], Thm 5.9) : *If we define the derived global sections functors on the open admissible subsets $U$ for the Whitney stratification of a topological space $X$, $U \mapsto \mathbb{R}_{\mathbb{E}_\infty}\Gamma(U,\underline{\mathbb{Q}}_U)$, where $\underline{\mathbb{Q}}_U$ denotes the constant sheaf on $U$, we get the* intersection complex sheaf *of $\mathbb{E}_\infty$-rings, which will be denoted:*

$$\mathcal{IC}^{\bullet} : U \mapsto \mathcal{IC}^{\bullet}(U, \underline{\mathbb{Q}}_U)$$

*satisfying the set of Deligne's axioms for sheaves of perverse $\mathbb{E}_\infty$-rings.*

REMARK 2.6.36 : Moreover, the category of perverse sheaves of $\mathbb{E}_\infty$-rings, as an abelian subcategory of the derived category $\mathscr{D}(\mathsf{QCoh}(X))$, is spanned by the following family of simple objects: $\{\mathcal{IC}_{\bullet}(\overline{S_\alpha}, \mathcal{L}_{S_\alpha})[\ell_\alpha]\}_\alpha$, where the $\mathcal{L}_{S_\alpha}$ are irreductible local systems on the closure $\overline{S_\alpha}$ of connected strata $S_\alpha$, and the $[\ell_\alpha]$ shift the complexes $\mathcal{IC}_{\bullet}(\overline{S_\alpha}, \mathcal{L}_{S_\alpha})$.

### 2.6.6 Intersection complexes and the decomposition theorem.

Let $f : X \to Y$ be a projective map of quasi-projective varieties.

THEOREM 2.6.37 ([66], Thm. 8.4.3) : *The right derived functor $\mathbb{R} f_*$ of the induced map $f_*$ sends the category of perverse sheaves on $X$ on the category of perverse sheaves on $Y$. Moreover, we have a decomposition into simple objects in the category of perverse sheaves on $Y$ of the pushforward $\mathbb{R} f_* \mathcal{F}$ for any perverse sheaf $\mathcal{F}$ such that:*

$$\mathbb{R} f_* \mathcal{F} \simeq \bigoplus_\alpha \mathcal{IC}^{\bullet}_{\mathcal{L}_{S_\alpha}}[\ell_\alpha]$$

*where the $\mathcal{L}_{S_\alpha}$ are irreductible local systems on the closure $\overline{S_\alpha}$ of connected strata $S_\alpha$ of $Y$, and the $[\ell_\alpha] \in \mathbb{Z}$ shift the complexes $\mathcal{IC}^{\bullet}_{\mathcal{L}_{S_\alpha}} := \mathcal{IC}^{\bullet}(\overline{S_\alpha}, \mathcal{L}_{S_\alpha})$.*

#### 2.6.6.1 There are enough projective perverse sheaves.

Finally, we have to check that there are enough projective perverse sheaves to cover the $\infty$-topos determined by the $\mathcal{G}_{\mathrm{PERF}^{\mathrm{PERV}}}$-geometry on $X$. This is indeed the case, thanks to [29], with respect to the condition that the space $X$ has finitely many strata and so the category of local systems on each of these is a category of finite-dimensional modules.

## 2.7. The positroid-étale $\infty$-topology on the TNN Grassmannian DM stack.

### 2.7.1 Open gluings.

Let us remind ourselves the étale $\infty$-topology constructed for $\infty$-topoi in section 2.2. We are going to construct an étale $\infty$-topology for the TNN (real) Grassmannian, which we will call: positroid-étale $\infty$-topology for the Deligne-Mumford stacks on the TNN Grassmannian.
First of all, we have a disjoint covering of the TNN (real) Grassmannian by (open) quasiaffine subschemes (over $\mathbb{C}$):

$$\coprod_{f \in \mathscr{B}(k,n)} \overset{\circ}{\Pi}_f \longrightarrow \mathrm{Gr}^{\mathbb{R}}_{\geqslant 0}(k,n) \subseteq \mathrm{Gr}^{\mathbb{C}}(k,n)$$

We also know thanks to proposition 2.2.30 that the essential image of the Yoneda embedding:

$$\begin{array}{ccc} \mathbf{Sch} & \xhookrightarrow{\text{fully faithful}} & \boldsymbol{\mathcal{S}}\mathbf{hv}^{\mathbf{\acute{e}t}} \\ X & \longmapsto & \big(R \mapsto X(R)\big) \end{array}$$

consists of the sheaves $\mathscr{F}$ in $\boldsymbol{\mathcal{S}}\mathbf{hv}^{\mathbf{\acute{e}t}}$, where $\boldsymbol{\mathcal{S}}\mathbf{hv}^{\mathbf{\acute{e}t}}$, defined in definition 2.2.21, is the category of étale $\infty$-sheaves of spaces on $\boldsymbol{\mathcal{CA}}\mathbf{lg}^{\heartsuit}_{\mathbb{Z}}$ (the classical rings in $\boldsymbol{\mathcal{R}}\mathbf{ings}$) as the full subcategory of the $\infty$-category of functors $\boldsymbol{\mathcal{F}}\mathbf{un}^{\infty}(\boldsymbol{\mathcal{CA}}\mathbf{lg}^{\heartsuit}_{\mathbb{Z}}, \boldsymbol{\mathcal{S}})$, satisfying the étale topology $\infty$-sheaf condition for $\boldsymbol{\mathcal{CA}}\mathbf{lg}^{\heartsuit}_{\mathbb{Z}}$. Let us recall it here:

$$\coprod_{f \in \mathscr{B}(k,n)} \overset{\circ}{\Pi}_f \xrightarrow[\text{surj. (on alg. closed points)}]{\text{representable}} \mathscr{F}$$

with each $\mathring{\Pi}_f \hookrightarrow \mathscr{F}$ is an open immersion, in the sense that in the following pullback square:

$$\begin{array}{ccc} \coprod_{f\in\mathscr{B}(k,n)} \mathring{\Pi}_f & \longrightarrow & \mathsf{Spec} R \\ \downarrow & \lrcorner & \downarrow \\ \mathring{\Pi}_f & \xrightarrow{f} & \mathscr{F} \end{array} \tag{2.7.1}$$

each $\mathring{\Pi}_f$ is quasiaffine with $\mathring{\Pi}_f \hookrightarrow \mathsf{Spec}\ R$ is an open immersion.

Thanks to definition 2.2.34, we can then state the following:

DEFINITION 2.7.2 : *We define and denote* $\mathfrak{Chr}^{\geqslant 0}_{\text{ét}}$, *the* category of Deligne-Mumford stacks on the TNN Grassmannian $\mathrm{Gr}^{\geqslant 0}_{k,n}$ *as the sheaves* $\mathscr{F}$ *in* $\boldsymbol{\mathcal{S}}\mathbf{hv}^{\text{ét}}$, *such that there exists, for each* $\mathscr{F} \in \boldsymbol{\mathcal{S}}\mathbf{hv}^{\text{ét}}$, *a surjective étale map* $\coprod_{f\in\mathscr{B}(k,n)} \mathring{\Pi}_f \longrightarrow \mathscr{F}$.

THEOREM 2.7.3 : *We call an* étale cover of the sheaf $\mathscr{F}$ *in* $\boldsymbol{\mathcal{S}}\mathbf{hv}^{\text{ét}}$, *a surjective étale map* $\coprod_{f\in\mathscr{B}(k,n)} \mathring{\Pi}_f \longrightarrow \mathscr{F}$. *These etale covers for all* $\mathscr{F}$ *in the category* $\mathfrak{Chr}^{\geqslant 0}_{\text{ét}}$ *of Deligne-Mumford stacks on the TNN Grassmannian* $\mathrm{Gr}^{\geqslant 0}_{k,n}$ *define a* Grothendieck $\infty$-topology, *which we name the* positroid-étale $\infty$-topology for the category $\mathfrak{Chr}^{\geqslant 0}_{\text{ét}}$ of Deligne-Mumford stacks on the TNN Grassmannian $\mathrm{Gr}^{\geqslant 0}_{k,n}$.

*Proof.* Indeed, it defines well a Grothendieck topology, since, for each Deligne-Mumford stack $\mathscr{F}$ on the TNN Grassmannian $\mathrm{Gr}_{k,n}$, that is for each sheaf $\mathscr{F}$ in $\mathfrak{Chr}^{\geqslant 0}_{\text{ét}}$, by definition 2.2.37, we can denote $\mathscr{F}^{\text{ét}}$ the full subcategory of the slice (comma) category $\boldsymbol{\mathcal{S}}\mathbf{hv}^{\text{ét}}/\mathscr{F}$, such that if $\mathscr{G}$ is an object of $\mathscr{F}^{\text{ét}}$, there exists a natural transformation $\mathscr{G} \xrightarrow{\text{ét}} \mathscr{F}$ which is étale (and then representable). This full subcategory $\mathscr{F}^{\text{ét}}$ is small, and can then be turned into a Grothendieck $\infty$-category, by making it an $\infty$-site with the étale topology, where an étale cover of any object in $\mathscr{F}^{\text{ét}}$ is a family of surjective étale maps (see theorem 2.7.3).
We can then consider the subcategory of sheaves $\mathcal{S}\mathbf{hv}(\mathscr{F}^{\text{ét}})$ in the presheaves category $\mathcal{F}\text{un}((\mathscr{F}^{\text{ét}})^{\text{op}}, \mathcal{S})$, satisfying the étale $\infty$-sheaf condition. The category $\mathscr{F}^{\text{ét}}$ being equipped in an $\infty$-site, we then have defined a Grothendieck $\infty$-topology for the whole category $\mathfrak{Gr}^{\geqslant 0}_{\text{ét}}$. □

This subcategory of sheaves $\mathcal{S}\mathbf{hv}(\mathscr{F}^{\text{ét}})$ in the presheaves category $\mathcal{F}\text{un}((\mathscr{F}^{\text{ét}})^{\text{op}}, \mathcal{S})$, satisfying the étale $\infty$-sheaf condition is exactly what we had called a prototypal $\infty$-topos, and proposition 2.2.39 gives us that the $\infty$-category $\mathcal{S}\mathbf{hv}(\mathscr{F}^{\text{ét}})$ associated to the Deligne-Mumford stack $\mathscr{F}$ in $\mathcal{S}\mathbf{hv}^{\textbf{ét}}$ is an $\infty$-topos.

We then have the following:

PROPOSITION 2.7.4 : *The category* $\mathfrak{Gr}^{\geqslant 0}_{\text{ét}}$ *of Deligne-Mumford stacks on the TNN Grassmannian* $\text{Gr}^{\geqslant 0}_{k,n}$ *is an* $\infty$*-topos.*

Let us now consider a piece of the covering $\mathring{\Pi}_f \longrightarrow \mathscr{F}$ on the Deligne-Mumford stack $\mathscr{F}$ for $f \in \mathscr{B}(k,n)$. We can equip $\mathfrak{Gr}^{\geqslant 0}_{\text{ét}}$ with a sheaf of rings $\mathscr{O}_{\mathfrak{Gr}^{\geqslant 0}_{\text{ét}}}$ such that $\mathscr{O}_{\mathfrak{Gr}^{\geqslant 0}_{\text{ét}}}(\mathring{\Pi}_f)$ is the ring of coordinates of $\mathring{\Pi}_f$. The (affinization of) these rings are cluster algebras thanks to [50].
The couple $\left(\mathfrak{Gr}^{\geqslant 0}_{\text{ét}},\ \mathscr{O}_{\mathfrak{Gr}^{\geqslant 0}_{\text{ét}}}\right)$ becomes a spectral Deligne-Mumford stack in the following way:
Following subsection 2.3.6, let us denote $\mathcal{M}\text{od}$ the $\infty$-category of pairs $(A, M)$ where $A$ is an $\mathbb{E}_\infty$-ring and $M$ is an $A$-module spectrum. We can then identify thanks to **??**, a $\mathcal{M}\text{od}$-valued sheaf on $\mathfrak{Gr}^{\geqslant 0}_{\text{ét}}$ with a pair $(\mathscr{O}_{\mathfrak{Gr}^{\geqslant 0}_{\text{ét}}}, \mathcal{IC}_{\mathfrak{Gr}^{\geqslant 0}_{\text{ét}}})$, where $\mathcal{IC}_{\mathfrak{Gr}^{\geqslant 0}_{\text{ét}}}$ is the intersection complex sheaf of $\mathbb{E}_\infty$-rings of theorem 2.6.35.

We then have the following:

DEFINITION 2.7.5 : *We define the* $\infty$-category of spectral Deligne-Mumford stacks on the $\infty$-topos of the TNN Grassmannian $\mathfrak{Gr}^{\geqslant 0}_{\text{ét}}$ *by:*

$$\mathfrak{G} := \left(\mathfrak{Gr}^{\geqslant 0}_{\text{ét}},\ \mathscr{O}_{\mathfrak{Gr}^{\geqslant 0}_{\text{ét}}}, \mathcal{IC}_{\mathfrak{Gr}^{\geqslant 0}_{\text{ét}}}\right)$$

In the langage of étale sheaves on étale spectra, slice topoi and effective covers (see subsubsection 2.2.3.1), where $R_f = \mathbb{C}[\overset{\circ}{\Pi}_f]$ (see Equation 2.7.6), we have the following description:

The ring of coordinate $R_f$ of $\overset{\circ}{\Pi}_f$ is in $\mathcal{CAlg}^{\heartsuit}_{\mathbb{Z}}$, with $R_f := \mathbb{C}[\overset{\circ}{\Pi}_f]$. It can then be viewed as an $\mathbb{E}_\infty$-ring spectrum in $\mathcal{CAlg}$ concentrated in degree 0.

Let us define now $\mathfrak{X} := \mathsf{Spét}\, R_f$ the $\infty$-topos of the étale spectrum of $R_f$. On the $\infty$-subcategory $\mathcal{CAlg}^{\text{ét}}_{R_f}$ (definition 2.2.50), we had defined the category of étale sheaves: $\mathcal{S}\mathbf{hv}^{\mathbf{\acute{e}t}}(\mathcal{CAlg}^{\text{ét}}_{R_f})$, denoted $\mathsf{Spét}\ R_f$, and we had called it the **étale spectrum**, as the subcategory of the presheaves category $\mathcal{F}\mathbf{un}((\mathcal{CAlg}^{\text{ét}}_{R_f})^{\text{op}}, \mathcal{S})$ with the étale $\infty$-sheaf condition (see definition 2.2.20):

$$\mathsf{Spét}\, R_f := \mathcal{S}\mathbf{hv}^{\text{ét}}_{\mathcal{S}}\,(\mathcal{CAlg}^{\text{ét}}_{R_f}) \tag{2.7.6}$$

(see Equation 2.7.6) We had obtained (see proposition 2.2.54) that the étale spectrum is a localization of the category of $\mathcal{S}$-valued presheaves $\mathcal{F}\mathbf{un}((\mathcal{CAlg}^{\text{ét}}_{R_f})^{\text{op}}, \mathcal{S})$ and that the $\infty$-category $\mathsf{Spét}\ R_f$ is an $\infty$-topos.

The $\infty$-topos $\mathsf{Spét}\, R_f := \mathcal{S}\mathbf{hv}^{\text{ét}}_{\mathcal{S}}\,(\mathcal{CAlg}^{\text{ét}}_{R_f})$ plays for the $\infty$-topos $\mathfrak{S}$ of spectral Deligne-Mumford stacks on the $\infty$-topos of the TNN Grassmannian $\mathrm{Gr}(k, n)$ the model role of the affine schemes in the theory of classical algebraic geometry.

We can then define the canonical sheaf $\mathscr{O}_{R_f}$ in $\mathcal{S}\mathbf{hv}_{\mathsf{Sp}}(\mathfrak{X})$ as the functor of points: $\mathscr{O}_{R_f}(S) = S$ for all $S \in \mathsf{Spét}\, R_f$ (see definition 2.2.62).

Following the definition of a sheaf on an $\infty$-topos, if $\mathscr{C}$ is an $\infty$-category, and $\mathfrak{X}$ an $\infty$-topos, the category of sheaves on the $\infty$-topos $\mathfrak{X}$, denoted $\mathcal{S}\mathbf{hv}_{\mathscr{C}}(\mathfrak{X})$, is the category of $\mathscr{C}$-valued presheaves on $\mathfrak{X}$ which preserves limits on $\mathfrak{X}^{\text{op}}$ to $\mathscr{C}$, where we write:

$$\mathcal{S}\mathbf{hv}_{\mathscr{C}}(\mathfrak{X}) := \mathcal{F}\mathbf{un}^{\lim}(\mathfrak{X}^{\text{op}}, \mathscr{C}) \tag{2.7.7}$$

Following definition 2.2.62, with $R_f = \mathbb{C}[\overset{\circ}{\Pi}_f]$ (see Equation 2.7.6) is in $\mathcal{CAlg}^{\heartsuit}_{\mathbb{Z}}$, and with $\mathfrak{X} := \mathsf{Spét}\, R_f$ the $\infty$-topos of the étale spectrum of $R_f$, we define the canonical sheaf $\mathscr{O}_{R_f}$ in $\boldsymbol{\mathcal{S}}\mathbf{hv}_{\mathsf{Sp}}(\mathfrak{X})$ as the functor of points: $\mathscr{O}_{R_f}(S) = S$ for all $S \in \mathsf{Sp\acute{e}t}\, R_f$. The embedding $\boldsymbol{\mathcal{S}}\mathbf{hv}_{\mathsf{Sp}}(\mathsf{Sp\acute{e}t}\, R_f) := \boldsymbol{\mathcal{F}}\mathbf{un}^{\lim}((\mathsf{Sp\acute{e}t}\, R_f)^{\mathrm{op}}, \mathsf{Sp}) \hookrightarrow \boldsymbol{\mathcal{F}}\mathbf{un}((\mathsf{Sp\acute{e}t}\, R_f)^{\mathrm{op}}, \mathsf{Sp})$ gives us that the canonical sheaf $\mathscr{O}_{R_f}$ in $\boldsymbol{\mathcal{S}}\mathbf{hv}_{\mathsf{Sp}}(\mathfrak{X})$ is just the forgetful functor of the étale topology $\infty$-sheaf condition.

We have that the pair $(\mathfrak{X}, \mathscr{O}_{\mathfrak{X}}) := (\mathsf{Sp\acute{e}t}\, R_f, \mathscr{O}_{R_f})$ is a spectrally ringed $\infty$-topos and gives us that if the terminal object of an $\infty$-topos $\mathfrak{X}$ is denoted $*$, one has an effective cover of $\mathfrak{X}$ since there exists a collection of objects $\mathfrak{U}_f \in \mathfrak{X}$ and we have an effective epimorphism from the disjoint union of these objects of $\mathfrak{X}$ : $\coprod \mathfrak{U}_f \to *$.

Finally, following definition 2.2.67, the spectrally ringed topos $\mathfrak{X}_f$ such that $\mathfrak{X}_f := (\mathfrak{X}_{/\mathfrak{U}_f}, \mathscr{O}_{\mathfrak{X}}\big|_{\mathfrak{X}_{/\mathfrak{U}_f}})$ where $\mathfrak{X}_{/\mathfrak{U}_f}$ is the slice $\infty$-topos $\mathfrak{X}$ by $\mathfrak{U}_f$ endowed with the canonical sheaf $\mathscr{O}_{\mathfrak{X}_{/\mathfrak{U}_f}}$ which is the restriction $\mathscr{O}_{\mathfrak{X}}\big|_{\mathfrak{X}_{/\mathfrak{U}_f}}$ to $\mathfrak{X}_{/\mathfrak{U}_f}$ of the canonical sheaf $\mathscr{O}_{\mathfrak{X}}$.

The definition 2.2.69 gives us in conclusion that $\mathfrak{X} = (\mathfrak{X}, \mathscr{O}_{\mathfrak{X}})$ for which the terminal object $*$ of $\mathfrak{X}$ has an effective cover $\{\coprod \mathfrak{U}_f \to *\}$, and such that each $\mathfrak{U}_f \in \mathfrak{X}$ of the effective cover is such that $\mathfrak{X}_f := (\mathfrak{X}_{/\mathfrak{U}_f}, \mathscr{O}_{\mathfrak{X}}\big|_{\mathfrak{X}_{/\mathfrak{U}_f}})$ is equivalent to some $(\mathsf{Sp\acute{e}t}\, R_f, \mathscr{O}_{R_f})$ for an $R_f \in \boldsymbol{\mathcal{CA}}\mathbf{lg}_{\mathbb{Z}}^{\heartsuit}$, and is then a spectrally ringed $\infty$-topos which is a (nonconnective) spectral Deligne-Mumford stack.

As $R_f$ is concentrated in degree 0, we have that the negative homotopy groups of $\mathscr{O}_{\mathfrak{X}}$ are zero: $\pi_i(\mathscr{O}_{\mathfrak{X}}) = 0$ for $i < 0$ and then definition 2.2.71 gives us the expected (and same) $\infty$-category $\mathfrak{Gr}_{\text{ét}}^{\geqslant 0}$ of Deligne-Mumford stacks on the TNN Grassmannian $\mathrm{Gr}_{\geqslant 0}^{\mathbb{R}}(k,n)$, now described by its effective cover $\mathfrak{X}_f := (\mathfrak{X}_{/\mathfrak{U}_f}, \mathscr{O}_{\mathfrak{X}}\big|_{\mathfrak{X}_{/\mathfrak{U}_f}})$, each piece being equivalent to some $(\mathsf{Sp\acute{e}t}\, R_f, \mathscr{O}_{R_f})$ for an $R_f \in \mathcal{CA}\mathrm{lg}_{\mathbb{Z}}^{\heartsuit}$, that is $R_f \in \mathcal{R}\mathrm{ings}$.

Following again subsection 2.3.6, we identify thanks to **??**, a $\mathcal{M}$od-valued sheaf on $\mathfrak{Gr}_{\text{ét}}^{\geqslant 0}$ with a pair $(\mathscr{O}_{\mathfrak{Gr}_{\text{ét}}^{\geqslant 0}}, \mathcal{IC}_{\mathfrak{Gr}_{\text{ét}}^{\geqslant 0}})$, where $\mathcal{IC}_{\mathfrak{Gr}_{\text{ét}}^{\geqslant 0}}$ is the intersection complex sheaf of $\mathbb{E}_\infty$-rings of theorem 2.6.35.

We have then the same following definition of the $\infty$-category of spectral Deligne-Mumford stacks on the $\infty$-topos of the TNN Grassmannian $\mathfrak{Gr}_{\text{ét}}^{\geqslant 0}$ :

$$\mathfrak{G} := \left(\mathfrak{Gr}_{\text{ét}}^{\geqslant 0}, \mathscr{O}_{\mathfrak{Gr}_{\text{ét}}^{\geqslant 0}}, \mathcal{IC}_{\mathfrak{Gr}_{\text{ét}}^{\geqslant 0}}\right)$$

### 2.7.2 The positive Grassmannian spectral algebraic space.

Let us consider now the closure $\Pi_f$ of each open positroid varieties $\mathring{\Pi}_f$ in the stratification of the positive Grassmannian.

We can consider the induced topology on this object, that is that of a closed substack of an open substack.

DEFINITION 2.7.8 ([103] Def 2.5.2.2) : *Let $\mathfrak{X}$ be an $\infty$-topos. And suppose we are given a diagram:*

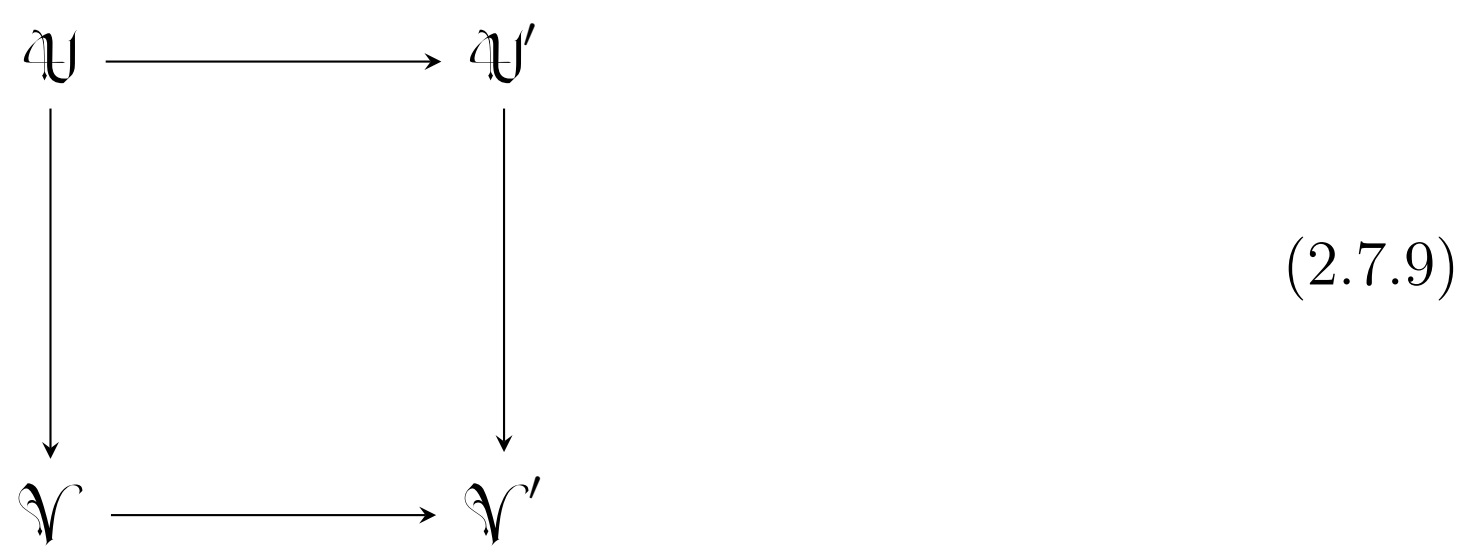

$$\begin{array}{ccc} \mathfrak{U} & \longrightarrow & \mathfrak{U}' \\ \downarrow & & \downarrow \\ \mathfrak{V} & \longrightarrow & \mathfrak{V}' \end{array} \tag{2.7.9}$$

*We call this diagram an* excision square *if one of the following equivalent statements happens:*

- *The diagram is both a pushout square and a pullback square and the map $f'$ is (-1)-truncated*
- *The diagram is a pushout square and the map $f$ is (-1)-truncated*
- *The diagram is a pullback square, the map $f'$ is (-1)-truncated, and if we let $i^* : \mathfrak{X}_{/\mathfrak{V}'} \longrightarrow \mathfrak{X}_{\mathfrak{V}'}/\mathfrak{V}$ denote the corresponding closed immersion, then $i^*\mathfrak{U}'$ is a final object of $\mathfrak{X}_{\mathfrak{V}'}/\mathfrak{V}$*

DEFINITION 2.7.10 ([103] Def 2.5.3.1) : *Let $(\mathfrak{X}, \mathscr{O}_{\mathfrak{X}})$ be a nonconnective spectral Deligne-Mumford stack. We name a* scallop decomposition *of $\mathfrak{X}$ a sequence of (-1)-truncated morphisms in $\mathfrak{X}$:*

$$\mathfrak{U}_0 \longrightarrow \mathfrak{U}_1 \longrightarrow \dots \mathfrak{U}_n$$

*satisfying the following conditions:*

- *the object $\mathfrak{U}_0 \in \mathfrak{X}$ is initial and the object $\mathfrak{U}_n$ is final*
- *for $1 \leqslant i \leqslant n$, there exists an excision square:*

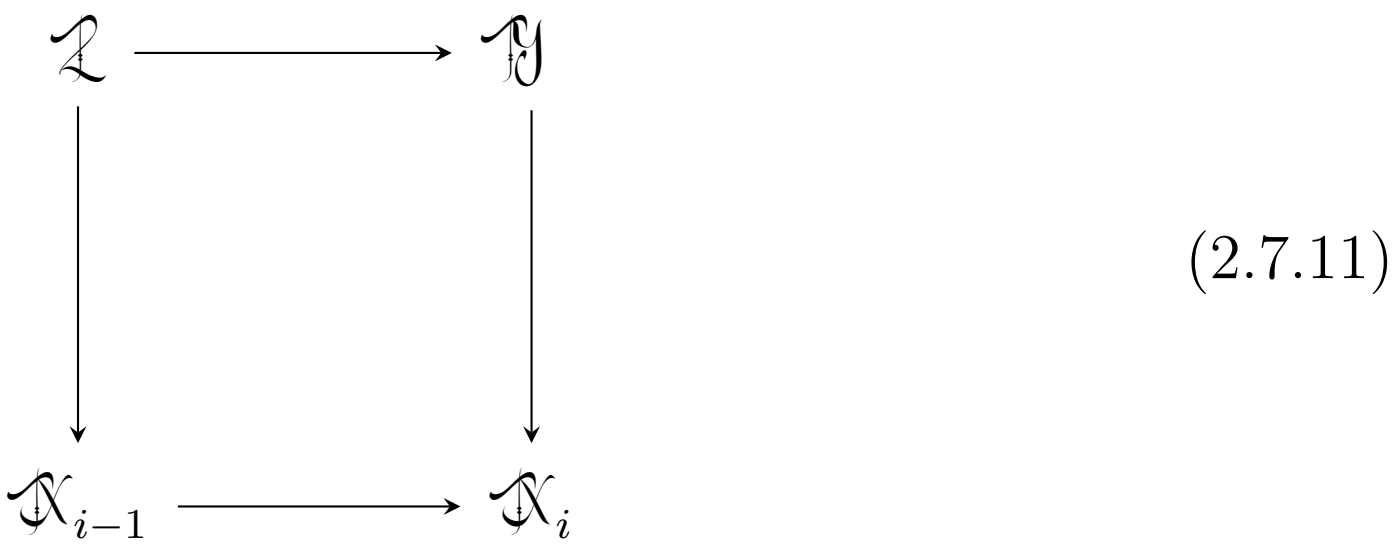

(2.7.11)

*where $\mathfrak{Y}$ is affine and $\mathfrak{Z}$ is quasi-compact*

REMARK 2.7.12 : Each $\mathfrak{U}_i$ determines an open substack $\mathsf{U}_i := (\mathfrak{X}_{\mathfrak{U}_i}, \mathscr{O}_{\mathfrak{X}_{|\mathfrak{U}_i}})$ of $\mathfrak{X}$. We then also refer to the sequence of open immersions:

$$\varnothing \simeq \mathsf{U}_0 \longrightarrow \mathsf{U}_1 \longrightarrow \ldots \mathsf{U}_n \simeq \mathfrak{X}$$

as a scallop decomposition.

EXAMPLE 2.7.13 : If $\mathfrak{X}$ is a quasi-affine nonconnective spectral Deligne-Mumford stack, then $\mathfrak{X}$ admits a scallop decomposition.

THEOREM 2.7.14 ([103] Thm 3.4.2.1) : *Let $(\mathfrak{X}, \mathscr{O}_{\mathfrak{X}})$ be a spectral Deligne-Mumford stack. Then $\mathfrak{X}$ admits a scallop decomposition if and only if it is a quasi-compact quasi-separated spectral algebraic space.*

Consider each open positroid varieties $\mathring{\Pi}_f$, that is an open substack.

We can consider the following process:

1. Take an $f_0 \in \mathscr{B}(k, n)$ and the corresponding open substack $\mathring{\Pi}_{f_0}$
2. Take the closure of $\mathring{\Pi}_{f_0}$
3. Glue along closed immersions of their frontier with all closures of open substacks $\mathring{\Pi}_f$ where $\mathring{\Pi}_f$ are next to $\mathring{\Pi}_f$ (this is possible since the stratification tells us that the closure of a cell does not belong to the closure of strictly more than 1 higher dimensional cell)

The result is open, and an open gluing of cells.

We can the reiterate the process and obtain many different scallop decompositions (it depends on the first choice), where the cells are interconnected, until covering the all TNN Grassmannian. All the combinatorics which is possible to compose/decompose the TNN Grassmannian could be studied. Nevertheless, we have the following:

THEOREM 2.7.15 : *Consider now one of the possible scallop decompositions of* $\mathfrak{G}$ *seen previously (trivial or not). We know by the preceeding theorem that having a scallop decomposition is equivalent to being a spectral algebraic space. We obtain that:*

$$\mathfrak{G} := \left(\mathfrak{Gr}^{\geqslant 0}_{\text{ét}}, \mathscr{O}_{\mathfrak{Gr}^{\geqslant 0}_{\text{ét}}}, \mathcal{IC}_{\mathfrak{Gr}^{\geqslant 0}_{\text{ét}}}\right)$$

*is a spectral algebraic space.*

## 2.8.The infinitesimal site and the De Rham stack.

Being in the field of algebraic geometry, we do not have the notion of "*infinitesimal neighborhood*". However, Grothendieck defined the "*infinitesimal site*" which comes to solution this problem.

DEFINITION 2.8.1 : *Let* $\mathscr{X}$ *be a scheme over* $\mathscr{S}$. *We call* infinitesimal site *of* $\mathscr{X}$ *over* $\mathscr{S}$, *denoted* $\mathrm{Inf}\,(\mathscr{X}/\mathscr{S})$, *the category such that:*

- *the objects are the closed immersions* $\mathscr{U} \subseteq \mathscr{T}$, *where* $\mathscr{U}$ *is an open of* $\mathscr{X}$, *and where the ideal defining this immersion is nilpotent*

- *the morphisms between* $\mathscr{U} \subseteq \mathscr{T}$ *and* $\mathscr{U}' \subseteq \mathscr{T}'$ *are those from* $\mathscr{T}$ *to* $\mathscr{T}'$ *inducing* *an open immersion from $\mathscr{U}$ into $\mathscr{U}'$*

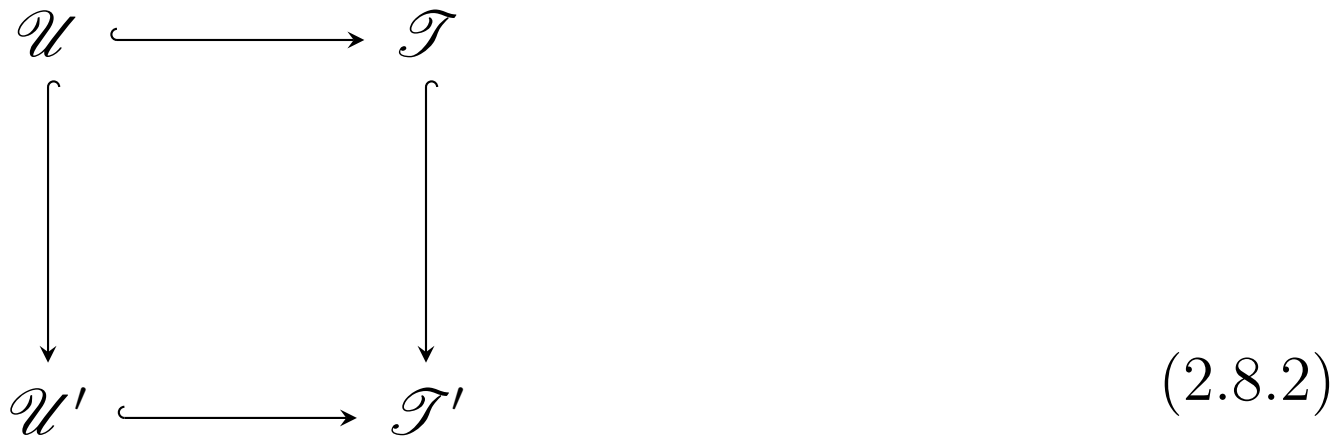

$$\begin{array}{ccc} \mathscr{U} & \hookrightarrow & \mathscr{T} \\ \downarrow & & \downarrow \\ \mathscr{U}' & \hookrightarrow & \mathscr{T}' \end{array} \tag{2.8.2}$$

*endowed with the structure of site by stating that the coverings of an object $\mathscr{U} \subseteq \mathscr{T}$ are the families $(\mathscr{U}_i \subseteq \mathscr{T}_i)_{i\in I}$ where $(\mathscr{T}_i)_{i\in I}$ is a Zariski open covering of $\mathscr{T}$. We have a canonical sheaf of rings on the infinitesimal site, sending $\mathscr{U} \subseteq \mathscr{T}$ to $\mathscr{O}(\mathscr{T})$. This sheaf is denoted $\mathscr{O}_{(\mathscr{X}/\mathscr{S})_{\text{inf}}}$ or simply $\mathscr{O}$ if there is no ambiguity.*

Sending to the construction of the cotangent complex in subsection 2.1.2 for the definition of the algebraic De Rham cohomology from Kähler differentials in the smooth case, we have the following:

THEOREM 2.8.3 (Grothendieck '66) : *Let $\mathscr{X}$ be a smooth scheme over $\mathscr{S}$ whose characteristic is* 0*. Then we have a canonic isomorphism:*

$$H^*((\mathscr{X}/\mathscr{S})_{\text{inf}}, \mathscr{O}_{(\mathscr{X}/\mathscr{S})_{\text{inf}}}) \simeq H^*_{\text{DRh}}(\mathscr{X}/\mathscr{S})$$

REMARK 2.8.4 : *The goal of all the subsequent efforts is to define an algebraic de Rham cohomology that allows us to remove the smoothness assumption and thus also handle cases where the scheme, the stack, etc., are singular. The singularities will, of course, remain within the limits set by the constructions to come, and this is what will effectively allow us to extend the current framework to broader horizons :* A brave new positive Schubert geometry.

We will not reproduce the demonstration of the theorem by Grothendieck to point the step where the smothness hypothesis can be removed. Nonetheless, lots of concepts of the demonstration will be extremely useful for our following work, in particular about formal moduli problems and deformation theory on one hand and about cohesiveness and synthetic differential geometry on the other hand. So let us make a few remarks and constructions we will use later.

REMARK 2.8.5 : There is another definition of the sheaf of Kähler differentials, which is the one that is of interest for the following construction: If $\mathscr{I}$ is the ideal of the immersion of the diagonal in the fiber product $\mathscr{X} \times_{\mathscr{S}} \mathscr{X}$, then the sheaf of Kähler differentials is the pullback by the diagonal of the quotient sheaf $\mathscr{I}/\mathscr{I}^2$. This gives also a third definition which is that this sheaf of Kähler differentials is the pullback by the diagonal of the definition ideal of the immersion of $\mathscr{X}$ in $\Delta_2^2$, which is the 2-order infinitesimal neibourghood in $\mathscr{X} \times_{\mathscr{S}} \mathscr{X}$. This will be the point we use to define De Rham cohomology using the infinitesimal neighbourhoods at all orders.

REMARK 2.8.6 : Grothendieck's demonstration did not use the infinitesimal site, despite its highly interest, its useful character and the inventive notion it reveals. Indeed, he used the named "stratifying site", denoted $\mathrm{Strat}(\mathscr{X}/\mathscr{S})$, defined such that we only take in the infinitesimal site the full subcategory made of the objects $\mathscr{U} \subseteq \mathscr{T}$ for which there locally exists a retraction $\mathscr{T} \longrightarrow \mathscr{X}$ with the induced topology. In the smooth case over $\mathscr{S}$, the two sites are two equivalent categories, since we have the following "formal smooth criterium":

- If $\mathscr{X}$ is a smooth $\mathscr{S}$-scheme, then for all commutative diagram

[formal smooth criterium]
$$\begin{array}{ccc} \mathscr{R} = \mathsf{Spec}\ A/I & \longrightarrow & \mathscr{X} \\ \downarrow & \nearrow & \downarrow \\ \mathscr{T} = \mathsf{Spec}\ A & \longrightarrow & \mathscr{S} \end{array} \tag{2.8.7}$$

where $A$ is a ring, $I$ is a nilpotent ideal of $A$, there exists a morphism $\mathscr{T} \longrightarrow \mathscr{X}$ lifting the base map $\mathscr{T} \longrightarrow \mathscr{S}$.

If we take $\mathscr{R} = \mathscr{U}$ in the infinitesimal site, with $\mathscr{U} \subseteq \mathscr{T}$, then the locally existence of a retraction is automatic, so that the two sites are equivalent. However, in the general case, the infinitesimal site will give the "*good cohomology*", without the restriction of smoothness, that is what we expect later.

The idea of connexion comes from differential geometry and can be adapted to the algebraico-geometric world. Given a path between two points in the neighbourhood of each other reciprocally, the idea is to "*transport*" a tangent vector of this common neibourhood along this path. We usually use the Kähler differentials in order to define a connexion, but, as they are forbidden in our construction, we have to use another way.

DEFINITION 2.8.8 : *Let $\mathscr{X}$ be an $\mathscr{S}$-scheme. For an $n \in \mathbb{N}^*$, we denote $\Delta_2^n$ the $n$-th infinitesimal neighbourhood in the diagonal $\mathscr{X} \times_{\mathscr{S}} \mathscr{X}$, $\Delta_3^n$ the $n$-th infinitesimal neighbourhood in the triple product $\mathscr{X} \times_{\mathscr{S}} \mathscr{X} \times_{\mathscr{S}} \mathscr{X}$, etc. and we consider the following projections:*

$$\mathscr{X} = \Delta_1^n \quad \underset{\longleftarrow p_2^n \text{———}}{\overset{\longleftarrow p_1^n \text{———}}{}} \quad \Delta_2^n \quad \begin{matrix} \longleftarrow p_{12}^n \text{———} \\ \longleftarrow p_{23}^n \text{———} \\ \longleftarrow p_{12}^n \text{———} \end{matrix} \quad \Delta_3^n \tag{2.8.9}$$

*Let $M$ be an $\mathscr{O}_{\mathscr{X}}$−module. We call $n$-connexion over $M$ the data of an isomorphism $\varphi : (p_1^n)^* M \simeq (p_2^n)^* M$ such that:*

$$(p_{23}^n)^*(\varphi) \circ (p_{12}^n)^*(\varphi) = (p_{13}^n)^*(\varphi) \tag{2.8.10}$$

*A stratification over a module $M$ is the data of a family of $n$-connexions for all $n \in \mathbb{N}^*$, which are compatible between themselves through the inclusion of the $n$-th infinitesimal neighbourhood into the $(n+1)$-th infinitesimal neighbourhood.*

REMARK 2.8.11 : The data of a $n$-connexion identifies the fibers of the module $M$ in two $n$-th order infinitesimally closed points $x$ and $y$. The extra isomorphisms garantee the compositions of the proximities whatever we take as the infinitesimal paths between $x$ and $y$.

REMARK 2.8.12 : The isomorphism $\varphi : (p_1^n)^* M \simeq (p_2^n)^* M$ is a priori an isomorphism of sheaves on $\Delta_2^n$, but it is in fact an equality of sheaves on the same underlying space of $\mathscr{X}$ since the underlying topological space of $\Delta_2^n$ is exactly the underlying space of $\mathscr{X}$, and that is for all $n \in \mathbb{N}^*$. The sheaves live, in a certain manner, in the same space. Seen as endofunctors of $\mathscr{O}_{\mathscr{X}}-$modules, $(p_1^n)^*$ and $(p_2^n)^*$ corrrespond to the tensor product by $(\mathscr{O}_{\mathscr{X}} \otimes_{\mathscr{O}_{\mathscr{S}}} \mathscr{O}_{\mathscr{X}})/\mathscr{I}^n$, where $\mathscr{I}$ is the definition ideal of the imersion of the diagonal in $\mathscr{X} \times_{\mathscr{S}} \mathscr{X}$. The structure of $\mathscr{O}_{\mathscr{X}}-$module lives in the multiplication on the right or on the left of the tensor product.

REMARK 2.8.13 : The data of a stratified module $M$ over the $\mathscr{S}$-scheme $\mathscr{X}$ allows the construction of a sheaf $\widetilde{M}$ on the stratifying site Strat$(\mathscr{X}/\mathscr{S})$. Let $\mathscr{U} \subseteq \mathscr{T}$ be an open of $\mathscr{X}$ and an infinitesimal thickening admitting a retraction $r : \mathscr{T} \longrightarrow \mathscr{X}$. We want $\widetilde{M}(\mathscr{U} \subseteq \mathscr{T}) = r^* M$. An important fact justified by all this machinery is that the all construction is independent of the retraction. The existence of a retraction will exactly allow to identify $r^* M$ and a possibly other $r'^* M$ for an other rtraction $r'$. We omit the quite heavy demonstration which is without interest for us here.

REMARK 2.8.14 : The remainder of the proof of the theorem now consists in constructing an adjunction between the following categories:

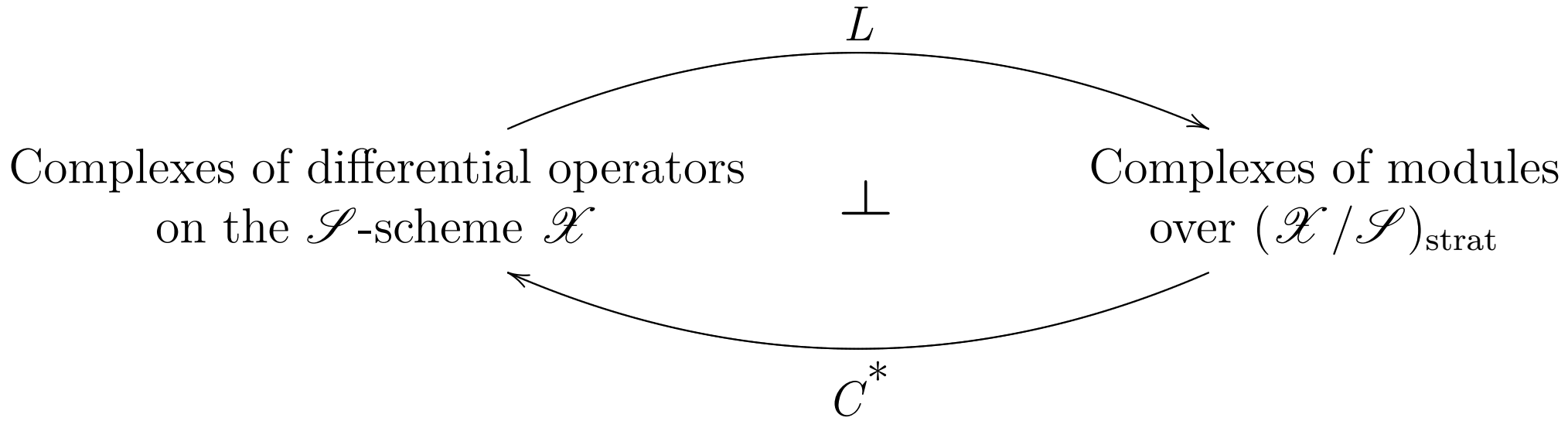

where the composition $C^* \circ L$ is the identity as endofunctor in the derived category of complexes of differential operators on the $\mathscr{S}$-scheme $\mathscr{X}$.

REMARK 2.8.15 : The functor $L$ associate a complexe of stratified (pro-)modules to a complexe of differential operators: if $M$ is an $\mathscr{O}_{\mathscr{X}}$-module then $(p_1^n)_* \mathscr{O}_{\Delta_2^n} \otimes_{\mathscr{O}_{\mathscr{X}}} M$ admits an $n$-connexion. We define $L(M)$ as the topological projective limit:

$$L(M) := \varprojlim((p_1^n)_* \mathscr{O}_{\Delta_2^n} \otimes_{\mathscr{O}_{\mathscr{X}}} M)_{n\in\mathbb{N}^*} \tag{2.8.16}$$

where the family $((p_1^n)_* \mathscr{O}_{\Delta_2^n} \otimes_{\mathscr{O}_{\mathscr{X}}} M)_{n\in\mathbb{N}^*}$ is endowed with a stratification. A sheaf $\mathscr{F} = L(M)$ on the stratifying site is named a crystal, which consists in the following property:

- For all morphisms $(\mathscr{U} \hookrightarrow \mathscr{T}) \longrightarrow (\mathscr{U}' \hookrightarrow \mathscr{T}')$ in the stratifying site, we have: $\mathscr{F}(\mathscr{U} \hookrightarrow \mathscr{T}) = f^*\mathscr{F}(\mathscr{U}' \hookrightarrow \mathscr{T}')$

Indeed, the sheaf on the stratifying site associated to a stratifyied module is defined as the pullback by retractions. Since if $r : \mathscr{T}' \longrightarrow \mathscr{X}$ is a retraction, then $r' \circ f$ is a retraction of $\mathscr{T}$, which is the property of being a crystal.

REMARK 2.8.17 : We define the functor $C^*$ if $\mathscr{F}$ is a module on the stratifying site, by the complex whose each $m$-th term is defined by the following:

$$C^m(\mathscr{F}) := \varprojlim \mathscr{F}(\mathscr{X} \hookrightarrow \Delta_m^n) \tag{2.8.18}$$

We have the following lemma by Grothendieck:

LEMMA 2.8.19 : *Let $\mathscr{F}$ be a sheaf of modules on the stratifying site $(\mathscr{X}/\mathscr{S})_{\mathrm{strat}}$. Let $\mathscr{F}$ be a crystal. Then we have the canonical isomorphism:*

$$H^*((\mathscr{X}/\mathscr{S})_{\mathrm{strat}}, \mathscr{F}) \simeq H^*(\mathscr{X}_{\mathrm{Zar}}, C^*(\mathscr{F})) \tag{2.8.20}$$

Now, the algebraic De Rham complex $(\Omega^*_{\mathscr{X}/\mathscr{S}})$ is a complex of differential operators. Applying the functor $L$ on it, we have, thanks to preceeding lemma:

$$H^*((\mathscr{X}/\mathscr{S})_{\mathrm{strat}}, L(\Omega^*_{\mathscr{X}/\mathscr{S}})) \simeq H^*(\mathscr{X}_{\mathrm{Zar}}, C^*(L(\Omega^*_{\mathscr{X}/\mathscr{S}}))) \tag{2.8.21}$$

On the other hand, we announce in introduction of this construction that the composition $C^* \circ L$ is the identity in the derived category of complexes of differential operators, that is the complex $C^*(L(\Omega^*_{\mathscr{X}/\mathscr{S}}))$ is quasi-isomorphic to $\Omega^*_{\mathscr{X}/\mathscr{S}}$. This gives the isomorphism:

$$H^*((\mathscr{X}/\mathscr{S})_{\text{strat}}, L(\Omega^*_{\mathscr{X}/\mathscr{S}})) \simeq H^*(\mathscr{X}_{\text{Zar}}, \Omega^*_{\mathscr{X}/\mathscr{S}}) \tag{2.8.22}$$

Since we always have a canonical morphism : $\mathscr{O}_{\text{strat}} \longrightarrow L(\Omega^*_{\mathscr{X}/\mathscr{S}})$, the demonstration of the theorem is then reduced to the following lemma:

LEMMA 2.8.23 (Grothendieck '66, Stratifying Poincaré lemma) : *If $\mathscr{X}$ is an $\mathscr{S}$-scheme, with characteristic 0, then the morphism $\mathscr{O}_{\text{strat}} \longrightarrow L(\Omega^*_{\mathscr{X}/\mathscr{S}})$ is a quasi-isomorphism.*

The characteristic zero plays a role in the fact that we need that the only polynomials whose derivation is null are constant polynomials. Except this condition which does not concern us here and speaks first to number theoretists for *crystalline cohomology*, the door is wide open for the De Rham cohomology in the singular cases, which we will also call crystalline cohomology since the corner stone of the construction is Grothendieck's notion of crystals.

We can then conclude this section by the construction we expected : the De Rham stack

DEFINITION 2.8.24 : *Let $\mathscr{X}$ be an $\mathscr{S}$-scheme, not necessarily smooth. We call* De Rham stack *over the $\mathscr{S}$-scheme $\mathscr{X}$, and we denote it $\mathscr{X}^{\text{dRh}}_{/\mathscr{S}}$, the stack induced by the functor:*

$$\begin{cases} \mathscr{S}\textbf{-Sch} & \longrightarrow & \boldsymbol{\mathcal{S}}\textbf{et} \\ \mathscr{Y} & \longmapsto & \mathscr{X}(\mathscr{Y}_{\text{red}}) \end{cases}$$

*where we denote $\mathscr{Y}_{\text{red}}$ the reduced $\mathscr{S}$-scheme of $\mathscr{Y}$.*

The stack $\mathscr{X}^{\text{dRh}}_{/\mathscr{S}}$ is then the fibered category of couples $(\mathscr{Y}, t)$ where $\mathscr{Y}$ is an $\mathscr{S}$-scheme and $t$ an $\mathscr{Y}_{\text{red}}$-point of $\mathscr{X}$ above $\mathscr{S}$.

We can then prove the following:

PROPOSITION 2.8.25 : *Let $\mathscr{X}$ be an $\mathscr{S}$-scheme, not necessarily smooth. The category of quasi-coherent sheaves on the De Rham stack $\mathscr{X}^{\text{dRh}}_{/\mathscr{S}}$ is equivalent to the full subcategory of quasi-coherent sheaves on the infinitesimal site Inf($\mathscr{X}/\mathscr{S}$) consisting in the sheaves $\mathscr{F}$ statisfying:*

- *For all morphisms* $f : (\mathscr{U} \hookrightarrow \mathscr{T}) \xrightarrow{f} (\mathscr{U}' \hookrightarrow \mathscr{T}')$ *in the infinitesimal site, the induced morphism* $f^*\mathscr{F}(\mathscr{U}' \hookrightarrow \mathscr{T}') \simeq \mathscr{F}(\mathscr{U} \hookrightarrow \mathscr{T})$ *is an isomorphism.*

This proposition says exactly the expected following statements:

[De Rham stack statement] *The category of quasi-coherent sheaves on the De Rham stack is equivalent to the category of quasi-coherent sheaves on the infinitesimal site which are isomorphic to a crystal.*

[De Rham stack crystals monadic generators] *The category of crystals over the infinitesimal site being monadic, it suffices to look at crystals, which are generators of the category of quasi-coherent sheaves on the initial scheme.*

[Smoothness inconsistency] *No criterium of smoothness is needed when we work with the infinitesimal site.*

## 2.9.Homotopy type theory, adjunctions and monads.

### 2.9.1 Category Theory: Yoneda, Adjunctions and Monads.

THEOREM 2.9.1 (Yoneda lemma, [121] Thm 2.2.4) : *Let* $\mathcal{C}$ *be a locally small category, and let* $x$ *be an object of the category* $\mathcal{C}$*. Let* $\mathscr{F}$ *be any functor from* $\mathcal{C}$ *to the category* $\mathcal{S}et$ *of sets. We then have a bijection between the collection of natural transformations from the (representable) functor* $h_x := \mathrm{Hom}_{\mathcal{C}}(x, -)$ *to the functor* $\mathscr{F}$ *and the set* $\mathscr{F}x$*:*

$$\begin{cases} \mathrm{Nat}(h_x, \mathscr{F}) := \mathrm{Nat}(\mathrm{Hom}_{\mathcal{C}}(x,-), \mathscr{F}) & \xrightarrow{\sim} & \mathscr{F}x \\ \alpha : \mathrm{Nat}(h_x, \mathscr{F}) := \mathrm{Hom}_{\mathcal{C}}(x,-) \Rightarrow \mathscr{F} & \longmapsto & \alpha_x(\mathrm{Id}_x) \in \mathcal{S}\mathrm{et} \end{cases} \tag{2.9.2}$$

*The bijection is natural in both* $x$ *of* $\mathcal{C}$ *and* $\mathscr{F}$ *from* $\mathcal{C}$ *to* $\mathcal{S}et$*. In particular, we obtain that this collection of natural transformations is a set, even though the category* $\mathcal{C}$ *is locally small but not necessarily small.*

REMARK 2.9.3 : Any functor $\mathscr{F} : \mathcal{C} \longrightarrow \mathcal{S}\mathrm{et}$ from a small category $\mathcal{C}$ to the category $\mathcal{S}$et of sets is equivalent to an inductive limit of representable functors: $\mathscr{F} \simeq \varprojlim_{x \in \mathcal{C}} h_x$.

This remark will play a role in *accessibility* for topoï and $\infty$-topoï (the $\infty$-Yoneda).

DEFINITION 2.9.4 (Adjunctions I, [121], Def 4.1.1) : *Let $\mathcal{C}$ and $\mathcal{D}$ be a pair of categories linked by a pair of functors $\mathscr{F}$ and $\mathscr{G}$. An* adjunction *between the two functors means that they enjoy the following relationship: For each $c \in \mathcal{C}$ and $d \in \mathcal{D}$, we have an isomorphism at the level of* Hom *which is natural in both variables $c$ and $d$ :*

$$\mathrm{Hom}_{\mathcal{D}}(\mathscr{F}c, d) \simeq \mathrm{Hom}_{\mathcal{C}}(c, \mathscr{G}d) \tag{2.9.5}$$

*We say here that $\mathscr{F}$ is* left adjoint *to $\mathscr{G}$ and $\mathscr{G}$ is* right adjoint *to $\mathscr{F}$ and we denote:*

$$\mathcal{C} \underset{\mathscr{F}}{\overset{\mathscr{G}}{\underset{\longrightarrow}{\overset{\longleftarrow}{\top}}}} \mathcal{D} \tag{2.9.6}$$

*The naturality statement in both variables and the isomorphism of Equation 2.9.17 mean exactly that the vertical double arrow in the following diagram is a natural isomorphism between the two bifunctors:*

$$\mathcal{C}^{\mathrm{op}} \times \mathcal{D} \; \underset{\mathrm{Hom}_{\mathcal{C}}(-, \mathscr{G}-)}{\overset{\mathrm{Hom}_{\mathcal{D}}(\mathscr{F}-, -)}{\Downarrow \simeq}} \; \mathcal{S}\mathrm{et} \tag{2.9.7}$$

REMARK 2.9.8 : It is right now noteworthy that a functor can have at the same time a left adjoint and a right adjoint. The case will be treated in detail in the case of discrete and indiscrete topologies for the following adjunctions to the forgetful functor between $\mathcal{T}$op, the category of topological spaces, and $\mathcal{S}$et, the category of sets. A discussion about the involved topologies will happen deeper at this time:

$$\mathcal{T}\mathrm{op} \; \underset{\text{top. discrete}(-)}{\overset{\text{top. indiscrete}(-)}{\overset{\top}{\underset{\top}{\xrightarrow{\text{forget}(-)}}}}} \; \mathcal{S}\mathrm{et} \tag{2.9.9}$$

The left adjoint "*top. discrete*$(-)$" constructs a topological space from a set $S$ in such a way that continuous maps from this space to an another space $T$ correspond naturally and bijectively to functions from $S$ to "*forget*$(T)$" :

$$\mathrm{Hom}_{\mathcal{T}\mathrm{op}}(\textit{top. discrete}(S), T) \simeq \mathrm{Hom}_{\mathcal{S}\mathrm{et}}(S, \textit{forget}(T)) \tag{2.9.10}$$

Similarly, the right adjoint "*top. indiscrete*$(-)$" constructs a topological space from a set $S$ in such a way that continuous maps from $T$ to this space correspond naturally and bijectively to functions from "*forget*$(T)$" to $S$ :

$$\mathrm{Hom}_{\mathcal{S}\mathrm{et}}(\textit{forget}(T), S) \simeq \mathrm{Hom}_{\mathcal{T}\mathrm{op}}(T, \textit{top. indiscrete}(S)) \tag{2.9.11}$$

This is sum up in the biadunction diagram Equation 2.9.9.

REMARK 2.9.12 : Another biadjunction phenomenon arises in the following situation, which appears later in relation to the concept of "*truth*" in an ($\infty$-)topos, and more generally in the sections concerning "*logics*" and "*elementary topoi*".
Indeed, we say that a function is **propositional** if it is a function:

$$P : X \longrightarrow \Omega = \{\textit{false}, \textit{true}\} \tag{2.9.13}$$

which we interpret as declaring, for each $x \in X$, wether $P(x)$ is *true* or *false*. The set $\Omega^X$ is then the set of propositional functions on $X$. The set $\Omega$ is given the partial order "*false* $\leqslant$ *true*", from which $\Omega^X$ inherits a pointwise-defined order: $P \leqslant Q$ if and only if $P(x) \leqslant Q(x)$ for all $x \in X$, which is the case if and only if $P \Rightarrow Q$. We then define the logical operations of universal and existential quantification by the functors $\forall_X, \exists_X : \Omega^X \rightrightarrows \Omega$ where $\forall_X P = \textit{true}$ if and only if $P(x) = \textit{true}$ for all $x \in X$, and $\exists_X P = \textit{true}$ if and only if there exists $x \in X$ with $P(x) = \textit{true}$. Considering the functor $\Delta_X : \Omega \longrightarrow \Omega^X$, defined such that $\Delta_X(\textit{false}) = \{x \mapsto \textit{false}\ \} \in \Omega^X$ and $\Delta_X(\textit{true}) = \{x \mapsto \textit{true}\ \} \in \Omega^X$ (see [17]), one have also the following triple of adjoints:

$$\Omega \underset{\exists_X(-)}{\overset{\forall_X(-)}{\leftleftarrows}} \xrightarrow{\Delta(-)} \Omega^X \quad (\top,\ \top) \tag{2.9.14}$$

Let us now define the concept of adjunction in a second way. In the adjunction:

$$\mathcal{C} \underset{\mathscr{F}}{\overset{\mathscr{G}}{\underset{\longrightarrow}{\overset{\longleftarrow}{\top}}}} \mathcal{D} \tag{2.9.15}$$

and the isomorphisms, natural in both variables $c \in \mathcal{C}$ and $d \in \mathcal{D}$ :

$$\mathrm{Hom}_{\mathcal{D}}(\mathscr{F}c, d) \simeq \mathrm{Hom}_{\mathcal{C}}(c, \mathscr{G}d) \tag{2.9.16}$$

let us fix $c \in \mathcal{C}$. This means exactly that the object $\mathscr{F}c \in \mathcal{D}$ represents the functor $\mathrm{Hom}_{\mathcal{C}}(c, \mathscr{G}-) : \mathcal{D} \longrightarrow \mathcal{S}\mathrm{et}$. By the Yoneda lemma (see Equation 2.9.2), the natural isomorphism:

$$\mathrm{Hom}_{\mathcal{D}}(\mathscr{F}c, -) \simeq \mathrm{Hom}_{\mathcal{C}}(c, \mathscr{G}-) \tag{2.9.17}$$

is determined by an element of $\mathrm{Hom}_{\mathcal{C}}(c, \mathscr{G}\mathscr{F}c)$, denoted $\eta_c$. Thanks to the naturality in $c \in \mathcal{C}$, let us define the natural transformation $\eta : \mathbf{1}_{\mathcal{C}} \Rightarrow \mathscr{G}\mathscr{F}$ in this way, and let us call it the unit of the adjunction. Dually, let us define the natural transformation $\epsilon : \mathscr{F}\mathscr{G} \Rightarrow \mathbf{1}_{\mathcal{D}}$, and let us call it the counit of the adjunction.

Remark 2.9.18 : If $\mathscr{F}$ and $\mathscr{G}$ were reciprocal isomorphisms between the catogories $\mathcal{C}$ and $\mathcal{D}$, then the unit $\eta$ and the counit $\epsilon$ would be identities $\eta : \mathbf{1}_{\mathcal{C}} = \mathbf{1}_{\mathcal{C}}$ and $\epsilon : \mathbf{1}_{\mathcal{D}} = \mathbf{1}_{\mathcal{D}}$, and if $\mathscr{F}$ and $\mathscr{G}$ were reciprocal equivalences, then the unit $\eta$ and the counit $\epsilon$ would be natural isomorphisms $\eta : \mathbf{1}_{\mathcal{C}} \simeq \mathbf{1}_{\mathcal{C}}$ and $\epsilon : \mathbf{1}_{\mathcal{D}} \simeq \mathbf{1}_{\mathcal{D}}$. But in general, the adjunction between $\mathscr{F}$ and $\mathscr{G}$ defines a unit $\eta : \mathbf{1}_{\mathcal{C}} \Rightarrow \mathscr{G}\mathscr{F}$ and a counit $\epsilon : \mathscr{F}\mathscr{G} \Rightarrow \mathbf{1}_{\mathcal{D}}$ which are any natural transformations.

However, from any natural transformations $\eta : \mathbf{1}_{\mathcal{C}} \Rightarrow \mathscr{G}\mathscr{F}$ and $\epsilon : \mathscr{F}\mathscr{G} \Rightarrow \mathbf{1}_{\mathcal{D}}$, we can define an adjunction $\mathscr{F} \dashv \mathscr{G}$, and there is even a bijective correspondence between the pair $(\mathscr{F}, \mathscr{G})$ of an adjunction $\mathscr{F} \dashv \mathscr{G}$ and the pair $(\eta, \epsilon)$ of the unit and the counit of this adjunction, which leads to the following :

Definition 2.9.19 (Adjunctions II, [121], Def 4.2.5) : *An* adjunction *consists of an opposing pair of functors* $\mathscr{F} : \mathcal{C} \rightleftarrows \mathcal{D} : \mathscr{G}$, *together with natural transformations* $\eta$ :

$\mathbf{1}_{\mathcal{C}} \Rightarrow \mathscr{G}\mathscr{F}$ *and* $\epsilon : \mathscr{F}\mathscr{G} \Rightarrow \mathbf{1}_{\mathcal{D}}$ *that satisfies the following triangle identities:*

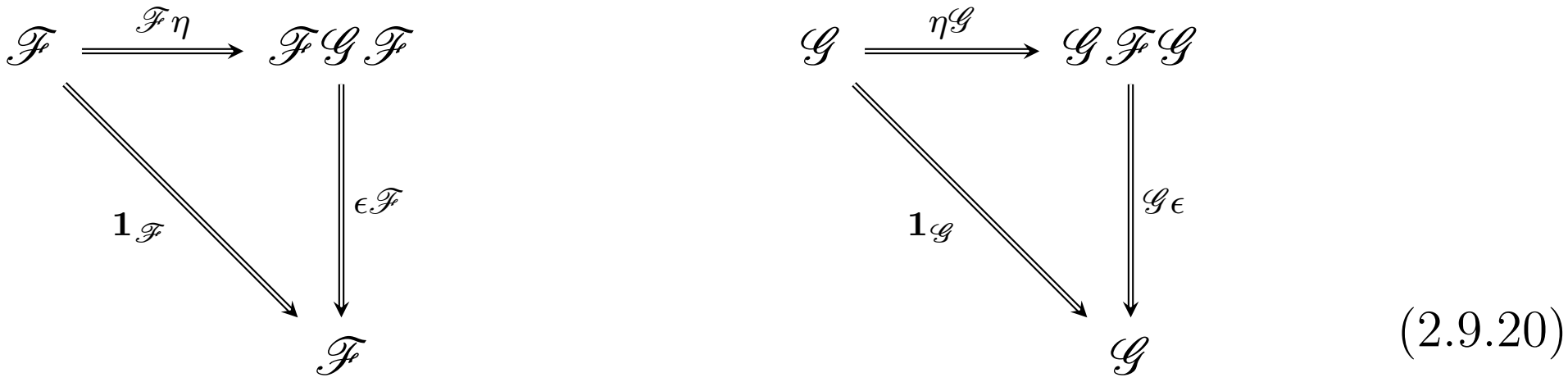

(2.9.20)

#### 2.9.1.1 Adjunctions, limits and colimits.

DEFINITION 2.9.21 ([121], Def 3.1.1) : *For any object $c$ in the category $\mathcal{C}$ and any category (of indexes) $\mathsf{J}$, the constant functor $c : \mathsf{J} \longrightarrow \mathcal{C}$ sends every object of $\mathsf{J}$ to $c$ and every morphism in $\mathsf{J}$ to the identity morphism $\mathbf{1}_c$. The constant functors define an embedding $\Delta : \mathcal{C} \longrightarrow \mathcal{C}^{\mathsf{J}}$ that sends an object $c$ to the constant functor at $c$ and a morphism $f : c \longrightarrow c'$ to the constant natural transformation, in which each component is defined to be the morphism $f$.*

PROPOSITION 2.9.22 ([121], Prop 4.5.1) : *A category $\mathcal{C}$ admits all limits of diagrams indexed by a small category $\mathsf{J}$ if and only if the constant diagram functor $\Delta : \mathcal{C} \longrightarrow \mathcal{C}^{\mathsf{J}}$ admits a right adjoint, and admits all colimits of $\mathsf{J}$-indexed diagrams if and only if $\Delta$ admits a left adjoint:*

$$\mathcal{C} \underset{\lim(-)}{\overset{\operatorname{colim}(-)}{\rightleftarrows}} \mathcal{C}^{\mathsf{J}} \tag{2.9.23}$$

PROPOSITION 2.9.24 ([121] Thm 4.5.2 and Thm 4.5.3) : *Right adjoints preserve limits. Left adjoints preserve colimits.*

COROLLARY 2.9.25 ([121] Cor 4.5.11) : *For any adjoint functors between abelian categories, the left adjoint is right exact and the right adjoint is left exact.*

DEFINITION 2.9.26 ([121] Def 4.5.12) : *A* reflective subcategory *of a category $\mathcal{C}$ is a full subcategory $\mathcal{D}$ so that the inclusion admits a left adjoint, called the localization:*

$$\mathcal{C} \underset{\mathscr{F}}{\overset{\text{localization}}{\underset{\hookrightarrow}{\overset{\longleftarrow}{\perp}}}} \mathcal{D} \tag{2.9.27}$$

DEFINITION 2.9.28 ((Pre)Sheaf on a topological space, [121] Def 3.3.4) : *Let $X$ be a topological space and write $\mathcal{O}(X)$ for the poset of open subsets, ordered by inclusion. An $\mathsf{J}$-indexed family of open subsets $\mathscr{U}_i \subseteq \mathscr{U}$ is said to* cover *$\mathscr{U}$ if the full diagram comprised of the sets $\mathscr{U}_i$ and the inclusions of their pairwise intersections $\mathscr{U}_i \cap \mathscr{U}_j$ has colimit $\mathscr{U}$. A* presheaf *$\mathscr{F}$ is a contravariant functor: $\mathcal{O}(X)^{\mathrm{op}} \longrightarrow \mathcal{S}\mathrm{et}$ from the opposite category $\mathcal{O}(X)^{\mathrm{op}}$ to the category of Sets. The presheaf $\mathscr{F}$ is* a sheaf *if it preserves these colimits, sending them to limits in the category $\mathcal{S}\mathrm{et}$ of sets. That means that the hypothesis to be a sheaf is the following:*

*For any open cover $\{\mathscr{U}_i,\, i \in \mathsf{J}\}$, of $\mathscr{U}$, the following sequence is an equalizer diagram:*

$$0 \hookrightarrow \mathscr{F}(\mathscr{U}) \xhookrightarrow{\mathscr{F}(\mathscr{U}_i \hookrightarrow \mathscr{U})} \prod_{i \in \mathsf{J}} \mathscr{F}(\mathscr{U}_i) \underset{\mathscr{F}(\mathscr{U}_i \cap \mathscr{U}_j \hookrightarrow \mathscr{U}_j)\circ \pi_j}{\overset{\mathscr{F}(\mathscr{U}_i \cap \mathscr{U}_j \hookrightarrow \mathscr{U}_i)\circ \pi_i}{\rightrightarrows}} \prod_{i,j \in \mathsf{J}} \mathscr{F}(\mathscr{U}_i \cap \mathscr{U}_j) \tag{2.9.29}$$

REMARK 2.9.30 : In this case, the left adjoint: localization is called here: sheafification.

PROPOSITION 2.9.31 ([121] Def 4.5.15) : *If $\mathcal{D} \hookrightarrow \mathcal{C}$ is a reflective subcategory, then:*

- *The inclusion $\mathcal{D} \hookrightarrow \mathcal{C}$ creates all limits that the category $\mathcal{C}$ admits.*
- *The category $\mathcal{D}$ has all colimits that the category $\mathcal{C}$ admits, formed by applying the localization to the colimit in the category $\mathcal{C}$.*

EXAMPLE 2.9.32 ((homotopy, nerve and simplicial sets)) : We have an adjunction between the category of small categories Cat and the category of simplicial sets $\mathcal{S}\mathrm{et}^{\Delta^{\mathrm{op}}}$ (it is a presheaves category), on which we will return later in the subsection treating about homotopy functor, nerve functor and simplicial sets in subsection 2.1.1: *The combinatorial category of simplicial objects*:

$$\mathsf{Cat} \underset{\mathsf{Nerve}}{\overset{\mathsf{Homotopy}}{\underset{\hookrightarrow}{\overset{\longleftarrow}{\perp}}}} \mathcal{S}\mathrm{et}^{\Delta^{\mathrm{op}}} \tag{2.9.33}$$

Here $\Delta \subset \mathsf{Cat}$ is the full subcategory whose objects are the finite non-empty ordinals, in this context denoted by $[0], [1], [2], \ldots$, and whose morphisms are order-preserving functions (functors) between them. The embedding $\mathsf{Nerve} : \mathsf{Cat} \hookrightarrow \mathcal{S}\mathrm{et}^{\Delta^{\mathrm{op}}}$ carries a small category $\mathcal{C}$ to its $\mathsf{nerve}$ : $\mathsf{Nerve}(\mathcal{C}) : \Delta^{\mathrm{op}} \longrightarrow \mathcal{S}\mathrm{et}$ sends $[\mathrm{n}] = 0 \to 1 \to \cdots \to n$ to the set of presheaves $\mathsf{Nerve}(\mathcal{C})(n) := \mathsf{Cat}([n], \mathcal{C})$. The left adjoint $\mathsf{Homotopy} : \mathcal{S}\mathrm{et}^{\Delta^{\mathrm{op}}} \longrightarrow \mathsf{Cat}$ sends a simplicial set to its homotopy category. Restricting to the objects $[0]$, $[1] \in \Delta$, a simplicial set $X : \Delta^{\mathrm{op}} \longrightarrow \mathcal{S}\mathrm{et}$ has an underlying reflexive directed graph $X_0 \begin{smallmatrix}\longrightarrow\\ \longleftarrow\\ \longrightarrow\end{smallmatrix} X_1$. The homotopy category $\mathsf{Homotopy}(X)$ is a quotient of the free category generated by this reflexive directed graph modulo relations that arise from elements of the set $X_2$. In particular, the counit defines an isomorphism $\mathsf{Homotopy}(\mathsf{Nerve}(\mathcal{C})) \simeq \mathcal{C}$ for any category $\mathcal{C}$, proving that the inclusion is full and faithful.

We can then conclude this section with the following result:

COROLLARY 2.9.34 ([121] Def 4.5.16) : *The category* $\mathsf{Cat}$ *of small categories is complete and cocomplete (admits finite limits and finite colimits).*

*Proof.* The category of presheaves $\mathcal{S}\mathrm{et}^{\Delta^{\mathrm{op}}}$ is complete and cocomplete, with limits and colimits defined objectwise in the category $\mathcal{S}\mathrm{et}$. The reflective subcategory $\mathsf{Cat}$ inherits these limits, defined objectwise in the category $\mathcal{S}\mathrm{et}$, and also these colimits, defined by applying the homotopy category functor to the colimit in the category of presheaves $\mathcal{S}\mathrm{et}^{\Delta^{\mathrm{op}}}$. □

#### 2.9.1.2 Bifunctors, two-variable adjunctions and cartesian closeness.

DEFINITION 2.9.35 ([121] Def 4.3.7) : *A triple of bifunctors:*

$$\mathcal{A} \times \mathcal{B} \xrightarrow{\mathscr{F}} \mathcal{C}\,,\ \mathcal{A}^{\mathrm{op}} \times \mathcal{C} \xrightarrow{\mathscr{G}} \mathcal{B}\,,\ \mathcal{B}^{\mathrm{op}} \times \mathcal{C} \xrightarrow{\mathscr{H}} \mathcal{A}$$

*equipped with natural isomorphisms:*

$$\mathrm{Hom}_{\mathcal{C}}(\mathscr{F}(a,b),c) \simeq \mathrm{Hom}_{\mathcal{B}}(b,\mathscr{G}(a,c)) \simeq \mathrm{Hom}_{\mathcal{A}}(a,\mathscr{H}(b,c))$$

*defines a* two-variable adjunction*.*

*Particularly, when* $\mathscr{F} : \mathcal{C}\times\mathcal{C} \longrightarrow \mathcal{C}$ *defines some sort of monoidal product, its pointwise-defined right adjoints* $\mathscr{G}$ *and* $\mathscr{H}$ *are called its* left and right closures*, respectively. When these are isomorphic, the bifunctor* $\mathscr{F}$ *is called* closed*.*

DEFINITION 2.9.36 ((cartesian closeness) [121], Def 4.3.9) : *A* cartesian closed category *is a category* $\mathcal{C}$ *with finite products in which the product bifunctor* $\mathcal{C}\times\mathcal{C} \longrightarrow \mathcal{C}$ *is closed.*

EXAMPLE 2.9.37 : The categories $\mathcal{S}$et, Cat, and $\mathfrak{Fun}(\mathcal{C},\mathcal{S}\mathrm{et})$ are cartesian closed.

REMARK 2.9.38 : Let us note that the category of topological spaces is not cartesian closed. However, we construct in the following example of adjunction, the Steerod solution to this problem. This convenient category of spaces will be of great importance in the theory of stable homotopy in general, and in our solutions in particular.

EXAMPLE 2.9.39 : Let us define the subcategory cgHaus in the category of Hausdorff spaces of compactly generated Hausdorff spaces such that we only pick the Hausdorff spaces $X$ with the property that any subset $A \subset X$ that intersects each compact subset $K \subset X$ in a closed subset $A \cap K$ is itself closed in $X$. The inclusion cgHaus $\hookrightarrow$ Haus has a right adjoint Haus $\longrightarrow$ cgHaus called $k$-ification: the space $k(X)$ refines the topology on $X$ by adding to the collection of closed sets those subsets $A \subset X$ whose intersections with all compact subsets $K$ are closed.

Firts, the presence of this adjoint implies that cgHaus is complete and cocomplete. Moreover, the product in cgHaus is the $k$-ification of the product in Haus, which is preserved in the category Top of topological spaces. The pointwise right adjoint is given by the following construction of function spaces $\mathrm{Map}_{\mathsf{cgHaus}}(X,Y)$: the underlying set of $\mathrm{Map}_{\mathsf{cgHaus}}(X,Y)$ is the set of continuous maps $X \longrightarrow Y$, and the topology is the $k$-ification of the compact-open topology.

We then have constructed a subcategory of topological spaces, which is complete, co-complete, cartesian closed, and sufficiently large so as to contain the CW complexes. If we denote $S^1$ the classical circle, then the two-variable adjunction specialize to the following (1-variable) adjunction:

$$\mathsf{cgHaus} \underset{\mathrm{Map}_{\mathsf{cgHaus}}(S^1,-)}{\overset{S^1\times -}{\underset{\longleftarrow}{\overset{\longrightarrow}{\bot}}}} \mathcal{S}\mathrm{et}^{\Delta^{\mathrm{op}}} \tag{2.9.40}$$

The space $\mathrm{Map}_{\mathsf{cgHaus}}(S^1, X)$ is the free loop space on $X$; points in $\mathrm{Map}_{\mathsf{cgHaus}}(S^1, X)$ are loops in the space $X$. The slice category of cgHaus under the singleton space $*$ defines a convenient category $\mathsf{cgHaus}_*$ of based topological spaces.

The category $\mathsf{cgHaus}_*$ admits the two-variable functors:

$$\mathsf{cgHaus}_* \times \mathsf{cgHaus}_* \xrightarrow{\wedge} \mathsf{cgHaus}_*$$

$$\mathsf{cgHaus}_*^{\mathrm{op}} \times \mathsf{cgHaus}_* \xrightarrow{\mathrm{Map}_*} \mathsf{cgHaus}_*$$

$$\mathsf{cgHaus}_*^{\mathrm{op}} \times \mathsf{cgHaus}_* \xrightarrow{\mathrm{Map}_*} \mathsf{cgHaus}_*$$

where $\mathrm{Map}_*((X,x),(Y,y))$ denotes the based space of basepoint-preserving continuous functions from $(X,x)$ to $(Y,y)$.

The bifunctor $\wedge$ is called smash product. The space $\mathrm{Map}_*((S^1,pt),(X,x)) =: \Omega X$ is the based loop space on $(X,x)$. The functor $\Omega$ has a left adjoint $\Sigma X := S^1 \wedge X$, called the reduced suspension on the based space $(X,x)$.

We can then write the adjunction between loop and suspension in $\mathsf{cgHaus}_*$:

$$\mathsf{cgHaus} \underset{\Omega}{\overset{\Sigma}{\underset{\longleftarrow}{\overset{\longrightarrow}{\bot}}}} \mathcal{S}\mathrm{et}^{\Delta^{\mathrm{op}}} \tag{2.9.41}$$

$$\mathsf{cgHaus}(\Sigma X, Y) \simeq \mathsf{cgHaus}(X, \Omega Y) \tag{2.9.42}$$

#### 2.9.1.3 Existence of adjoints, accessibility and representability.

We have now what we could call "*technical*" theorems for pursuing our goal, about topoï, factorization systems, structured spaces, and existence of adjoints for (differential) cohesive topoï, even if these results have their powerful interest in themselves.

THEOREM 2.9.43 ((General Adjoint Functor Theorem) [121] Thm 4.6.3) : *Let* $\mathfrak{F} : \mathfrak{X} \to \mathcal{S}$ *be a continuous functor (which respects limits) whose domain is locally small and complete. Suppose that* $\mathfrak{F}$ *satisfies the following solution set condition:*

- *For every* $S \in \mathcal{S}$*, there exists a set of morphisms* $\Phi_s : \{f_i : s \longrightarrow \mathfrak{F}a_i\}$ *so that any* $f : s \to \mathfrak{F}a$ *factors through some* $f_i \in \Phi s$ *along a morphism* $a_i \to a$ *in* $\mathfrak{X}$*.*

*Then the functor* $\mathfrak{F}$ *admits a left adjoint.*

DEFINITION 2.9.44 (((Co)separating set) [121] Def 4.6.7) : *A* separating set *for a category* $\mathcal{C}$ *is a set* $\Phi$ *of objects that can distinguish between distinct parallel morphisms in the following sense: given* $f, g : x \rightrightarrows y$*, if* $f \neq g$ *then there exists some* $h : c \to x$ *with* $c \in \Phi$ *so that* $fh \neq gh$*. A* coseparating set *in* $\mathcal{C}$ *is a separating set in* $\mathcal{C}^{\mathrm{op}}$*.*

We will return to the definition of a subobject later, and also to closed and open immersions, much longer about Grothendieck contexts, which is a corner stone of our construction. However, let us give a very simple definition for the moment:

DEFINITION 2.9.45 ((Subobject) [121] Def 4.6.8) : *A* subobject *of an object* $c \in \mathcal{C}$ *is a monomorphism* $c' \rightarrowtail c$ *with codomain* $c$*.*

DEFINITION 2.9.46 ((Intersection of a family of subobjects) [121] Def 4.6.9) : *The* intersection of a family of subobjects *of an object* $c$ *in the category* $\mathcal{C}$ *is the limit of the diagram of monomorphisms with codomain* $c$*.*

THEOREM 2.9.47 ((Special Adjoint Functor Theorem) [121] Thm 4.6.10) : *Let us consider* $\mathfrak{F} : \mathfrak{X} \longrightarrow \mathcal{S}$ *a continuous functor whose domain is complete and whose domain and codomain are locally small. Furthermore, if the category* $\mathcal{S}$ *has a small coseparating set and every collection of subobjects of a fixed object in* $\mathfrak{X}$ *admits an intersection, then*

$\mathfrak{F}$ *admits a left adjoint.*

COROLLARY 2.9.48 ((Some completeness implies cocompleteness) [121] Cor 4.6.13) : *Suppose the category $\mathcal{C}$ is locally small, complete, has a small coseparating set, and has the property that every collection of subobjects of a fixed object has an intersection. Then the category $\mathcal{C}$ is also cocomplete.*

COROLLARY 2.9.49 ((Some completeness implies representability) [121] Cor 4.6.14) : *Suppose the category $\mathcal{C}$ is locally small, complete, has a small coseparating set, and has the property that every collection of subobjects of a fixed object has an intersection. Then any continuous functor $\mathfrak{F} : \mathfrak{X} \to \mathcal{S}$ is representable.*

THEOREM 2.9.50 ((Freyd's Representability Theorem) [121] Cor 4.6.15) : *Let us consider $\mathfrak{F} : \mathfrak{X} \to \mathcal{S}\mathrm{et}$ a continuous functor and suppose that $\mathfrak{X}$ is complete and locally small. If the functor $\mathfrak{F}$ satisfies the solution set condition:*

- *There exists a set $\Phi$ of objects of $\mathfrak{X}$ so that for any $x \in \mathfrak{X}$ and any element $s \in \mathfrak{F}x$, there exists an $y \in \Phi$, an element $t \in \mathfrak{F}y$, and a morphism $f : y \to x$ so that $\mathfrak{F}(f(t)) = s$*

*then the functor $\mathfrak{F} : \mathfrak{X} \to \mathcal{S}\mathrm{et}$ is representable.*

DEFINITION 2.9.51 (((locally) presentability and accessibility) [121] Def 4.6.16) : *Let $\kappa$ be a regular cardinal. A locally small category $\mathcal{C}$ is* locally $\kappa$-presentable *if it is cocomplete and if it has a set of objects $S$ so that:*

- *Every object in $\mathcal{C}$ can be written as a colimit of a diagram valued in the subcategory spanned by the objects in $S$*
- *For each object $s \in S$, the functor $\mathrm{Hom}_{\mathcal{C}}(s, ) : \mathcal{C} \longrightarrow \mathcal{S}\mathrm{et}$ preserves $\kappa$-filtered colimits.*

*A functor between locally $\kappa$-presentable categories is* accessible *if it preserves $\kappa$-filtered colimits.*

Between locally presentable categories, meaning categories that are locally $\kappa$-presentable for some $\kappa$, the adjoint functor theorem takes a particularly appealing form.

THEOREM 2.9.52 ((Left and right adjointness and (co)continuity) [121] Thm 4.6.17) : *A functor* $\mathfrak{F} : \mathfrak{X} \longrightarrow \mathfrak{Y}$ *between locally presentable categories:*

- *admits a right adjoint if and only if it is cocontinuous*
- *admits a left adjoint if and only if it is continuous and accessible*

#### 2.9.1.4 Monads, comonads and algebras.

DEFINITION 2.9.53 (Monad, [121] Def 5.1.1) : *We call a* monad *on a category* $\mathcal{C}$ *the following data:*

- *an endofunctor* $T : \mathcal{C} \longrightarrow \mathcal{C}$
- *a unit natural transformation* $\eta : \mathbf{1}_{\mathcal{C}} \Rightarrow T$
- *a multiplication natural transformation* $\mu : T^2 \Rightarrow T$*, so that the following diagrams commute in* $\mathrm{Hom}_{\mathsf{CAT}}(\mathcal{C}, \mathcal{C})$ *:*

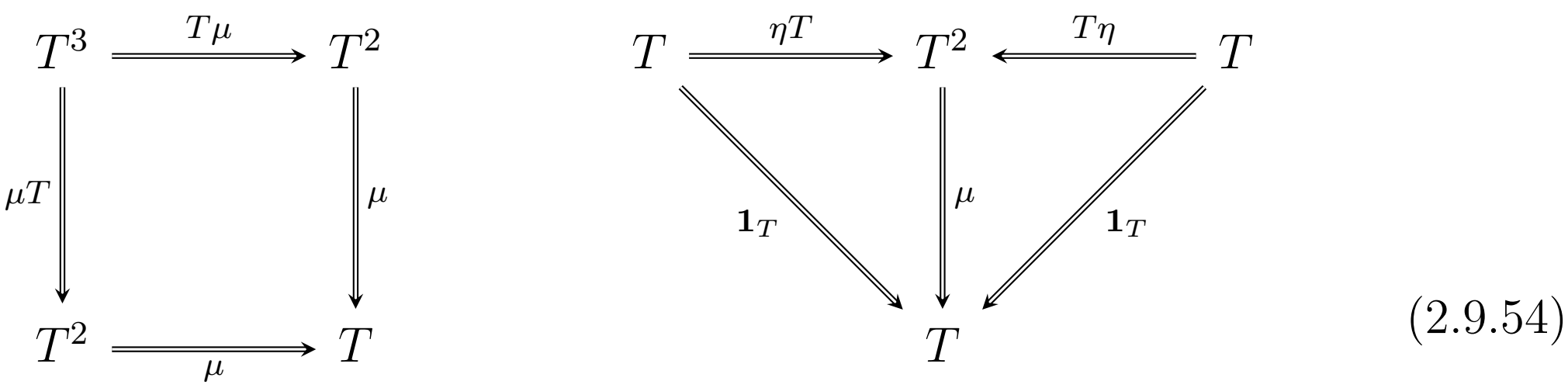

(2.9.54)

REMARK 2.9.55 : The category $\mathrm{Hom}_{\mathsf{CAT}}(\mathcal{C}, \mathcal{C})$ of endofunctors on $\mathcal{C}$ is a strict monoidal category, where the coherence diagrams are made of natural isomorphisms which are identities. A monad on $\mathcal{C}$ is a monoid in this monoidal category, where the binary functor $\mathrm{Hom}_{\mathsf{CAT}}(\mathcal{C},\mathcal{C}) \times \mathrm{Hom}_{\mathsf{CAT}}(\mathcal{C},\mathcal{C}) \longrightarrow \mathrm{Hom}_{\mathsf{CAT}}(\mathcal{C},\mathcal{C})$ is composition, and the unit object is the identity functor $\mathbf{1}_{\mathcal{C}} \in \mathrm{Hom}_{\mathsf{CAT}}(\mathcal{C},\mathcal{C})$.

proposition 2.9.56 : *Any adjunction:*

$$\mathcal{C} \underset{\mathscr{F}}{\overset{\mathscr{G}}{\rightleftarrows}} \mathcal{D} \quad (\top) \qquad \eta : \mathbf{1}_{\mathcal{C}} \Rightarrow \mathscr{G}\mathscr{F} \qquad \epsilon : \mathscr{F}\mathscr{G} \Rightarrow \mathbf{1}_{\mathcal{C}} \tag{2.9.57}$$

*gives rise to a monad on the category $\mathcal{C}$, which is the domain of the left adjoint, with:*

- *the endofunctor $T : \mathcal{C} \longrightarrow \mathcal{C}$ is $T := \mathscr{G}\mathscr{F}$*
- *the unit natural transformation $\eta : \mathbf{1}_{\mathcal{C}} \Rightarrow T$ is $\eta : \mathbf{1}_{\mathcal{C}} \Rightarrow \mathscr{G}\mathscr{F}$*
- *the multiplication natural transformation $\mu : T^2 \Rightarrow T$ comes from the functor: $\mathscr{G}\epsilon\mathscr{F} : \mathscr{G}\mathscr{F}\mathscr{G}\mathscr{F} \Rightarrow \mathscr{G}\mathscr{F}$ in the following whiskered diagram and square diagram which commute in* $\mathrm{Hom}_{\mathsf{CAT}}(\mathcal{C},\mathcal{C})$ *:*

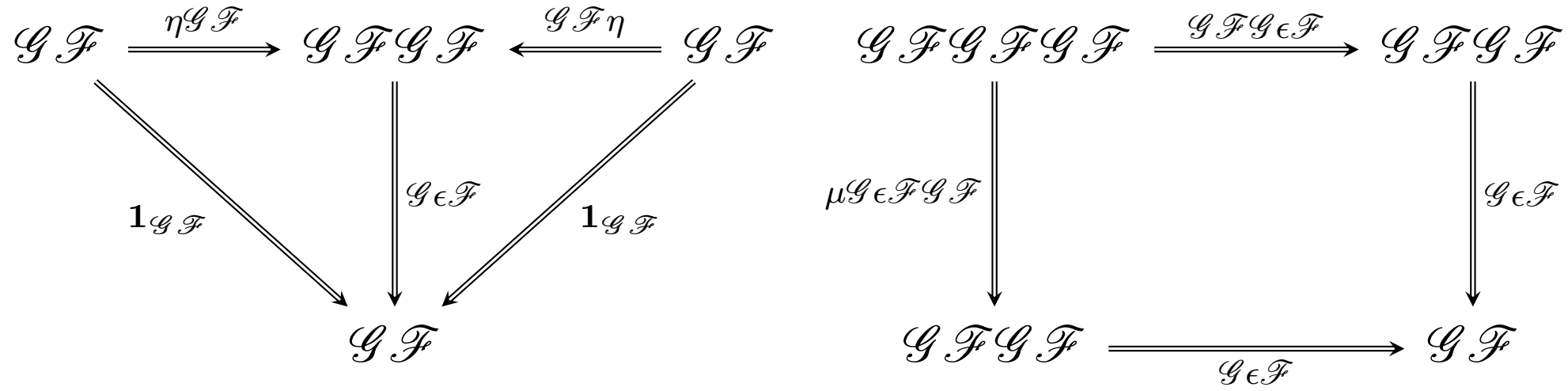

(2.9.58)

definition 2.9.59 (Comonad, [121] Def 5.1.6) : *We call a* **comonad** *on a category $\mathcal{C}$ a monad on the category $\mathcal{C}^{\mathrm{op}}$, that is the following data:*

- *an endofunctor $K : \mathcal{C} \longrightarrow \mathcal{C}$*
- *a counit natural transformation $\epsilon : K \Rightarrow \mathbf{1}_{\mathcal{C}}$*

- *a comultiplication natural transformation $\delta : K \Rightarrow K^2$, so that the following diagram commute in* $\mathrm{Hom}_{\mathsf{CAT}}(\mathcal{C}, \mathcal{C})$ *:*

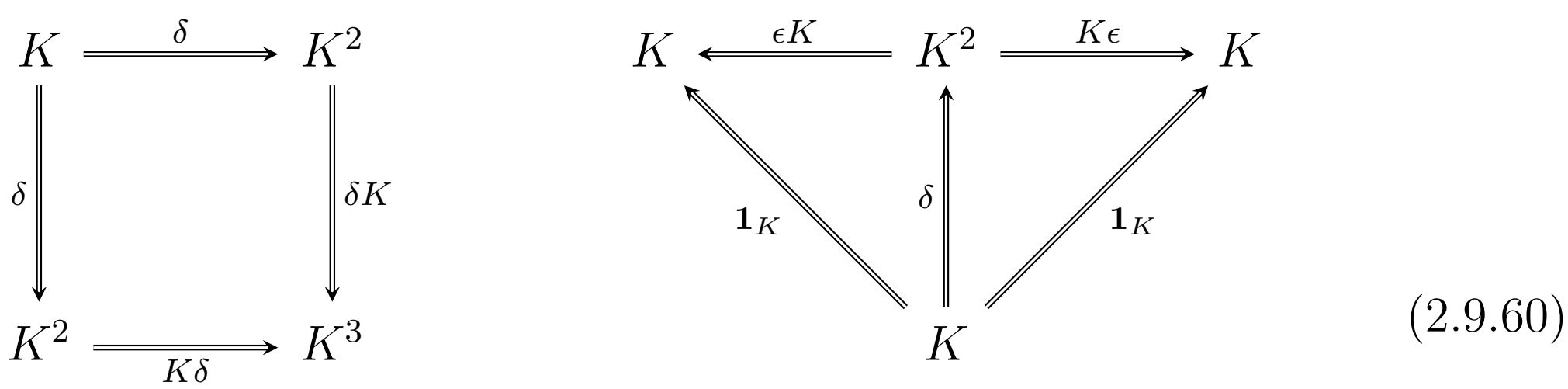

(2.9.60)

REMARK 2.9.61 : Dually, a comonad is a comonoid in the category of endofunctors of $\mathcal{C}$, and any adjunction induces a comonad on the domain of its right adjoint.

DEFINITION 2.9.62 : *We call an* **idempotent monad** *the following data:*

- *The multiplication $\mu : T^2 \Rightarrow T$ is a natural isomorphism (hence, the appellation "*idempotent*")*

- *Each component of $\mu : T^2 \Rightarrow T$ is a monomorphism*

- *The natural transformations $\eta T$, $T\eta : T \Rightarrow T^2$ are equal*

REMARK 2.9.63 : An adjunction associated to a reflective subcategory (see **??**) of the category $\mathcal{C}$ induces an idempotent monad on the category $\mathcal{C}$.

A question may be posed at this stage : Do all monads come from an adjunction ?

The answer is: Yes, but with the necessity to be careful : There are lots of adjunctions for a single monad. In that sense, we can consider that monads are "*generating poles*" in the sense of monopoles in electromagnetism of physics.

Two adjonctions play a particular role, between the Kleisli category, usually denoted $C_T$ and the Eilenberg-Moore category, usually denoted $C^T$ and also called "*category of $T$-algebras*".

Here are two construction models :

- Let $A$ be an affine space on a $\Bbbk$-vector space $V$:

$$\mathcal{S}\text{et} \underset{\textit{Aff. Build.}_{\Bbbk}}{\overset{\mathscr{F}\text{orget}}{\underset{\longrightarrow}{\overset{\longleftarrow}{\top}}}} \mathsf{Aff}_{\Bbbk} \tag{2.9.64}$$

  It is an adjointness derived from the monade $T := (\textit{Aff. Build.}_{\Bbbk}, \eta, \mu)$.

  We construct the $T$-algebra in the monadic sense, here of $\Bbbk$-affine space $A$, thanks to the monade $T := (\textit{Aff. Build.}_{\Bbbk} : \mathcal{S}\text{et} \longrightarrow \mathcal{S}\text{et}, \eta, \mu)$, where we choose $A$ initially a set without structure which we equip with the additive homomorphisme $V \times A \longrightarrow A$ which give the structure of $T$-algebra in the monadic sense, that is of affine space with an action by translation of the $\Bbbk$-espace vectoriel $V$. The monad $T := (\textit{Aff. Build.}_{\Bbbk}, \eta, \mu)$ in then the "*plan of construction*" of any affine space, in the classical sense of affine spaces with action by translation of a $\Bbbk$-vector on the elements of a set which initially have no structure. The monad $T := (\textit{Aff. Build.}_{\Bbbk} : \mathcal{S}\text{et} \longrightarrow \mathcal{S}\text{et}, \eta, \mu)$ captures the essence of the concept of affine space.

- Let us considet a $R$-module $M$ with $R$ a unitary ring:

$$\mathsf{Ab} \underset{R \otimes_{\mathbb{Z}} -}{\overset{\mathscr{F}\text{orget}}{\underset{\longrightarrow}{\overset{\longleftarrow}{\top}}}} \mathsf{Mod}_R \tag{2.9.65}$$

  We have an adjunction derived from the monad $T := (R \otimes_{\mathbb{Z}} - : \mathsf{Ab} \longrightarrow \mathsf{Ab}, \eta, \mu)$

  We construct the structure of $T$-algebra in the monadic sense on $M$ initially considered as an abelian group which we equip with an homomorphisme $R \otimes_{\mathbb{Z}}$ $M \longrightarrow M$ which gives the structure of $T$-algebra in the monadic sense. The monade $T := (R \otimes_{\mathbb{Z}} - : \mathsf{Ab} \longrightarrow \mathsf{Ab}, \eta, \mu)$ is then the "*plan of construction*" of any $R$-module, in the classical sense of $R$-modules with multiplication by a scalar of $R$ of the elements of an abelian group $M$. The monad $T := (R \otimes_{\mathbb{Z}} - : \mathsf{Ab} \longrightarrow \mathsf{Ab}, \eta, \mu)$ capture the essence of module.

PROPOSITION 2.9.66 ([121] Prop 5.2.12) : *Kleisli category is initial in the category* $\mathsf{Adj}_T$ *and Eilenberg-Moore is final in the category* $\mathsf{Adj}_T$. *That is, for any adjunction*

$$\mathcal{C} \underset{\mathscr{F}}{\overset{\mathscr{G}}{\underset{\longrightarrow}{\overset{\longleftarrow}{\top}}}} \mathcal{D} \qquad \eta : \mathbf{1}_{\mathcal{C}} \Rightarrow \mathscr{G}\mathscr{F} \qquad \epsilon : \mathscr{F}\mathscr{G} \Rightarrow \mathbf{1}_{\mathcal{C}} \tag{2.9.67}$$

*inducing the monad* $(T, \eta, \mu)$ *on* $\mathcal{C}$, *there exists unique functors*

$$\begin{array}{ccccc} \mathcal{C}_T & \overset{J}{\underset{\exists!}{\dashrightarrow}} & \mathcal{D} & \overset{K}{\underset{\exists!}{\dashrightarrow}} & \mathcal{C}^T \\ & {}_{\mathscr{G}_T}\searrow\nwarrow^{\mathscr{F}_T} \; \dashv & \mathscr{G}\downarrow \vdash \uparrow\mathscr{F} & \vdash \; {}^{\mathscr{G}^T}\swarrow\nearrow_{\mathscr{F}^T} & \\ & & \mathcal{C} & & \end{array} \tag{2.9.68}$$

*commuting with the left and right adjoints.*

DEFINITION 2.9.69 : *The unique functor of proposition 2.9.66 from the Kleisli category for any monad to the Eilenberg-Moore category, which commutes with the free functor* $\mathscr{F}^T$ *and the forgetful functor* $\mathscr{G}^T$ *from the underlying category* $\mathcal{C}$ *is called the* canonical comparison functor.

DEFINITION 2.9.70 ([121], Def. 5.3.1.) :

- *An adjunction* $\mathcal{C} \underset{\mathscr{F}}{\overset{\mathscr{G}}{\underset{\longrightarrow}{\overset{\longleftarrow}{\top}}}} \mathcal{D}$ *is* monadic *if the canonical comparison functor* $K$ *of proposition 2.9.66 from* $\mathcal{D}$ *to the category of* $T$*-algebras* $\mathcal{C}^T$ *for the induced monad* $T$ *on* $\mathcal{C}$ *defines an equivalence of categories.*

- *A functor* $\mathscr{G} : \mathcal{D} \longrightarrow \mathcal{C}$ *is* monadic *if it admits a left adjoint that defines a monadic adjunction.*
- *A functor* $\mathscr{G} : \mathcal{D} \longrightarrow \mathcal{C}$ *is* strictly monadic *if it is monadic and the comparison functor defines an isomorphism of categories.*

PROPOSITION 2.9.71 ([121] Prop. 5.3.3) : *Inclusions of reflexive subcategories* $\mathcal{C} \underset{i}{\overset{L}{\underset{\hookrightarrow}{\overset{\longleftarrow}{\top}}}} \mathcal{D}$ *with reflector (or localisation)* $L$ *define monadic functors from* $\mathcal{D}$ *to* $\mathcal{C}$*. The induced endofunctor* $L : \mathcal{D} \longrightarrow \mathcal{D}$ *defines a monad on* $\mathcal{D}$ *such that* $L^2\mathcal{D} \simeq L\mathcal{D}$*, that is an idempotent monad.*

DEFINITION 2.9.72 ([121], Def. 5.4.5) : *A* split coequalizer *diagram consists of maps:*

$$A \underset{g}{\overset{f}{\rightrightarrows}} B \xrightarrow{h} C, \quad t : B \to A, \quad s : C \to B \tag{2.9.73}$$

*such that :* $h \circ f = h \circ g$, $h \circ s = \mathbf{1}_C$, *and* $f \circ t = s \circ h$.

EXAMPLE 2.9.74 : For any algebra $(A, \alpha : TA \longrightarrow A)$ in the Eilenberg-Moore category $\mathcal{C}^T$ induced by a monad $(T, \eta, \mu)$ on $\mathcal{C}$, the diagram

$$T^2A \underset{\mu_A}{\overset{T\alpha}{\rightrightarrows}} TA \xrightarrow{\alpha} C, \quad \eta_{TA} : TA \to T^2A, \quad \eta_A : C \to TA \tag{2.9.75}$$

is a split coequalizer diagram in $\mathcal{C}$.

DEFINITION 2.9.76 : *Given a functor* $\mathscr{G} : \mathcal{D} \longrightarrow \mathcal{C}$

- *A* $\mathscr{G}$-split coequalizer *is a parallel pair* $f, g : A \rightrightarrows B$ *in* $\mathcal{D}$ *together with an extension of the pair* $\mathscr{G}f, \mathscr{G}g : \mathscr{G}A \rightrightarrows \mathscr{G}B$ *to a split coequalizer diagram in* $\mathcal{C}$*:*

$$\mathscr{G}A \underset{\mathscr{G}g}{\overset{\mathscr{G}f}{\rightrightarrows}} \mathscr{G}B \xrightarrow{h} C, \quad t : \mathscr{G}B \to \mathscr{G}A, \quad s : C \to \mathscr{G}B \tag{2.9.77}$$

- $\mathscr{G}$ creates coequalizers of $\mathscr{G}$-split pairs *if any $\mathscr{G}$-split coqualizer admits a coequalizer in $\mathcal{D}$ whose image under $\mathscr{G}$ is isomorphic to the coequalizer $\mathscr{G}h$ of $\mathscr{G}f$ and $\mathscr{G}g$ underlying the given $\mathscr{G}$-split coequalizer diagram in $\mathcal{C}$, and if any such diagram in $\mathcal{D}$ is a coequalizer*

- $\mathscr{G}$ strictly creates coequalizers of $\mathscr{G}$-split pairs *if any $\mathscr{G}$-split coqualizer admits a unique lift to a coequalizer in $\mathcal{D}$ for the given parallel pair.*

PROPOSITION 2.9.78 : *For any monad $(T, \eta, \mu)$ acting on a category $\mathcal{C}$, the monadic forgetful functor $\mathscr{G}^T : \mathcal{C}^T \longrightarrow \mathcal{C}$ strictly creates coequalizers of $\mathscr{G}^T$-split pairs.*

THEOREM 2.9.79 (Monadicity [21],[121] Thm 5.1.1) : *A right adjoint functor $\mathscr{G} : \mathcal{D} \longrightarrow \mathcal{C}$ is monadic if and only if it creates coequalizers of $\mathscr{G}$-split pairs.*

By recalling the right hand part of the diagram of proposition 2.9.66:

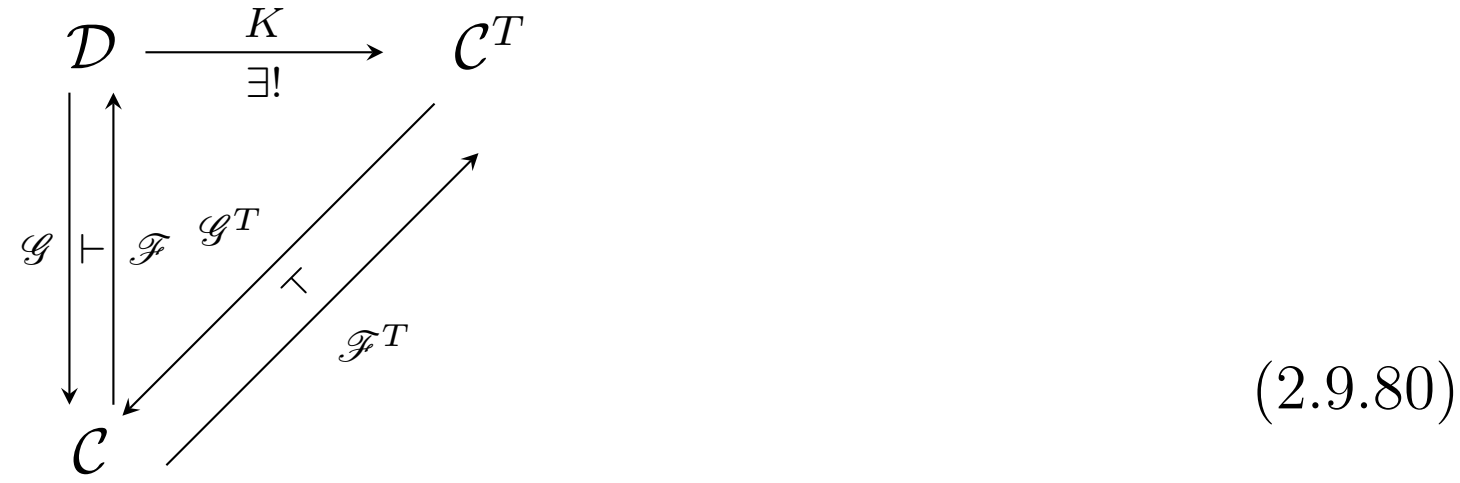

(2.9.80)

theorem 2.9.79 means exactly the following equivalence:

- $K$ is an equivalence (resp. isomorphism) of categories

- $\mathscr{G}$ creates (resp. strictly creates) coequalizers of $\mathscr{G}$-split pairs

By choosing for $\mathscr{G}$ the forgetful functor and its left adjoint the free construction, we get that the categories of monoids, groups, abelian groups, rings, commutative rings, unital commutative rings, modules over a ring, $\Bbbk$-vector spaces, $\Bbbk$-affine spaces, lattices, and pointed sets are monadic.

The category of unital commutative rings being monadic, the opposite category, that of schemes, is comonadic.

More generally, a common unification can be written as follows:

DEFINITION 2.9.81 : *A functor is said* finitary *if it preserves filtered colimits.*

In particular, a monad $T : \mathcal{C} \longrightarrow \mathcal{C}$ is finitary if it preserves filtered colimits in $\mathcal{C}$. If a right adjoint is finitary, then so is its monad because its left adjoint preserves all colimits.

DEFINITION 2.9.82 : *A category $\mathcal{A}$ is said to be a* category of models for an algebraic theory *if there is a finitary monadic functor $\mathcal{G} : \mathcal{A} \longrightarrow \mathcal{S}\mathrm{et}$.*

Any category of models for an algebraic theory is locally finitely presentable ([2], Thm 2.78).

Although the monad is not finitary, we have the monadicity of the forgetful functor between the category of compact Hausdorff spaces and the category of Sets.

THEOREM 2.9.83 ([112]) : *The contravariant power set functor $P : \mathcal{S}\mathrm{et}^{\mathrm{op}} \longrightarrow \mathcal{S}\mathrm{et}$ is monadic.*

THEOREM 2.9.84 ([112], Beck-Chevalley condition in [104] for subobject lattices but in [82] for base change adjoint functors) : *Let $\mathcal{E}$ be an elementary topos. The contravariant power set functor $P : \mathcal{E}^{\mathrm{op}} \longrightarrow \mathcal{E}$ is monadic.*

The condition on split coequalizers allows us to conclude with the following theorems.

THEOREM 2.9.85 : *A monadic functor $\mathcal{G} : \mathcal{C} \longrightarrow \mathcal{D}$ creates any limits that $\mathcal{D}$ has, and it creates any colimits that $\mathcal{D}$ has and the monad and its square preserve.*

*Proof.* It suffices to prove it for the forgetful functor from $\mathcal{C}^T \longrightarrow \mathcal{C}$. □

COROLLARY 2.9.86 : *The inclusion of a reflective subcategory creates all limits. In particular, a reflective subcategory of a complete category is complete.*

COROLLARY 2.9.87 : *Any category that is monadic over* $\mathcal{S}et$ *is complete, with limits created by the monadic forgetful functor.*

COROLLARY 2.9.88 : *The category* $\mathcal{S}et$ *is cocomplete. Any category of models for an algebraic theory is cocomplete (for the colimits preserved by the monad).*

THEOREM 2.9.89 : *If* $T : \mathcal{C} \longrightarrow \mathcal{C}$ *is a finitary monad on a complete and cocomplete, locally small category* $\mathcal{C}$*, then the Eilenberg-Moore category* $\mathcal{C}^T$ *of* $T$*-algebras is also complete and cocomplete.*

In particular, the category $\mathcal{S}et$ is complete and cocomplete.

COROLLARY 2.9.90 : *Any category of models for an algebraic theory is cocomplete*

That is without the restriction of the colimits preserved by the monad.

THEOREM 2.9.91 : *Even if the monad for compact Hausdorff spaces is not finitary, the category* cHaus *is cocomplete, and is a reflective subcategory of* $\mathcal{T}op$*, the category of topological spaces.*

PROPOSITION 2.9.92 (Adunctions as Kan extension) : *If* $F \dashv G$ *is an adjunction with unit* $\eta : \mathbf{1} \Rightarrow \mathscr{G}\mathscr{F}$ *and counit* $\epsilon : \mathscr{F}\mathscr{G} \Rightarrow \mathbf{1}$*, then* $(G, \eta)$ *is a left Kan extension of the identity functor along* $\mathscr{F}$ *and* $(\mathscr{F}, \epsilon)$ *is a right Kan extension of the identity functor along* $\mathscr{G}$*. Conversely, if* $(G, \eta \Rightarrow \mathscr{G}\mathscr{F})$ *is a left Kan extension of the identity along* $\mathscr{F}$ *and if* $\mathscr{F}$ *preserves this Kan exetnsion, then we have an adjointness* $\mathscr{F} \dashv \mathscr{G}$ *with unit* $\eta$*.*

### 2.9.2 Dependent homotopy-types and locally cartesian closed ∞-categories.

DEFINITION 2.9.93 : *A category* $\mathcal{C}$ *is called* cartesian closed *if it has cartesian products* $X \times Y$ *of all objects* $X, Y \in \mathcal{C}$ *and if there is for each* $X \in \mathcal{C}$ *a mapping space functor* $[X, -] : \mathcal{C} \longrightarrow \mathcal{C}$*, characterized by the fact that there is a bijection at the level of Hom sets:*

$$\mathrm{Hom}_{\mathcal{C}}(X \times Y, Z) \simeq \mathrm{Hom}_{\mathcal{C}}(X, [Y, Z]) \simeq \mathrm{Hom}_{\mathcal{C}}(X, Z^Y) \tag{2.9.94}$$

*which is natural in $X, Y, Z$. We say that the category $\mathcal{C}$ has* exponential objects *and we note that having cartesian products means in particular having a terminal object (for the empty product). The category $\mathcal{C}$ is called* locally cartesian closed *if for each object $X \in \mathcal{C}$, the slice category $\mathcal{C}_{/X}$ is cartesian closed.*

REMARK 2.9.95 : We will see below that the main examples of locally cartesian closed categories which interest us are the topoi. However, it is useful to interpret local cartesian closure in terms of base change thanks to the the following:

PROPOSITION 2.9.96 (See Beck-Chevalley condition in [82]) : *If the category $\mathcal{C}$ is locally cartesian closed, then for any morphism $f : X \longrightarrow Y$ in $\mathcal{C}$, there exists an adjoint triple of functors between the slice categories over $X$ and $Y$, called* base change *functors, such that the first diagram is a pullback diagram along $f$ for all $C \in \mathcal{C}$ and $c \in \mathcal{C}_Y$, and $f_!$ and $f_*$ are left and right adjoint to the pullback functor $f^*$ in the second diagram.*

$$\begin{array}{ccc} C & \overset{\mathbf{1}_c}{=\!=\!=} & C \\ {\scriptstyle f^*(c)}\downarrow & \lrcorner & \downarrow{\scriptstyle c} \\ X & \underset{f}{\longrightarrow} & Y \end{array} \qquad\qquad \begin{array}{ccc} & \xrightarrow{\quad f_! \quad} & \\ & \perp & \\ \mathcal{C}_X & \xleftarrow{\quad f^* \quad} & \mathcal{C}_Y \\ & \perp & \\ & \xrightarrow{\quad f_* \quad} & \end{array} \tag{2.9.97}$$

*Conversely, if a category $\mathcal{C}$ has pullbacks and has for every morphism $f : X \longrightarrow Y$ in $\mathcal{C}$ a left and right adjoint $f_!$ and $f_*$ in the second diagram to the pullback functor $f^*$ from the first diagram, then $\mathcal{C}$ is a locally cartesian closed category.*

DEFINITION 2.9.98 (Dependent type theory/Categorical semantics [105]) : *If $\mathcal{C}$ is a locally cartesian closed category, one says that its internal logic is a* dependent type theory *(or it provides* categorical semantics *for dependent type theory) if:*

- *the objects of $\mathcal{C}$ are called* types
- *the objects of $\mathcal{C}_{/\Gamma}$ are called the* types in context $\Gamma$ *(or the* types dependent in the context $\Gamma$*), and we denote:*

$$\Gamma \vdash X : \text{Type}$$

- *a morphism* $\{*\} \longrightarrow X$ *from the terminal objet into any objet* $X$ *in the slice category* $\mathcal{C}_{/\Gamma}$ *is called* a term of type $X$ in context $\Gamma$ *and we denote:*

$$\Gamma \vdash x : X$$

*or more explicitly*

$$a : \Gamma \vdash x(a) : X(a)$$

- *given a morphism* $f : \Gamma_1 \longrightarrow \Gamma_2$ *in* $\mathcal{C}$ *with its induced base change adjoint triple of functors between slice categories:*

$$\mathcal{C}_{/\Gamma_1} \underset{f_*}{\overset{f_!}{\underset{\perp}{\overset{\perp}{\longleftarrow f^* \longrightarrow}}}} \mathcal{C}_{/\Gamma_2} \qquad (2.9.99)$$

*then*

  - *given a morphism* $\{*\} \longrightarrow X$ *from the terminal objet into any objet* $X$ *in the slice category* $\mathcal{C}_{/\Gamma_2}$, *hence a terme* $\Gamma_2 \vdash x : X$, *then its pullback by* $f^*$ *is denoted by* substitution of variables

$$a : \Gamma_1 \vdash x(f(a)) : X(f(a))$$

  - *given an object* $X \in \mathcal{C}_{/\Gamma_1}$, *its image* $f_!(X) \in \mathcal{C}_{/\Gamma_2}$ *is called the* dependent sum *of* $X$ *along* $f$ *and is denoted as*

$$\Gamma_2 \vdash \sum_f X \,:\, \text{Type}$$

  - *given an object* $X \in \mathcal{C}_{/\Gamma_1}$, *its image* $f_*(X) \in \mathcal{C}_{/\Gamma_2}$ *is called the* dependent product *of* $X$ *along* $f$ *and is denoted as*

$$\Gamma_2 \vdash \prod_f X \,:\, \text{Type}$$

- *the universal property of the ajdoints* $(f_! \vdash f^* \vdash f_*)$ *translates to evident rules for introducing and for transforming terms of these dependent sum/product types, called* term induction *and* term elimination *rules.*

When this synthactic translation is properly formalized, it yields an equivalent description of locally cartesian closed categories:

PROPOSITION 2.9.100 : *There is an equivalence of* 2*-categories between locally cartesian closed categories and dependent type theories.*

REMARK 2.9.101 : Given an object $X \in \mathcal{C}_{/\Gamma}$, its diagonal $X \longrightarrow X \times X$ regarded as an object of $\mathcal{C}_{(/\Gamma \times X \times X})$ serves as the identity type of $X$ denoted

$$\Gamma, (x_1, x_2) : X \times X \vdash (x_1 = x_2) : \text{Type}$$

REMARK 2.9.102 : Namely given two terms $x_1, x_2 : X$, then a term $\Gamma \vdash p : (x_1 = x_2)$ is as a morphism in $\mathcal{C}$ an element on the diagonal of $X \times X$ and in the type theory is a proof of equality of $x_1$ and $x_2$. If there is such a proof of equality then it is unique, since the diagonal is always a monomorphism.

REMARK 2.9.103 : But consider now the case that $\mathcal{C}$ in addition carries the structure of a model category. Then there is for each $X$ a path space object $X^I \longrightarrow X \times X$. Using this as the categorical semantics of identity types, instead of the plain diagonal $X \longrightarrow X \times X$, means to make identity behave instead like if there were then possibily many equivalences between two terms of a given type, and many equivalences between equivalences and so on, that is higher equivalence homotopies. If $\mathcal{C}$ is moreover right proper as a model category and such that its cofibrations are precisely its monomorphisms, then there exists a variant of the dependent type theory reflecting these homotopy-theoretic identity-types. This is what is called homotopy type theory. At the same time, such a model category is a presentation for the homotopy-theoretic analogy of a locally cartesian closed category, which is then called a locally cartesian closed $(\infty,1)$-category (see A.3 in [98]).

PROPOSITION 2.9.104 ([60] for the first appearance, [18] for the refinement to $\infty$-toposes, [34, 138] for technical details, and [128] for our context) : *Up to equivalence, the internal type theory of a locally cartesian closed* $(\infty,1)$*-category is* homotopy type theory *(without necessary univalence) and conversely homotopy type theory (without necessary univalence) has* categorical semantics in locally cartesian closed $(\infty,1)$-category.

Urs Schreiber has [129] (Schreiber's Draft Prop. 4.1.2) an "internal" perspective of homotopy type theory and to the question of univalence for our setting. This will not be developped further in this thesis, but however, we already have at this stage the triple adjunction in a locally cartesian closed category (for example an $\infty$-topos), in which its internal logic provides categorical semantics for dependent type theory. This will be developed concretely in our setting in the following sections, without pushing the theory of homotopy type theory and the univalence axiom. A further exploration is thereby open for a future work.

## 2.10.Differential (infinitesimal) cohesion.

### 2.10.1 Six-functor formalism and Grothendieck context.

Let $\mathscr{C}$ and $\mathscr{D}$ be two complete, cocomplete, small and locally small$\infty$-categories.

DEFINITION 2.10.1 : *A functor* $\mathscr{F} : \mathscr{C} \longrightarrow \mathscr{D}$ *is said* conservative *if* $\mathscr{G} : A \to B$ *is a morphism in* $\mathscr{C}$ *such that* $\mathscr{F}(\mathscr{G})$ *is an isomorphism in* $\mathscr{D}$*, then* $\mathscr{G}$ *is an isomorphism in* $\mathscr{C}$*.*

For $f : \mathfrak{X} \longrightarrow \mathfrak{Y}$ in $\mathscr{C}$ we have a natural adjunction:

$$\mathcal{S}\mathbf{hv}_{\mathscr{D}}(\mathfrak{Y}) \underset{f_*}{\overset{f^*}{\rightleftarrows}} \mathcal{S}\mathbf{hv}_{\mathscr{D}}(\mathfrak{X}) \qquad (2.10.2)$$

If we take $Y = \{*\}$ the point, we have $f_* : \mathcal{S}\mathbf{hv}_{\mathscr{D}}(\mathfrak{X}) \longrightarrow \mathcal{S}\mathbf{hv}_{\mathscr{D}}(*) \simeq \mathbb{Z}-\mathcal{M}\mathbf{od}$ which coincides with the global sections functor. This means that the right derived functor coincides with the sheaf cohomology:

$$\mathbb{R}^n f_*(\mathscr{F}) \simeq H^n(\mathfrak{X}, \mathscr{F})$$

Another functor, named "*direct image with compact support*" is defined such that:

$$\Gamma(\mathfrak{X}_i, f_! \mathscr{F}) := \{s \in \Gamma(\mathfrak{X}_i, f_* \mathscr{F}) = \Gamma(f^{-1}(\mathfrak{X}_i, \mathscr{F}) \mid s \text{ has compact support}\}$$

Agian for $Y = \{*\}$ the point, we obtain the cohomology with compact support:

$$\mathbb{R}^n f_!(\mathscr{F}) = H^n_c(\mathfrak{X}, \mathscr{F})$$

REMARK 2.10.3 : The functor $f_!$ does not admit an adjoint in general.

Let us now turn to derived categories of sheaves. We obtain an adjoint:

$$\mathscr{D}(\boldsymbol{\mathcal{S}}\mathbf{hv}_{\mathscr{D}}(\mathfrak{X})) \underset{f^!}{\overset{\mathbb{R}f_!}{\underset{\longleftarrow}{\overset{\longrightarrow}{\perp}}}} \mathscr{D}(\boldsymbol{\mathcal{S}}\mathbf{hv}_{\mathscr{D}}(\mathfrak{X})) \tag{2.10.4}$$

We call $f^!$ the "*exceptional inverse image*" and $f_!$ the "*exceptional direct image image*".

DEFINITION 2.10.5 (Six-functor formalism) : *Together with the functors* "left derived tensor product" *and* "internal hom of sheaves" *and by dropping the indicators for derived functors and derived categories, we obtain* "our"– "six-functor formalism" *with the following set of adjunctions, where we know that in our situation, $f_! = f_*$ if $f$ is proper :*

$$(\otimes, \underline{\mathrm{Hom}}) \tag{2.10.6}$$

$$(f^* \dashv f_* \dashv f^!) : \ \boldsymbol{\mathcal{S}}\mathbf{hv}_{\mathscr{D}}(\mathfrak{Y}) \ \overset{f^*}{\longrightarrow} \underset{f^!}{\overset{f_! = f_*}{\longleftarrow}} \overset{}{\longrightarrow} \boldsymbol{\mathcal{S}}\mathbf{hv}_{\mathscr{D}}(\mathfrak{X}) \tag{2.10.7}$$

DEFINITION 2.10.8 (Grothendieck context) : *The construction and the situation described above is what we call a* Grothendieck context.

### 2.10.2 Cohesive topos.

We will enhance the situation by adjoining an additional adjoint to this Grothendieck context, thereby achieving our initial goal: to be able to measure volumes using De Rham's volume form as an integral along thickened infinitesimals pieces.

DEFINITION 2.10.9 : *We call and denote* $\mathcal{E}$ *a* cohesif 1-topos *on a 1-topos* $\mathcal{S}$, *a topos* $\mathcal{E}$ *equipped with a* geometric morphism *with the topos* $\mathcal{S}$*:*

$$(f^* \dashv f_*): \ \mathcal{E} \underset{f_*}{\overset{f^*}{\leftrightarrows}} \mathcal{S} \tag{2.10.10}$$

*such that the topos* $\mathcal{E}$ *is a topos* strongly connected *and a* local topos*.*
*To be more precise, this means that:*

- *the topos* $\mathcal{E}$ *is a* locally connected topos*, that is that beyond* $f^*$ *there exists another adjunction on the left* $(f_! \dashv f^*)$
- *the topos* $\mathcal{E}$ *is a* connected topos*, that is that this new adjoint on the left* $f_!$ *preserves the terminal object, ou equivalently* $f^*$ *is fully faithful*
- *the topos* $\mathcal{E}$ *is a* strongly connected topos*, that is that* $f_!$ *preserves finite products*
- *the topos* $\mathcal{E}$ *is a* local topos*, that is that, more that to be connected, there exists beyond* $f_*$ *another adjonction on the right* $(f_* \dashv f^!)$

*This is resumed in a quadruple adjunction of adjoint* $\infty$*-functors of the form:*

$$(f_! \dashv f^* \dashv f_* \dashv f^!): \ \mathcal{E} \begin{array}{c} \xrightarrow{f_!} \\ \xleftarrow{f^*} \\ \xrightarrow{f_*} \\ \xleftarrow{f^!} \end{array} \mathcal{S} \tag{2.10.11}$$

*where* $f^*$ *is full and faithful and preserves finite products. The functor* $f^!$ *is then also full and faithful.*

More classicaly in the case of cohésion, we name the quadruple adjonction in the following way:

$$(\Pi_0 \dashv \mathrm{Disc} \dashv \Gamma \dashv \mathrm{coDisc}) : \ \mathcal{E} \underset{\mathrm{coDisc}}{\overset{\Pi_0}{\underset{\Gamma}{\overset{\mathrm{Disc}}{\rightleftarrows}}}} \mathcal{S} \tag{2.10.12}$$

The quadruple adjonction canonically induces a triple adjunction of endofoncteurs of $\mathcal{E}$:

$$(\int \dashv \flat \dashv \sharp) : \ \mathcal{E} \underset{\Gamma}{\overset{\Pi_0}{\underset{\mathrm{Disc}}{\rightleftarrows}}} \mathcal{S} \underset{\mathrm{coDisc}}{\overset{\mathrm{Disc}}{\underset{\Gamma}{\rightleftarrows}}} \mathcal{E} \tag{2.10.13}$$

où $\int = \mathrm{Disc} \circ \Pi_0$ , $\flat = \mathrm{Disc} \circ \Gamma$ , $\sharp = \mathrm{coDisc} \circ \Gamma$.

These three adjoint endofoncteurs being idempotent (co)monads on the topos $\mathcal{E}$, they are what we call modalities in the theory of types (in internal logic) of $\mathcal{E}$. We call these modalities:

$$\text{shape modality } \int \ \dashv \text{ flat modality } \flat \ \dashv \text{ sharp modality } \sharp$$

DEFINITION 2.10.14 : *We say that a cohesive topos has* Aufhebung of becoming *if the sharp modality* $\sharp$ *preserves the initial object:*

$$\sharp \varnothing \simeq \varnothing$$

REMARK 2.10.15 : In the axiom of Aufhebung of becoming, the extra exactness condition on the shape modality, which says in particular that the shape modality $\int$ preserves the terminal object:

$$\int * \simeq *$$

is a dual companion (co-Aufhebung) to Aufhebung of becoming, which leads us to consider a variant of the axioms of cohesion. The following variant of cohesion is:

- There is an adjoint triple of idempotent (co-)monads $(\int \dashv \flat \dashv \sharp)$

- such that the sharp modality $\sharp$ satisfies Aufhebung and the shape modality $\int$ satisfies co-Aufhebung.

DEFINITION 2.10.16 : *In a cohesive topos, the adjunction between the shape modality* $\int$ *and the flat modality* $\flat$ *gives a natural transformation:*

$$\flat\mathfrak{X} \xrightarrow{\flat\text{-counit}} \mathfrak{X} \xrightarrow{\int\text{-unit}} \int\mathfrak{X}$$

*This transformation may be called the transformation* "from points to their pieces"*.*
*If this is an epimorphism for all* $\mathfrak{X} \in \mathfrak{H}$ *a cohesive topos, we say that:*

pieces have points *or that* the Nullstellensatz is verified*.*

DEFINITION 2.10.17 : *In a cohesive topos* $\mathfrak{H}$*, the adjunction* $(f^* \dashv f_*)$*,* "from points to pieces" *gives a natural transformation:*

$$\flat\mathfrak{X} \xrightarrow{\flat\text{-counit}} \mathfrak{X} \xrightarrow{\int\text{-unit}} \int\mathfrak{X}$$

*We saw that this transformation may be called the transformation* "from points to pieces" *or* "points-to-pieces transform"*. But if we consider the map*

$$f_*\mathfrak{X} \longrightarrow f_!\mathfrak{X} := (f_*\mathfrak{X} \longrightarrow f_*f^*f_!\mathfrak{X} \longrightarrow f_!\mathfrak{X})$$

*we obtain by applying backwards* $f^*$ *and postcomposing with the unit of the adjunction* $(f^* \dashv f_*)$*, an action equivalent to just applying* $f^*$*, since the idempotency of* $\flat$ *gives an isomorphism which is the counit restricted to the discrete object* $f^*f_!\mathfrak{X}$*. Therefore, the points-to-pieces transformation and its adjunct are related by:*

$$(\flat\mathfrak{X} \xrightarrow{\flat\text{-counit}} \mathfrak{X} \xrightarrow{\int\text{-unit}} \int\mathfrak{X}) = f^*(f_*\mathfrak{X} \longrightarrow f_!\mathfrak{X})$$

*This means that, since* $f^*$ *is a full and faithful left and right adjoint, the* points-to-pieces transform *is an epi/iso/mono precisely if its adjunct, the* pieces-to-points transform*:* $f_*\mathfrak{X} \longrightarrow f_!\mathfrak{X}$ *is, respectively.*
*If this adjunct:*

$$\flat \mathcal{X} \xrightarrow{\epsilon^\flat} \mathcal{X} \xrightarrow{\eta^\sharp} \sharp \mathcal{X}$$

*is a monomorphism, we say that* discrete objects are concrete.

Here are some immediate consequencies of axioms:

PROPOSITION 2.10.18 : *The axioms* "pieces have points"*, i.e.*

- $\flat \mathcal{X} \longrightarrow \mathcal{X} \longrightarrow \int X$ *is an epimorphism for all* $\mathbf{X} \in \mathfrak{H}$

*and* "discrete objects are concrete"*, i.e.*

- $\flat \mathcal{X} \xrightarrow{\eta^\sharp_{\flat \mathcal{X}}} \sharp \flat \mathcal{X}$ *is a monomorphism for all* $\mathbf{X} \in \mathfrak{H}$

*are equivalent.*

REMARK 2.10.19 : For a cohesive 1-topos, if the pieces-to-points transform is an epimorphism then there is Aufhebung of the initial opposition $\varnothing \dashv *$ in that $\sharp \varnothing \simeq \varnothing$.

PROPOSITION 2.10.20 (JOHNST 11, 1.6) : *A sheaf topos that:*

- *is locally connected and connected*
- *satisfies* "pieces have points"
  *also is:*
- *local*
- *hyperconnected*
- *strongly locally connected*
- *cohesive*

PROPOSITION 2.10.21 (JOHNST 11, Thm 3.4) : *For a sheaf topos, the condition that it:*

- *is locally connected*
- *is connected*
- *satisfies* "pieces have points"
  *is equivalent to the condition that it:*
- *is locally connected*
- *is hyperconnected*
- *is local*

### 2.10.3 Remaks on the denomination of theories.

#### 2.10.3.1 Modality of elasticity.

The triple adjoints of (co)monads is also called $\Re(\mathfrak{X}) \dashv \Im(\mathfrak{X}) \dashv \&(\mathfrak{X})$ in Schreiber's Draft ([129]). In particular, between the two versions of Schreiber's work on differential cohesion, $\mathbf{H}_{\mathfrak{th}} := \Pi_{\mathrm{inf}}(\mathbf{H}) = \Im(\mathbf{H})$ denotes the same object: *the thickening* $\mathbf{H}_{\mathfrak{th}}$ *of the initial $\infty$-topos*" $\mathbf{H}$. The typography and the symbols just differ.

#### 2.10.3.2 Shape theory for $\infty$-topoï (Lurie HTT Chap 7).

Shape theory in Lurie's work HTT Chap 7 ([98] Chap 7) has a strong connexion with Schreiber's work of his shape theory of an $\infty$-topos $\int \mathbf{H}$.

#### 2.10.3.3 Lurie's structured and fractured topos.

We do not make the distinction in J. Lurie's work in SAG ([103]) about fractured topoï and structured spaces of DAG V ([91]), DAG VII ([92]) and all other DAG. Indeed, the theory of structured spaces and $\mathcal{G}$-geometries seems deeper in DAG's, and this is this one we used.

#### 2.10.3.4 Scholze's condensed anima.

The work of P. Scholze on condensed anima, in the setting of Hausdorff spaces and Stone duality, Frames and Locales, seems not to correspond to our, since Zariski topology is not Hausdorff. Howevever, since we introduced differential cohesion, a strong link is to make between the two theories. Unfortunetely, condensed anima and fractured (or structured) topoï do not verify the same set of adjunctions and condensed anima become fractured topoï after a costly work which has not been chosen to be made in this thesis.

### 2.10.4 Differential (infinitesimal) cohesion, Artin-Lurie representability theorem, and Formal moduli problems.

Let $\mathcal{CR}\mathrm{ing}_\infty^{\mathrm{cn}}$ be the $\infty$-category of connective $\mathbb{E}_\infty$-rings.

DEFINITION 2.10.22 ((Formal moduli problem: DAG XIV, Def 2.1.1)) : *An $\infty$-functor:*

$$\mathfrak{X} : \mathcal{CR}\mathrm{ing}_\infty^{\mathrm{cn}} \longrightarrow \mathcal{S}$$

*that is a presheaf on $(\mathcal{CR}\mathrm{ing}_\infty^{\mathrm{cn}})^{\mathrm{op}}$, is called "*cohesive*" or a "*formal moduli problem*" if $\mathfrak{X}$ sends the base ring to a contractible space, and if $\mathfrak{X}$ sends pullbacks of morphisms which are surjective on $\pi_0$ to pullbacks.*

REMARK 2.10.23 ((Cohesive ($\infty$,1)-presheaf on $\mathbb{E}_\infty$-rings)) : Infinitesimal cohesion (in the sens of J. Lurie) is the defining property that makes such a functor a formal moduli problem, hence equivalently an $L_\infty$-algebra. It is also one of the characteristics of a Deligne-Mumford stack, due to the Artin-Lurie representability theorem (it is one of the hypothesis to be representable by a Deligne-Mumford stack), hence it is satisfied by those functors which are infinitesimal approximations to geometric $\infty$-stacks.

### 2.10.5 Idea of the differential (infinitesimal) cohesion.

The cohesion axioms imply the existence of infinitesimal spaces. However, they do not provide ex nihilo in all cases an explicit synthetic notion of infinitesimal extension. We must introduce an additional axiom on a cohesive $\infty$-topos, which then implies a proper intrinsic notion of synthetic differential extension, compatible with cohesion.

As we already constructed crystalline cohomology which intrinsically gives the notion of synthetic differential extension, we are then equipped to what will refer now to "*differential cohesion*", to which we now have full access for our goal.

In a cohesive $\infty$-topos with differential cohesion, there are then well-defined notions of "*formal smoothness*", and "*De Rham spaces*".

### 2.10.6 Differential (infinitesimal) cohesion.

It is an additional structure to a cohesive $\infty$-topos, which encodes a refinement of the corresponding notion of cohesion named "*differential cohesion*".
We consider the inclusion $\mathfrak{H} \overset{i}{\hookrightarrow} \mathfrak{H}_{\mathfrak{th}}$ of $\infty$-topoi that treat objects of $\mathfrak{H}_{\mathfrak{th}}$ as cohesive infinitesimal neighborhoods of objects of $\mathfrak{H}$.

Definition 2.10.24 : *We call and denote* $\mathfrak{H}_{\mathfrak{th}}$ *an* infinitesimal cohesif neighbourhood *of a cohesif* $\infty$*-topos* $\mathfrak{H}$ *if* $\mathfrak{H}_{\mathfrak{th}}$ *is equipped with a quaduple set of adjoint* $\infty$*functors of the form:*

$$(i_! \dashv i^* \dashv i_* \dashv i^!) : \mathfrak{H} \underset{i^!}{\overset{i_!}{\underset{\leftarrow}{\overset{\hookrightarrow}{\underset{i_*}{\overset{i^*}{\underset{\hookrightarrow}{\overset{\leftarrow}{}}}}}}}} \mathfrak{H}_{\mathfrak{th}} \tag{2.10.25}$$

*where* $i_!$ *is full and faithful and preserves finite products. The functor* $i_*$ *is then fully faithfull too.*
*Recoprocally, we will say that the data of Equation 2.10.25 equip the cohesif* $\infty$*-topos* $\mathfrak{H}$ *of a* differential (infinitesimal) cohesion.

### 2.10.7 The thickening monad.

We then have four adjoint functors with three adjoint (co)monades $\Re/\Im/\&$ with their (co)units such that:

Definition 2.10.26 : *Let us define the adjoint triple of adjoint (∞,1)-functors corresponding to the adjoint quadruple:*

$$(\Re \dashv \Im \dashv \&) = (i_! \, i^* \dashv i_* \, i^* \dashv i_* \, i^!) : \mathbf{H}_{\mathfrak{th}} \longrightarrow \mathbf{H}_{\mathfrak{th}}$$

*We say that*

- $\Re$ *is the* reduction modality
- $\Im$ *is the* infinitesimal shape modality
- $\&$ *is the* infinitesimal flat modality

*An object in the full sub-∞-category*

- *of* $\Re$ *is called a* reduced object
- *of* $\Im$ *is called a* coreduced object

Definition 2.10.27 : *For* $\mathfrak{X} \in \mathbf{H}$*, we say that:*

- $\Im(\mathfrak{X})$ *is the* infinitesimal path ∞-groupoid *of* $\mathfrak{X}$ *(also called the* De Rham space *or* De Rham stack *of* $\mathfrak{X}$*) and the* $(i^* \dashv i_*)$*-unit* $\mathfrak{X} \longrightarrow \Im(\mathfrak{X})$ *is the* constant infinitesimal path inclu
- $\Re(\mathfrak{X})$ *is the* reduced cohesive ∞-groupoid underlying $\mathfrak{X}$ *and the* $(i^* \dashv i_*)$*-unit* $\Re(\mathfrak{X}) \longrightarrow \mathfrak{X}$ *is the* inclusion of the reduced part *of* $\mathfrak{X}$

Definition 2.10.28 : *For any* $\mathfrak{X} \in \mathbf{H}$*, we denote* Jet *the base change:*

$$\mathrm{Jet} : \mathbf{H}_{/\mathfrak{X}} \underset{i_*}{\overset{i^*}{\rightleftarrows}} \mathbf{H}_{/\Im(\mathfrak{X})}$$

*For* $(\mathfrak{X} \longrightarrow \mathfrak{Y}) \in \mathbf{H}_{/\mathfrak{Y}}$ *, we call Jet*$(\mathfrak{X}) \longrightarrow \Im(\mathfrak{Y})$ *the* jet bundle *of* $(\mathfrak{X} \longrightarrow \mathfrak{Y})$*.*

DEFINITION 2.10.29 : *We say an object* $\mathfrak{X} \in \mathbf{H}$ *is* formally smooth *if the constant infinitesimal path inclusion* $\mathfrak{X} \longrightarrow \Im(\mathfrak{X})$ *is an effective epimorphism.*

*We say that* $f$ *is a* formally unramified morphism *if this is a (-1)-truncated morphism. More generally,* $f$ *is an order-*$k$ formally unramified morphisms *for* $(2) \leqslant k \leqslant \infty$ *if this is a* $k$*-truncated morphism.*

*A morphism* $f : \mathfrak{X} \longrightarrow \mathfrak{Y}$ *in* $\mathbf{H}_{\mathfrak{th}}$ *is* formally étale *if it is* $\Im$*-closed, hence if its* $\Im$*-unit naturality square:*

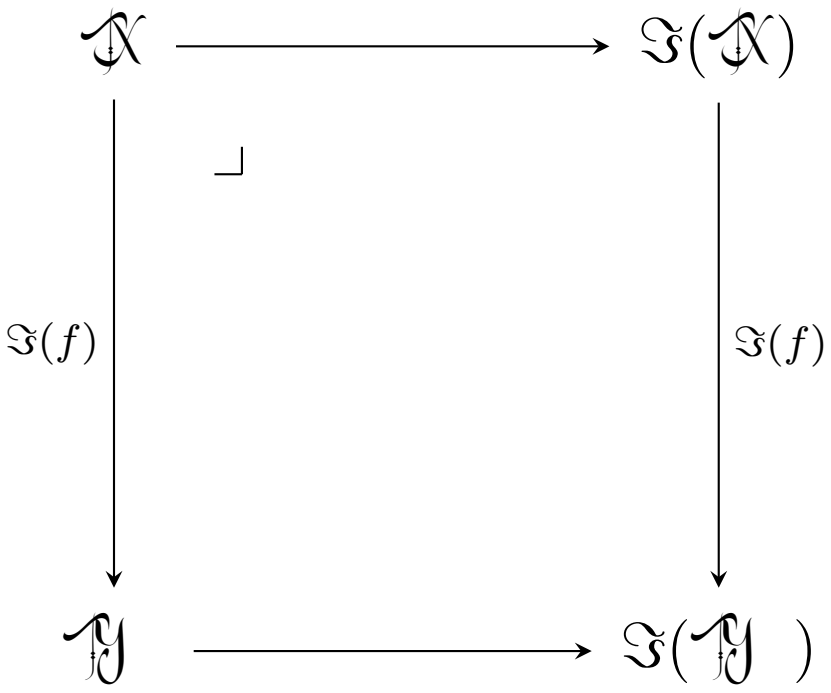

*is an* $(\infty,1)$*-pullback.*

The properties of formally etale morphisms (see Rem. 5.3.31 in Schreiber's Draft [129]) correspond to the axioms on the admissible maps in a geometry for structured $(\infty,1)$-topoï (see definition 2.3.18 or Def. 1.2.1 in [91]). This means that, in this cohesive context, $\mathbf{H}_{\mathfrak{th}}$is a $\mathcal{G}^{\mathrm{ad}}$-structure topos on $\mathbf{H}$where $\mathcal{G}^{\mathrm{ad}}$ is spanned by formally étale morphisms.

The thickening $\Im : \mathbf{H}_{\mathfrak{th}} \longrightarrow \mathbf{H}_{\mathfrak{th}}$ constructed from the monad of the adjunction $(\Re \dashv \Im)$ is directly linked to the formal infinitesimal $\Delta^n_k$ in synthetic differential geometry and crystalline cohomology by this $\mathcal{G}^{\mathrm{ad}}$-structure given by formally etale morphisms.

The nilcompleteness, the integrability ans the infinitesimally cohesiveness in the sense of J. Lurie's in Lurie–Artin's representation theorem is therefore automatically satisfied.

The complex cotangent $\mathbb{L}_{\mathfrak{H}_{\mathfrak{th}}/\mathfrak{H}}$ is 0 since the immersion is étale.

Let us remark that $\mathbb{L}_{\mathfrak{Y}/\mathfrak{X}} = 0$ if $\mathfrak{X} \longrightarrow \mathfrak{Y}$ étale, but $\mathfrak{H}_{\mathfrak{th}} \longrightarrow \mathfrak{H}$ in the reverse direction is not étale.

That is: $\mathfrak{X}_{\mathrm{dRh}} \hookrightarrow \mathfrak{X}$ is étale as an immersion but $\mathfrak{X} \to \mathfrak{X}_{\mathrm{dRh}}$ is not étale, its a projection/reduction.

The crystalline cohomologie constructed for non-smooth underlying varieties, in particular for the Grassmannian (which is smooth), for the positroid varieties (which are smooth taken one by one), but also their gluing along closed immersions which we already constructed and which is not smooth, can then be applied.

We then get the following conclusion:

THEOREM 2.10.30 : *Considering the differential cohesive structured spectral Deligne-Mumford $\infty$-stack:*

$$\mathfrak{H}_{\mathfrak{th}} := \mathrm{Gr}(k, n)_{\mathrm{th}}$$

*is an $\mathcal{IC}^{\bullet}$-structured spectral algebraic space which represents the initial formal moduli problem definition 2.11.4.Its crystalline cohomology is:*

$$H^*_{\mathrm{dRh}}(\mathfrak{H}_{\mathfrak{th}}) := H^*_{\mathrm{dRh}}(\mathrm{Gr}(k, n)_{\mathrm{th}}, \mathbf{1}_{\otimes})$$

*The symetric monoïdality of $\mathsf{QCoh}(\mathfrak{H}_{\mathfrak{th}})$ and of $\mathsf{QStk}(\mathfrak{H}_{\mathfrak{th}})$ of (compact perfect) dualizable generator $\mathscr{O}_{\mathfrak{H}_{\mathfrak{th}}}$ allows us to apply Tannaka duality to recover the inital algebraic space, which then has a volume given by the exterior power of maximal rank of the De Rham derivative (of rank 1), that is the volume form on a differential compact space, that is:*

$$\mathrm{Volume}(\mathfrak{H}_{\mathfrak{th}}) = \int \mathrm{d}\Omega^n(\mathfrak{H}_{\mathfrak{th}}) = \mathrm{Volume}(\mathfrak{H})$$

*since the cohomology is insensible to the reduction.*

REMARK 2.10.31 : The $\mathcal{IC}^{\bullet}$-structured spectral algebraic space $\mathfrak{H}_{\mathfrak{th}}$ :=$\mathrm{Gr}(k,n)_{\mathrm{th}}$ is autodual for Tannaka reconstruction, as we will see it in the second part concerning the non-autodual Amplituhedron.

## 2.11.Conclusion of Part I: The Grassmannian formal moduli problem.

Here is the definition of a Lurie's formal moduli problem:

DEFINITION 2.11.1 : *We say that a $\mathbb{E}_\infty$-algebra over the ring $\Bbbk$ is augmented if it is an $\mathbb{E}_\infty$-algebra over $\Bbbk$ viewed as a $\mathbb{E}_\infty$-ring. Augmented algebras are denoted $\mathcal{CA}\mathbf{lg}_{\Bbbk}^{aug}$.*

DEFINITION 2.11.2 : *If $A \in \mathcal{CA}\mathbf{lg}_{\Bbbk}^{aug}$ is connective, $\pi_* A$ is a finite dimensional vector space over $\Bbbk$, and $\pi_0 A$ is a local ring, we say that $A$ is* Artinian. *Artinian algebras are denoted $\mathcal{CA}\mathbf{lg}_{\Bbbk}^{art}$.*

REMARK 2.11.3 : We may call these algebras *locally* artinian since it is a *derived* version of artineness.

DEFINITION 2.11.4 : *Let $\mathscr{F} : \mathcal{CA}\mathbf{lg}_{\Bbbk}^{\mathrm{art}} \longrightarrow \mathcal{S}$ be a functor. We say that $\mathscr{F}$ is a* formal moduli problem *if it satisfies the following conditions:*

- *The space $\mathscr{F}(\Bbbk)$ is contractible*

- *For every pullback square*

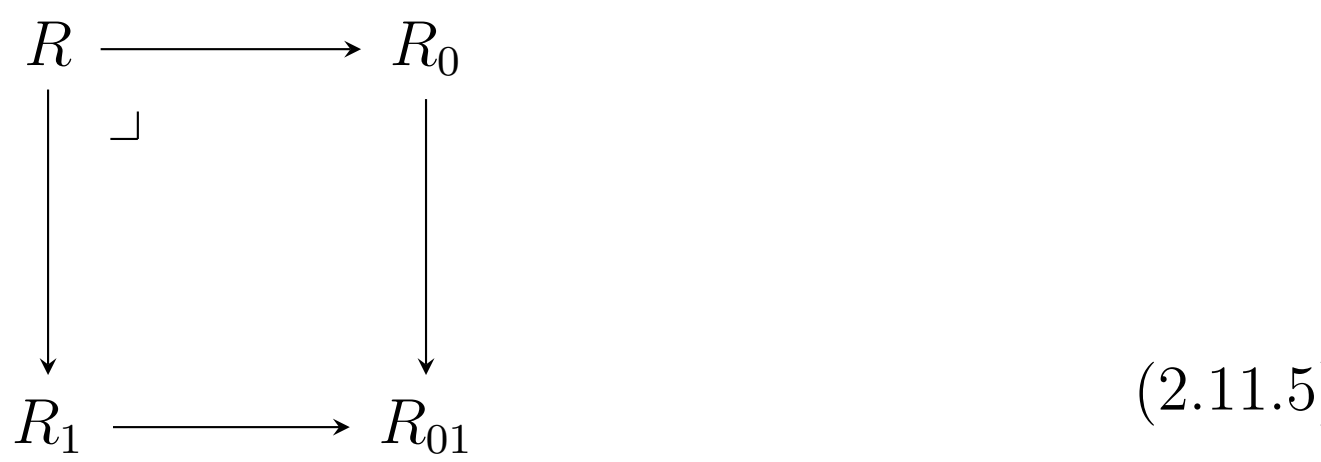

(2.11.5)

*in* $\mathcal{CA}\mathbf{lg}_{\Bbbk}^{\mathrm{art}}$, *if both the maps* $\pi_0 R_0 \longrightarrow \pi_0 R_{01} \longleftarrow \pi_0 R_1$ *are surjective, then the diagram of spaces*

$$\begin{array}{ccc} \mathscr{F}(R) & \longrightarrow & \mathscr{F}(R_0) \\ \downarrow & \lrcorner & \downarrow \\ \mathscr{F}(R_1) & \longrightarrow & \mathscr{F}(R_{01}) \end{array} \tag{2.11.6}$$

*is also a pullback square.*

We can interpret the formal moduli problem $\mathscr{F}$ in particular if the pullback diagrams Equation 2.11.5 are exact gluings along closed immersions over common intersections. This comes back to the existence of the diagram Equation 2.11.6 exactly constructed along the exact image gluings along the images of these closed immersions involved in the initial diagram Equation 2.11.5.

We have proved in the preceeding section that this formal moduli problem gives a representative for the functor $\mathscr{F}$, that is the $\mathcal{IC}^{\bullet}$-structured spectral algebraic space $\mathfrak{H}_{\mathfrak{th}}$ :=$\mathrm{Gr}(k,n)_{\mathrm{th}}$, which is differential cohesive, with a volume form well-defined by cristalline cohomology, with enough perfect perverse sheaves (see subsubsection 2.6.6.1) for the $\mathcal{IC}^{\bullet}$-structured spectral algebraic space.

Conclusion: These are the *brave new objects, the thickened positive Grassmannian and its thickened positroid varieties, which represent this differential cohesive spectral formal moduli problem as a new brave Schubert geometry. In othen words, this is the spectral De Rham stacks of these objects which are now of interest for us, equipped with their canonical sheaf of* $\mathbb{E}_\infty$*-rings, and their structure sheaf of* $\mathbb{E}_\infty$*-modules given by the perverse sheaves and their intersection complexes on them.*

A NEW POSITIVE SCHUBERT GEOMETRY

$$\mathfrak{I}(\maltese)_{\mathrm{dRh}} \ := \ \mathfrak{I}\big(\mathfrak{Gr}^{\geqslant 0}_{\text{ét}}, \mathscr{O}_{\mathfrak{Gr}^{\geqslant 0}_{\text{ét}}}, \mathcal{IC}^{\bullet}_{\mathfrak{Gr}^{\geqslant 0}_{\text{ét}}}\big)_{/\,\mathrm{dRh}}$$

To recall the whole construction in a few words, $\mathcal{G}^{\mathrm{perv}}_{\mathrm{perf}}$-structured $(\mathscr{O}_{\mathfrak{G}}, \mathscr{F}_{\mathfrak{G}})-$module geometry writen in a Luries's formal moduli problem $\mathscr{F}$ structures the objects which represent this formal moduli functor. This one is represented by a spectral stably $t$-structured $(\mathsf{QStk}(\maltese)^{\geqslant 0}, \mathsf{QStk}(\maltese)^{\leqslant 0})$ quasi-coherent Deligne-Mumford stack $\mathsf{QStk}(\maltese)$ in $\mathsf{Groth}_\infty \subseteq \mathsf{QCoh}^{\mathrm{Pst}}_{\mathrm{comp}}(\maltese)$. This spectral algebraic space is infinitesimally thickened cohesive in Grothendieck modal adjunction contexts $\mathfrak{R}(\maltese) \dashv \mathfrak{I}(\maltese) \dashv \&(\maltese)$. Tannaka 0-affineness and 1-affineness of $\mathcal{BCFW}_{n,4}$ positroid $\mathring{\Pi}_f$ (which we will see soon their perverslely twisted $H^2_{\text{ét}}(\maltese, \mathcal{IC}^{\bullet\times} \mathit{cohomology})$ allow us to recover modal infinitesimal shapes of differential cohesive De Rham stacks $\mathfrak{I}(\maltese)_{\mathrm{dRh}}$, that we describe now the new brave objects of Schubert geometry.

This description is a bit long, but it reflects the whole construction, even if a few details will be explained below about twisted cohomology. Let us now see how this apparent complexity can bring us some non-trivial result : *the existence and the volume of the Dual Amplituhedron.*

# Chapter 3

# PART II: A non-trivial application : The Dual Amplituhedron

In our case of the Amplituhedron $\mathscr{A}_{n,k,4}$, which is a tiling made of exact étale gluings:

$$\mathscr{A}_{n,k,4} = \coprod_{f \in \mathscr{B}(k,n)} Z_{\overset{\circ}{\Pi}_f} := \coprod_{f \in \mathscr{B}(k,n)} \widetilde{Z}\,(\overset{\circ}{\Pi}_f)$$

of $\mathcal{BCFW}_{n,4}$ cells $Z_{\overset{\circ}{\Pi}_f}$, images themselves of $\mathsf{BCFW}_{n,4}$ (closed) positroid varieties $\overset{\circ}{\Pi}_f$ of the TNN Grassmannian $\mathrm{Gr}_{k,n}^{\geqslant 0}$ glued along smaller positroid varieties $\overset{\circ}{\Pi}_g$ of this TNN Grassmannian $\mathrm{Gr}_{k,n}^{\geqslant 0}$ in the étale topology, the existence of our formal moduli problem $\mathscr{F}$ means here exactly that the image made of the $\mathcal{BCFW}_{n,4}$ image cells of the Amplituhedron $\mathscr{A}_{n,k,4}$ respects exactly the image of the Schubert stratification in the initial $\mathsf{BCFW}$ (closed) positroid varieties.

☆ ☆ ☆

As the Amplituhedron $\mathscr{A}_{n,k,4}$, as a spectral algebraic space, represents the formal moduli problem $\mathscr{F}$, finally, the enhancement of the amplituhedron $\mathscr{A}_{n,k,4}$ as an image object by $\widetilde{Z}$ of the TNN Grassmannian $\mathrm{Gr}_{k,n}^{\geqslant 0}$ into a functor:

$$\mathfrak{Z} : \mathfrak{G} := \left(\mathfrak{Gr}^{\geqslant 0}_{\text{ét}}, \mathscr{O}_{\mathfrak{Gr}^{\geqslant 0}_{\text{ét}}}, \mathcal{IC}_{\mathfrak{Gr}^{\geqslant 0}_{\text{ét}}}\right) \longrightarrow \left(\mathcal{A}^{\text{ét}}_{n,k,4}, \mathscr{O}_{\mathcal{A}^{\text{ét}}_{n,k,4}}, \mathcal{IC}_{\mathcal{A}^{\text{ét}}_{n,k,4}}\right) =: \mathfrak{A} \tag{3.0.1}$$

between spectral Deligne-Mumford stacks (in this case between spectral algebraic spaces) permits to boil down the construction of a dual to the Amplituhedron, that is a dual to the object Amplituhedron as a spectral algebraic space, to demonstrate that the dual formal moduli problem $\mathscr{F}^{\star}$ (if exists and is represented) is indeed represented by what we will then be allowed to call the *dual Amplituhedron* $\mathscr{A}^{\star}_{n,k,4}$.

We will then have to prove, as a last step, that this dual Amplituhedron $\mathscr{A}^{\star}_{n,k,4}$ is a modal 1-localic spatial compact infinitesimally thickened differential cohesive corporeal De Rham space, and therefore has a volume which can be calculated, in order to allow scattering amplitudes computations.

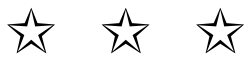

Having to indicate the list of properties that we will need to prove, the references we will use in the literature, and the notation we will employ, we can observe in the following detailed table of contents the increasing complexity of the concepts involved, due to the progressive categorification of objects and phenomena, and not simply to the transition from Grothendieck's (1-)topos to $\infty$-topos.

Indeed, the problem of finding a dual to something can already be observed in the problem of finding a dual to a vector space, say of finite dimension if we take this simplest example (here, the Amplituedron is compact). Indeed, let us take a finite dimensional $\Bbbk$-vector space $V$, its dual is then easy to define, and there is up to isomorphism only one object that fits: $V^{\vee} \simeq \mathrm{Hom}(V, \Bbbk)$. But, speaking categorically, this is the monoidal dual of the object $V$ in the symmetric monoidal category of finite-dimensional vector spaces $\mathcal{V}\mathrm{ect}^{fin}$, equipped with the well-known tensor product $\otimes$ and the unit element $\mathbf{1}_{\Bbbk}$, that is the base field $\Bbbk$. In fact, we know of no other dual than the monoidal dual in a monoidal supercategory (see Milne) with or without a unit element. We will do nothing else here. The objects are simply more complicated to put into a manipulable form, namely dualizable categories of modules and categories of quasi-coherent sheaves on spectral algebraic stacks, since the objects involved are highly not smooth but rather particularly singular.

## 3.1. The Amplituhedron coaffine spectral Deligne-Mumford stack.

### 3.1.1 The Amplituhedron map.

DEFINITION 3.1.1 : *If we fix $Z$ as a $(k+m)\times n$ real matrix in $\mathbb{R}^{(k+m)\times n}$ with all maximal minors are positive, where $k+m \leqslant n$ (and especially here $m=4$ for our purpose), then the matrix $Z$ induces a map:*

$$\widetilde{Z} : \mathrm{Gr}_{k,n}^{\geqslant 0} \longrightarrow \mathrm{Gr}_{k,k+m}$$

*defined by*

$$\widetilde{Z}(\langle v_1, \dots, v_k\rangle) := \langle Z(v_1), \dots Z(v_k)\rangle$$

*where $\langle v_1, \dots, v_k\rangle$ represents an element of $\mathrm{Gr}_{k,n}^{\geqslant 0}$ with $k$ basis vectors. We then define the* (tree) Amplituhedron $\mathcal{A}_{n,k,m}(Z)$ *as the image* $\widetilde{Z}(\mathrm{Gr}_{k,n}^{\geqslant 0})$ *of* $\mathrm{Gr}_{k,n}^{\geqslant 0}$ *inside* $\mathrm{Gr}_{k,k+m}$.

REMARK 3.1.2 : The case of immediate relevance to physics is $m = 4$. In this case, scattering amplitudes can be expressed as an integral over the dual Amplituhedron. More precisely, the initial conjecture of N. Arkani-Hamed and J. Trnka for this computation ([13, 12, 10, 14]) was that a certain collection of $4k$-dimensional cells of $\mathrm{Gr}_{k,k+m}^{\geqslant 0}$ would provide a "*triangulation*" of the Amplituhedron $\mathcal{A}_{n,k,m}(Z)$ which would allow the computation. This conjecture was realized thanks to the so-called "*standard BCFW recurrence*" ([37]) which provided a way to compute scattering amplitudes. In particular, the internal open positroid varieties and their walls in $\mathrm{Gr}_{k,k+m}^{\geqslant 0}$ which were "*selected*" to obtain a good tiling of the Amplituhedron were called standard BCFW cells, denoted $\mathsf{BCFW}_{n,k}$, and are ruled by the following definitions and theorems.
The first definition and theorem deal with triangulation and BCFW cells, namely injectivity, separation, surjectivity, and internal walls of standard BCFW cells involved in the triangulation of the standard BCFW tiling of the Amplituhedron. The second theorem ([39]) will be used to extend the standard BCFW recurrence to general

Amplituhedron tiling. The third one ([38]) will give the "*canonical form*" for the Amplituhedron, which was the expected object coming from the work of T. Lam ([9]) to bring the subject to a satisfactory conclusion in order to manage the calculation of scattering amplitudes.

Nevertheless, as the work of N. Arkani-Hamed, J. Trnka, T. Lam and others relies on the positivity of geometries, it is still expected that this calculation is one of the volume of an object defined in the real world corresponding to the "*dual*" of this famous *Amplituhedron*. And this definition of a "*dual Amplituhedron*" continued to resist combinatorial and algebro-geometric assaults. Let me present here the points of the solutions which are already obtained, and then, in the rest of this thesis, an attempt at an algebro-geometric solution to the existence and the volume of this dual Amplituhedron.

DEFINITION 3.1.3 ([37], Def 1.1) : *We call a* triangulation *of the Amplituhedron $\mathcal{A}_{n,k,m}(Z)$ a collection $\{S_{\overset{\circ}{\Pi}_f}, f \in \mathscr{B} \subset \mathscr{B}(k,n)\}$ of km-dimensional open positroid cells of $\mathrm{Gr}_{k,n}^{\geqslant 0}$, such that for every $Z \in \mathrm{Mat}_{n\times(k+m)}^{>0}$ :*

- Injectivity*: $S_{\overset{\circ}{\Pi}_f} \mapsto \widetilde{Z}.S_{\overset{\circ}{\Pi}_f}$ is an injective map for every $f \in \mathscr{B}(k,n)$*
- Separation *: $\widetilde{Z}.S_{\overset{\circ}{\Pi}_f}$ and $\widetilde{Z}.S_{\overset{\circ}{\Pi}_{f'}}$ are disjoint for every $\overset{\circ}{\Pi}_f \neq \overset{\circ}{\Pi}_{f'}$*
- Surjectivity *: $\bigcup_{f\in\mathscr{B}(k,n)} \widetilde{Z}.S_{\overset{\circ}{\Pi}_f}$ is an open dense subset of $\mathcal{A}_{n,k,m}(Z)$*

DEFINITION 3.1.4 ([37]) : *The collection $\{S_{\overset{\circ}{\Pi}_f}, f \in \mathscr{B}(k,n)\}$ corresponding to a triangulation is named* BCFW cells *and denoted* $\mathsf{BCFW}_{n,k}$*. The cardinality of this collection is the Narayana number: $\frac{1}{3}.\binom{n-3}{k+1}.\binom{n-3}{k}$ cells.*

THEOREM 3.1.5 ([37], Thm 1.3, 1.4 & 1.5) : *There exists a* standard *triangulation and* standard BCFW cells*, denoted* $\mathsf{BCFW}_{n,4}$*, such that the application $\widetilde{Z}$ is a smooth submersion and:*

- *The map $S_{\overset{\circ}{\Pi}_f} \mapsto \widetilde{Z}.S_{\overset{\circ}{\Pi}_f}$ is injective for every cell $S_{\overset{\circ}{\Pi}_f} \in \mathsf{BCFW}_{n,4}$*
- *The images $\widetilde{Z}.S_{\overset{\circ}{\Pi}_f}$ and $\widetilde{Z}.S_{\overset{\circ}{\Pi}_{f'}}$ are disjoint for every two different cells $S_{\overset{\circ}{\Pi}_f}$, $S_{\overset{\circ}{\Pi}_{f'}}$ in* $\mathsf{BCFW}_{n,4}$

- *The union of* $\widetilde{Z}.S_{\Pi_f}^{\circ}$ *over* $f \in \mathsf{BCFW}_{n,4}$ *is an open dense subset of* $\mathcal{A}_{n,k,4}(Z)$.

THEOREM 3.1.6 ([37],Proof of Thm. 1.5, Prop. 5.5, Prop. 8.1) : *Since the application* $\widetilde{Z}$ *is a smooth submersion, we also have, since* $\mathrm{Gr}_{4,n}^{\geqslant 0}$ *is compact and connected, that the amplituhedron is also compact and connected. Since* $\mathrm{Gr}_{4,n}^{>0}$ *is an open submanifold of* $\mathrm{Gr}_{4,n}$ *whose closure is* $\mathrm{Gr}_{4,n}^{\geqslant 0}$*, also* $\widetilde{Z}.(\mathrm{Gr}_{4,n}^{>0})$ *is an open submanifold of* $\mathrm{Gr}_{4,k+4}$ *whose closure is* $\mathcal{A}_{n,k,4}(Z)$.

REMARK 3.1.7 : Following [37],Proof of Thm. 1.5, Prop. 5.5, and Prop. 8.1, we can say the following facts:

- The Amplituhedron map is injective on standard BCFW cells in $\mathsf{BCFW}_{n,4}$.
- By [37] Prop. 8.1, it is also injective for codimension one boundary strata which do not map to the boundary. Call such codimension one strata internal boundaries.
- Moreover, the inverse procedure, as described in [37] Prop. 5.5 and Prop. 8.1, shows that on the image of the standard BCFW cells in $\mathsf{BCFW}_{n,4}$ and the internal boundaries, the inverse map is also continuous.
- Since continuous maps map a compact set to a compact set, the amplituhedron map maps the closure of any positroid cell, which is compact, to the closure of the image of the cell.
- Since the amplituhedron map is smooth, positroid cells of dimension $d \leqslant 4k$ map to subspaces of dimension at most $d$.
- Thus, for each BCFW cell the amplituhedron map is a homeomorphism of manifolds with boundary from the union of $S_{\Pi_f}^{\circ}$ and its internal boundaries to their image.
- The other boundaries either map to the boundary of the amplituhedron or into subspaces of codimension at least 2.

- The image of any codimension one stratum of a standard BCFW cell which does not map to the boundary, is also the image of a codimension one boundary of another standard BCFW cell.

- The amplituhedron map is injective on these internal boundaries.

Following L. Williams's work, and that of her collaborators, we can extend the notion of "*standard*" BCFW cells to "*general*" BCFW cells. Indeed, standard BCFW cells are a special case of (general) BCFW cells $S_{\mathfrak{r}}$, which are certain cells in $\mathrm{Gr}_{4,n}^{\geqslant 0}$ of dimension $4k$. We can construct general BCFW cells recursively, and obtain the following:

THEOREM 3.1.8 ([39] Thm. General BCFW cells give tiles) : *The amplituhedron map is injective on each general BCFW cell. That is, the closure* $\widetilde{Z}_{\mathfrak{r}} := \widetilde{Z}(S_{\mathfrak{r}})$ *of the image of a general BCFW cell* $S_{\mathfrak{r}}$ *is a tile, that is the closure of the image of a* $4k$*-dimensional cell of* $\mathrm{Gr}_{4,n}^{\geqslant 0}$ *on which* $\widetilde{Z}$ *is injective, which we refer to as a* general BCFW tile*, which form a tiling collection of* $\mathcal{A}_{n,k,4}(Z)$ *denoted* $\mathcal{BCFW}_{n,4}$.

The key point for the definition of the dual Amplituhedron is to refer to T. Lam's work ([9, 71, 70]), and by removing all the graphical elements corresponding to physics, which are not relevant here, despite their real importance in the final calculation of scattering amplitudes in physics, let us write with T. Lam [9] Conjecture 35 in Section 7.4.3, then Section 7.4.4 Formula 7.208 and Section 7.4.8 Formula 7.277 (see also [41] Section 2 Formula 2.1 and Formula 3.22 and the discussion which follows about the dual Grassmannian and the exterior poles of the integrand):

$$\Omega_{\mathscr{A}_{n,k,4}}(\underline{x}) \;=\; \mathrm{Vol}\,(\mathscr{A}_{n,k,4}^{\star})\,\mathrm{d}^{k,k+4}\,\underline{x}$$

where:

- $\mathscr{A}_{n,k,4}$ denotes the Amplituhedron as physico-mathematical object which correspond exactly to the purely mathematical object $\mathcal{A}_{n,k,4}(Z)$, and $\mathscr{A}_{n,k,4}^{\star}$ denotes the dual Amplituhedron

- $\mathrm{Vol}(-)$ denotes the volume of a real mathematical object in $\mathrm{Gr}_{k,k+4}(\mathbb{R})$
- "$\mathrm{d}^{k,k+4}$" denotes the classical $(k, k+4)$-dimensional De Rham volume differential
- $\underline{x}$ denotes a mute variable inside the Amplituhedron $\mathscr{A}_{n,k,4}$ bordered on its boundaries by logarithmic singularities of the meromorphic form $\Omega_{\mathscr{A}_{n,k,4}}$.

REMARK 3.1.9 : The role of the exteriority of the logarithmic poles with respect to the Amplituhedron is essential for the correct functioning of a calculation based on triangulation.

The point of this article is then to construct and determine what is:

$$\text{the } \textsf{dual Amplituhedron:} \qquad \mathscr{A}_{n,k,4}^{\star}$$

in order to determine:

$$\text{the } \textsf{volume of the dual Amplituhedron:} \qquad \mathrm{Vol}(\mathscr{A}_{n,k,4}^{\star})$$

which is our real purpose, instead of the dual Amplituhedron itself.
Indeed, in order to find the canonical form of the amplituhedron, one method is to tile $\mathcal{A}_{n,k,4}(Z)$ and sum over the canonical forms of the tiles.

DEFINITION 3.1.10 ([38] Def. 7.10 (Candidate canonical form of a tile)) : *Let $Z_{\overset{\circ}{\Pi}_f}$ be a tile of $\mathcal{A}_{n,k,4}(Z)$. As the amplituhedron map $\widetilde{Z}$ is injective on the open tile $Z_{\overset{\circ}{\Pi}_f}$, we can define the inverse $Z^{-1}_{\overset{\circ}{\Pi}_f}(\overset{\circ}{\Pi}_f) \longrightarrow \overset{\circ}{\Pi}_f$. Then, let us consider the pullback of the canonical form of the positroid cell under $\widetilde{Z}^{-1}$:*

$$\widetilde{\Omega}\,(Z_{\overset{\circ}{\Pi}_f}) := (\widetilde{Z}^{-1})^{*}\,\Omega\Big(\overline{\overset{\circ}{\Pi}_f},\ \overline{\overset{\circ}{\Pi}_f}^{\geqslant 0}\Big) \tag{3.1.11}$$

DEFINITION 3.1.12 (General $\mathsf{BCFW}_{n,4}$ cells) : *We call* $\textsf{general BCFW cells}$*, and denote also $\mathsf{BCFW}_{n,4}$, without conflict with (standard) BCFW cells, the $4k$-dimensional cells in $\mathrm{Gr}_{4,n}^{\geqslant 0}$, which constitue the inverse image by $Z_{\overset{\circ}{\Pi}_f}$.*

DEFINITION 3.1.13 ($\mathcal{BCFW}_{n,4}$ tiling of the Amplituhedron) : *We call* general BCFW tile*, or* tiling of $\mathcal{A}_{n,k,4}(Z)$*, and denote* $\mathcal{BCFW}_{n,4}$ *with respect to the same name for the same object of theorem 3.1.8, the direct image by* $Z$ *of* $\mathsf{BCFW}_{n,4}$*, the general BCFW cells chosen to constitute the (general) tiling.*

THEOREM 3.1.14 ([38], Conj. 7.11 & Thm 7.14 (Tiles are positive geometries)) : *Let* $\{Z_{\overset{\circ}{\Pi}_f},\ f \in \mathsf{BCFW}_{n,4}\}$ *be a tiling of* $\mathcal{A}_{n,k,4}(Z)$*. Then* $(\mathrm{Gr}_{k,k+4}(\mathbb{C}), Z_{\overset{\circ}{\Pi}_f})$ *is a positive geometry and its canonical form* $\Omega\big(\mathrm{Gr}_{k,k+4}(\mathbb{C}), Z_{\overset{\circ}{\Pi}_f}\big)$ *is the previous canonical form* $\widetilde{\Omega}\,(Z_{\overset{\circ}{\Pi}_f})$ *in Equation 3.1.11 of definition 3.1.10.*

THEOREM 3.1.15 ([38], Conj. 7.12 & Thm 7.14 (Canonical form for tilings)) : *Let* $\{Z_{\overset{\circ}{\Pi}_f},\ f \in \mathsf{BCFW}_{n,4}\}$ *be a tiling of* $\mathcal{A}_{n,k,4}(Z)$*. Then the canonical form of the amplituhedron* $\mathcal{A}_{n,k,4}(Z)$ *is obtained such that:*

$$\Omega\big(\mathrm{Gr}_{k,k+4}(\mathbb{C}), \mathcal{A}_{n,k,4}(Z)\big) = \sum_{f\in\mathsf{BCFW}_{n,4}} \Omega\big(\mathrm{Gr}_{k,k+4}(\mathbb{C}), Z_{\overset{\circ}{\Pi}_f}\big) \tag{3.1.16}$$

*In particular, the right-hand side of the equation is independent of the tiling.*

REMARK 3.1.17 : We notice that in order to compute the canonical form of the Amplituhedron, and nevertheless the Amplituhedron is tiled by the finite collection $\{Z_{\overset{\circ}{\Pi}_f},\ f \in \mathsf{BCFW}_{n,4}\}$, even if the finite collection $\mathsf{BCFW}_{n,4}$ of standard BCFW cells (or of general BCFW cells) recovers the TNN Grassmannian “*points by points*” since the entire TNN Grassmannian $\mathrm{Gr}^{\geqslant 0}_{k,k+4}$ is sent on $\mathcal{A}_{n,k,4}(Z)$, the finite collection $\mathsf{BCFW}_{n,4}$ does not cover the TNN Grassmannian *with respect to the positroid cells of the TNN Grassmannian.* A simple argument about the dimension $4k$ of the positroid cells in $\mathsf{BCFW}_{n,4}$ justifies this point. But the Amplituhedron map $\widetilde{Z}$ being injective on each open tile $Z_{\overset{\circ}{\Pi}_f}$ in the finite collection $\{Z_{\overset{\circ}{\Pi}_f},\ f \in \mathsf{BCFW}_{n,4}\}$, the inverse map $Z^{-1}_{\overset{\circ}{\Pi}_f}$ sends $\mathcal{BCFW}_{n,4}$ to $\mathsf{BCFW}_{n,4}$, cell by cell, that is from a tile of $\mathcal{BCFW}_{n,4}$ in $\mathcal{A}_{n,k,4}(Z)$ to a positroid cell $\overset{\circ}{\Pi}_f$ in $\mathsf{BCFW}_{n,4}$ and we have a bijection:

$$\widetilde{Z} \;:\; \begin{cases} \mathsf{BCFW}_{n,4} & \overset{\sim}{\longrightarrow} \quad \mathcal{BCFW}_{n,4} \\ \quad \overset{\circ}{\Pi}_f & \longmapsto \quad Z_{\overset{\circ}{\Pi}_f} \end{cases}$$

with the constraints given in remark 3.1.7, which can be summarized in that the Amplituhedron is well-tiled by the images of the standard BCFW cells, in a suitable way for the boundaries of the cells, which are also positroid cells, but of decreasing dimension. Characterizations of this map $\widetilde{Z}$ and a categorification of it will follow below. It is the core of the application of spectral algebraic geometry to the Amplituhedron map.

### 3.1.2 The Amplituhedron functor.

The construction of the $\mathcal{BCFW}_{n,4}$-étale $\infty$-topology for the Amplituhedron is made exactly the same way as that of the positroid-étale $\infty$-topology for the TNN Grassmannian. This leads to the definition of a categorification of the Amplituhedron map $Z$ between the corresponding $\infty$-topoi, and finally to the categorified Amplituhedron functor between the categories of spectral Deligne-Mumford stacks on these $\infty$-topoi.

The construction of the $\mathcal{BCFW}_{n,4}$-étale $\infty$-topology is analogous with the previous section. We do not write it again, but we make a number of comments below about the $\infty$-sheaf condition.

Let us now categorize the Amplituhedron map Z. Let $\mathcal{C}$ be an $\infty$-category, and let $Z : \mathrm{Gr}_{k,n}^{\geqslant 0} \longrightarrow \mathcal{A}_{n,k,4}$ be the Amplituhedron map, which is a continuous map of topological spaces, if we view $\mathrm{Gr}_{k,n}^{\geqslant 0}$ and $\mathcal{A}_{n,k,4}$ as topological spaces. Let $\mathscr{F}: \mathcal{U}^{\mathrm{op}} \longrightarrow \mathcal{C}$ be a $\mathcal{C}$-valued presheaf on $\mathrm{Gr}_{k,n}^{\geqslant 0}$, where $\mathcal{U}^{\mathrm{op}}$ is the opposite category of open sets in $\mathrm{Gr}_{k,n}^{\geqslant 0}$. Then we can define another $\mathcal{C}$-valued presheaf $(Z_*\mathscr{F}) : \mathcal{U}^{\mathrm{op}} \longrightarrow \mathcal{C}$ , given on objects by the formula :

$$(Z_*\mathscr{F})(U) = \mathscr{F}(Z^{-1}U)$$

If $\mathscr{F}$ is a $\mathcal{C}$-valued sheaf on $\mathrm{Gr}_{k,n}^{\geqslant 0}$, then $Z_*\mathscr{F}$ is a $\mathcal{C}$-valued sheaf on $\mathcal{A}_{n,k,4}$. Moreover, the construction $\mathscr{F} \longrightarrow Z_*\mathscr{F}$ determines a functor $\mathscr{F}_* : \mathcal{S}\mathrm{hv}_{\mathcal{C}}(\mathrm{Gr}_{k,n}^{\geqslant 0}) \longrightarrow \mathcal{S}\mathrm{hv}_{\mathcal{C}}(\mathcal{A}_{n,k,4})$, given by the precomposition with the map of partially ordered sets $Z^{-1} : \mathcal{U}(\mathcal{A}_{n,k,4}) \longrightarrow \mathcal{U}(\mathrm{Gr}_{k,n}^{\geqslant 0})$.

We have defined a map, denoted $\mathscr{Z}_*$, which is a pushforward functor associated to the

Amplituhedron map $Z$ :

$$\mathscr{Z}_* \ : \mathcal{S}\mathrm{hv}_{\mathcal{C}}(\mathrm{Gr}_{k,n}^{\geqslant 0})^{\mathrm{op}} \longrightarrow \mathcal{S}\mathrm{hv}_{\mathcal{C}}(\mathcal{A}_{n,k,4})^{\mathrm{op}}$$

We can regard the construction $\mathrm{Gr}_{k,n}^{\geqslant 0} \longmapsto \mathcal{S}\mathrm{hv}_{\mathcal{C}}(\mathrm{Gr}_{k,n}^{\geqslant 0})^{\mathrm{op}}$ as a functor from the category $\mathcal{S}$ of topological spaces to the category of simplicial sets, which carries a continuous map $Z : \mathrm{Gr}_{k,n}^{\geqslant 0} \longrightarrow \mathcal{A}_{n,k,4}$ to the pushforward functor $\mathscr{Z}_*$.

Let us define a left exact left adjoint pullback functor $\mathscr{Z}^*$ to this pusforward functor:

$$\mathscr{Z}^* \ : \mathcal{S}\mathrm{hv}_{\mathcal{C}}(\mathcal{A}_{n,k,4})^{\mathrm{op}} \longrightarrow \mathcal{S}\mathrm{hv}_{\mathcal{C}}(\mathrm{Gr}_{k,n}^{\geqslant 0})^{\mathrm{op}}$$

If we take $\mathcal{C} = \mathcal{S}$, the $\infty$-category of spaces, and when restricting $\mathscr{Z}_*$ to étale sheaves on $\mathrm{Gr}_{k,n}^{\geqslant 0}$ satisfying the étale topology $\infty$-sheaf condition, we can choose a sheaf $\mathscr{F}$ in what is exactly defined in definition 2.7.2, that is: $\mathfrak{Gr}_{\text{ét}}^{\geqslant 0}$.

Let us define the category $\mathscr{A}_{n,k,4}^{\text{ét}}$ as the subcategory of $\mathcal{S}\mathrm{hv}_{\mathcal{S}}(\mathcal{A}_{n,k,4}))$ which is in the image of $\mathscr{Z}_* \ : \mathcal{S}\mathrm{hv}_{\mathcal{C}}(\mathrm{Gr}_{k,n}^{\geqslant 0})^{\mathrm{op}}$.

LEMMA 3.1.18 : *The sheaves of the category $\mathscr{A}_{n,k,4}^{\text{ét}}$ satisfy the étale $\infty$-sheaf condition in the presheaves category $\mathcal{F}\mathrm{un}((\mathscr{G}^{\text{ét}})^{\mathrm{op}}, \mathcal{S})$. The category $\mathscr{A}_{n,k,4}^{\text{ét}}$ is then an $\infty$-topos.*

We then have that each sheaf $\mathscr{G} := \mathscr{Z}_*(\mathscr{F}) \in \mathcal{S}\mathrm{hv}_{\mathcal{S}}(\mathcal{A}_{n,k,4}) =: \mathscr{A}_{n,k,4}^{\text{ét}}$, that is in the image of $\mathfrak{Gr}_{\text{ét}}^{\geqslant 0}$ by $\mathscr{Z}_*$ in the $\infty$-category $\mathscr{Shv}(\mathscr{F}^{\text{ét}})$ of $\infty$-sheaves, satisfies the étale $\infty$-sheaf condition in the presheaves category $\mathcal{F}\mathrm{un}((\mathscr{G}^{\text{ét}})^{\mathrm{op}}, \mathcal{S})$.

The image of $\mathscr{Z}_*(\mathfrak{Gr}_{\text{ét}}^{\geqslant 0})$ is an $\infty$-topos, which, thanks to the $\infty$-functoriality of **??**, is now sent in an $\infty$-topos, which gives the following, remembering the Equation 3.1.11 and the pullback $(\widetilde{Z}^{-1})^*$ :

DEFINITION 3.1.19 : *Since the general BCFW cells $\mathsf{BCFW}_{n,4}$ and the general tiling collection $\mathcal{BCFW}_{n,4}$ are in smooth bijection for their interiors, theirs internal walls, their external boundaries, and the injectivity, separation and surjectivity of the homeomor-* *phism* $\widetilde{Z}$ *and* $\widetilde{Z}^{-1}$*, we can define the* $\mathcal{BCFW}_{n,4}$-étale ∞-topology for the Amplituhedron *in the same way as previously, and define and denote* $\mathcal{A}^{\text{ét}}_{n,k,4}$*, the* ∞-category of Deligne-Mumford stacks on the ∞-topos of the Amplituhedron $\mathcal{A}_{n,k,4}$ *, endowed with this topology, by the same construction of that of the category of Deligne-Mumford stacks on the* ∞*-topos of the TNN Grassmannian endowed with the positroid-étale* ∞*-topology.*

REMARK 3.1.20 : We can notice that the sheaves $\mathscr{G}$ in $\mathscr{Shv}^{\text{ét}}$ which satisfy the ∞-sheaf condition are exactly those such that there exists, for each $\mathscr{G} \in \mathscr{Shv}^{\text{ét}}$, a collection $\mathscr{A}$ of bounded permutations $f$ in $\mathcal{B}(k,n)$, and a surjective étale map $\coprod_{f \in \mathcal{A}} \mathring{\Pi}_f \longrightarrow (\widetilde{Z}^{-1})^*(\mathscr{G})$, which is then called an étale cover of the sheaf $\mathscr{G}$ in $\mathscr{Shv}^{\text{ét}}$ and we have:

$$\mathcal{A}^{\text{ét}}_{n,k,4} := \mathscr{Z}_* \left(\mathfrak{Gr}^{\geqslant 0}_{\text{ét}}\right)$$

Let us define the category $\mathcal{A}^{\text{ét}}_{n,k,4}$ as the subcategory of $\mathcal{S}\text{hv}_{\mathcal{S}}(\mathcal{A}_{n,k,4}))$ which is in the image of $\mathscr{Z}_* : \mathcal{S}\text{hv}_{\mathcal{C}}(\text{Gr}^{\geqslant 0}_{k,n})^{\text{op}}$.

LEMMA 3.1.21 : *The sheaves of the category* $\mathcal{A}^{\text{ét}}_{n,k,4}$ *satisfy the étale* ∞*-sheaf condition in the presheaves category* $\mathcal{F}\text{un}((\mathscr{G}^{\text{ét}})^{\text{op}}, \mathcal{S})$*. The category* $\mathcal{A}^{\text{ét}}_{n,k,4}$ *is then an* ∞*-topos.*

We then have constructed a geometric morphism between topoi:

DEFINITION 3.1.22 : *We define the* Amplituhedron geometric morphism $\mathscr{Z}$ *between the* ∞*-topos of the TNN Grassmannian* $\mathfrak{Gr}^{\geqslant 0}_{\text{ét}}$ *and the* ∞*-topos of the Amplituhedron* $\mathcal{A}^{\text{ét}}_{n,k,4}$ *such that:*

$$\mathscr{Z}_* \;:\; \mathfrak{Gr}^{\geqslant 0}_{\text{ét}} \underset{\longrightarrow}{\overset{\text{left ex. \& left adj.}}{\overset{\longleftarrow}{\perp}}} \mathcal{A}^{\text{ét}}_{n,k,4} \;:\; \mathscr{Z}^* \tag{3.1.23}$$

LEMMA 3.1.24 : *The pullback functor* $\mathscr{Z}^*$ *is a left exact left adjoint functor to the pushforward functor* $\mathscr{Z}_*$*.*

THEOREM 3.1.25 (Consistency of the Amplituhedron functor) : *Provided that the simple checks confirm that the morphism is afffine, proper and semi-small, that it preserves Grothendieck contexts and the perfect perverse geometries that arre enough (see subsubsection 2.6.6.1), we define the* ∞*-geometric morphism, named* Amplituhedron functor, *denoted* $\mathfrak{Z}$*, between the spectral Deligne-Mumford stacks of the TNN Grassmannian* $\mathfrak{G}$ *and the spectral Deligne-Mumford stacks of the Amplituhedron, which are here both spectral algebraic spaces:*

$$\mathfrak{Z} : \mathfrak{G} := \left(\mathfrak{Gr}^{\geqslant 0}_{\text{ét}}, \mathscr{O}_{\mathfrak{Gr}^{\geqslant 0}_{\text{ét}}}, \mathcal{IC}_{\mathfrak{Gr}^{\geqslant 0}_{\text{ét}}}\right) \longrightarrow \left(\mathcal{A}^{\text{ét}}_{n,k,4}, \mathscr{O}_{\mathcal{A}^{\text{ét}}_{n,k,4}}, \mathcal{IC}_{\mathcal{A}^{\text{ét}}_{n,k,4}}\right) \tag{3.1.26}$$

DEFINITION 3.1.27 : *We define, and denote* $\mathfrak{A}$ *, the* Amplituhedron spectral algebraic space*:*

$$\mathfrak{A} := \left(\mathcal{A}^{\text{ét}}_{n,k,4}, \mathscr{O}_{\mathcal{A}^{\text{ét}}_{n,k,4}}, \mathcal{IC}_{\mathcal{A}^{\text{ét}}_{n,k,4}}\right) \tag{3.1.28}$$

## 3.2.The Amplituhedron functor between spectral Deligne-Mumford stacks.

### 3.2.1 The étale geometry on the Amplituhedron spectral DM stack.

REMARK 3.2.1 : Following [37],Proof of Thm. 1.5, Prop. 5.5, and Prop. 8.1, we can say the following facts:

- The Amplituhedron map is injective on standard BCFW cells in $\mathsf{BCFW}_{n,4}$.
- By [37] Prop. 8.1, it is also injective for codimension one boundary strata which do not map to the boundary. Call such codimension one strata internal boundaries.
- Moreover, the inverse procedure, as described in [37] Prop. 5.5 and Prop. 8.1, shows that on the image of the standard BCFW cells in $\mathsf{BCFW}_{n,4}$ and the internal boundaries, the inverse map is also continuous.
- Since continuous maps map a compact set to a compact set, the amplituhedron map maps the closure of any positroid cell, which is compact, to the closure of the image of the cell.
- Since the amplituhedron map is smooth, positroid cells of dimension $d \leqslant 4k$ map to subspaces of dimension at most $d$.

- Thus, for each BCFW cell the amplituhedron map is a homeomorphism of manifolds with boundary from the union of $S_{\Pi_f}^{\circ}$ and its internal boundaries to their image.
- The other boundaries either map to the boundary of the amplituhedron or into subspaces of codimension at least 2.
- The image of any codimension one stratum of a standard BCFW cell which does not map to the boundary, is also the image of a codimension one boundary of another standard BCFW cell.
- The amplituhedron map is injective on these internal boundaries.

### 3.2.2 The perfect perverse geometrie of the Amplituhedron.

To invite to the solution given in this thesis for the Dual Amplituhedron, let us state right now that the $\infty$-category of triples which are considered pursuing the global solution to the problem solved, is here the $\infty$-category of $(\mathfrak{X}, \mathscr{O}_{\mathfrak{X}}, \mathscr{F}_{\mathfrak{X}})$, where:

- $\mathfrak{X}$ are gluing of spectral Deligne-Mumford stacks along closed immersions (see [90]) of the images by $\mathscr{L}_*$ of 4k-dimensional closed positroid varieties, belonging to the BCFW tiling of the Amplituhedron (see [37, 38])
- $\mathscr{O}_{\mathfrak{X}}$ are the canonical sheaves of the underlying topological 0-truncated 1-localic 1-topoï of these spectral Deligne-Mumford stacks, which we view as a global algebraic space whose "*affine pieces*" we glue together are involved in as many Grothendieck contexts as there are tiles in the BCFW tiling
- $\mathscr{F}_{\mathfrak{X}}$ are perfect perverse stacks, that is certain sheaves of $\mathbb{E}_{\infty}$-ring spectra, on these spectral Deligne-Mumford stacks, and they are enough subsubsection 2.6.6.1.

The morphisms between these triples will be made of geometric morphisms of topoï (that is in particular with localization by left adjoints), verifying the following:

- the gluings of Deligne-Mumford stacks along closed immersions preserve scallop decomposition (see [103] and [90]) and structure of algebraic spaces (see [90]), applied here to the images by $\mathscr{Z}_*$ of 4k-dimensional closed positroid varieties belonging to the BCFW tiling of the Amplituhedron (see [38])

- the canonical sheaves are glued in the same manner for the underlying "*affine pieces*" of the considered algebraic spaces (see [90])

- the perfect perverse stacks behave well for gluing (see [103]) and [90]) and for natural transformations between algebraic spaces thanks to the decomposition theorem for derived functors, and they are enough subsubsection 2.6.6.1, since the Amplituhedron map is proper and semi-small (see [103]) and gives rise to an Amplituhedron functor at the level of spectral Deligne-Mumford stacks

## 3.3. The dual Amplituhedron coaffine spectral stack.

### 3.3.1 The dual of a spectral Deligne-Mumford stack.

Let us recall that we can identify a $\mathcal{M}\mathbf{od}$-valued sheaf on $\mathfrak{X}$ with a pair $(\mathscr{O}_{\mathfrak{X}}, \mathscr{F}_{\mathfrak{X}})$, where $\mathscr{O}_{\mathfrak{X}}$ is a sheaf of $\mathbb{E}_\infty$-rings on $\mathscr{O}_{\mathfrak{X}}$ and $\mathscr{F}_{\mathfrak{X}}$ is a $\mathscr{O}_{\mathfrak{X}}$-module. That is if we consider $(\mathfrak{X}, \mathscr{O}_{\mathfrak{X}}, \mathscr{F}_{\mathfrak{X}})$ a spectral Deligne-Mumford stack endowed with a module geometry represented by the attachment of a $\mathscr{O}_{\mathfrak{X}}$-module $\mathscr{F}_{\mathfrak{X}}$, we can take, for example, in the Equation 3.1.23:

$$\mathfrak{Z} : \mathfrak{G} := \left(\mathfrak{Gr}^{\geqslant 0}_{\text{ét}}, \mathscr{O}_{\mathfrak{Gr}^{\geqslant 0}_{\text{ét}}}, \mathcal{IC}_{\mathfrak{Gr}^{\geqslant 0}_{\text{ét}}}\right) \longrightarrow \left(\mathcal{A}^{\text{ét}}_{n,k,4}, \mathscr{O}_{\mathcal{A}^{\text{ét}}_{n,k,4}}, \mathcal{IC}_{\mathcal{A}^{\text{ét}}_{n,k,4}}\right) =: \mathcal{A}^{\text{ét}}_{n,k,4} \quad (3.3.1)$$

the following assignments:

- $\mathfrak{X} := \mathfrak{Gr}^{\geqslant 0}_{\text{ét}}$, $\mathscr{O}_{\mathfrak{X}} := \mathscr{O}_{\mathfrak{Gr}^{\geqslant 0}_{\text{ét}}}$, and $\mathscr{F}_{\mathfrak{X}} := \mathcal{IC}_{\mathfrak{Gr}^{\geqslant 0}_{\text{ét}}}$ for the TNN Grassmannian Deligne-Mumford stack

- $\mathfrak{X} := \mathcal{A}^{\text{ét}}_{n,k,4}$, $\mathscr{O}_{\mathfrak{X}} := \mathscr{O}_{\mathcal{A}^{\text{ét}}_{n,k,4}}$, and $\mathscr{F}_{\mathfrak{X}} := \mathcal{IC}_{\mathcal{A}^{\text{ét}}_{n,k,4}}$ for the Amplituhedron Deligne-Mumford stack

DEFINITION 3.3.2 : *We define the* dual spectral Deligne-Mumford stack $(\mathfrak{X}, \mathscr{O}_{\mathfrak{X}}, \mathscr{F}_{\mathfrak{X}})^{\star} := (\mathfrak{X}^{\star}, \mathscr{O}_{\mathfrak{X}}^{\star}, \mathscr{F}_{\mathfrak{X}}^{\star})$ *such that:*

- $\mathfrak{X}^{\star}$ *is constructed from definition 2.5.14 such that* $\mathsf{QCoh}(\mathfrak{X}^{\star})$ *is the monoidal dual of* $\mathsf{QCoh}(\mathfrak{X})$ *defined in definition 2.5.17, since the monoidal structure* $\mathsf{QStk}^{\mathrm{PSt}}(\mathfrak{X})$ *then allows us to take the monoidal dual* $\mathsf{QCoh}(\mathfrak{X})^{\vee}$ *of* $\mathsf{QCoh}(\mathfrak{X})$ *in* $\mathsf{QStk}^{\mathrm{PSt}}(\mathfrak{X})$
- $\mathscr{O}_{\mathfrak{X}}^{\star}$ *is the canonical sheaf of* $\mathfrak{X}^{\star}$ *constructed in definition 2.2.62*
- $\mathscr{F}_{\mathfrak{X}}^{\star}$ *is defined like in Verdier duality of definition 2.5.21*

$$\mathscr{F}_{\mathfrak{X}}^{\star} := \mathbb{R}_{\mathbb{E}_\infty}\Gamma(\bullet, \underline{\mathbb{Q}}_{\bullet})^{\star} := \mathcal{H}\mathrm{om}_{\mathcal{D}^{+}(\mathsf{Coh}(\bullet))}(\mathbb{R}\Gamma(\bullet, \mathscr{O}_{\bullet}), \omega_{\bullet})$$

*since following theorem 2.6.35 stating the intersection complex sheaf of* $\mathbb{E}_\infty$*-rings from the derived global sections functors,* $\mathscr{F}_{\mathfrak{X}} = \mathcal{IC}_{\mathfrak{X}}^{\bullet} : U \mapsto \mathbb{R}_{\mathbb{E}_\infty}\Gamma(U, \underline{\mathbb{Q}}_U)$*, where* $\underline{\mathbb{Q}}_U$ *denotes the constant sheaf on an open admissible subset* $U$ *for the stratification of the* $\infty$*-topos constructed from the underlying topological space of* $\mathfrak{X}$

#### 3.3.1.1 The dual TNN Grassmannian spectral stack.

In the case of the TNN Grassmannian, we have the following:

DEFINITION 3.3.3 : *The dual TNN Grassmannian endowed with its dual module geometry:*

$$\left(\mathfrak{Gr}^{\geqslant 0}_{\text{ét}}, \mathscr{O}_{\mathfrak{Gr}^{\geqslant 0}_{\text{ét}}}, \mathcal{IC}_{\mathfrak{Gr}^{\geqslant 0}_{\text{ét}}}\right)^{\star} := \left(\mathfrak{Gr}^{\geqslant 0\ \star}_{\text{ét}}, \mathscr{O}_{\mathfrak{Gr}^{\geqslant 0\ \star}_{\text{ét}}}, \mathcal{IC}_{\mathfrak{Gr}^{\geqslant 0\ \star}_{\text{ét}}}\right)$$

*defines, for short, the* dual Grassmannian Deligne-Mumford stack*:*

$$\left(\mathfrak{Gr}^{\geqslant 0\ \star}_{\text{ét}}, \mathcal{IC}_{\mathfrak{Gr}^{\geqslant 0\ \star}_{\text{ét}}}\right)$$

#### 3.3.1.2 The dual Amplituhedron spectral stack.

In the case of the Amplituhedron, we have the following:

DEFINITION 3.3.4 : *The dual Amplituhedron endowed with its dual module geometry:*

$$\left(\mathcal{A}^{\text{ét}}_{n,k,4}, \mathscr{O}_{\mathcal{A}^{\text{ét}}_{n,k,4}}, \mathcal{IC}_{\mathcal{A}^{\text{ét}}_{n,k,4}}\right)^{\star} := \left(\mathfrak{Gr}^{\geqslant 0\,\star}_{\text{ét}}, \mathscr{O}_{\mathfrak{Gr}^{\geqslant 0\,\star}_{\text{ét}}}, \mathcal{IC}_{\mathfrak{Gr}^{\geqslant 0\,\star}_{\text{ét}}}\right)$$

*defines, for short, the* dual Amplituhedron Deligne-Mumford stack*:*

$$\left(\mathcal{A}^{\text{ét}\,\star}_{n,k,4}, \overset{\star}{\mathscr{O}}_{\mathcal{A}^{\text{ét}}_{n,k,4}}, \mathcal{IC}_{\mathcal{A}^{\text{ét}\,\star}_{n,k,4}}\right)$$

## 3.4.Conclusion of Part II: The Dual Amplituhedron.

### 3.4.1 The existence of the dual Amplituhedron and its Schubert geometry.

On the same path as described in J. Lurie's "*Structured spaces*" [91] and "*Spectral schemes*" [92], we have associated a geometry to a spectrally Deligne-Mumford stack which suits well with Schubert Geometry. As well as in the case of Zariski topology, where J. Lurie associated $\mathcal{G}_{\text{ZAR}}$, and in the case of étale topology, where he associated $\mathcal{G}_{\text{ÉT}}$, here, we associate a module geometry $\mathcal{G}_{\text{PERV}}$ incarned in the complex intersection sheaf $\mathscr{F}_{\mathfrak{X}} := \mathcal{IC}^{\bullet}$.

But we saw that this intersection complex sheaf verifies the transfer of geometries. That is the functor Amplituhedron:

$$\mathfrak{Z} : \mathfrak{G} := \left(\mathfrak{Gr}^{\geqslant 0}_{\text{ét}}, \mathscr{O}_{\mathfrak{Gr}^{\geqslant 0}_{\text{ét}}}, \mathcal{IC}_{\mathfrak{Gr}^{\geqslant 0}_{\text{ét}}}\right) \longrightarrow \left(\mathcal{A}^{\text{ét}}_{n,k,4}, \mathscr{O}_{\mathcal{A}^{\text{ét}}_{n,k,4}}, \mathcal{IC}_{\mathcal{A}^{\text{ét}}_{n,k,4}}\right) =: \mathcal{A}^{\text{ét}}_{n,k,4} \quad (3.4.1)$$

sends the $\mathcal{IC}$-geometry $\mathcal{IC}_{\mathfrak{Gr}^{\geqslant 0}_{\text{ét}}}$ of $\mathfrak{Gr}^{\geqslant 0}_{\text{ét}}$ to the $\mathcal{IC}$-geometry $\mathcal{IC}_{\mathcal{A}^{\text{ét}}_{n,k,4}}$ of $\mathcal{A}^{\text{ét}}_{n,k,4}$.

As we saw in subsection 2.6.6 and subsection 3.2.1, let us recall the following:

Let $f : X \to Y$ be a projective map of quasi-projective varieties.

THEOREM 3.4.2 ([66], Thm. 8.4.3) : *The right derived functor $\mathbb{R} f_*$ of the induced map $f_*$ sends the category of perverse sheaves on $X$ on the category of perverse sheaves on $Y$. Moreover, we have a decomposition into simple objects in the category of perverse sheaves on $Y$ of the pushforward $\mathbb{R} f_* \mathcal{F}$ for any constant sheaf $\underline{\mathbb{C}}_X$ such that:*

$$\mathbb{R}\, f_* \underline{\mathbb{C}}_X \simeq \bigoplus_\alpha \mathcal{IC}^\bullet_{\mathcal{L}_{S_\alpha}}[\ell_\alpha]$$

*where the $\mathcal{L}_{S_\alpha}$ are irreductible local systems on the closure $\overline{S_\alpha}$ of connected strata $S_\alpha$ of $Y$, and the $[\ell_\alpha] \in \mathbb{Z}$ shift the complexes $\mathcal{IC}^\bullet_{\mathcal{L}_{S_\alpha}} := \mathcal{IC}^\bullet(\overline{S_\alpha}, \mathcal{L}_{S_\alpha})$.*

Here, each $S_\alpha$ is a $4k$-dimensional tile of the Amplituhedron tiling, contractible, that is : $\mathring{\Pi}_f$ so is $\overline{\mathring{\Pi}_f} =: \Pi_f$, and we have then each local system $\mathcal{L}_{\Pi_f} \simeq \mathbb{C}[m]$ for a certain shift $m \in \mathbb{N}$ .

So the sum can be written as semi-simple connective DG algebra:

$$\mathbb{R}\, \mathfrak{Z}_* \underline{\mathbb{C}}_{\Pi_f} \simeq \oplus_{n\in\mathbb{N}} \mathbb{C}[n]$$

We remark that $\mathbb{R}\, \mathfrak{Z}_* \underline{\mathbb{C}}_{\Pi_f}$ is in degree $\geqslant 0$ since the DG Algebra is connective ($\mathfrak{G}$ and $\mathfrak{Y}$ are coaffine so coconnective stacks) and the sum is finite since $\mathfrak{G}$ and $\Pi_f$ for each $f$ have finite étale cohomological dimension.

Following ([103] Chap 11, et Toën, Antieau, et Chang-Yeong Chough), we look for the group $H^2(\mathcal{A}^{\text{ét}}_{n,k,4}, \mathcal{IC}^\times_{\mathcal{A}^{\text{ét}}_{n,k,4}})$.

The inversibles of $\mathcal{IC}^\bullet_{\mathcal{A}^{\text{ét}}_{n,k,4}}$ as $\mathbb{E}_\infty$-ring spectra correspond to perfect (that is quasi-isomorphic to dualisable) perverses (thanks to Verdier duality) non-zero cochains.

The extended Brauer group of the Amplituhedron will then be non-trivial, as soon as there is at least one cell in the Amplituhedron tiling.

Moreover, for an open cell $\mathring{\Pi}_f$, which is an affine variety given by the cluster structure of (the affinization of) its ring of coordinates, and can thereby be written $\mathring{\Pi}_f = \mathbb{C}[x_1, x_2, \ldots, x_n]/\mathcal{I} =: \mathsf{Spec}\ \ \pi_0 R$, we have:

corollary 3.4.3 ([6] Cor 7.13) : *If $R$ is a connective $\mathbb{E}_\infty$-ring spectrum, then homotopy*

*groups of* $\mathbf{Br}(R)$ *are described by:*

$$\pi_k\mathbf{Br}(R) \simeq \begin{cases} H^1_{\text{ét}}(\mathsf{Spec}\ \pi_0 R, \mathbb{Z}) \times H^2_{\text{ét}}(\mathsf{Spec}\ \pi_0 R, \mathbb{G}_{\mathrm{m}}) & k=0 \\ H^0_{\text{ét}}(\mathsf{Spec}\ \pi_0 R, \mathbb{Z}) \times H^1_{\text{ét}}(\mathsf{Spec}\ \pi_0 R, \mathbb{G}_{\mathrm{m}}) & k=1 \\ \pi_0 R^{\times} & k=2 \\ \pi_{k-2} R & k \geqslant 3 \end{cases} \tag{3.4.4}$$

In our situation:

THEOREM 3.4.5 ([103] Construction 11.5.5.3) : *If* $(\mathring{\mathfrak{P}}_f, \mathscr{O}_{\mathring{\mathfrak{P}}_f}, \mathcal{IC}^{\bullet}_{\mathring{\mathfrak{P}}_f})$ *is the spectral Deligne-Mumford stack induced from* $\mathring{\Pi}_f$*, then the homotopy groups of the* Brauer sheafs $\underline{\mathscr{B}\mathrm{r}}^{\dagger}_{\mathring{\mathfrak{P}}_f}$ *are described by:*

$$\pi_k \underline{\mathscr{B}\mathrm{r}}^{\dagger}_{\mathring{\mathfrak{P}}_f} \simeq \begin{cases} 0 & k=0 \\ \underline{\mathbb{Z}} & k=1 \\ (\pi_0 \mathcal{IC}^{\bullet}_{\mathring{\mathfrak{P}}_f})^{\times} & k=2 \\ \pi_{k-2}\mathcal{IC}^{\bullet}_{\mathring{\mathfrak{P}}_f} & k \geqslant 3 \end{cases} \tag{3.4.6}$$

*where* $\underline{\mathbb{Z}}$ *is the constant sheaf associated to* $\mathbb{Z}$*.*

Thanks to [103] Example 11.5.5.5, for the spectral Deligne-Mumford stack $\mathring{\mathfrak{P}}_f$, we have:

$$\underline{\mathscr{B}\mathrm{r}}^{\dagger}_{\mathring{\mathfrak{P}}_f} \simeq H^2(\mathring{\mathfrak{P}}_f, \mathcal{IC}^{\bullet,\times}_{\mathring{\mathfrak{P}}_f}) \times H^1(\mathring{\mathfrak{P}}_f, \underline{\mathbb{Z}})$$

And since $\mathring{\mathfrak{P}}_f$ is a quasi-compact, quasi-separated, normal, locally Noetherian spectral algebraic space, then we obtain:

$$\underline{\mathscr{B}\mathrm{r}}^{\dagger}_{\mathring{\mathfrak{P}}_f} \simeq H^2(\mathring{\mathfrak{P}}_f, \mathcal{IC}^{\bullet,\times}_{\mathring{\mathfrak{P}}_f})$$

That is:

$$\underline{\mathscr{B}\mathrm{r}}^{\dagger}_{\mathring{\mathfrak{P}}_f} \simeq H^2(\mathring{\mathfrak{P}}_f, \mathcal{IC}^{\bullet,\times}_{\mathring{\mathfrak{P}}_f}) \simeq (\pi_0 \mathcal{IC}^{\bullet}_{\mathring{\mathfrak{P}}_f})^{\times}$$

But $H^2(\mathring{\mathfrak{P}}_f, \mathcal{IC}^{\times}) = c_1(\mathring{\mathfrak{P}}_f)$ for $k=2$, the first Chern class of the corresponding positroid variety, since the étale map gives an isomorphism of the Brauer group.

And this first Chern class of an (open) positroid variety is $\mathbb{Z}$, with one generator **1** since it is an affine variety given by the cluster structure of its (affinization of the) ring of coordinates.

We will remark that the twisted sheaf has -**1** as generator of its first Chern class.

*And, again, if one of the tiles flips over, then all the tiles flip over due to their overlap along the gluing of the closed immersions with common closed intersections.*

We can then conclude:

### 3.4.2 Main theorems of Part II.

#### 3.4.2.1 First version of the main theorems of Part II.

THEOREM 3.4.7 (Main theorem 1: Existence of the Dual Amplituhedron) : *Let $Z$ be a full rank $(k+4)\times n$ real matrix in $\mathbb{R}^{(k+4)\times n}$ with all maximal minors are positive, inducing a map:*

$$\widetilde{Z} : \mathrm{Gr}_{k,n}^{\geqslant 0} \longrightarrow \mathrm{Gr}_{k,k+m}$$

*defined by*

$$\widetilde{Z}\left(\langle v_1,\ldots,v_k\rangle\right) := \langle Z(v_1),\ldots Z(v_k)\rangle$$

*where $\langle v_1,\ldots,v_k\rangle$ represents with $k$ basis vectors, an element of $\mathrm{Gr}_{k,n}^{\geqslant 0}$, the real TNN Grassmannian, defines the so-called* (tree) Amplituhedron $\mathscr{A}_{n,k,4}$ *as the image $\widetilde{Z}\left(\mathrm{Gr}_{k,n}^{\geqslant 0}\right)$ of $\mathrm{Gr}_{k,n}^{\geqslant 0}$ inside $\mathrm{Gr}_{k,k+m}^{\geqslant 0}$.*
*Then we can construct a* unique non-trivial topological dual Amplituhedron $\mathscr{A}_{n,k,4}^{\star}$ *to this (tree) Amplituhedron $\mathscr{A}_{n,k,4}$ satisfying the following properties:*

- *the prestable $\infty$-category* $\mathsf{QCoh}(\mathscr{A}_{n,k,4}^{\star})$ *of quasi-coherent $\infty$-sheaves on $\mathscr{A}_{n,k,4}^{\star}$ is the monoidal dual of* $\mathsf{QCoh}(\mathscr{A}_{n,k,4})$ *for Lurie's tensor product of $\infty$-categories in the subcategory* $\mathsf{QStk}(\mathscr{A}_{n,k,4})$ *of the $\infty$-category of Grothendieck prestable tensor $\infty$-categories* $\mathsf{Groth}_{\infty}$*, where the $\infty$-subcategory* $\mathsf{QStk}(\mathscr{A}_{n,k,4})$ *is made of quasi-coherent $\infty$-stacks on $\mathscr{A}_{n,k,4}$ endowed with the étale $\infty$-topology*

- *this dual Amplituhedron $\mathscr{A}_{n,k,4}^{\star}$ has a canonical $\infty$-sheaf of $\mathbb{E}_{\infty}$-rings, $\mathscr{O}_{\mathscr{A}_{n,k,4}}^{\star}$, which is the Verdier dual to the canonical $\infty$-sheaf $\mathscr{O}_{\mathscr{A}_{n,k,4}}$ in the sense that $\mathscr{O}_{\mathscr{A}_{n,k,4}}^{\star} = \mathrm{Hom}_{\mathcal{D}^+(\mathsf{Coh}(\mathscr{A}_{n,k,4}))}(\mathscr{O}_{\mathscr{A}_{n,k,4}},\omega_{\mathscr{A}_{n,k,4}})$ where $\omega_{\mathscr{A}_{n,k,4}}$ is the dualizing complex in the $\infty$-category of $\mathbb{E}_{\infty}$-spectra.*

THEOREM 3.4.8 (Main theorem 2: Consistency of the $\mathcal{BCFW}$-structure cells of the topological Dual Amplituhedron) : *The unique non-trivial topological dual Amplituhedron*

*$\mathscr{A}_{n,k,4}^{\star}$ constructed in theorem 3.4.7 respects the Schubert cells of the (tree) Amplituhedron $\mathscr{A}_{n,k,4}$, that is it satisfies the following properties, such that the dual Amplituhedron admits a tiling in dual $\mathcal{BCFW}$ open positroid cells whose Zariski closure overlap themselves along other dual $\mathcal{BCFW}$ open positroid cells, that is:*

- $\widetilde{Z}(\Pi_f)^{\star} = \overline{\widetilde{Z}(\mathring{\Pi}_f^{\star})}$
- $\widetilde{Z}(\Pi_f)^{\star} \setminus \widetilde{Z}(\mathring{\Pi}_f)^{\star} = \coprod_{g \in \mathscr{B} \subset \mathscr{B}(k,n)} \widetilde{Z}(\mathring{\Pi}_g)^{\star}$
- *this $4k$-dimensional $\mathcal{BCFW}$ tiling is a perfect partition* $\mathscr{A}_{n,k,4}^{\star} = \coprod_{f \in \mathscr{B}(k,n)} \widetilde{Z}(\mathring{\Pi}_f)^{\star}$

*Another characterization is that the $\infty$-topos constructed from the topological dual Amplituhedron $\mathscr{A}_{n,k,4}^{\star}$ has a $\mathcal{G}^{\star}$-structure in the sense of Lurie ([91]) where $\mathcal{G}^{\star}$ is the $\infty$-topos spanned by the perfect perverse complexes of $\infty$-sheaves ([103] Chap 9) on $\mathscr{A}_{n,k,4}^{\star}$, that is Deligne's intersection complexes of $\infty$-sheaves with respect to the initial Whitney stratification, which are dualizable as compact objects, and preserved as perfect objects by the decomposition theorem for perverse sheaves, as the Amplituhedron map $Z$ is small and proper.*
*This means that the dual Amplituhedron $\mathscr{A}_{n,k,4}^{\star}$ admits a $\mathcal{G}^{\star}$-structure, where $\mathcal{G}^{\star}$ is the $\infty$-topos spanned by the perfect perverse complexes of $\infty$-sheaves in the heart $\mathsf{QCoh}(\mathscr{A}_{n,k,4}^{\star})^{\heartsuit}$ of the $\infty$-category $\mathsf{QCoh}(\mathscr{A}_{n,k,4}^{\star})$, which admits a triangulated $t$-structure preserved by duality.*
*Endly, a last identical characterization is that the Schubert geometry of the Amplituhedron is stable with respect to taking the dual $\mathscr{A}_{n,k,4}^{\star}$ of the Amplituhedron.*

theorem 3.4.9 (Main theorem 3: Volume of the Dual Amplituhedron) : *Recall theorem 2.10.30: theorem Considering the differential cohesive structured spectral Deligne-Mumford $\infty$-stack:*

$$\mathfrak{H}_{\mathfrak{th}} := \mathrm{Gr}(k,n)_{\mathrm{th}}$$

*is an $\mathcal{IC}^{\bullet}$-structured spectral algebraic space which represents the initial formal moduli problem definition 2.11.4.Its crystalline cohomology is:*

$$H^*_{\mathrm{dRh}}(\mathfrak{H}_{\mathfrak{th}}) := H^*_{\mathrm{dRh}}(\mathrm{Gr}(k,n)_{\mathrm{th}}, \mathbf{1}_{\otimes})$$

*The symetric monoïdality of* $\mathsf{QCoh}(\mathfrak{H}_{\mathfrak{th}})$ *and of* $\mathsf{QStk}(\mathfrak{H}_{\mathfrak{th}})$ *of (compact perfect) dualizable generator* $\mathscr{O}_{\mathfrak{H}_{\mathfrak{th}}}$ *allows us to apply Tannaka duality to recover the inital algebraic space, which then has a volume given by the exterior power of maximal rank of the De Rham derivative (of rank 1), that is the volume form on a differential compact space, that is:*

$$\mathrm{Volume}(\mathfrak{H}_{\mathfrak{th}}) = \int \mathrm{d}\Omega^n(\mathfrak{H}_{\mathfrak{th}}) = \mathrm{Volume}(\mathfrak{H})$$

*since the cohomology is insensible to the reduction. theorem We apply then this theorem to the dual Amplituhedron* $\mathscr{A}^{\star}_{n,k,4}$*:*

$$\mathrm{Volume}(\mathscr{A}^{\star}_{n,k,4}) = \int \mathrm{d}\Omega^n(\mathscr{A}^{\star}_{n,k,4}) = \mathrm{Volume}(\mathscr{A}^{\star}_{n,k,4})$$

#### 3.4.2.2 Second version of main theorems of Part II.

THEOREM 3.4.10 (First main theorem of Part II: The Dual Amplituhedron) : *We have the following assertions:*

- *The Amplituhedron functor sends the dual TNN Grassmannian spectral* $\mathcal{IC}$*-algebraic space on the dual Amplituhedron spectral* $\mathcal{IC}_{\mathcal{A}^{\text{ét}\,\star}_{n,k,4}}$*-algebraic space.*
- *The image is not empty since Lurie's extended Brauer-Grothendieck group of the dual Amplituhedron is isomorphic to Lurie's extended Brauer-Grothendieck group of the dual TNN Grassmannian:*

$$H^2_{\text{ét}}(\mathcal{A}^{\text{ét}\,\star}_{n,k,4}, \mathcal{IC}^{\times}_{\mathcal{A}^{\text{ét}\,\star}_{n,k,4}}) \simeq H^2_{\text{ét}}(\mathfrak{Gr}^{\geqslant 0\,\star}_{\text{ét}}, \mathcal{IC}^{\times}_{\mathfrak{Gr}^{\geqslant 0\,\star}_{\text{ét}}}) \simeq c_1(\mathfrak{Gr}^{\geqslant 0\,\star}_{\text{ét}}) \simeq c_1(\mathfrak{Gr}^{\geqslant 0}_{\text{ét}}) \simeq \mathbb{Z}$$

*since* $\mathfrak{Gr}^{\geqslant 0}_{\text{ét}}$ *is autodual. The first Chern class is spanned by one generator, that is the product of the first Chern classes for each BCFW cell in the étale (dual) TNN Grassmannian, since if one of the tiles flips over, then all the tiles flip over* *due to their overlap along the gluing of the closed immersions with common closed intersections. The dual twist is here given by the swap between the generators* **1** *and* -**1**.

remark 3.4.11 : The previous theorem 3.4.10 means exactly:

- The duality functor commutes with the Amplituhedron functor
- We constructed exactly one dual.

theorem 3.4.12 (Second main theorem of Part II: The Schubert Geometry of the Dual Amplituhedron) : *We have the following assertion:*

- *The Lurie's $\mathcal{IC}$-geometry of the dual TNN Grassmannian spectral $\mathcal{IC}$-algebraic space is sent on the dual Lurie's $\mathcal{IC}$-geometry of the dual Amplituhedron spectral $\mathcal{IC}$-algebraic space by a transformation of geometries between perfect perverse admissible geometries on the corresponding spectral algebraic spaces.*

remark 3.4.13 : The previous theorem 3.4.12 means exactly:

- The duality functor respects the Schubert geometry of the TNN Grassmannian since Lurie's geometry of the corresponding dual Amplituhedron spectral algebraic space is the image of Lurie's geometry of the corresponding coaffine dual TNN Grassmannian spectral algebraic space by a transformation of Lurie's geometries at the level of perfect perverse $\mathcal{IC}$-structures.

*Proof.*

- The decomposition theorem for perverse sheaves sends local systems on the TNN Grassmannian, which reflects the Schubert geometry, to semi-simple intersection complexes of local systems on subvarieties compatible with the Whitney BCFW étale stratification of the dual Amplituhedron.
- The Amplituhedron functor sends perverse sheaves on perverse sheaves.

- These semi-simple Deligne's intersection complexes of local systems on subvarieties of the dual Amplituhedron are perfect spectral Deligne-Mumford stacks, which allow us to use Lurie's framework of spectral algebraic geometry.

□

THEOREM 3.4.14 (Third main theorem of Part II: The Volume of the Dual Amplituhedron) : *Recall theorem 2.10.30: theorem Considering the differential cohesive structured spectral Deligne-Mumford ∞-stack:*

$$\mathfrak{H}_{\mathfrak{th}} := \mathrm{Gr}(k,n)_{\mathrm{th}}$$

*is an $\mathcal{IC}^{\bullet}$-structured spectral algebraic space which represents the initial formal moduli problem definition 2.11.4.Its crystalline cohomology is:*

$$H^*_{\mathrm{dRh}}(\mathfrak{H}_{\mathfrak{th}}) := H^*_{\mathrm{dRh}}(\mathrm{Gr}(k,n)_{\mathrm{th}}, \mathbf{1}_{\otimes})$$

*The symetric monoïdality of* $\mathsf{QCoh}(\mathfrak{H}_{\mathfrak{th}})$ *and of* $\mathsf{QStk}(\mathfrak{H}_{\mathfrak{th}})$ *of (compact perfect) dualizable generator* $\mathscr{O}_{\mathfrak{H}_{\mathfrak{th}}}$ *allows us to apply Tannaka duality to recover the inital algebraic space, which then has a volume given by the exterior power of maximal rank of the De Rham derivative (of rank 1), that is the volume form on a differential compact space, that is:*

$$\mathrm{Volume}(\mathfrak{H}_{\mathfrak{th}}) = \int \mathrm{d}\Omega^n(\mathfrak{H}_{\mathfrak{th}}) = \mathrm{Volume}(\mathfrak{H})$$

*since the cohomology is insensible to the reduction. theorem We apply then this theorem to the dual Amplituhedron* $\mathscr{A}^{\star}_{n,k,4}$*:*

$$\mathrm{Volume}(\mathscr{A}^{\star}_{n,k,4}) = \int \mathrm{d}\Omega^n(\mathscr{A}^{\star}_{n,k,4}) = \mathrm{Volume}(\mathscr{A}^{\star}_{n,k,4})$$

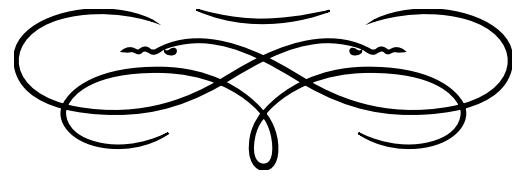